\documentclass[english,11pt]{amsart}
\usepackage[latin1]{inputenc}
\usepackage{color}
\usepackage{amsmath,amsthm,amssymb, babel}
\usepackage{bbm,graphicx}
\usepackage{caption}

\usepackage{hyperref}
\usepackage[notref,notcite]{}

\numberwithin{equation}{section}

\def\RR{{\mathbb R} }

\def\ve{{\varepsilon}}

 \title[A new mathematical   model for cell motility  with nonlocal repulsion from saturated areas]
 {A new mathematical   model for cell motility  with nonlocal repulsion from saturated areas}
 \author[C. Giambiagi Ferrari]{Carlo Giambiagi Ferrari}
 \author[F. Guill\'en-Gonz\'alez]{Francisco Guill\'en-Gonz\'alez}
\author[M. P\'erez-Llanos]{Mayte P\'erez-Llanos}
\author[A. Su\'arez]{Antonio Su\'arez}

\address{Carlo Giambiagi Ferrari,Francisco Guill\'en-Gonz\'alez, Mayte P\'erez-Llanos and Antonio Su\'arez \hfill\break\indent
Departamento de Ecuaciones Diferenciales y An\'alisis Num\'erico
\hfill\break\indent  Facultad de Matem\'aticas, U. de Sevilla,  C. Tarfia, s/n, 41012 Sevilla, Spain }
\email{{\tt cgf@gmail.com, guillen@us.es, mpperez@us.es, suarez@us.es }}


\begin{document}

\maketitle
\begin{abstract} The main purpose of  this work is the mathematical  modelling of large  populations of  cells under different deterministic interactions among themselves, in balance with naturally random diffusion. We focus on cell-cell adhesion mechanisms for a single population confined to an isolated domain.  Our most relevant contribution is to derive a  mathematical model including a nonlocal saturation coefficient as part of an appropriate nonlocal drift term,   including repulsion effects, depending on the level of saturation of  the area. For this purpose, we use two  discrete approaches  taking into account different perspectives: Eulerian and Lagrangian reference systems.

\end{abstract}

\section{Introduction}

Cellular movements are critical for a wide number of biological processes, including normal embryogenesis, tissue formation, wound healing and even defense against infection, among others. It is also an important factor in diseases such as cancer metastasis and birth defects.  It is well known that  cell interaction involves attraction and repulsion. Cells need to be sufficiently close to attract each other. At the same time, they also exhibit certain  random movement. The overall evolution is a result of the balance between the diffusion originated from the Brownian motion and the transport due to deterministic interactions. For instance, see  \cite{Berg} and \cite{Bressloff} for a general description on cell motion.  Cells might move attracted by other cells.However, as an area becomes too crowded,  the attraction weakens as other cells tend  to  less populated areas and a repulsion effect appears. This way  high accumulations at narrow zones or single points are prevented from occurring.

In the last decade,  introducing non-local terms in the classical PDEs governing many of the biological processes is an issue that is increasingly employed,  see for example \cite{kavallaris}. The process of gaining information about presence
or absence of other species in the environment is intrinsically non-local, be it higher level species or cells. Animals explore their surroundings  to find prey, avoid predators or aggregate in colonies, herds  or swarms, (we refer to e.g. \cite{carrillo-animal}, \cite{cucker-and-smale}, \cite{self-propelled}, \cite{swarm-nolocal}). This non-local sensing, which in animals is due mainly to using smell, hearing or sight, it also occurs at cellular level, through the extension of long thin protrusions, see \cite{Teddy}.

In mathematical modelling on the population level, non-local integro-differential operators allow to capture the drift motion caused by interactions between cells more accurately than local differential operators. This is due to their ability to take into account the effect of the surrounding  environment to describe what happens at certain point, in contrast to local differential operators. A recent review on this topic is \cite{painter-hillen-potts}. In this respect one pioneer reference is the work \cite{APS},
where a non-local equation was adopted for the first time in the context of cell-cell adhesion modelling. There, the velocity of the drift due to adhesion is defined so that its motion is described by an integral operator. Moreover, the authors point out the  ability of the model to replicate fundamental behaviour associated with cell-cell adhesion  and
the active sorting process of two or more cell types from a randomly distributed mixture. This model inspired a large number of subsequent contributions, including \cite{CMSTT}, \cite{MT}, \cite{painter-protusions},  \cite{apli-APS}, see also the review \cite{compisuru} and references therein.


According to the previous short overview on mathematical modelling of cell-to-cell interactions, non-locality in continuum models often arises
as an integral term inside a derivative. Consequently, many of the mathematical techniques developed to deal with classical PDEs (such us comparison and monotonicity methods, spectral theory, Lyapunov fuctionals, etc) need to be adapted or even redefined appropriately to apply them within this non-local framework, we refer for instance to \cite{hillen-painter-winkler}, \cite{libro-buten}, \cite{no-apli}.  However, it is not always possible to find suitable arguments based on classical theories, which gives rise to new challenging directions of work within the analysis of PDEs. As a counterpart, in  single-population models of aggregation, the nonlocal structure of the advective term
is fundamental to avoiding blow-up phenomena, that is, unbounded solutions at finite time, and ensuring global existence,  \cite{no-apli}, \cite{coron},   \cite{Franchesco}, \cite{Giunta}  or \cite{painter-esmais}, among others.


\

The primary focus of this work is to describe large populations of  cells: how do they evolve in confinement when  interaction among them is driven by attraction and avoidance of too saturated areas. Moreover,  we determine adequate microscopic rules  specifically regulating different mechanisms of attraction or repulsion,  evaluating the effect of each parameter or term in those relations and  their impact on the underlying macroscopic behavior.

\

The content of this paper is organized as follows: We begin by summarizing in Section \ref{resumen}, some existing literature that is a starting point of our work. In Section~\ref{euler} we propose several discrete models inspired by some of the references mentioned in  Section \ref{resumen}, thereby adopting the Eulerian perspective. Namely, defining an initial (discrete) density of $N$ cells,  we  determine  its changing rate at every time step, as a balance between  loss and  gain terms,  resulting from different deterministic rules in combination with random motion. This analysis will relate our models as $N\to \infty$ with certain PDEs containing nonlocal terms. In Section~\ref{lagra} this study will be complemented with a Lagrangian formulation, where we focus on the motion of individual cells, according to the mentioned rules. This leads to a system of $N$ SDEs describing the cell trajectories. As it could be expected, in the hydrodynamic limit as $N\to \infty$ this system recovers  nonlocal PDEs. However, both points of view are useful. On the one side, Lagrangian perspective allows us to follow every single cell, which is not possible when dealing with densities in the eulerian prespective. On the other side however, the trajectory simulations correspond to a reduced number of particles, thus  saturation effects are better described by the PDEs eulerian approach. At this point, the discrete  eulerian modelling could be viewed as a numerical discretization of such PDEs. Nevertheless,  due to the intricate nature of the problem we use appropriate  semi-explicit upwind schemes instead, to obtain reliable results when representing numerical simulations in Section \ref{section:nume}. Several observations for different interactions will be illustrated by a variety of examples, based on both, Lagrangian trajectories and PDEs discretizations. We conclude the paper with the main conclusions of our analysis in Section \ref{conclusiones}.

\

\section{Overview of related cell-cell adhesion models} \label{resumen}

We begin briefly reviewing the Armstrong Painter Sherratt $1D$ model (APS model), from \cite{APS}, where the density of a population of cells  $u=u(x,t)\ge 0$ evolves in a spatial $1D$ domain $x\in \Omega=(-L,L)$ during a time interval $t\in [0,+\infty)$  according to the nonlocal PDE
\begin{equation}\label{eq:APS}
u_t(x,t)=u_{xx}(x,t)-\gamma\,\partial_x\Big(u(x,t)K(u)(x,t)\Big),
\end{equation}
where
\begin{equation}\label{gradiente-nl}
K(u)(x,t)=\int_{0}^R 
 (g(u(x+r,t))-g(u(x-r,t)) ) \omega(r)dr.
\end{equation}

Here, $\gamma>0$ is a parameter representing the strength of the adhesion with respect to the diffusion. The function 
    $g(u(x+r,t))$ and $g(u(x-r,t))$ describes the dependence of the local force acting on the cell at position $x$ due to the cells at position $x+r$ or $x-r$.
    The weight $\omega(r)$ modifies the previous force according to the distance $r$ between the cells. The total force acting on the cell at position  $x$ is given by $\gamma \, K(u)(x,t)$. This global force accounts only for the action of the cells at distance from $x$ less than the sensing radius $R>0$. Namely, $\omega(r)=0$ for $r>R$ is imposed. 

The main novelty of \cite{APS} resides in the proposed adhesion term $\partial_x(u(x,t)K(u)(x,t))$, where $K(u)(x,t)$ is the drift nonlocal term given in \eqref{gradiente-nl}. The authors simplify the model  by taking $g(u)=u$ and $\omega\equiv 1$ in $[0,R]$, though warned that a non-constant function  would be more realistic. In this case, the integral term $K(u)(x,t)$ moves the cells towards more populated areas, giving raise to strong aggregation in several peaks, separated by $2R$ distance. In these examples periodic boundary conditions were considered by convenience. Further simulations are presented in \cite{APS}, taking $\omega \equiv 1$ in $[0,R]$ and $g(u)=u(1-u/K)_+$, which is  truncated at the crowding  capacity $K>0$. This modification ensures  that cells travel to moderately populated areas, softening the above mentioned peaks. In addition \cite{APS} also contains a study of  several kinds of interactions between two  populations of different type of cells in the two dimensional setting, which is out of the scope of the present work.

 Murakawa and Togashi found in \cite{MT} that, when extended to two cell types, the APS model  was unable to replicate properly certain cell sorting
patterns observed experimentally. Adopting a logistic function for $g(u)$, they proposed  to replace the linear laplacian diffusion by a porous media diffusion, obtaining a model that captures  more accurately the behaviour the authors expected, according to experimental results. Note that already  in the  case of a single population, this nonlinear porous media diffusion introduces intrinsically a repulsion effect from more crowded areas to empty spaces, that is stronger than in  the Laplacian case, where diffusivity is uniform across the domain.
 To avoid  cells drifting to medium populated areas, as it happens in \cite{MT}, the following modification of the drift velocity was proposed in \cite{CMSTT} 
$$
K(u)(x,t)=(1-u(x,t))\int_{0}^R \Big(u(x+r,t)-u(x-r,t)\Big) \omega(r)dr,
$$
where the crowding capacity  has been rescaled to $K=1$. In fact, the variable $u$ in this model can be interpreted as the volume fraction of cells, since it doesn't exceed one if the initial values are below one as well.
Solutions of the model in \cite{CMSTT} behave differently to those of the APS model. On the one hand,  even though cells are directed  to most populated areas  by the non-local integral term, the velocity of their movement is multiplied by the factor $(1-u(x,t))$. This velocity is slower as the density of cells reaches gradually higher values at point $x$. Therefore, this coefficient slows down cells motion according to the local crowdedness of the model at each $x$. On the other hand, the nonlinear diffusion of porous media type also considered in  \cite{CMSTT} introduces a repulsion effect to further avoid accumulating in high crowded areas.

We notice that  repulsion due to saturation is absent in the original APS model, while in   \cite{CMSTT} and \cite{MT} it is driven  by the nonlinear diffusion of porous media type, and not by the drift term. Some remarks in \cite{CMSTT} admit that  the weight $\omega(r)$ needs to take positive (aggregation) and negative (repulsion) values, depending on the distance among the interacting cells. See also \cite{carrillo2018}, where  an intercellular interaction kernel introduces change of sign of parabolic shape. Therefore, this kernel defines short, middle and long range interactions.  However, it is not taken into account that the saturation of the area should also lead to repulsion. The main aims of this paper is to explore the dynamics including repulsion by saturation as part of the drift term. In this respect the choice of the weight will contribute to regulate the effect of repulsion or aggregation.

%
%
%

\

\

Although \cite{APS}, \cite{CMSTT} and much of the existing literature on mathematical modelling of cell-cell adhesion consider periodic boundary conditions, with the consequent simplification of the simulation codes, this approach is just appropriate whenever the size of the domain can be considered infinitely large with respect to  cells size. Nevertheless, it is also relevant the study of single-cell migrating in isolation medium. The work \cite{butten-bounded} introduces several models in bounded domains considering either adhesive, or repulsive or even neutral boundary effects. Moreover, some numerical simulations are included illustrating the different boundary effects. 


In the above description, we have limited ourselves to summarizing the models that have mainly inspired our work. Nevertheless, as we already mentioned, there exists a wide variety of contributions based on the APS model due to Armstrong, Painter and Sherrat, see the book \cite{libro-buten} and references therein.

Summarizing, the main novelty of our modelling is that it captures more accurately the repulsion from saturated areas by a nonlocal factor included in an appropriate nonlocal drift term  in isolated domains.

\section{Euler perspective}\label{euler}

In this Section we develope various discrete models under the Euler perspective. We use the approach similar to \cite{PBS} (see also \cite{AS}  for a similar modelling approach to a biological application). For simplicity, we work in a one-dimensional space  and  consider a domain $\Omega=(-L,L)$. The extension to higher dimensions is straightforward, but would cause high computational cost in the subsequent simulations.

We construct models in which the unknown variable is the proportion of cells at certain position, combining infinitesimal motion directed by a nonlocal signal depending on the nonlocal term given in \eqref{gradiente-nl} with random diffusive processes.

\subsection{Generic balance} \label{Sec:3}

Consider a population of $N\in\mathbb{N}$ cells that we enumerate from $1$ to $N$. At time $t \ge 0$ we assign a real number $x_i(t)\in \overline\Omega=[-L,L]$,  representing the position  of the cell number $i\in\{1,\cdots,N\}$.  The cell 's position can only change due to a diffusion process and  different nonlocal interactions with other cells, that are at a distance smaller than a sensing radius $R>0$. 

In order to obtain master equations, let us subdivide the interval $[-L,L]$ in a family of intervals $\{ I_j \}, j\in \{1,\cdots,M\}$ of equal length $h=2L/M$.  In the sequel, it will be convenient to declare as $y_j$ the midpoint of the interval $I_j$. For $t\geq0$, we denote by 
$$s_j(t) = \frac{\# \{i : x_i(t) \in I_j \}}{ N}$$
the proportion of cells with position  within  the interval $I_j$ at time $t$. Obviously, $\sum_{j=1}^Ms_j(t)=1$ for all $t\geq0$.

As in \cite{PBS}, we want to derive a system of $M$ recursive discrete equations based on the assumption that the proportion of cells at the interval $I_j$ at stage $t+\Delta t$, namely $s_j(t+\Delta t)$, equals the proportion of cells at $I_j$ in the previous stage $s_j(t)$, plus certain gain and loss terms, representing the changes as cells enter and leave $I_j$, respectively, under the action of intercell interactions and diffusion. Obviously, the magnitude of these gain and loss terms depend on how long we consider the time step $\Delta t$ to be. At this point we stress an important difference with respect to \cite{PBS}.  In our model, the gain and loss terms result from two different effects, interactions and diffusion, which do not occur necessarily with the same frequency. The detailed analysis below shows that indeed, each of these processes  has its own scale.
Let us fix the characteristic time step $\Delta t$.  Then, 
\begin{equation}\label{general}s_j (t + \Delta t) = s_j ( t) +\Delta t\left[Q^h_{\mbox{diff}}\bigg(G_{\mbox{diff}}(j,t)- L_{\mbox{diff}}(j,t)\bigg)  + Q^h_{\mbox{int}}\bigg( G_{\mbox{int}}(j,t) -L_{\mbox{int}}(j,t)\bigg)\right],
\end{equation}
for $j\in \{1,\cdots, M\}$, where $G_{\mbox{diff}}(j,t)$ and $G_{\mbox{int}}(j,t)$ stand for probability gain rates while $L_{\mbox{diff}}(j,t)$ and $L_{\mbox{int}}(j,t)$ represent  loss rates, corresponding to diffusion and interaction effects, respectively.  The coefficients $Q^h_{\mbox{diff}}$ and $Q^h_{\mbox{int}}$ account for the frequency of diffusion and interaction processes, respectively, per unit of time.
The key to assemble correctly these two mechanisms  of motion is precisely the determination of the exact dependence on $h$ of coefficients $Q^h_{\mbox{diff}}$ and $Q^h_{\mbox{int}}$.

The study of the resulting system of $M$ equations (for $M$ large) is easier  if considering its reformulation for probability density. To this end, we introduce  a piecewise constant in space function $u_h$, 
representing a discrete probability density of the cell population, namely
\begin{equation}\label{ucont}
u_h(y,t): [-L,L] \times[0 ,\infty) \to \mathbb R_0^+\; \mbox{ such that  } (u_h(\cdot,t))|_{I_j}=s_j(t)/h .
\end{equation}
 In particular, $s_j ( t) = \int_{I_j} u_h (y, t) dy$. 

From now on we focus on the system of equations (\ref{general}) and rewrite  it in terms of the  function $u_h$.  Assuming that there is a limit function $u(x,t)$ such that $u_h\to u$ as $h\to 0$ in some sense  and also that $(u_h(\cdot,t)|_{I_j})$ is bounded as $h$ goes to zero, our aim is to determine (formally) which are the different  PDEs satisfied by $u$, depending on the model under consideration. Notice that for example, $u_h$ could be regarded as an appropriate numerical approximation of $u$.

 \subsection{Motion due to aggregation through nonlocal terms \cite{APS}.}  
 \label{Sec:APS}
 We begin our study with the simplest case of the APS model given in \cite{APS}, for which the direction of the transport term is determined by  nonlocal terms depending on  the density of cells. The action of the surrounding cells  could be  neglected    beyond certain threshold, $R>0$, known as the {\it sensing cell radius}. Moreover, this action could also decrease with respect to  the distance between the interacting cells, trough a  weight $w=w(r)$ with 
\begin{equation}\label{peso}w:[0,+\infty)\to \mathbb{R}_0^+  \mbox{ such that $w(r)>0$ if $r\in[0,R)$ whereas $w(r)=0$ if $r\geq R$.}
\end{equation}
 Interaction forces between cells belonging to  $I_i$ and $I_j$, respectively, will be weighted by $w(|y_j-y_i|)$ being $y_k$ the midpoint of each $I_k$.
 This model corresponds to the following  gain and loss terms:  the proportion of cells arriving to $I_j$ and leaving $I_j$ are, respectively
 \begin{align*}
  G_{\mbox{int}}(j,t)&=s_{j+1}(t) P_{I_{j+1}\to I_{j}}(t)+s_{j-1}(t)P_{I_{j-1}\to I_{j}}(t),\\
 L_{\mbox{int}}(j,t)&=s_{j}(t) \Big[P_{I_j\to I_{j-1}}(t)+P_{I_j\to I_{j+1}}(t)\Big].
 \end{align*}
 
 Notice that space motion is considered to be local, 
 because the arriving cells to interval $I_j$ come just from the adjacent positions, $I_{j-1}$, $I_{j+1}$. On the other hand, the cells leaving $I_j$ can just move to wether $I_{j-1}$ or $I_{j+1}$.
    
 In contrast, the transition probabilities are nonlocal. They depend on the neighbouring  configuration of the population inside  the sensing radius $R$ from interval $I_j$. Namely, related to the gain term
 \begin{equation}\label{prob-left}
P_{I_{j+1}\to I_{j}}(t)=\sum_{i=j-r}^{j}s_i(t)w(y_{j+1}-y_i),\quad\quad P_{I_{j-1}\to I_{j}}(t)=\sum_{i=j}^{j+r}s_i(t)w(y_i-y_{j-1}),
\end{equation}
and similarly, regarding the loss term, 
$$
 P_{I_j\to I_{j-1}}(t)=\sum_{i=j-1-r}^{j-1}s_i(t)w(y_j-y_i)\quad\quad P_{I_j\to I_{j+1}}(t)=\sum_{i=j+1}^{j+1+r}s_i(t)w(y_i-y_j).
 $$
Here, the natural number $r$ is in fact $r=r(h)$, such that $r(h)=[R/h]$, where $[x]$ denotes the integer part of $x$. 
 To shorten the notation, we introduce  the following terms 
$$
 \hbox{Left}(j,t):=s_{j+1}(t)P_{I_{j+1}\to I_{j}}(t) - s_{j}(t) P_{I_{j}\to I_{j-1}}(t), 
 $$
 $$
 \hbox{Right}(j,t):=s_{j-1}(t) P_{I_{j-1}\to I_{j}}(t) - s_{j}(t) P_{I_{j}\to I_{j+1}}(t),
 $$
hence
 \begin{equation}\label{L-R}
G_{\mbox{int}}(j,t) -L_{\mbox{int}}(j,t) = \hbox{Left}(j,t) +  \hbox{Right}(j,t).
\end{equation}

Note that, by the definition (\ref{ucont}) and assuming that $w|_{[0,R]}\in C^2([0,R])$, by the midpoint quadrature rule, and recalling that $ u_h(y_i,t)$ is bounded with respect to $h$,
one has 
 \begin{equation} \label{ucontcarrillo}
s_i(t)w(|y_{j+1}-y_i|)=\int_{I_i}u_h(y,t)w(|y_{j+1}-y|)dy+ O(h^3),
\end{equation}
 hence, from \eqref{prob-left},
 \begin{align*}\label{grad}
 P_{I_{j+1}\to I_{j}}(t)
 =\int_{\bigcup\limits_{i=j-r}^jI_i }  u_h(y,t)w(y_{j+1}-y)dy+O(h^2). 
 \end{align*}
 Similarly, one has 
 $$
 P_{I_{j-1}\to I_{j}}(t)
=\int_{\;\;\bigcup\limits_{i=j}^{j+r}I_i }u_h(y,t)w(y-y_{j-1})dy+O(h^2).
 $$
 At this point, we define
the following operators acting on a generic  function $v$:
$$
\mathcal{G}^-v(x,t)=   \int_{x-R}^{x}v(y,t)w(x-y)dy,
 \quad
 \mathcal{G}^+v(x,t)=  \int_{x}^{x+R}v(y,t)w(y-x)dy,
$$
and its discrete version acting on a generic piecewise constant function $v_h$ on intervals $I_j$, for $j\in \{1,\cdots, M\}$:
 \begin{equation}\label{grad}\begin{array}{ll}
 \displaystyle\mathcal{G}_h^-v_h(y_{j+1},t)&:= \displaystyle\int_{\bigcup\limits_{i=j-r}^jI_i }  v_h(y,t)w(y_{j+1}-y)dy,\\
  \displaystyle\mathcal{G}_h^+v(y_{j-1},t)&:= \displaystyle\int_{\;\;\bigcup\limits_{i=j}^{j+r}I_i }v_h(y,t)w(y-y_{j-1})dy. 
 \end{array}
 \end{equation}
One arrives then, recalling that $u_h(y_j,t)$ is bounded,  at 
\begin{align*}
 \hbox{Left}(j,t)
 &=\bigg[h \, u_h(y_{j+1},t)\bigg] \bigg[\mathcal{G}_h^-u_h(y_{j+1},t)+O(h ^2)\bigg]-\bigg[h \, u_h(y_{j},t)\bigg]\bigg[\mathcal{G}_h^-u_h(y_{j},t)+O(h^2)\bigg]\\
 &= h^2 \left[\frac{u_h(y_{j+1},t)\mathcal{G}_h^-u_h(y_{j+1},t) - u_h(y_{j},t)\mathcal{G}_h^-u_h(y_{j},t)}{h} \right]+ O(h^3).
 \end{align*}
 Analogously,
  \begin{align*}
\hbox{Right}(j,t)
 =-h^2  \left[\frac{u_h(y_{j},t)\mathcal{G}_h^+u_h(y_{j},t) - u_h(y_{j-1},t)\mathcal{G}_h^+u_h(y_{j-1},t)}{h}\right]+ O(h^3).
 \end{align*}
 
 We combine the left and right terms and denote forward and backward differences with respect to $y$ variable by
 \begin{equation}\label{deriv}
 \partial_h^+v(y,t):=\frac{v(y+h,t)-v(y,t)}h, \;\;\; \mbox{ and }\;\;\; \partial_h^-v(y,t):=\frac{v(y,t)-v(y-h,t)}h.
 \end{equation}
 Then, according to \eqref{L-R}, we have shown that
 \begin{equation}\label{nlg}
 G_{\mbox{int}}(t,j)-L_{\mbox{int}}(t,j)= -h^2 \mathcal T^h_{aps}(y_j,t) + O(h^3) ,
 \end{equation}
 where 
 \begin{align}\label{aps:inte}
 \mathcal T^h_{aps}(y_j,t):&=\partial_h^-\bigg[u_h(y_j,t)\mathcal G_h^+u_h(y_j,t)\bigg]-\partial_h^+\bigg[u_h(y_j,t)\mathcal G_h^-u_h(y_j,t)\bigg].
 \end{align}
Let $x\in[-L,L]$ be an arbitrary point. For each  uniform partition of length  $h$, we can  choose an index $j=j(h)$ such that $y_j\to x$ as $h\to 0$. In conclusion, assuming that $u_h\to u$ as $h\to 0$,  the following approximation can be (formally) deduced
 $$
  \lim_{h\to 0}\mathcal T^h_{aps}(y_j,t)=\frac{\partial}{\partial x}\bigg( u(x,t)(\nabla_{NL}u)(x,t)\bigg),
 $$
 where the nonlocal term  $\nabla_{NL}u$ is defined as
 $$
 (\nabla_{NL}u)(x,t):= \mathcal{G}^+u(x,t)-\mathcal{G}^-u(x,t)=
  \int_{x}^{x+R}u(y,t)w(y-x)dy
  -\int_{x-R}^{x}u(y,t)w(x-y)dy,
 $$
which can be rewritten as 
\begin{equation}\label{NLG}
 (\nabla_{NL}u)(x,t)= 
  \int_{0}^{R}(u(x+r,t)-u(x-r,t))w(r)dr.
   \end{equation}

As observed in Section 2, in \cite{APS}  it is considered as part of the drift  the nonlocal term $\nabla_{NL}g(u)$  of the truncated logistic function $g(u)=u(1-u)_+$, see \eqref{gradiente-nl}.  This choice allows to avoid the strong aggregation observed in the basic APS model, though it does not fit  the point of view of our modeling.  Note that in case $g(u)=u(1-u)_+$, it is not possible to  write the motion as a balance between gain and loss terms.  As reviewed earlier in Section 2, in \cite{CMSTT}  the appropriateness of the nonlinear term inside the integral is discussed. Among other issues, in \cite{CMSTT} the authors  did not agree with the fact that cells are attracted by moderately populated areas. Accordingly, to soften the potential strong aggregation they  proposed  instead a modification of the APS model, where they kept $g(u)=u$ but multiplied the nonlocal term with a coefficient that ensures saturation. Inspired by this proposal, in the next section we derive  a discrete model  that in the macroscopic limit yields the model from  \cite{CMSTT}.

 \

 \subsection{Motion due to aggregation through nonlocal term with local saturation coefficient   \cite{CMSTT}.}   \label{Sec:3.3}
As before, the dynamics is determined by the result of the gain-loss balance
 $$
 G_{\mbox{int}}(t,j)-L_{\mbox{int}}(t,j)=
 s_{j+1}(t)P_{I_{j+1}\to I_{j}}(t) - s_{j}(t) P_{I_{j}\to I_{j-1}} 
 +s_{j-1}(t) P_{I_{j-1}\to I_{j}}(t) - s_{j}(t) P_{I_{j}\to I_{j+1}}, 
 $$
where now the transition probabilities include a coefficient 
 $(1-u_h(y_j,t))$, such that 
 \begin{align*}
 P_{I_{j+1}\to I_{j}}(t)&= \left(1-u_h(y_{j+1},t)\right)\sum_{i=j-r}^{j}s_i(t)w(y_{j+1}-y_i)
 \\ P_{I_{j-1}\to I_{j}}(t)&= \left(1-u_h(y_{j-1},t)\right)\sum_{i=j}^{j+r}s_i(t)w(y_i-y_{j-1}).
 \end{align*}

 Notice that in \cite{CMSTT} the variable $u(x,t)$ is regarded to be bounded by one, representing the proportion of cells at position $x$ at time $t$. Therefore, now we will assume that $0\le u_h(y_j,t)\le 1$.
Again, ${u_h}$ could be seen as a discretization of $u$ by an appropriate numerical scheme.

We argue as in the previous case.
By (\ref{ucontcarrillo}), 
%
the above probabilities are rewritten as

\begin{equation*}
 P_{I_{j+1}\to I_{j}}(t)
=(1-{u_h}(y_{j+1},t))\bigg[\mathcal{G}^-{u_h}(y_{j+1},t)+O(h^2)\bigg],
  \end{equation*}
  and similarly,
   \begin{equation*}
 P_{I_{j-1}\to I_{j}}(t)
=(1-{u_h}(y_{j-1},t))\bigg[\mathcal{G}^+{u_h}(y_{j-1},t)+O(h^2)\bigg],
 \end{equation*}
 
 with the operators $\mathcal{G}_h^-, \mathcal{G}_h^+$ specified in (\ref{grad}).
 Then, using (\ref{ucontcarrillo}) we can write the terms
 
 \begin{align*}
 \hbox{Left}(j,t)& 
= h ^2 \partial^+_h\bigg[ {u_h}(y_{j},t)(1- {u_h}(y_{j},t))\mathcal{G}^-{u_h}(y_{j},t)\bigg]+ O(h ^3),
\\[0.3cm]
  \hbox{Right}(j,t)&
  =h^2\partial^-_h\bigg[ {u_h}(y_{j},t)(1- {u_h}(y_{j},t))\mathcal{G}^+{u_h}(y_{j},t)\bigg]+ O(h^3),
 \end{align*}
 
 with the forward and backward differences with respect to $y$ variable defined in (\ref{deriv}).

 This proves that
 \begin{align}\label{carrillo}
  G_{\mbox{int}}(t,j)-L_{\mbox{int}}(t,j)=  - h^2 \mathcal  T^h_{sat}(y_j,t)+O(h^3),
 \end{align}
 where,

 \begin{align}\label{carrillo:inte}
 \mathcal T^h_{sat}(y_j,t):&=\partial_h^-\bigg[u_h(y_j,t)(1-u_h(y_j,t))\mathcal G_h^+u_h(y_j,t)\bigg]-\partial_h^+\bigg[u_h(y_j,t)(1-u_h(y_j,t))\mathcal G_h^-u_h(y_j,t)\bigg].
 \end{align}
 
 Passing again to the limit formally, assuming that $  {u_h}\to u$ as $h\to 0$, one arrives at
 $$\lim_{h\to0} \mathcal T^h_{sat}(y_j,t)=\frac{\partial}{\partial x}\bigg( u(x,t)(1- u(x,t))(\nabla_{NL}  u)(x,t)\bigg),$$
 with the nonlocal term $\nabla_{NL}u$  given in (\ref{NLG}). 

 \subsection{Motion due to aggregation/repulsion through nonlocal term with nonlocal saturation coefficient.}\label{nuestro}  One of our main objectives in this paper  is to include the effect of repulsion as part of the cells interaction, via a nonlocal coefficient. With this aim in mind,  let us consider now a nonlocal saturation coefficient in the transition probabilities for  cells positioned at $I_j$ as follows:
  \begin{equation}\label{satu-nolocal}
 1-\frac1K \sum_{i=j-r}^{j+r} s_i(t)=1-\frac1K \mathcal U_h(y_j,t),
 \end{equation}
 being
 $$
 \mathcal U_h(y_j,t)=\int_{\;\;\bigcup\limits_{i=j-r}^{j+r}I_i }u_h(y,t)dy.
 $$
 Here, the parameter $K\in (0,1)$ represents a threshold  of nonlocal type  for the crowding capacity in the interval $[x_j-rh,x_j+rh]$ for each $x_j\in [-L,L]$ (recall that $r=r(h)$, such that $r(h)=[R/h]$). 
We could also ponder the influence of neighboring cells through a non-constant weight $\widehat w = \widehat w (r)$, namely changing  \eqref{satu-nolocal} by
$$
 1-\frac1K \sum_{i=j-r}^{j+r}
s_i(t)\widehat w (|y_j-y_i|).
$$
To make the presentation    of the model  simpler, we take $\widehat w \equiv 1$. Indeed, as Figures \ref{fig:NL_K06pesolineal_histogramas} and \ref{fig:NL_K06pesolineal_driftsatu}  show, the choice $\widehat w(r)=(R-r)/R$  allows more aggregation for the same value of that crowding capacity. Indeed high accumulations arise for $\widehat w$ linear decreasing, in contrast to Figures \ref{fig:NL_K06_histogramas} and  \ref{fig:NL_K06_driftsatu} where $\widehat w$ is constant equal to $1$.

\

The nonlocal coefficient (\ref{satu-nolocal}) has the following effect on cell motion: a cell positioned at certain $I_j$ moves according to aggregation influences  whenever the number of cells around its position does not exceed the crowding capacity $K$. Above this threshold, the cell begins to move towards less populated areas as a result of repulsion. 
%
%

\

As before, cells arriving to $I_j$ and leaving from $I_j$ are given by the gain-loss balance:
$$
 G_{\mbox{int}}(t,j)-L_{\mbox{int}}(t,j)=
 s_{j+1}(t)P_{I_{j+1}\to I_{j}}(t) - s_{j}(t) P_{I_{j}(t)\to I_{j-1}} 
 +s_{j-1}(t) P_{I_{j-1}\to I_{j}}(t) - s_{j}(t) P_{I_{j}(t)\to I_{j+1}} 
 $$
 where the transition probabilities are driven now by two effects. Cells moving from $I_{j+1}$ towards $I_j$ (resp. from $I_{j-1}$ towards $I_{j}$) engage in two processes: aggregation produced by cells at left hand side inside the sensing radius (resp. at right hand side) and saturation produced by the existing cells at right hand side (resp. at left hand side). To be more precise, denoting as $(g)_-=\min\{g,0\}$ and $(g)_+=\max\{g,0\}$ (hence $g=(g)_+ + (g)_-$), we will split the nonlocal saturation coefficient (\ref{satu-nolocal}) into its positive and negative parts representing the conditional transition probabilities due to attraction or repulsion effects, respectively. Accordingly,   we define 
 \begin{align*}
 P_{I_{j+1}\to I_{j}}(t)&=-\left(1-\frac1K \mathcal U_h(y_{j+1},t)
 \right)_{-} \, 
 \sum_{i=j+2}^{j+2+r}s_i(t)w(y_i-y_{j+1}) \quad\mbox{(repulsion term)}
 \\&\quad 
 +\left(1-\frac1K \mathcal U_h(y_{j+1},t)
 \right)_+ \, 
 \sum_{i=j-r}^{j}s_i(t)w(y_{j+1}-y_i)\quad \mbox{(aggregation term)}
 \\&=: P_{I_{j+1}\to I_{j}}(t)_{(1)} + P_{I_{j+1}\to I_{j}}(t)_{(2)} \, ,
 \\ P_{I_{j-1}\to I_{j}}(t)&=
 -\left(1-\frac1K \mathcal U_h(y_{j-1},t)
 \right)_- \,
 \sum_{i=j-2-r}^{j-2}s_i(t)w(y_{j-1}-y_i) \quad\mbox{(repulsion term)}
 \\&\quad 
 +\left(1-\frac1K \mathcal U_h(y_{j-1},t)
 \right)_+ \,
 \sum_{i=j}^{j+r}s_i(t)w(y_i-y_{j-1})\quad \mbox{(aggregation term)}
 \\&=: P_{I_{j-1}\to I_{j}}(t)_{(1)} + P_{I_{j-1}\to I_{j}}(t)_{(2)}  \, .
 \end{align*}

 Now we proceed rewriting these probabilities in terms of  the function $u_h$  in (\ref{ucont}). For example, let us consider  terms $P_{I_{j+1}\to I_{j}}(t)_{(1)}$ and $P_{I_{j}\to I_{j-1}}(t)_{(1)}$, corresponding to gain and loss  by repulsion effects. 
 Indeed, we can express
\begin{align*}s_{j+1}(t)&P_{I_{j+1}\to I_{j}}(t)_{(1)}\\
&=- h\, u_h(y_{j+1},t) \left(1-\frac1K \mathcal U_h(y_{j+1},t))\right)_- \bigg[\mathcal G_h^+u_h(y_{j+1},t)+O(h^2)\bigg]\\
&=-h\bigg[u_h(y_{j+1},t)\left(1-\frac1K \mathcal U_h(y_{j+1},t)\right)_-\mathcal G_h^+u_h(y_{j+1},t)\bigg] +O(h^3),
\end{align*}
with $\mathcal G_h^+$ specified in  (\ref{grad}).  Similarly, one has that
$$
s_j(t)P_{I_{j}\to I_{j-1}}(t)_{(1)}=-h\bigg[u_h(y_{j},t)\left(1-\frac1K \mathcal U_h(y_{j},t)\right)_-\mathcal G_h^+u_h(y_{j},t)\bigg] +O(h^3) .
$$

In conclusion,
 \begin{align*}
 \hbox{Left}_{(1)}(j,t) &:=s_{j+1}(t) P_{I_{j+1}\to I_{j}}(t)_{(1)}-s_j(t)P_{I_{j}\to I_{j-1}}(t)_{(1)}\\[0.3cm]
 &=-h^2\left[\frac{u_h(y_{j+1},t)(1-\frac1K \mathcal U_h(y_{j+1},t))_-\mathcal G_h^+u_h(y_{j+1},t)-u_h(y_{j},t)(1-\frac1K \mathcal U_h(y_{j},t))_-\mathcal G_h^+u(y_{j},t)}{h} \right] +O(h^3) \\
 &=-h^2\partial^+_h\bigg[ u_h(y_{j},t)\left(1-\frac1K \mathcal U_h(y_{j},t)\right)_-\mathcal G_h^+(y_{j},t)\bigg] +O(h^3),
 \end{align*}
where we are using the notation in (\ref{deriv}) for the forward derivative.

 Similar arguments for the rest of the terms yield 
\begin{equation}\label{carlo}
 G_{\mbox{int}}(t,j)-L_{\mbox{int}}(t,j)\sim -h^2 \mathcal T^h_{NLsat}(y_j,t)+O(h^3),
 \end{equation}
 where
$$
 \mathcal T^h_{NLsat}(y_j,t):=\partial_h^-\bigg[u_h(y_j,t)\left(1-\frac1K \mathcal U_h(y_{j},t)\right)\mathcal G_h^+u_h(y_j,t)\bigg]-\partial_h^+\bigg[u_h(y_j,t)\left(1-\frac1K \mathcal U_h(y_{j},t)\right)\mathcal G_h^-u_h(y_j,t)\bigg].
 $$

 Passing to the limit (formally), we get 
 $$
 \lim_{h\to 0}\mathcal T^h_{NLsat}(y_j,t)=\frac{\partial}{\partial x}\left(u(x,t)\left(1-\frac1K \int_{x-R}^{x+R}u(y,t)dy\right)\nabla_{NL}u(x,t)\right),
 $$
 where $\nabla_{NL}u$ is specified in (\ref{NLG}).

 \subsubsection{Boundary effects.}  For simplicity, in the  descriptions of gain and loss terms given in Sections \ref{Sec:APS}, \ref{Sec:3.3} and \ref{nuestro}, we have not taken into account the effects of the boundary, so that our derivations are only valid for $x \in (-L+R,L-R)$. We are going to consider isolating boundary conditions, which means that particles 
  can neither enter nor leave the domain. This fact has an influence only for the closest intervals  to the boundary. 
   In particular,   if capacity is large, cells at the boundary are more likely to be  pushed towards the interior by the attraction of other cells, while  outside there are not cells, see figures in ${\S 5.2}$ and ${\S 5.3}$.
    Under strong repulsion the situation is completely different: cells could be pushed towards the boundary, where they then  form high accumulation peaks due to its impenetrability, 
    as examples show in ${\S 5.4}$ and ${\S 5.5}$.  
    
 
\

\subsection{Dynamics driven by the interplay between diffusion and nonlocal drift}
In this Section we study  the combined effect of   diffusion with  aggregation/repulsion effects and derive the corresponding PDEs. 
With this aim in mind, let us  first  specify the gain and loss terms due to diffusion.  Following the standard approach   that the cells diffuse randomly and, as before, the jump of a cell is of length $h$, the proportion of cells arriving to $I_j$ due to diffusion comes from adjacent intervals, with respective probabilities $P_{I_{j+1}\to I_{j}}=p/2=P_{I_{j-1}\to I_{j}}$ and probability $1-p$ of remaining at $I_j$. Namely,
 $$
 G_{\mbox{diff}}(t,j)=s_{j+1}(t) P_{I_{j+1}\to I_{j}}+s_{j-1} P_{I_{j-1}\to I_{j}}=\frac p2 s_{j+1}(t)+\frac p2 s_{j-1}(t).
 $$
 
 On the other hand, also by diffusion the cells in the interval $I_j$ will leave it to travel to either $I_{j-1}$ or $I_{j+1}$ with respective probabilities $P_{I_j\to I_{j-1}}=p/2=P_{I_j\to I_{j+1}}$. Thus, the proportion of cells leaving $I_j$ is
 $$
 L_{\mbox{diff}}(t,j)= s_j(t)(P_{I_j\to I_{j-1}}+P_{I_j\to I_{j+1}})=p\, s_{j}(t).
 $$
Consequently, 
\begin{align}\label{lapla2}
 G_{\mbox{diff}}(t,j)- L_{\mbox{diff}}(t,j)&=\frac p2 s_{j+1}(t)+\frac p2 s_{j-1}(t)-p\, s_j(t)\nonumber\\&
 = h\frac p2\Big[(u_h(y_{j+1},t)-u_h(y_j,t))+(u_h(y_{j-1},t)-u_h(y_j,t))\Big]\\ \nonumber
 & ={h^3} \frac p2 \Big[\Delta_h u_h(y_j)\Big],
\end{align}
where 
$$
\Delta_h v(x)=\frac{v(x+h)-2v(x)+v(x-h)}{h^2}
$$
is the discrete Laplacian.
Then (\ref{lapla2}) can be rewritten as
\begin{align}\label{lapla}
 G_{\mbox{diff}}(t,j)- L_{\mbox{diff}}(t,j)
 =h^3\mathcal T^h_{diff}(y_j,t),
\end{align}
with
$$
\lim_{h\to 0} \mathcal T^h_{diff}(y_j,t) =\frac p2\Delta u(x,t).
$$

Then, combining \eqref{general} and  \eqref{lapla} with the different interaction terms,  we deduce the general relation
\begin{equation}\label{allefects}
h\left[\frac{u_h(y_j,t+\Delta t)-u_h(y_j,t)}{\Delta t } 
\right]= Q^h_{\mbox{diff}} h^3\mathcal {T}^h_{diff}(y_j,t)
-Q^h_{\mbox{int}} h^2  ( \mathcal T^h(y_j,t) + O(h)),
\end{equation}
where $\mathcal {T}^h$ is defined either as (\ref{nlg}), (\ref{carrillo}) or \ref{carlo}.


In order to preserve all of the effects in the limit PDEs, 
%
the following hypotheses on the frequencies must be imposed:
\begin{equation}\label{cociente-diff-constrain}
\lim_{h\to 0}h^2Q^h_{\mbox{diff}}=\alpha_{\mbox{diff}},
\end{equation}
and
\begin{equation}\label{cociente-int-constrain}
\lim_{h\to 0} hQ^h_{\mbox{int}}=\alpha_{\mbox{int}}.
\end{equation}
At this stage, it arises the importance of the correct choice of the frequencies $Q^h_{\mbox{diff}}$ and $Q^h_{\mbox{int}}$, with respect to $h$.  See \cite{AS}, where $h$ dependent coefficients are taken.

After appropriate time scaling in (\ref{allefects}), passing formally to the limit as $\Delta t,h\to 0$,  we arrive to the general PDE given  in (\ref{eq:APS}), for each $x\in(-L,L)$ and $t>0$. Namely,
\begin{equation}\label{PDE:APS}
u_t(x,t)= u_{xx}(x,t)
-\gamma \, \partial_x \bigg( u(x,t)K(u)(x,t)\bigg),
\end{equation}
 where $\gamma=2\alpha_{\mbox{int}}/(p\alpha_{\mbox{diff}})$ is a positive parameter related to frequencies of interactions and diffusion. The nonlocal drift term $K(u)$ is given by
\begin{itemize}
\item[3.5.1] Super aggregation APS   model \cite{APS}: $$K(u)=\nabla_{NL}u.$$
\item[3.5.2] Local saturation model \cite{CMSTT}:  $$K(u)=(1-u)\nabla_{NL}u.$$
\item[3.5.3] Nonlocal saturation model: $$K(u)(x,t)=\left(1-\frac1K\int_{x-R}^{x+R}u(y,t) \mathbbm{1}_{\Omega}(y) dy\right)\nabla_{NL}u(x,t),$$
where $\mathbbm{1}_{\Omega}$ the characteristic function in $\Omega$.
\end{itemize} 
Including the effect of the bounded domain $\Omega$, the nonlocal term $\nabla_{NL}u$ takes the form
\begin{equation}\label{NLGGeneral}
(\nabla_{NL}u)(x,t)=
 \int_{x}^{x+R}u(y,t)w(y-x)\mathbbm{1}_{\Omega}(y)dy 
-\int_{x-R}^{x}u(y,t)w(x-y)\mathbbm{1}_{\Omega}(y)dy.
\end{equation}
Moreover, for all of these models, the no-flux boundary conditions read as
 \begin{equation}\label{boundary}
-u_x(x,t)+\gamma \, u(x,t) K(u)(x,t) =0, \quad x=-L, L.
\end{equation}

 \section{Lagrangian perspective}\label{lagra}
 
 While in Section 3 we modelled cell densities distributed among intervals, here we choose a Lagrangian perspective, following the motion of individual cells. This optics is particularly useful when the number of cells in the population is not very large. 
 In this way we expect  
 to extract 
 additional features from cell interactions. Thus each perspective reveals complementary information with respect to the other while, on the other hand, they not only give rise identical PDEs, but also they are closely related at the discrete level, as subsequent simulations will reveal.

Now, agent based models are derived where the unknown variable is the position of the cells at certain time $t$. As in Section 3 we start by considering a single population of $N$ cells randomly distributed in  $\Omega= (-L,L)$, at initial positions $x_i(0)$, $i=1,\cdots,N$.
The following system of stochastic differential equations describes how the position $x_i(t)$ of the $i$-th cell evolves over  time: 
\begin{equation}\label{estocastica}
dx_i(t)=\frac 1 N \sum_{j\neq i}\mathcal{F}(x_i(t), x_j(t), S_a(x_i))dt+\sqrt{2}\, \ve \, B_t^i,\quad i=1,\cdots,N,\; t>0.
\end{equation}

 Here,  $\mathcal{F}(x_i,x_j, S_a(x_i))$ represents the action of  particle $j$  on  particle $i$, that is, the model evolves with pairwise interactions. In case of the APS aggregative interactions, described in $\S$\ref{Sec:APS}, cells experiment  attractive pairwise forces towards the cells within their sensing area. Indeed,  $\mathcal{F}(x_i,x_j, S_a(x_i))=F(|x_i-x_j|)$, with $|\cdot|$ denoting the Euclidean distance.   However, the attraction that cell  $j$  exerts on  cell  $i$ could be weakened  if the neighbouring area to  position $x_i$ is too crowded, or it could become even repulsive. This fact is taken into account by the  term  $S_a(x_i)$, which for the interactions treated in $\S$\ref{Sec:3.3} and $\S$\ref{nuestro}  will be an appropriate nonlocal saturation or repulsive coefficient of the individual forces acting on particle $i$. Namely,  these pairwise forces are of the form 
 \begin{equation}\label{fuerzas}
 \mathcal{F}(x_i,x_j, S_a(x_i))=S_a(x_i)F(|x_i-x_j|).
  \end{equation}
  The terms  $(B_t^i)$, $i=1,\cdots,N$, $t\geq 0$ denote independent standard Brownian motions on $\RR$. To ensure then that the particles remain inside the domain, the trajectories given by (\ref{estocastica}) must be truncated such that $x_i(t)\in [-L,L]$, for all $t\geq 0$.

Systems as (\ref{estocastica}) are classic in many contexts in physics, biology, economics, sociology, etc.,  where it is necessary to analyze the motion of a large number of particles described their positions, combining deterministic and stochastic processes. The existing literature in this field is vast and  increasing over time due to its great applicability. For a general overview on this topic and specific applications to biological processes, such as population dynamics including logistic growth and repulsion effects, tumoral cell growth, applications to finance and market analysis,  we refer to the  book \cite{Libro_capasso} and references therein.  For further references in stochastic equations with applications in swarming we refer to the Chapter \cite{CFTV}, and for those related to fluid mechanics, granular flows, traffic models, or large communications networks \cite{bellomo_book}, and references included in these compilations.

In the  so-called hydrodynamic limit (as $N\to\infty$) the independent Brownian motions produce diffusion, while the average sums  of (\ref{estocastica})  lead to the different  drift terms corresponding to each of the models, depending on the choice of the saturation coefficient $S_a$. The continuous density of cells in the limit $N\to\infty$, denoted by $u$, will be a weak solution to  mean field deterministic PDEs of the form
$$
u_t(x,t)=\ve^2 u_{xx}-\partial_x\Big(u(x,t)K(u)(x,t)\Big),
$$
where $K(u)(x,t)$ takes into account the balance of the forces exerted on the particles placed on $x$. After rescaling in time, we can rewrite it as equation (\ref{eq:APS}).

Our aim now is to recover the PDE (\ref{PDE:APS}) for the different choices of $K(u)$. According to the nature of the interacting forces, the drift term $K(u)$ can yield always aggregation, or  it could consider saturation or even repulsion effects. This fact will be determined in the SDEs system by the saturation coefficient, $S_a$.

The rigourous derivation of the limit  PDE from a discrete system of $N$ stochastic differential equations as $N\to\infty$ is beyond the scope of this paper.
Nevertheless we include  a short presentation following the course notes collected in  \cite{B}. The empirical measure of the system of cells is given by
 $$
 \mu_t^N=\frac1N\sum_{i=1}^N\delta_{x_i(t)},
 $$
where  $\delta_x$ is the Dirac mass on $x\in\RR$. 

Let $\varphi:[-L,L]\to\RR$ be an observable function and  $x_t=x(t)\in [-L,L]$ a trajectory of a cell. How does $\varphi(x_t)$ evolve in time according  to $\mu^N_t$ ? Recalling \^{I}to's formula  
and system \eqref{estocastica}, one has
\begin{align}\label{xi}
\langle \varphi,\delta_{x_i(t)}\rangle-\langle \varphi,\delta_{x_i(0)}\rangle&=\varphi(x_i(t))-\varphi(x_i(0))=\sqrt{2}\ve\int_0^t  \frac{d \varphi}{d x}(x_i(s))dB^i_s\nonumber\\&+\int_0^t \left[S_a(x_i(s))\frac 1{N}\sum_{j\neq i}F(x_i(s),x_j(s))\frac{d\varphi}{d x}(x_i(s))+\ve^2
\frac{d^2\varphi}{d x^2}(x_i(s))\right] ds.
\end{align}
The nonlocal saturation/repulsion coefficient for particle $i$ at time $t$ has the form
\begin{equation}\label{satu-cof}
S_a(x_i(t))=1-\frac \lambda{2a K}\langle \widehat w(|z-x_i|)\mathbbm{1}_{|z- x_i(t)| < a},\mu_t^N\rangle, \; 0<a\leq L,\; 0<K\leq 1, \mbox{ and }\lambda\in\{0,1\}. 
\end{equation}
If $\lambda=0$ the forces take into account only  aggregation effects, while saturation occurs when $\lambda=1$. Saturation strength will vary according to the value  $0<K\leq 1$, representing the crowding capacity. The parameter $0<a\leq L$  delimits the region around the $i-$th cell where such effect is taken into account. Finally, the weight $\widehat w=\widehat w(r)$ moderates the influence of the surrounding cells according to the distance from the position $x_i(t)$. 
At this stage, averaging equation \eqref{xi} in $(x_i)_{i=1} ^N$  yields 
\begin{align}\label{superequation}
\langle \varphi,\mu_t^N\rangle-\langle \varphi,\mu_0^N\rangle&=\int_0^t\frac 1N\sum_{i=1}^N \left[\left(1-\frac \lambda{2a K}\langle \widehat w(|z-x_i|) \mathbbm{1}_{|z- x_i(t)| < a},\mu_t^N\rangle\right)\frac 1{N}\sum_{j\neq i}F(x_i(s),x_j(s))\frac{d\varphi}{d x}(x_i(s))\right]\,ds\nonumber\\&+\int_0^t\frac 1N\sum_{i=1}^N \ve^2
\frac{d^2\varphi}{d x^2}(x_i(s)\,ds+\sqrt{2}\ve\int_0^t \frac 1N\sum_{i=1}^N \frac{d \varphi}{d x}(x_i(s))dB^i_s\nonumber\\
&=\int_0^t\int_{-L}^L\left(1-\frac\lambda{2aK}\int_{-L}^{L} \widehat w(|z-x|)\mathbbm{1}_{| z-x| < a}\:\mu_s^N(z)\,dz\right)\left(\int_{-L}^L F(x,y) d\mu_s^N(y)dy\right)\frac{d\varphi}{d x}(x)\,dx\,ds\\&+\int_0^t\int_{-L}^L\ve^2
\frac{d^2\varphi}{d x^2}(x)d\mu_s^N(x)\,ds+\sqrt{2}\ve\int_0^t \frac 1N\sum_{i=1}^N \frac{d \varphi}{d x}(x_i(s))dB^i_s.\nonumber
\end{align}
It is well known that by the law of large numbers it holds that
$$
\sqrt{2}\ve\int_0^t \frac 1N\sum_{i=1}^N \frac{d \varphi}{d x_i}(x_i(s))dB^i_s\to 0,\quad \mbox{ as }\,N\to\infty.
$$
Appropriate hypotheses  ensure the existence of a measure $\mu_t$ absolutely continuous with respect to the spatial variable, such that $\mu_t^N\to \mu_t$ as $N\to\infty$. We identify formally the limiting PDEs for each of the cases.

\subsection{Strong aggregation (APS) model \cite{APS}.}\label{Model:APS}

\

Let us begin our analysis with the Lagrangian formulation of the APS model,  which exhibits strong aggregation,
by taking  in the coefficient  \eqref{satu-cof} of the pairwise forces \eqref{fuerzas} the parameter $\lambda=0$ and defining
\begin{equation}\label{fuerzasAPS}
F(x_i,x_j) = 
\gamma\,
w (|x_i-x_j|)\, {\rm sign} (x_j-x_i)
\end{equation}
where $w:[0,+\infty)\to \mathbb{R}_+$ is a nonnegative weight function with $w(r)>0$ for $r\in [0,R)$ and $w(r)=0$ for $r\ge R$. For instance $w(r)$ could be  linearly decreasing with respect to the distance $r=|x_i-x_j|$. Then, in the limit as $N\to\infty$
\begin{align*}\label{forma-debil}
\int_{-L}^L \varphi d\mu_t-\int_{-L}^L\varphi d\mu_0&=\int_0^t\left[\int_{-L}^L \int_{-L}^L  \gamma w (|x-y|)\, {\rm sign} (y-x)d\mu_s(y)\nabla_x\,\varphi d\mu_s(x)+\ve^2
\Delta \varphi(x)d\mu_s(x)\right]ds,
\end{align*}
which is the weak form of equation \eqref{PDE:APS} with drift term according to model \S{3.5.1} assuming that $\mu_t(A)=\int_A u(x,t)dx$.


Starting from random distributed positions,  at  first stages, when the cells are still  all spread, it is improbable that attraction plays a role far from the boundary.  However, when the moving cell has a position $x_i(t)$   close to the boundary, it experiments an attraction  towards the interior where the rest of cells are. 
 
On the other hand, as the density of cells is no longer equally distributed, if these attraction forces are dominant with respect to diffusion, then it might yield to very contractive and highly oscillate dynamics. 
In the spirit of avoiding it, we will take the value of $\lambda=1$ in the saturation coefficient \eqref{satu-cof} and deduce aggregation/saturation models.

\subsection{Nonlocal saturation coefficient}\label{Model:nonlocal}

 We propose to consider a nonlocal saturation coefficient \eqref{satu-cof} with $\lambda=1$ and choosing $a=R$  the sensing radius and the forces \eqref{fuerzasAPS}. Namely,
 $$
 \mathcal{F}(x_i,x_j, S_R(x_i))=\left(1-\frac 1{2R K}\langle \widehat w(|z-x_i|)\mathbbm{1}_{|z- x_i(t)| < R},\mu_t^N\rangle\right)F(x_i,x_j),
 $$
 for some value of the crowding capacity $0<K\leq1$ and  some weight $\widehat w=\widehat w(r)$,   linear decreasing with the distance or a positive constant. As numerical simulations will show, choosing $\widehat w$ linearly decreasing is related to take higher crowding capacity $K$, permitting more aggregation, as it can be appreciated in  Figures \ref{fig:NL_K06pesolineal_histogramas} and \ref{fig:NL_K06pesolineal_driftsatu}.  Hence, in the sequel $\widehat w=1$. Taking $N\to\infty$ in  \eqref{superequation} yields the PDE \eqref{PDE:APS} with drift term corresponding to \S{3.5.3}.

  \subsection{Local saturation coefficient \cite{CMSTT}.}\label{Model:Carrillo}
  
  \
  
 Recall that trajectories follow \eqref{superequation} according to the nonlocal saturation coefficient \eqref{satu-cof} with $\lambda=1$ and pairwise interactions  $F(x_i,x_j)$ specified  by \eqref{fuerzasAPS}. In the limit $N\to\infty$ we wish to obtain  the PDE given in \eqref{PDE:APS}, corresponding to model \S{3.5.2}, thus the nonlocal coefficient \eqref{satu-cof} must  converge to the local term $(1-u(x,t))$. For this purpose, we  take $K=1$ and  $a=a(N)$ such that $a(N)\to 0$ as $N\to\infty$. With this choice
 $$
 S_a(x_i(t))=1-\frac 1{2a }\langle\mathbbm{1}_{|z- x_i(t)| < a},\mu_t^N\rangle=1-\frac 1{2a }\int_{x_i(t)-a}^{x_i(t)+a}d\mu_t^N(z),
 $$
hence the coefficient $S_a(x_i(t))$ approaches to the local PDE coefficient $ 1-u$ as $N\to\infty$.

\section{Numerical Simulations}\label{section:nume}

 We conclude this work by performing some numerical simulations
 to illustrate how the interaction effects of the different models under consideration   act on the diffusive random movements of the cells. We compare  the corresponding asymptotic  behaviour  obtained for each case, when varying the parameters involved. These numerical experiments will also reveal the coherence between both perspectives, Eulerian and Lagragian, when approaching the corresponding mean field PDEs of each model.
 
\subsection{Numerical methods}
 Notice that the general equation (\ref{allefects}) could be thought as an  explicit Euler time scheme, and centred spatial finite difference  approximation  of the diffusion and interaction terms, of the resulting PDE. 
  However, such a scheme does not generally give rise to an appropriate numerical method. Indeed, in the case of  interaction dominance, strong aggregation and repulsion could eventually destabilize the numerical results. 
 
 Willing to avoid this issue, we propose  a  semi-implicit linear scheme in time combined with an upwind finite volume scheme (of Godunov type) in space, a well known technique widely used when simulating propagating phenomena, see for instance \cite{saito}. Roughly speaking, upwind schemes attempt to discretize the  equation by using finite differences  biased in the direction determined by the sign of the characteristic speed. This technique results a very robust and efficient numerical approach. The interaction term, $\nabla\cdot (uK(u))$ for different choices of $K(u)$ in $\S 3.5.1$, $\S 3.5.2$, $\S 3.5.3$,  is linearized by approximating the nonlocal term $K(u)$ explicitly. This way it is possible to keep the sparse structure of the approximating matrices.


Regarding the Lagrangian perspective, recall that the motion of any cell $i$ adopts the following structure

$$d x_i(t) = \mathcal{F}(\vec{x}(t), S_a(x_i(t))) d t + \varepsilon^2 d B_i(t).$$
The independent Brownian motions are driven by a normal distribution for all $i=1,\cdots,N$, precisely $B_i(t)= \mathcal{N}(0,0.1)$,  and the diffusion coefficient is considered to be constant  $\varepsilon^2 =0.4$. 

We approximate the deterministic part of the trajectories  using Euler explicit scheme  and sample the normal distribution to reproduce the Brownian motion. Namely,

\begin{equation}\label{finite-disc}x_i ( t + \Delta t) = x_i (t) + \Delta t [ \mathcal{F}(\vec{x}(t),S_a(x_i(t))+  \varepsilon^2 \mathcal{N}(0,0.1)]. 
\end{equation}
The choice of saturation coefficient $S_a$, specified in \eqref{satu-cof}, will vary in order to depict the three different models studied in the previous section.

In all of the simulations we can appreciate the consistency between both perspectives, Eulerian and Lagrangian. From the set of trajectories of the SDEs system, we compute the drift terms of the different models, $\S 3.5.1$, $\S 3.5.2$, $\S 3.5.3$, and their corresponding   saturation coefficients in the cases $\S 3.5.2$, $\S 3.5.3$, according to a fix mesh in the spatial variable. When obtaining the drifts $\S 3.5.1$ and $\S 3.5.3$ from the SDEs system, we take a mesh of size $\hat h>0$ independent of the number of cells. However, to approximate the local coefficient in $\S 3.5.2$ by the nonlocal coefficient \eqref{satu-cof}, is more delicate than in the other two cases. In this particular approach we need $\hat h$ to depend on the number of cells $N$. Precisely, we take $\hat h=2/\sqrt N$.

{\bf Common parameters:} 
 All of the subsequent simulations are performed in the domain $\overline\Omega=[-1,1]$, choosing the sensing radius to be $R=0.5$. The interactions among cells  are weighted by (\ref{peso}), which in most of the cases will be taken inversely proportional to their distance as $w(r)=(R-r)/R$ and just in some cases being constant $w\equiv 1$, precisely in Figures \ref{fig:NL_K1_cw_histogramas}, \ref{fig:NL_K1_driftsatu}, \ref{fig:NL_K04_cw_histogramas} and \ref{fig:NL_K04_cw_driftsatu}. On the other hand, in  Figures \ref{fig:NL_K06pesolineal_histogramas} and \ref{fig:NL_K06pesolineal_driftsatu}  we investigate the impact of considering a weight $\hat w(r)=(R-r)/R$ in the nonlocal saturation coefficient in $\S 3.5.3$ and  in the SDEs system~ $\S$\ref{Model:nonlocal} in contrast to $\hat w\equiv 1$, as in the rest of the examples related to this model.
 
 \
  
  In case of Eulerian simulations we use a uniform spatial mesh of size $h=2\cdot10^{-3}$ and  time step $\Delta t=10^{-4}$. The number of time iterations will depend on the example under consideration, needing to be higher as the saturation is stronger. It is also essential to calibrate the frequency of occurrence of both mechanisms, diffusion and interactions. This is recorded in the parameter $\gamma$ accounting for the relative frequency between them. Indeed, recall that $\gamma=2\alpha_{\mbox{int}}/(p\alpha_{\mbox{diff}})$, with $\alpha_{\mbox{int}}$ and $\alpha_{\mbox{diff}}$ specified in \eqref{cociente-diff-constrain} and \eqref{cociente-int-constrain}, respectively and $p\in (0,1)$ the probability of moving randomly at one time step $\Delta t$. We take the value $\gamma=10^3$.  
 
 \
 
 In the Lagrangian perspective, for cases \ref{Model:APS} and \ref{Model:nonlocal} we start with $N\sim 300$, use a time step $\Delta t=0.01$ and  a mesh $\hat h=0.01$ to depict the drift and saturation coefficients. As we mentioned, model~\ref{Model:Carrillo} is more delicate. We start with $N\sim 3000$ and reduce the time step $\Delta t=0.001$, which increases the number of iterations. In addition, a mesh of size $\hat h=2/\sqrt{N}$ is needed to obtain a representation of the drift terms  and saturation coefficients  as similar as possible to $\S 3.5.2$  and histograms resembling the solutions of the PDEs.
 
 \ 

In the sequel we include several examples illustrating all of the models for different values of the parameter $K\in(0,1]$  regulating the crowding capacity, combined with diverse choices of weights $w$ and $\hat w$. We  compare the pictures obtained for Eulerian and for Lagrangian perspectives. In all of these examples, the initial configuration of the cells is uniformly distributed, either along the whole domain $[-L,L]$ or, as in Figures \ref{fig:NL_K06concentrado_histogramas}, \ref{fig:NL_K06concentrado_driftsatu}, \ref{fig:NL_K02concentrado_histogramas}, \ref{fig:NL_K02concentrado_driftsatu} concentrated in some $[-\delta,\delta]$. This option corresponds to take as initial datum a perturbation of an appropriated constant, either $1/(2L)$ or $1/(2\delta)$ in the Eulerian perspective.

\

\

{\bf Numerical simulations description:} We include several examples for the three  models treated in this work considering different intensities for repulsion. Precisely, absence of repulsion ($K=1$), weak repulsion ($K=0.6$), intermediate repulsion ($K=0.4$) and strong repulsion ($K=0.2$).

\

{\it Absence of repulsion:} Figures  \ref{fig:APS_Histogramas} and \ref{fig:APS_Drift} are devoted to illustrate the APS model, for which high accumulations peaks are formed. In contrast, even for high values of $K$, notice that the  local saturation model prevents the formation of these strong aggregations. Figures \ref{fig:Carrillo_K1_histogramas} and \ref{fig:Carrillo_K1_satu} correspond to $K=1$. As we mentioned, the approach of the local saturation coefficient $1-u$ by the discrete saturation coefficients $S_a$ given in Section $\S 4.3$ needs high values of $N$. We take $N=3.500$ and $a(N)=1/\sqrt{N}$, (hence $a(N)\to 0$, as $N\to\infty$). The  time step is reduced to $\Delta t=0.001$. However, we cannot avoid that the saturation coefficient $S_a$ of model $\S 4.3$ eventually becomes negative, as it is depicted in Figure~\ref{fig:Carrillo_K1_satu}. Nevertheless, the correspondence of both dynamics, Eulerian and Lagrangian is quite coherent, as it can be appreciated in Figures  \ref{fig:Carrillo_K1_histogramas} and \ref{fig:Carrillo_K1_drift}. In the absence of repulsion, the nonlocal saturation model avoids as well the formation of accumulations peaks, but only considering a constant weight $w\equiv1$ in the nonlocal term, see Figures \ref{fig:NL_K1_cw_histogramas} and \ref{fig:NL_K1_driftsatu}. For a linear decreasing weight $w$ and $K=1$ Figures \ref{fig:NL_K1_histogramas} and \ref{fig:NL_K1_driftsatu}  show strong aggregations.

\

{\it Weak repulsion:} As repulsion slightly   grows, $K=0.6$, particles tend to spread more distributed along the domain. This dissemination is more evident  and uniform in case of the local saturation coefficient model,  see Figures \ref{fig:Carrillo__K07_histogramas}  and \ref{fig:Carrillo__K07_Drift}. In contrast, considering the nonlocal saturation model in \S 4.2, cells move more aggregated and towards the interior of the domain, see Figures \ref{fig:NL_K06_histogramas} and \ref{fig:NL_K06_driftsatu}, and accumulate at the boundary of the ``aggregation tubes". The non-locality of the saturation coefficient allows higher accumulations on such boundaries.  In this range of weak repulsion we have also depicted the case of a constant weight $\widehat w\equiv 1$. This choice of weight in the coefficient reduces the saturation, hence accused peaks are formed, as Figures \ref{fig:NL_K06pesolineal_histogramas}  and \ref{fig:NL_K06pesolineal_driftsatu} reveal. We conclude this subsection by investigating the evolution of a uniformly concentrated initial data in the interval $[-0.3,0.3]$ when $K=0.6$. The dynamics reduces approximately to the interval $[-0.2,0.2]$, with accumulations at the boundaries. It can be observed in Figures \ref{fig:NL_K06concentrado_histogramas} and \ref{fig:NL_K06concentrado_driftsatu}.

 \

{\it Intermediate repulsion:} Increasing repulsion, ($K=0.4$),  the local saturation model $\S 4.3$ spreads the cells along the whole domain. In contrast to the choice $K=0.6$, if $K=0.4$ cells  accumulate  at the boundaries. This fact can be observed in  Figures \ref{fig:Carrillo_K04_histogramas} and \ref{fig:Carrillo_K04_Drift}. Regarding nonlocal saturation model $\S 4.2$, the saturation effect acts in a larger set  than in the case $K=0.6$. This fact disseminates the cells wider than with less repulsion and, as before, accumulated at the boundaries of the aggregation tubes,  see  Figures \ref{fig:NL_K04_histogramas} and \ref{fig:NL_K04_driftsatu}.  Choosing $\widehat w\equiv 1$ in the nonlocal coefficient introduces  repulsion, given that all of the existing cells inside the sensing radius have the same importance. Consequently, cells are pushed towards the boundaries, though there exit some aggregations in the interior, see Figures \ref{fig:NL_K04_cw_histogramas} and \ref{fig:NL_K04_cw_driftsatu}.

\

{\it High repulsion:} In the local saturation model $\S 4.3$, reducing the crowding capacity up to $K=0.2$ produces no significant change to the case $K=0.4$, hence we do not include the simulations. With respect to nonlocal saturation model $\S 4.2$, the cells spread now along the whole domain, Figures  \ref{fig:NL_K02_histogramas} and \ref{fig:NL_K02_driftsatu}, even when starting with concentrated data, Figures \ref{fig:NL_K02concentrado_histogramas} and \ref{fig:NL_K02concentrado_driftsatu}. 

\newpage

\subsection{Absence of repulsion $K=1$}\label{k=1}

\

\begin{minipage}{\textwidth}
\centering
\includegraphics[scale=0.22]{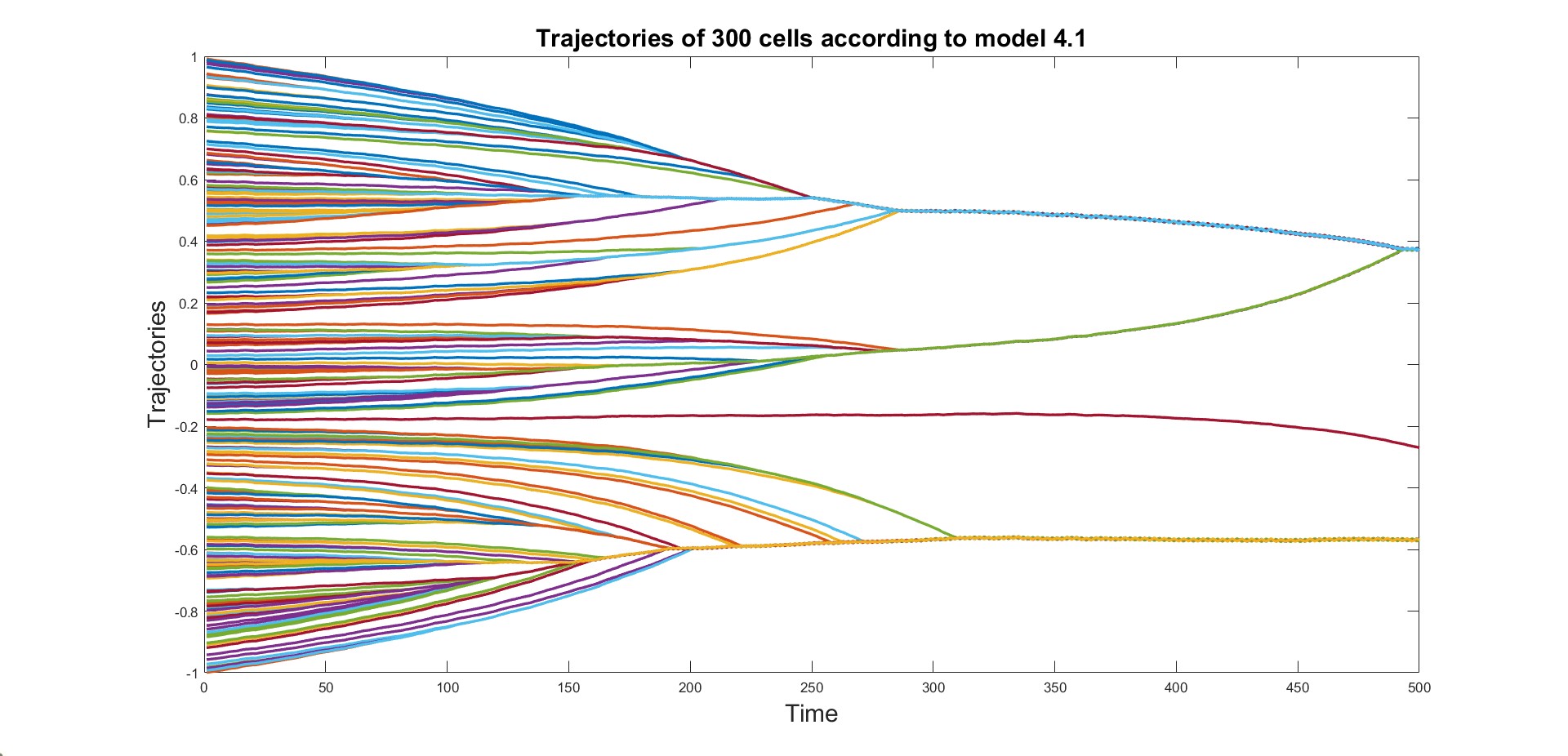}\\

\includegraphics[scale=0.28]{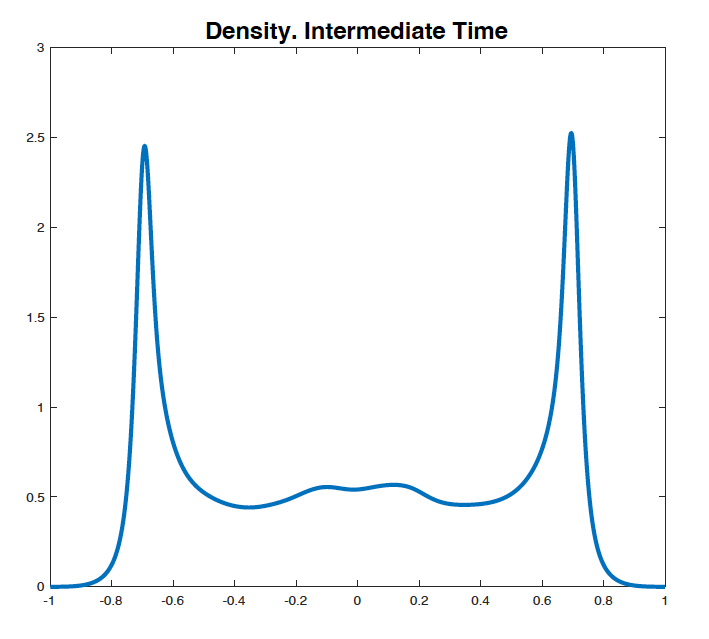}\;\includegraphics[scale=0.28]{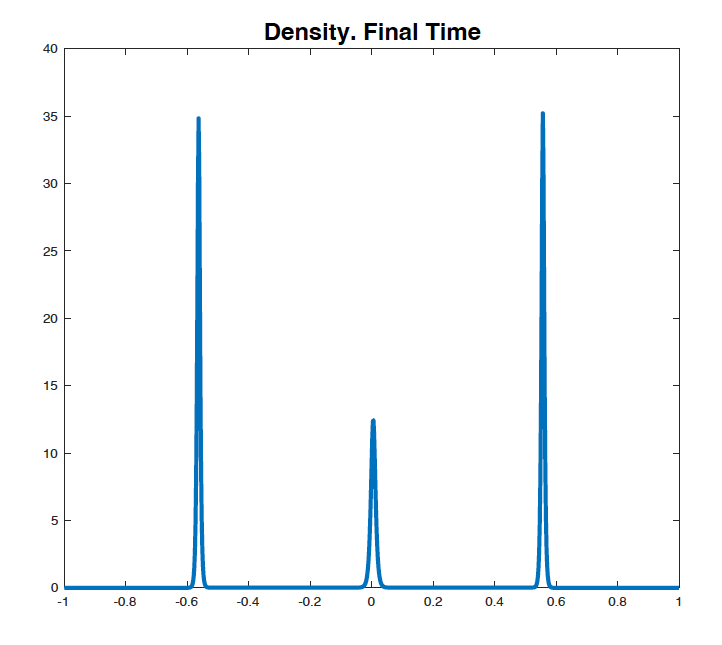}\\
\includegraphics[scale=0.28]{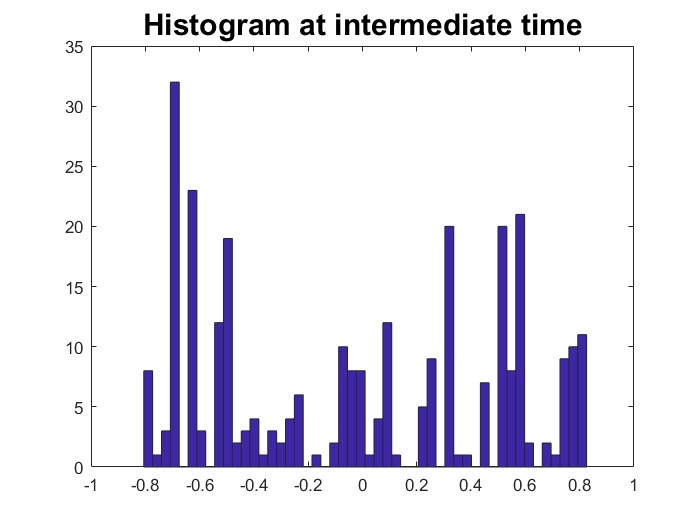}\!\includegraphics[scale=0.28]{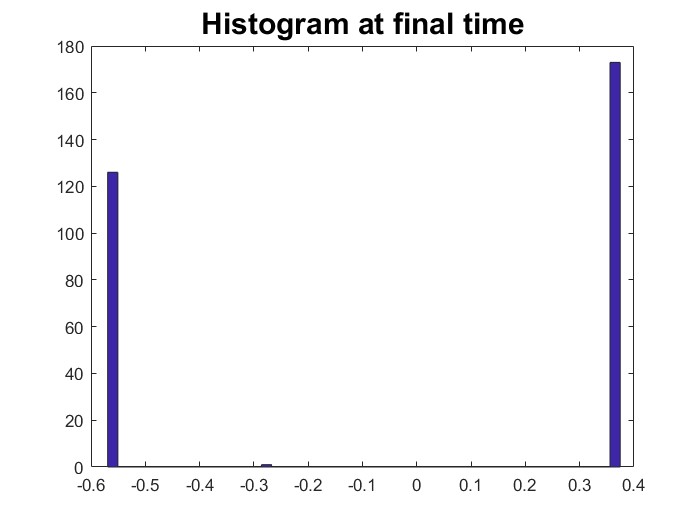}
\captionof{figure}{APS model. First row trajectories with $N=300$ cells, second row densities of the PDE model (\S 3.5.1) and third one histograms of SDEs (\S 4.1). Both at intermediate time on the left and at final time on the right.} 
\label{fig:APS_Histogramas}
\end{minipage}

 \begin{minipage}{\textwidth}
\centering
\includegraphics[scale=0.3]{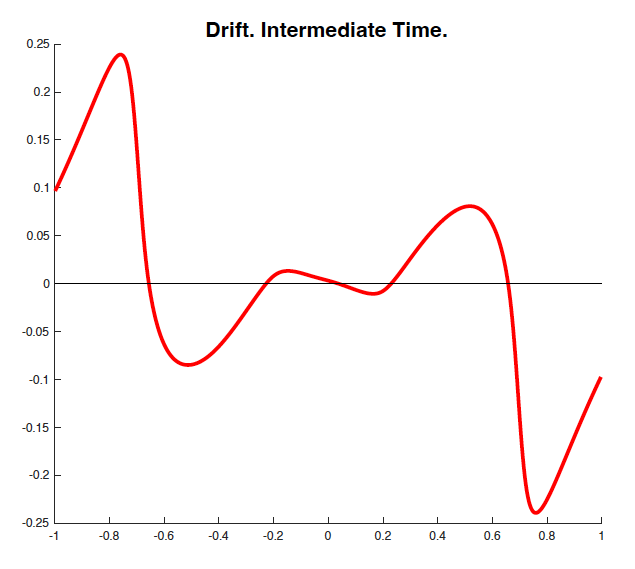}\,\,\,\;\;\includegraphics[scale=0.3]{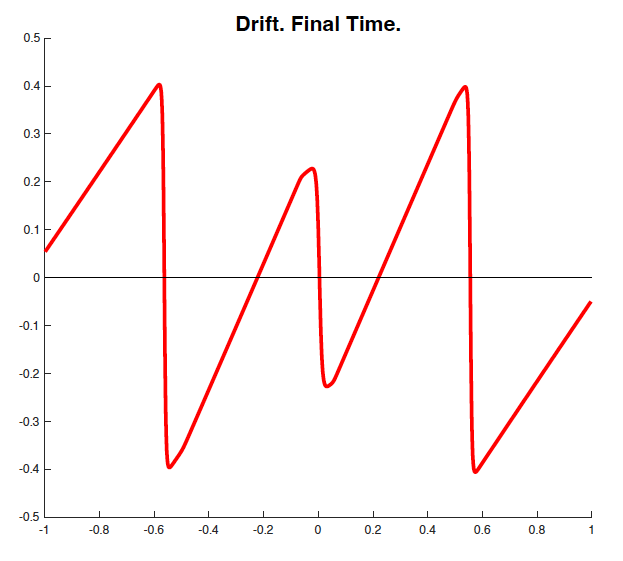}\\

\includegraphics[scale=0.3]{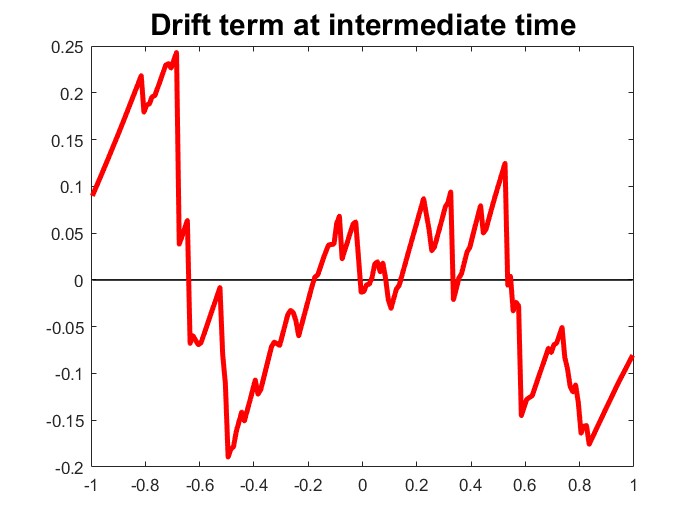}\!\!\includegraphics[scale=0.3]{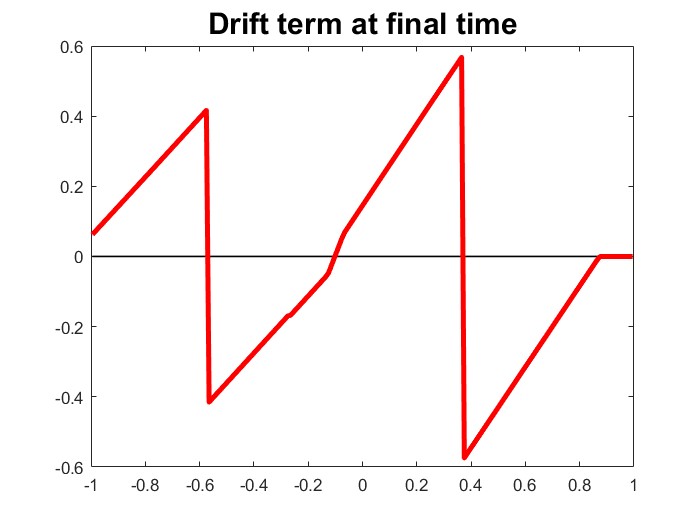}
\captionof{figure}{Drifts for APS model. First row PDE model (\S 3.5.1) and second row SDEs (\S 4.1). Both at intermediate time on the left and at final time on the right.} 
\label{fig:APS_Drift}
\end{minipage}

\

 \begin{minipage}{\textwidth}
\centering
\includegraphics[scale=0.24]{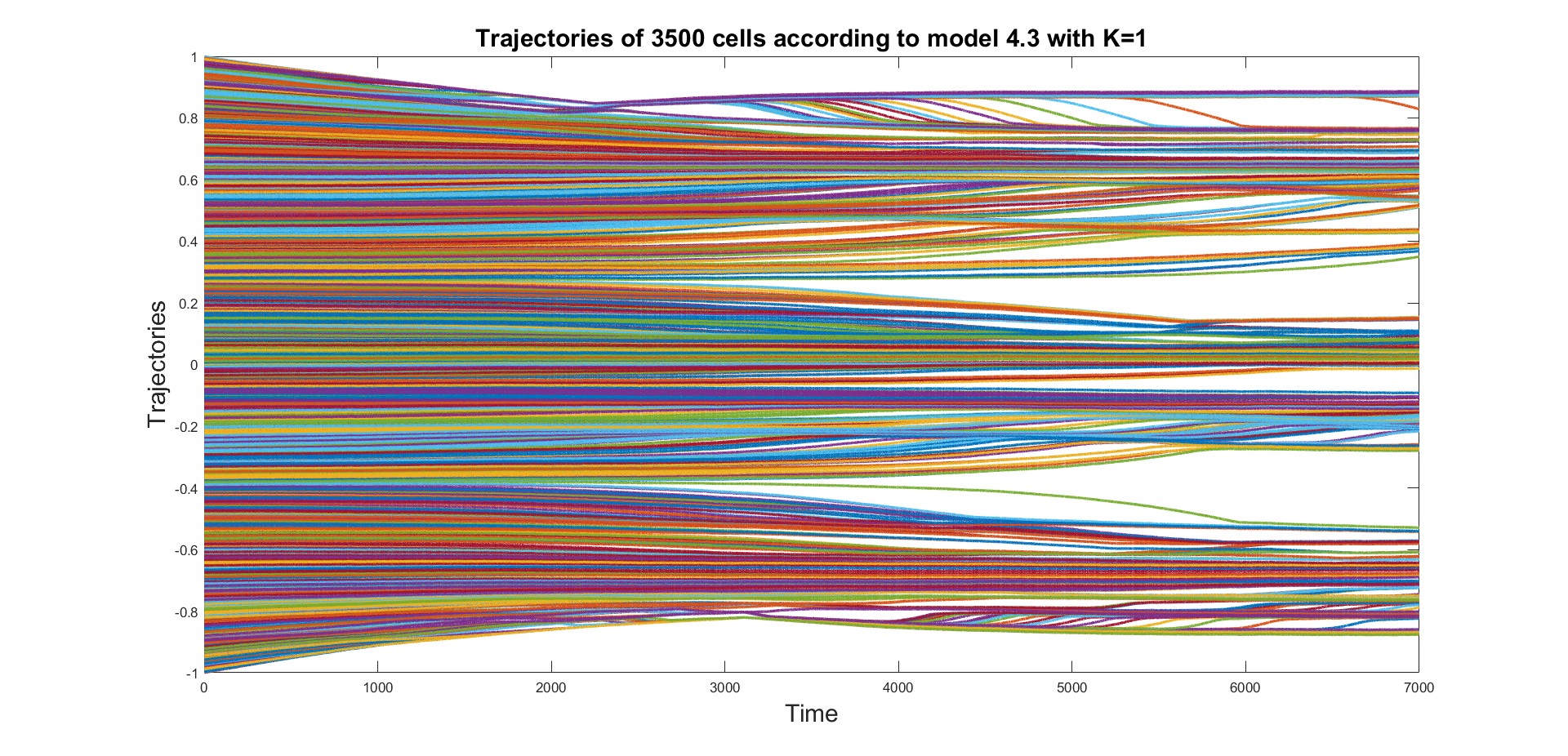}\\

\includegraphics[scale=0.19]{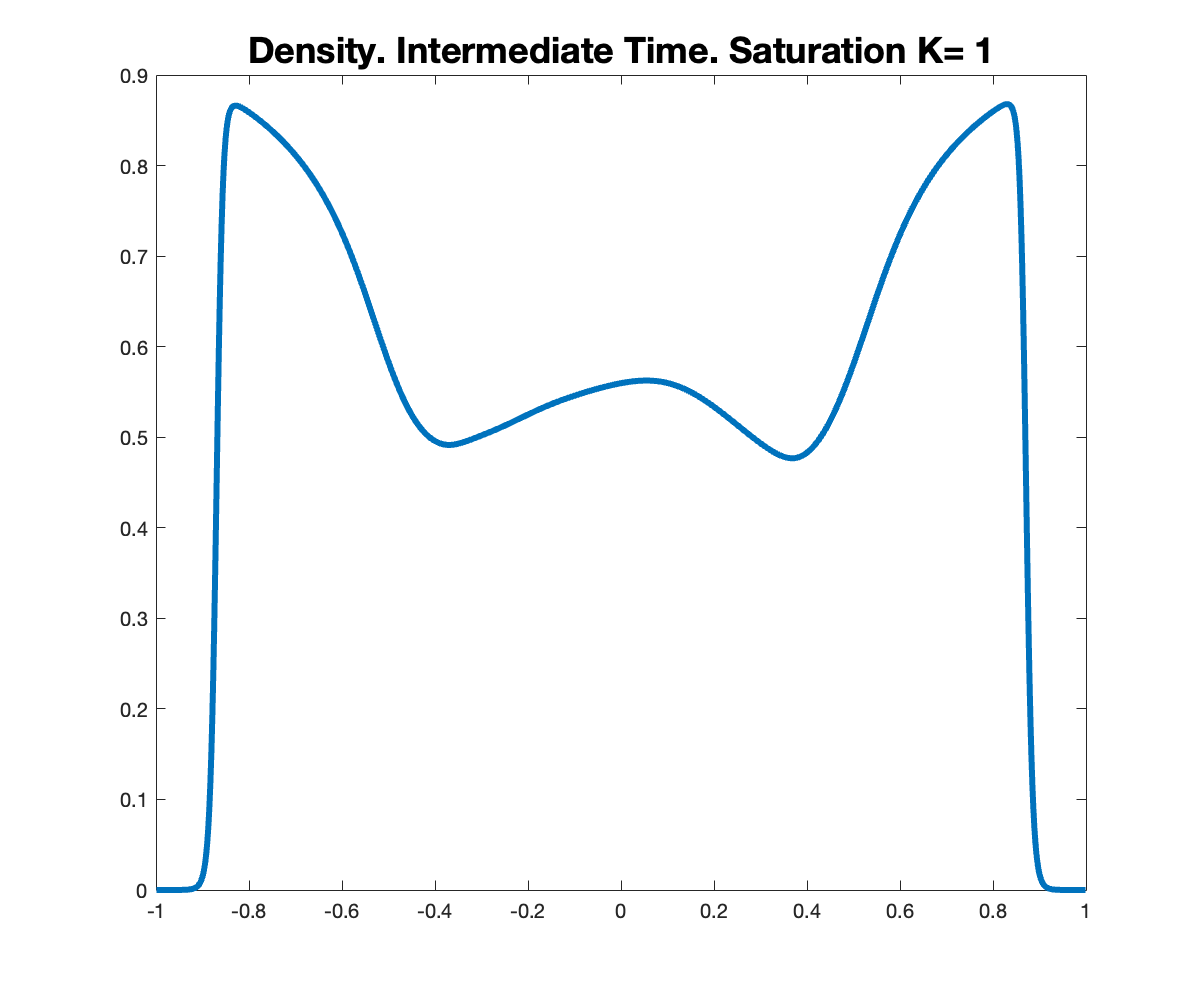}\!\!\includegraphics[scale=0.19]{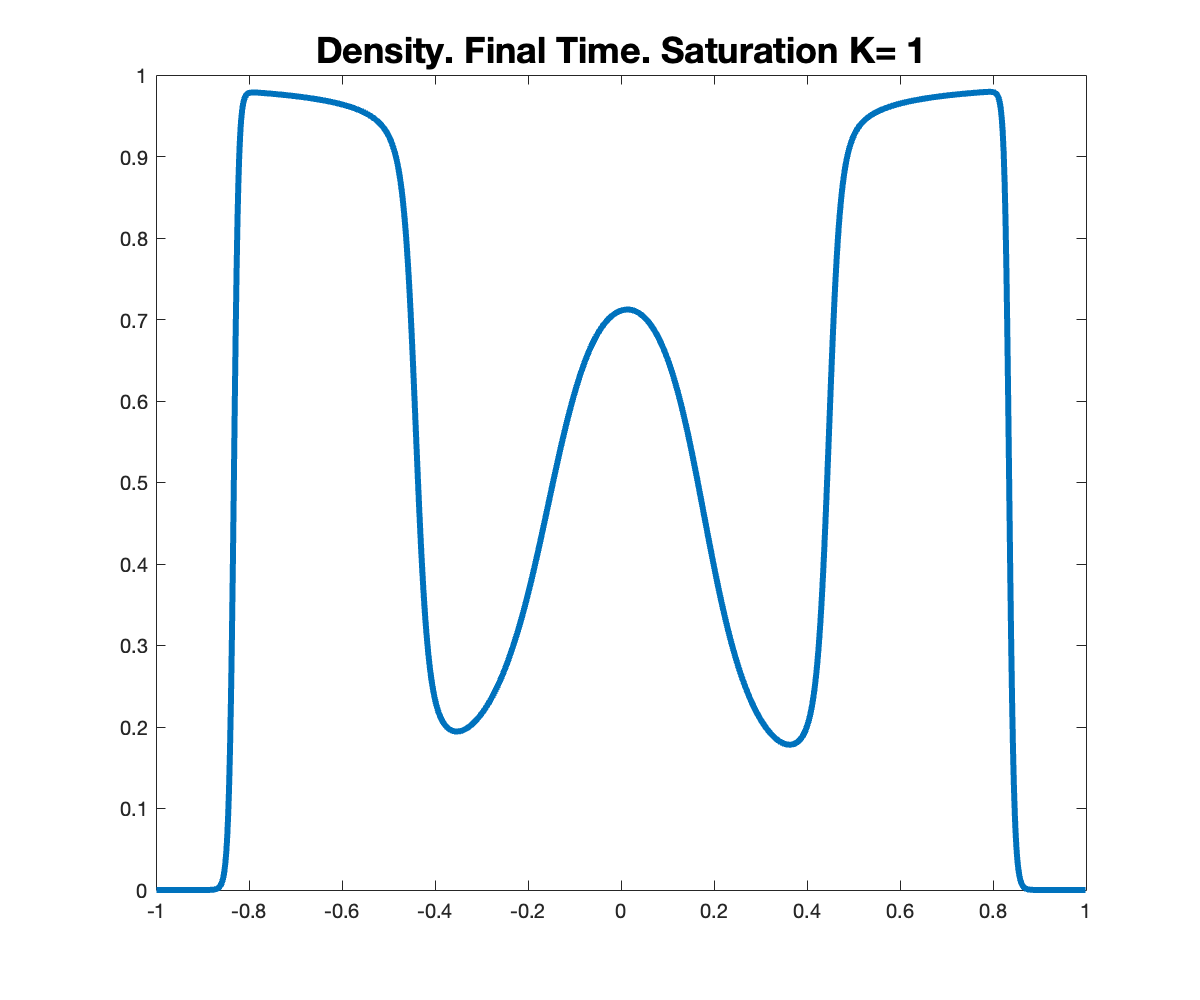}\\
\includegraphics[scale=0.32]{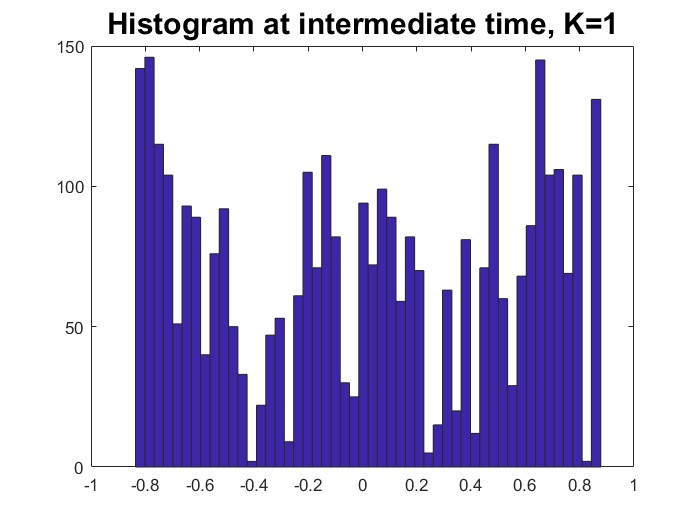}\!\!\includegraphics[scale=0.32]{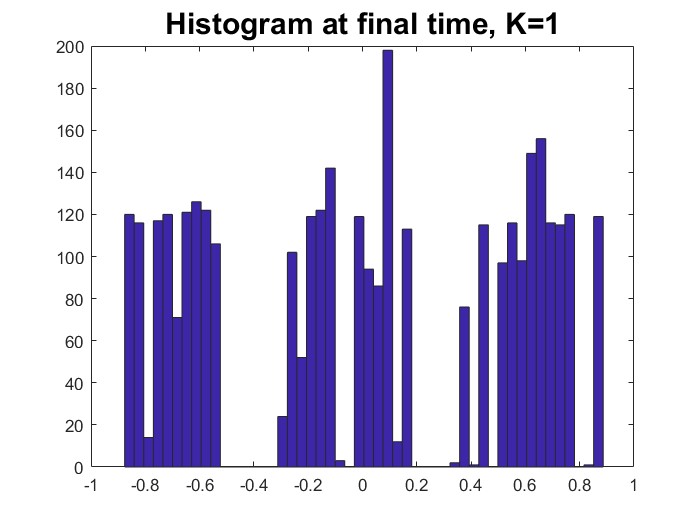}
\captionof{figure}{Local saturation model, $K=1$. First row trajectories with $N=3500$ cells, second row densities of the PDE model (\S 3.5.2) and third one histograms of SDEs (\S 4.3). Both at intermediate time on the left and at final time on the right.
} 
\label{fig:Carrillo_K1_histogramas}
\end{minipage}

\

 \begin{minipage}{\textwidth}
\centering
\includegraphics[scale=0.19]{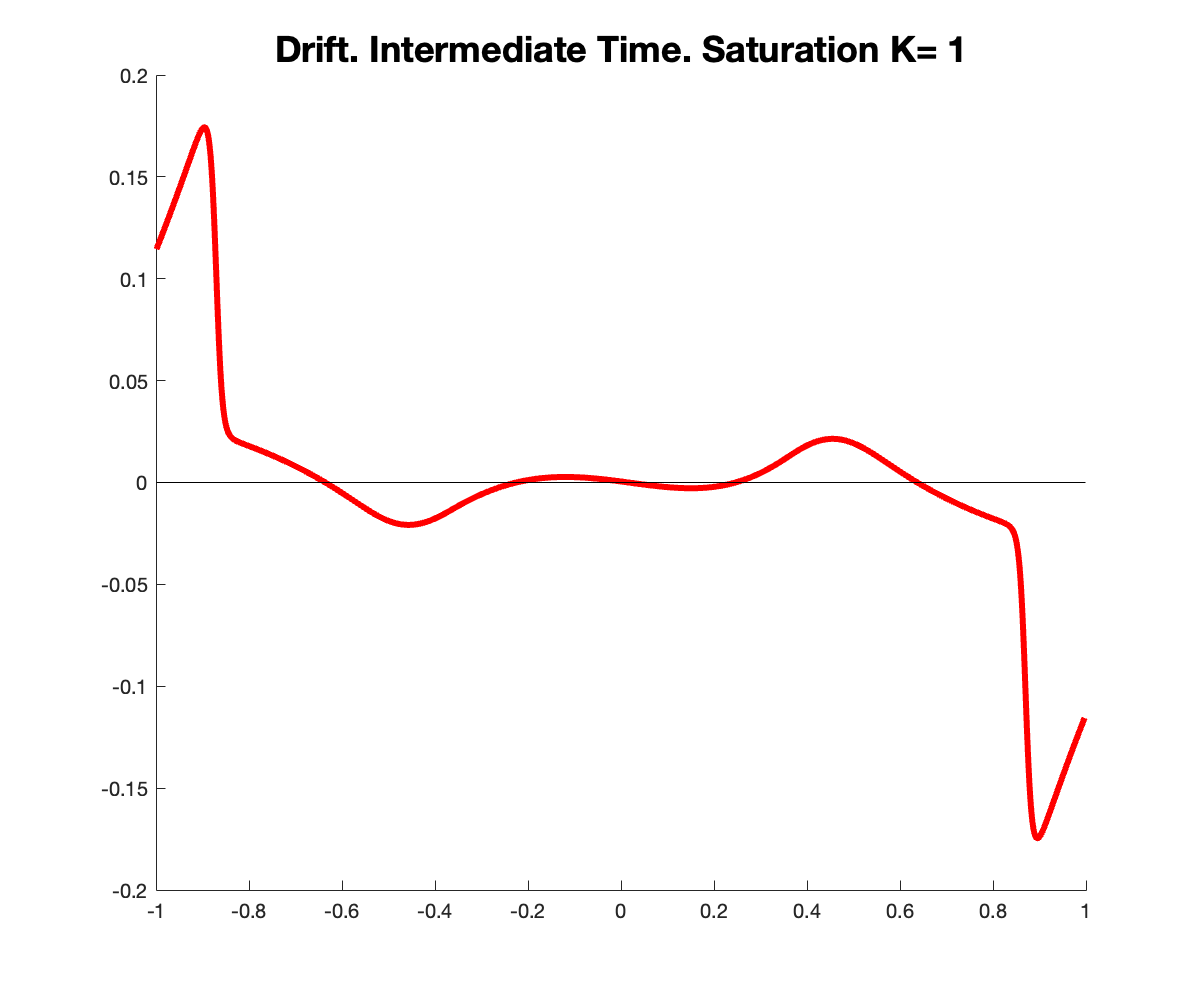}\!\!\!\!\!\!\!\!\!\includegraphics[scale=0.19]{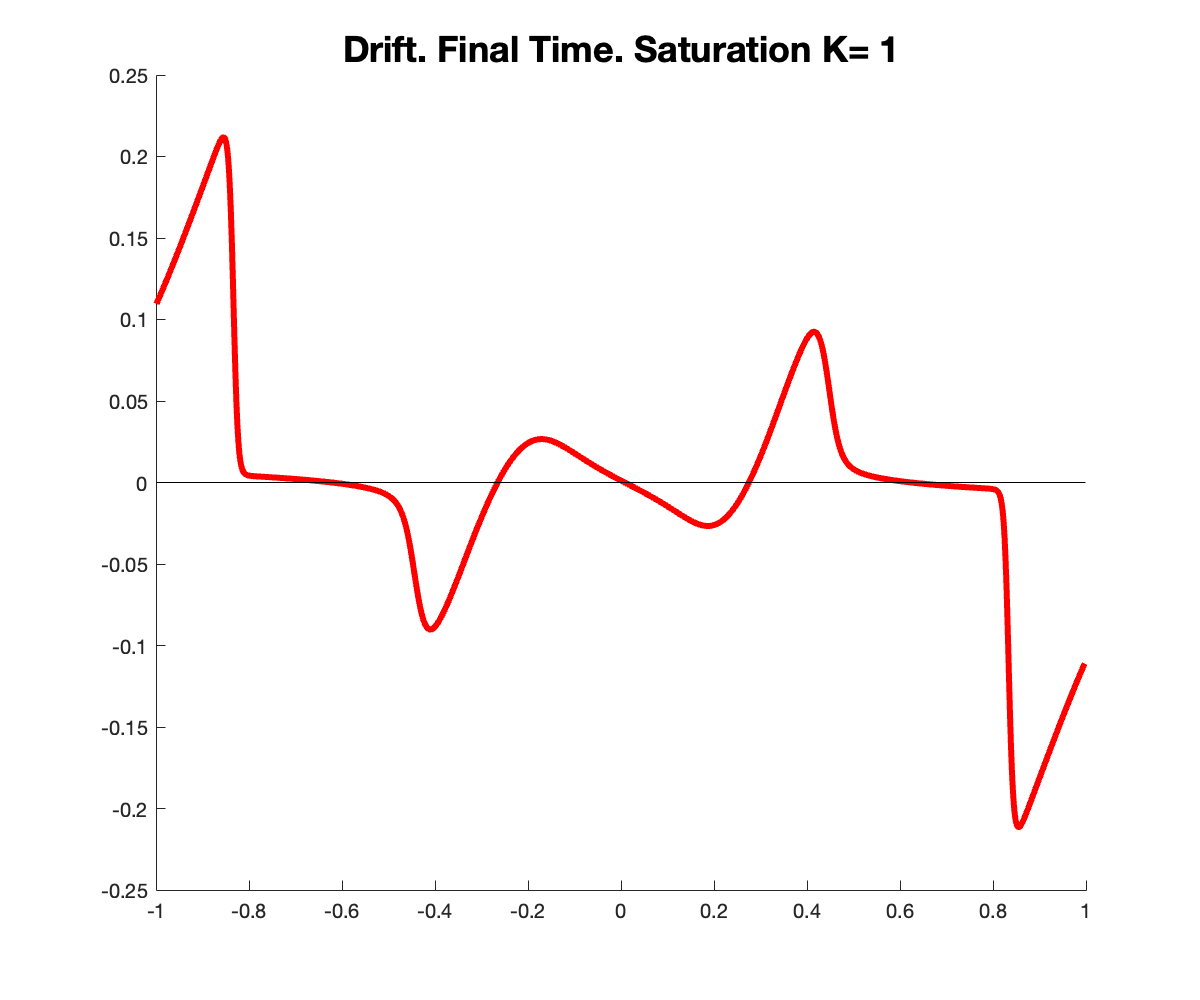}\\

\includegraphics[scale=0.3]{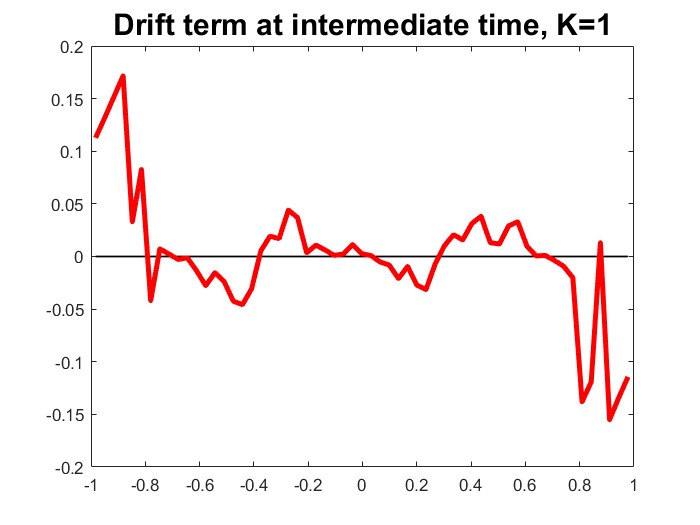}\includegraphics[scale=0.3]{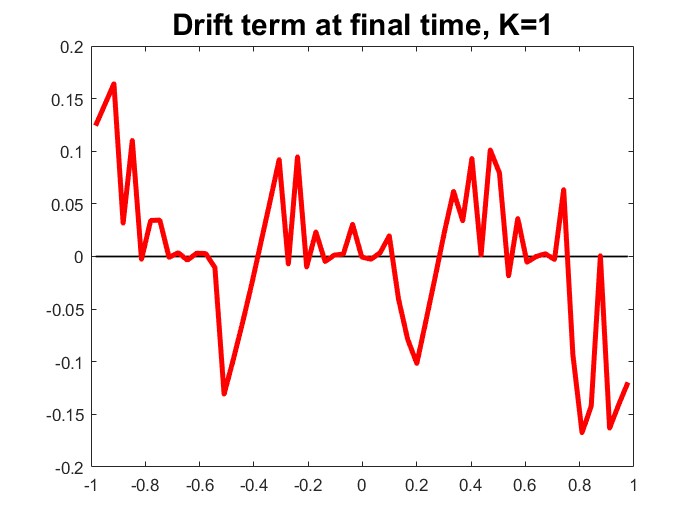}
\captionof{figure}{Local saturation model, $K=1$. First row trajectories with $N=3500$ cells, second row densities of the PDE model (\S 3.5.2) and third one histograms of SDEs (\S 4.3). Both at intermediate time on the left and at final time on the right.
} 
\label{fig:Carrillo_K1_drift}
\end{minipage}

\

\begin{minipage}{\textwidth}
\centering
\includegraphics[scale=0.3]{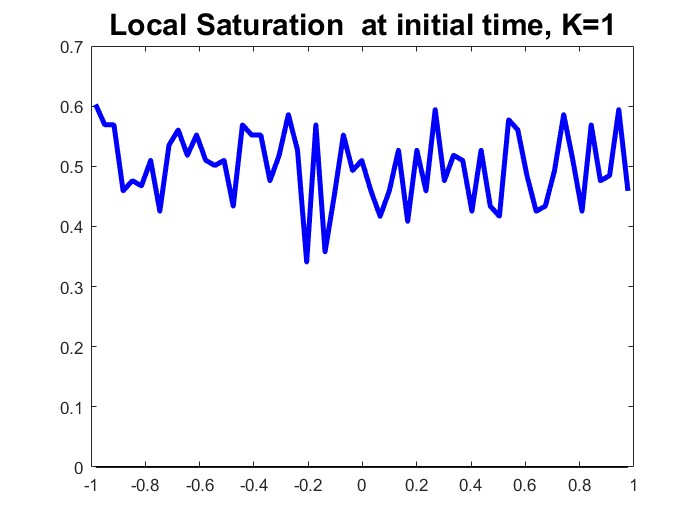} \includegraphics[scale=0.3]{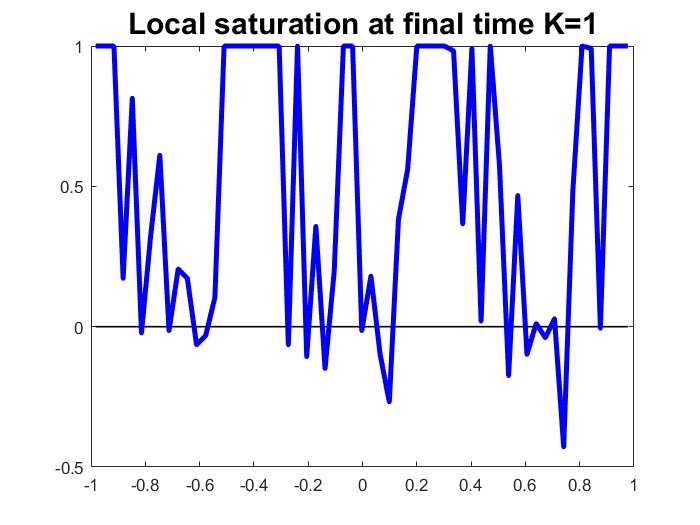}

\captionof{figure}{Local saturation coefficients for SDEs (\S 4.3), $K=1$. Initial time on the left and  final time on the right.
} 
\label{fig:Carrillo_K1_satu}
\end{minipage}


     \begin{minipage}{\textwidth}
\centering
\includegraphics[scale=0.24]{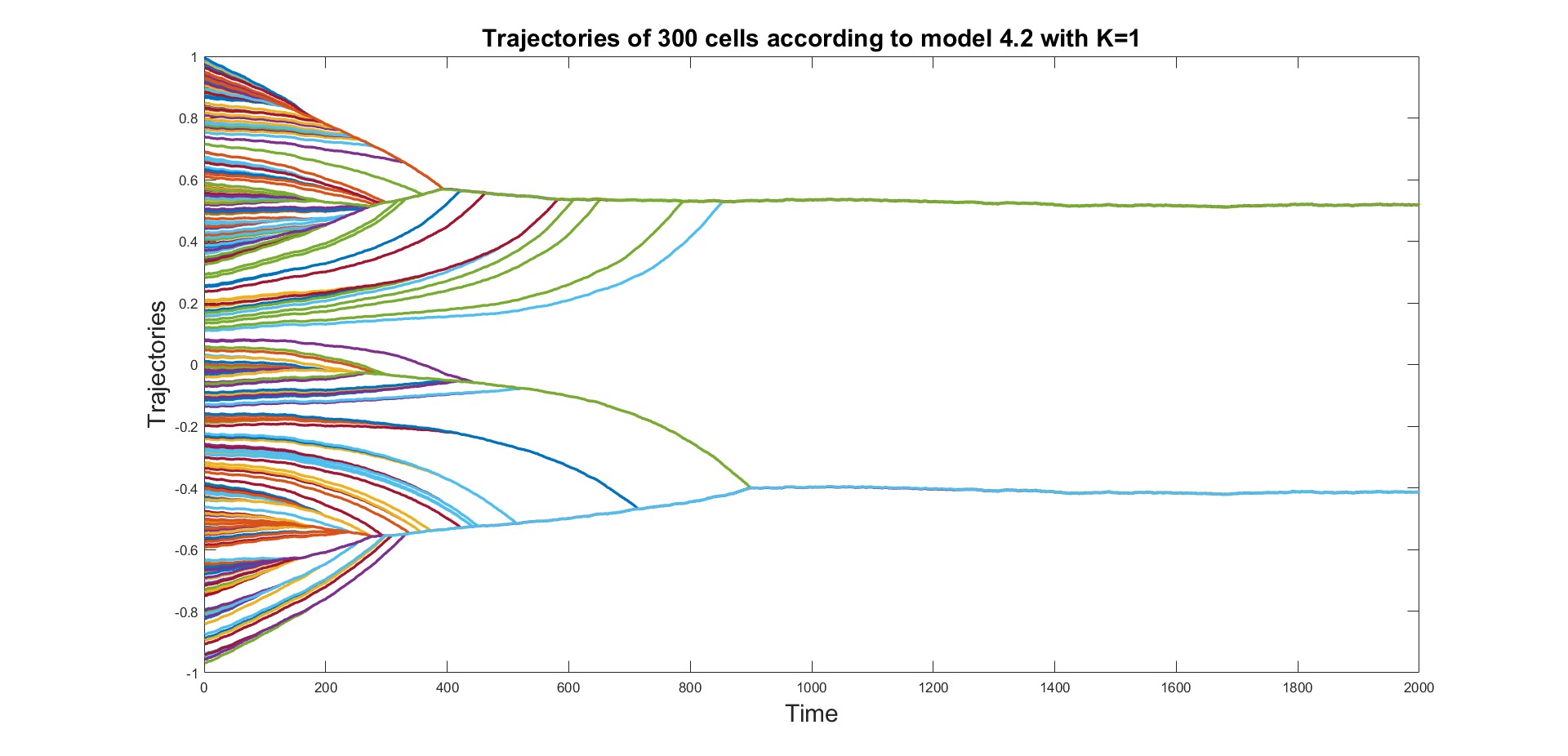}\\

\includegraphics[scale=0.18]{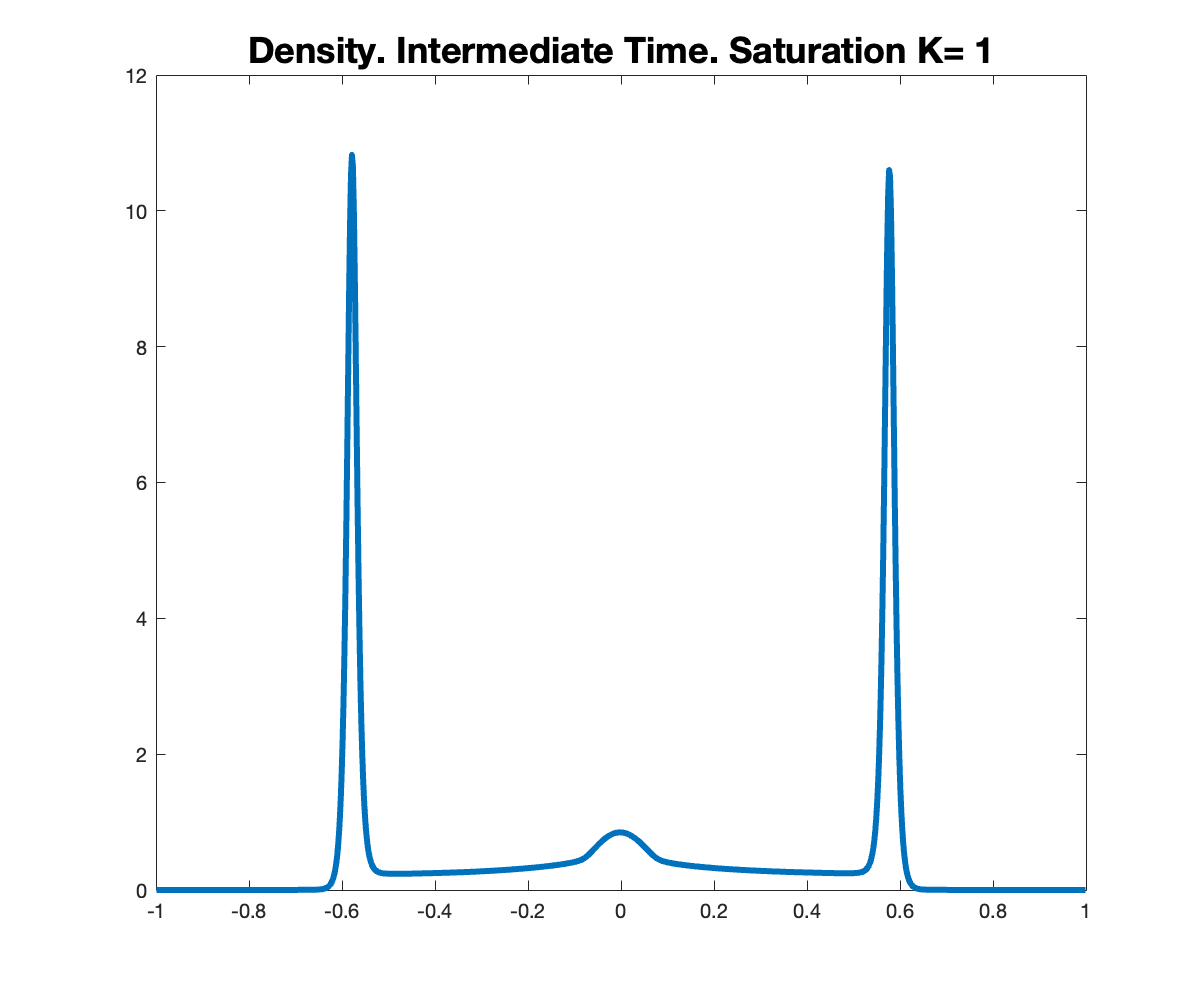}\includegraphics[scale=0.18]{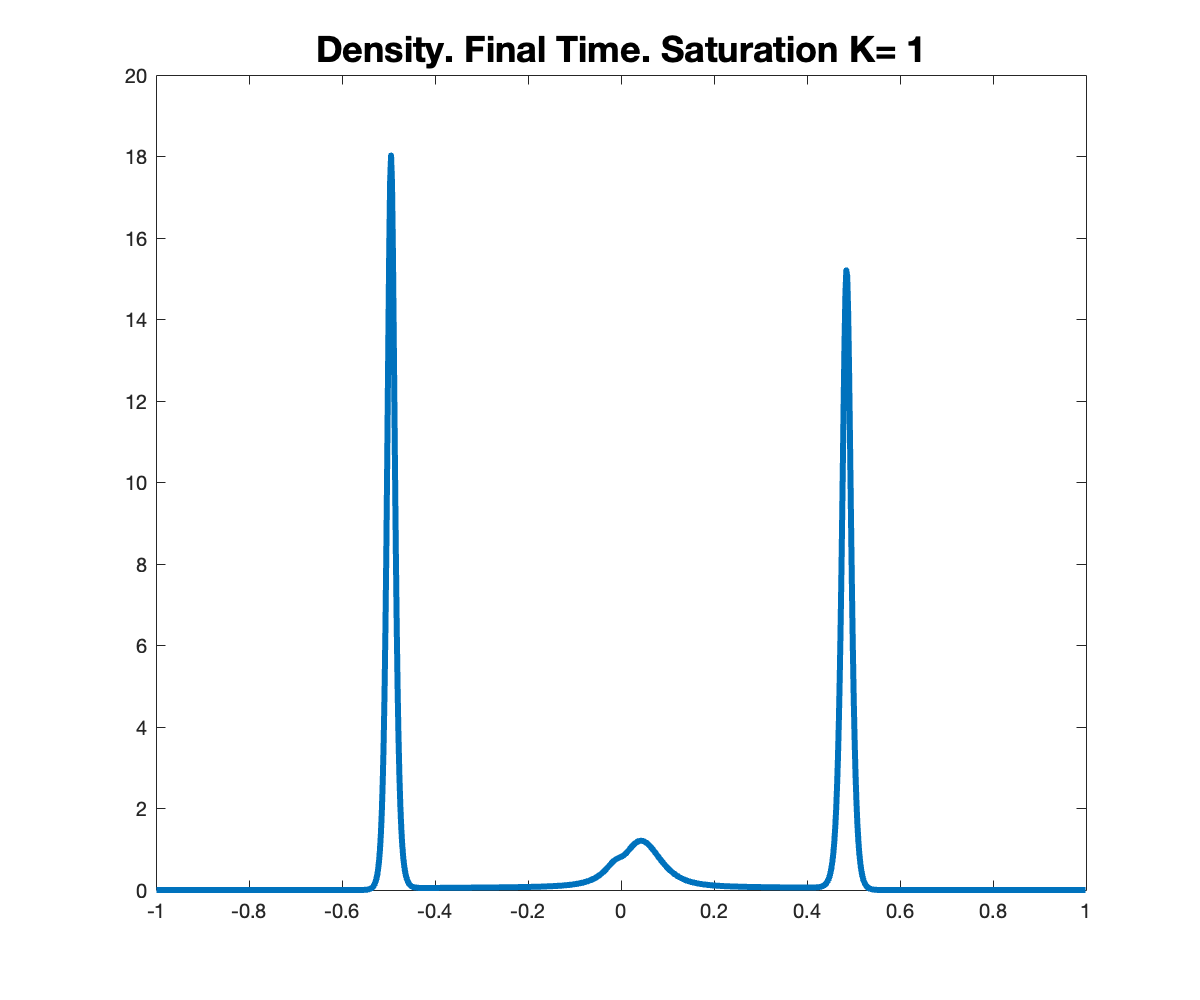}\\
\includegraphics[scale=0.3]{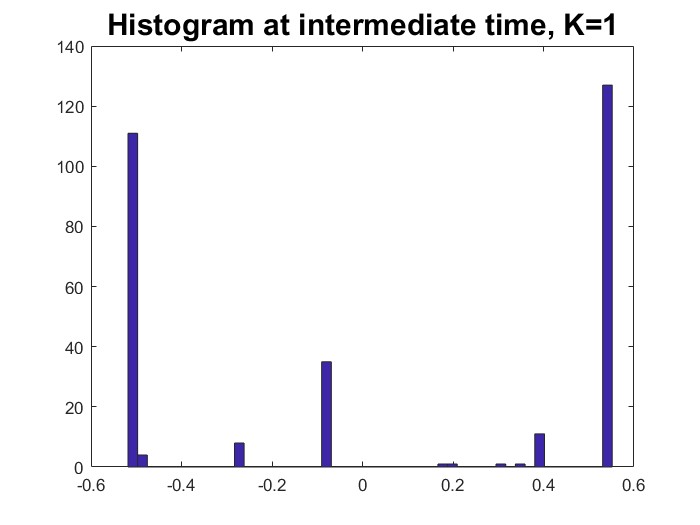}\!\!\includegraphics[scale=0.3]{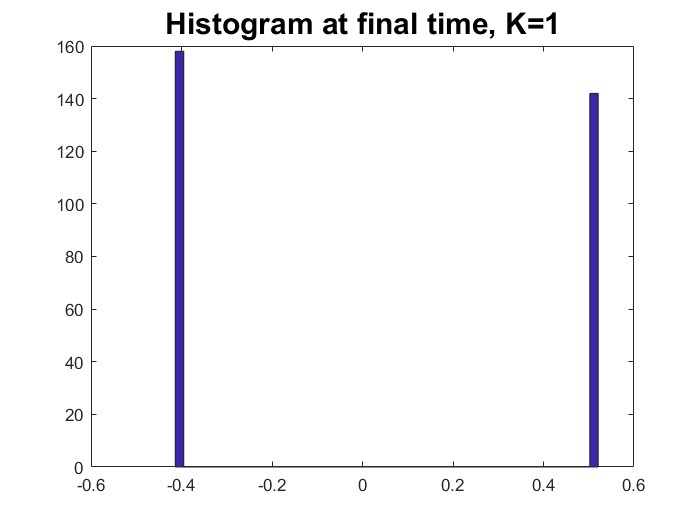}
\captionof{figure}{Nonlocal saturation model, $K=1$. First row trajectories with $N=300$ cells, second row densities of the PDE model (\S 3.5.3) and third one histograms of SDEs (\S 4.2). Both at intermediate time on the left and at final time on the right.} 
\label{fig:NL_K1_histogramas}
\end{minipage}

 \begin{minipage}{\textwidth}
\centering
\includegraphics[scale=0.15]{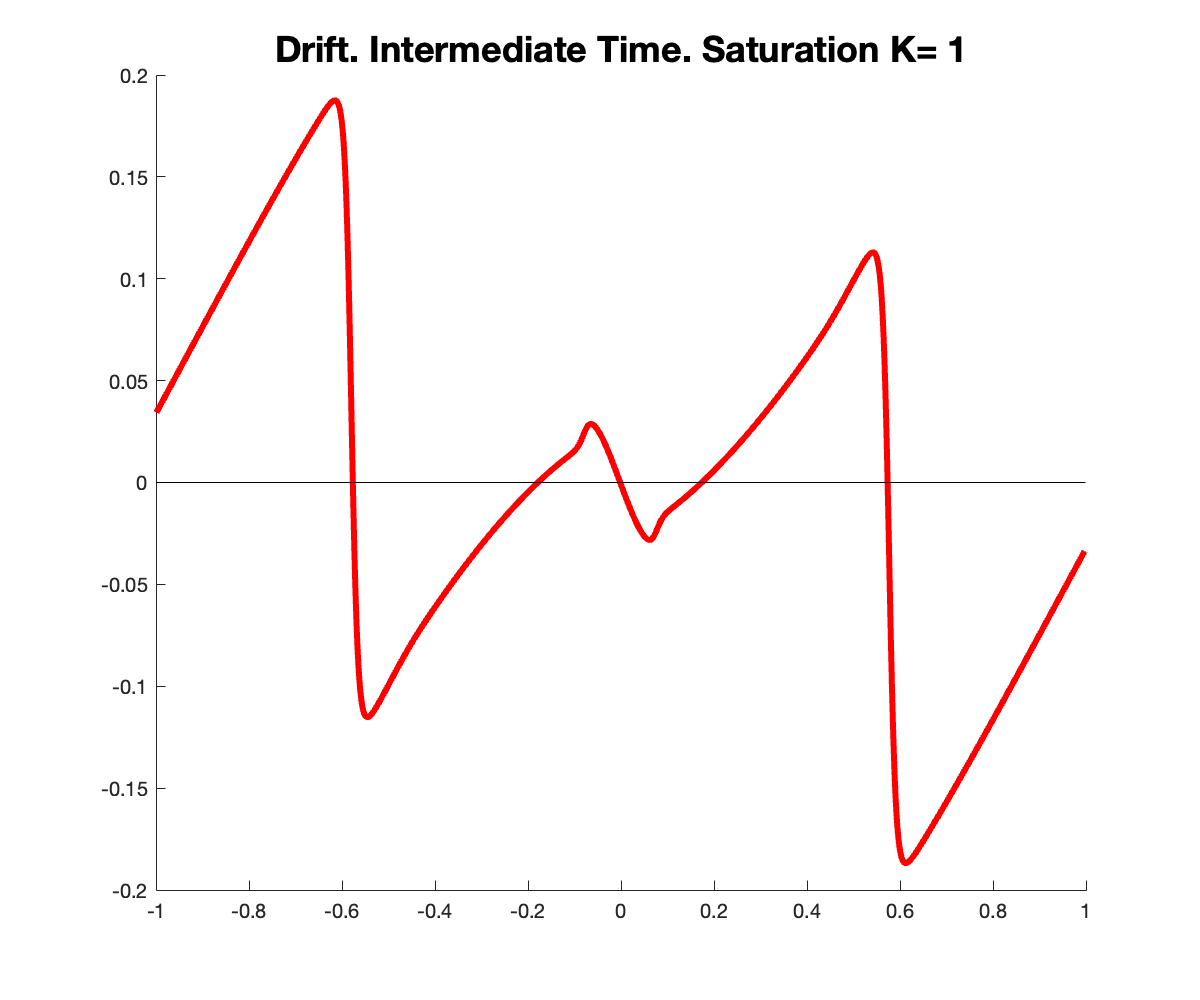}\!\!\!\!\includegraphics[scale=0.15]{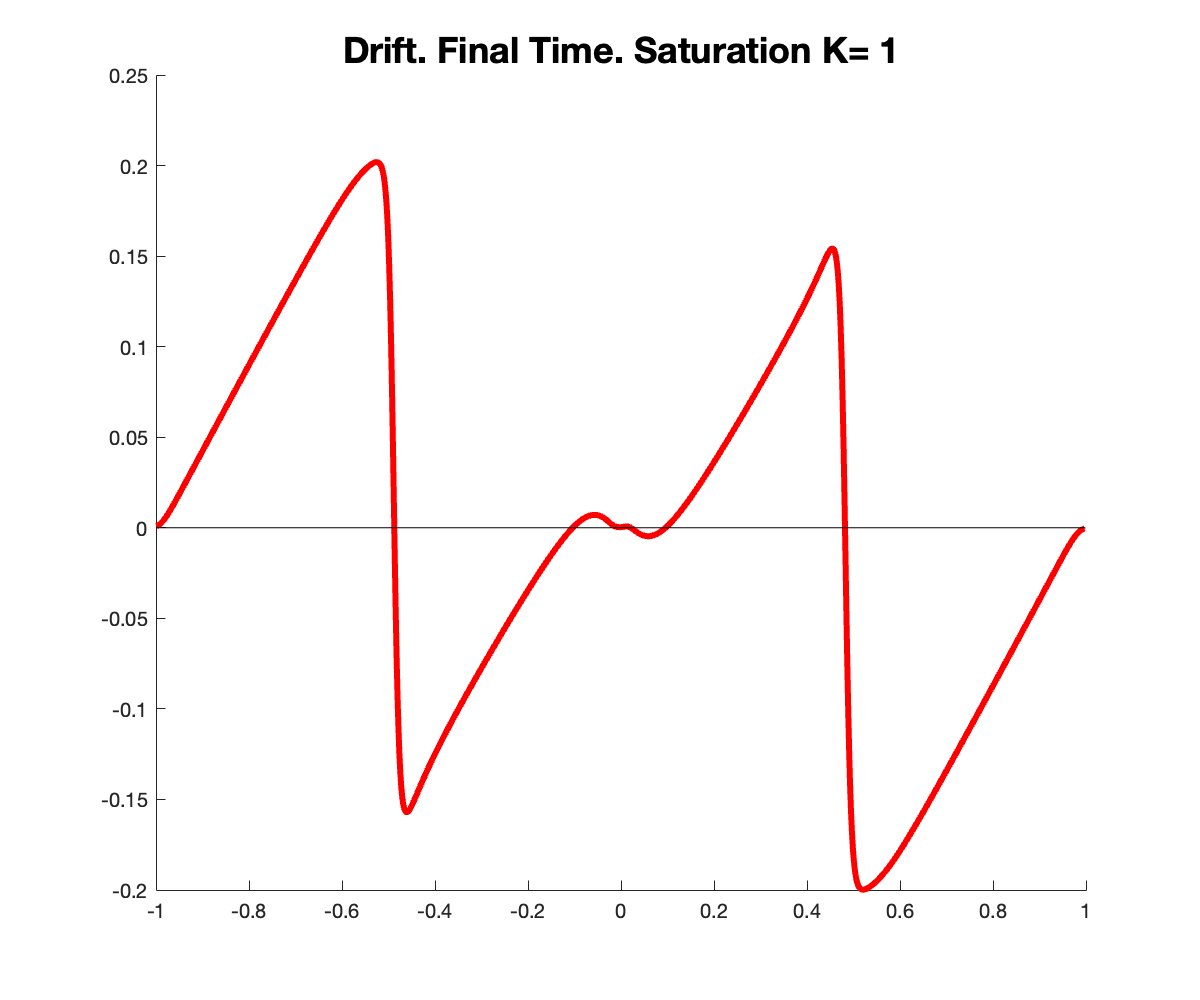}\\

\includegraphics[scale=0.26]{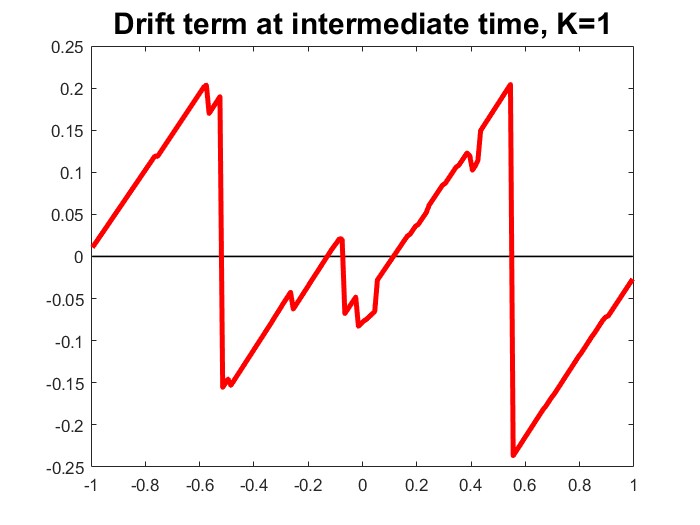}\!\!\includegraphics[scale=0.26]{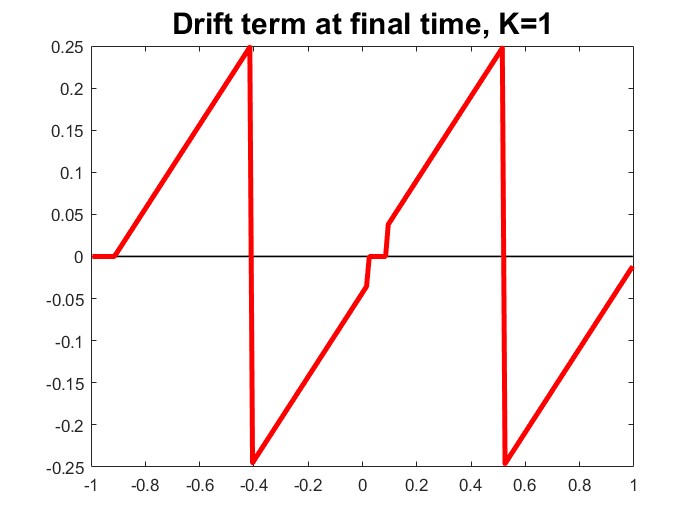}
\includegraphics[scale=0.15]{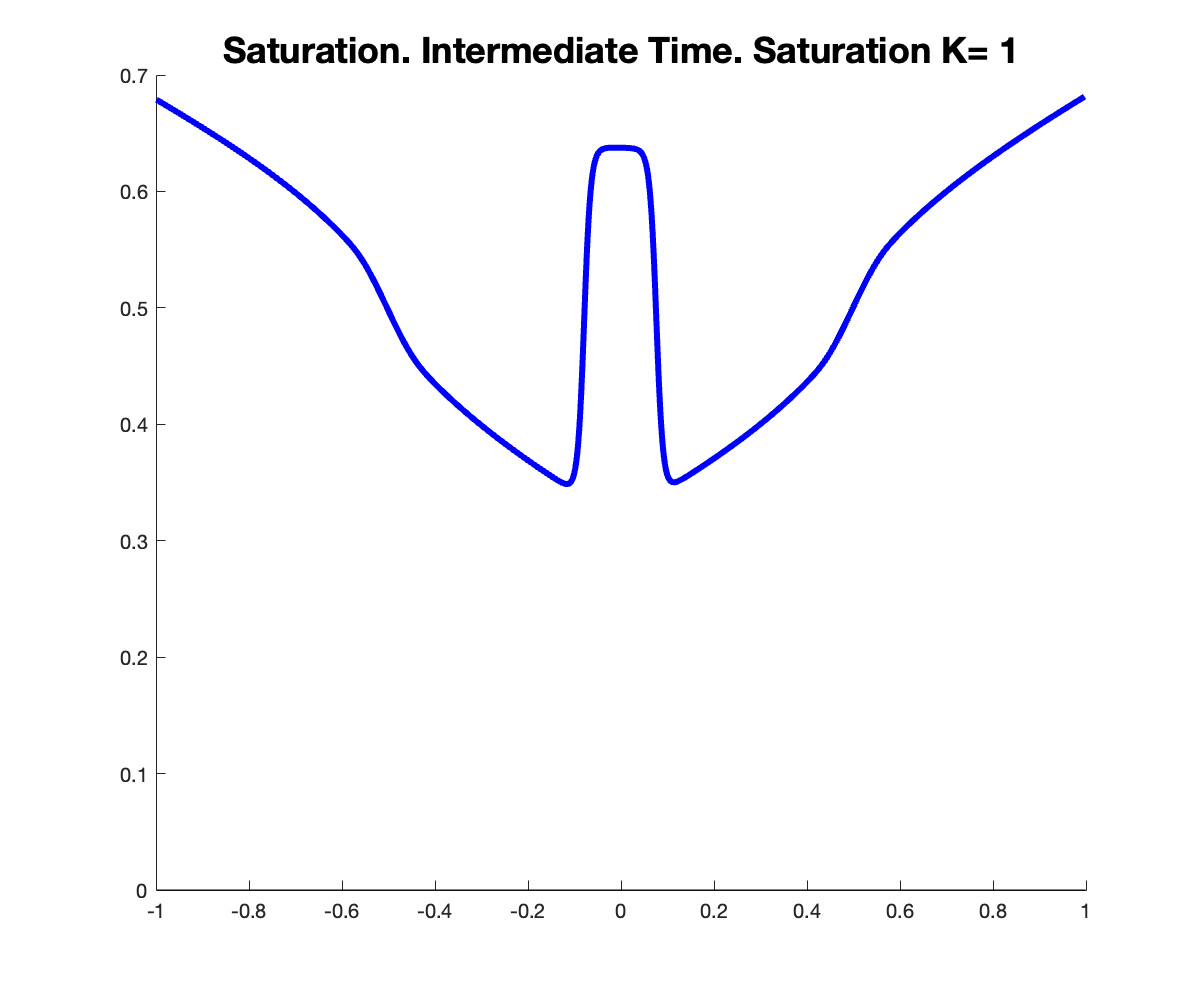}\!\!\!\!\includegraphics[scale=0.15]{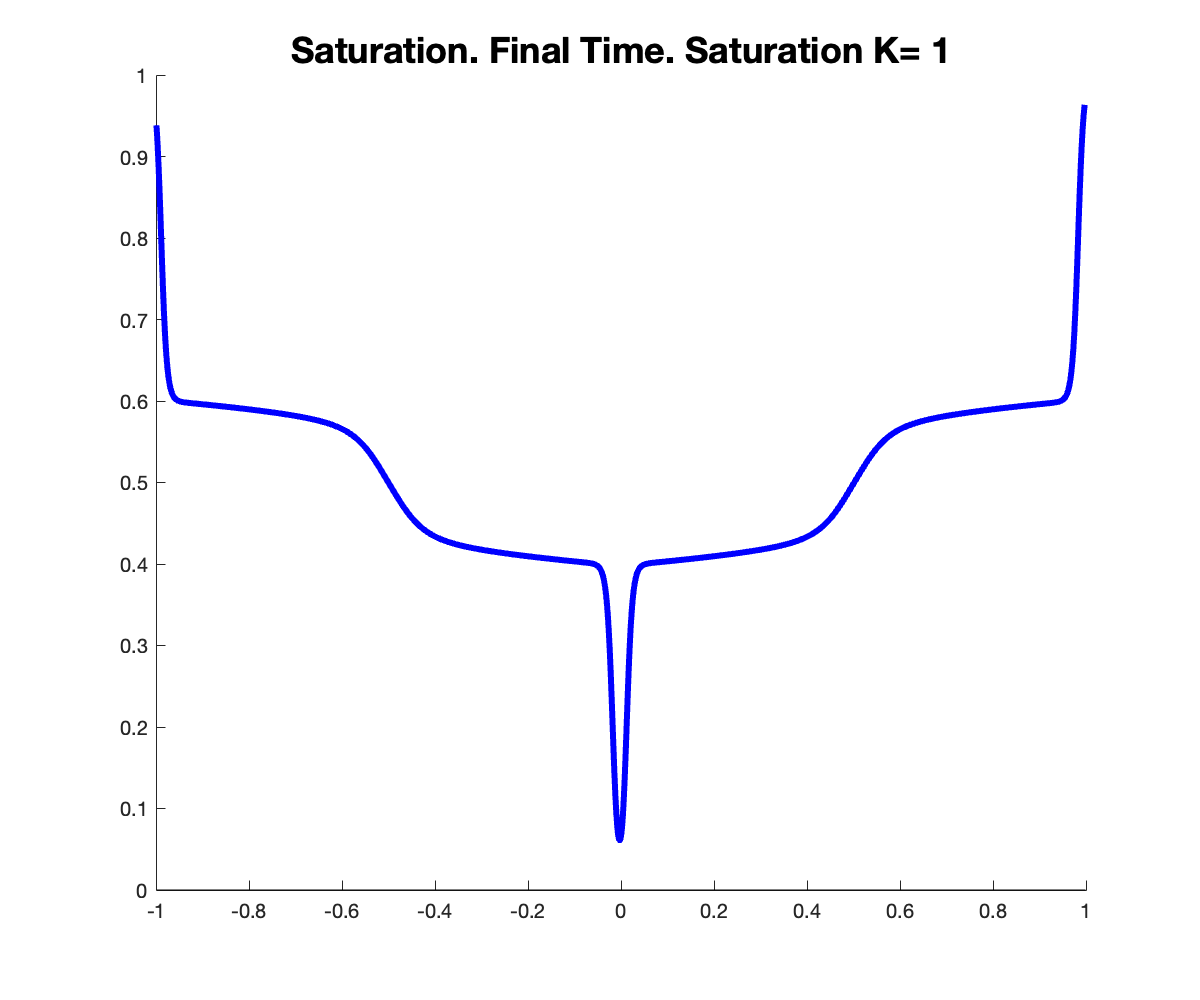}\\

\includegraphics[scale=0.26]{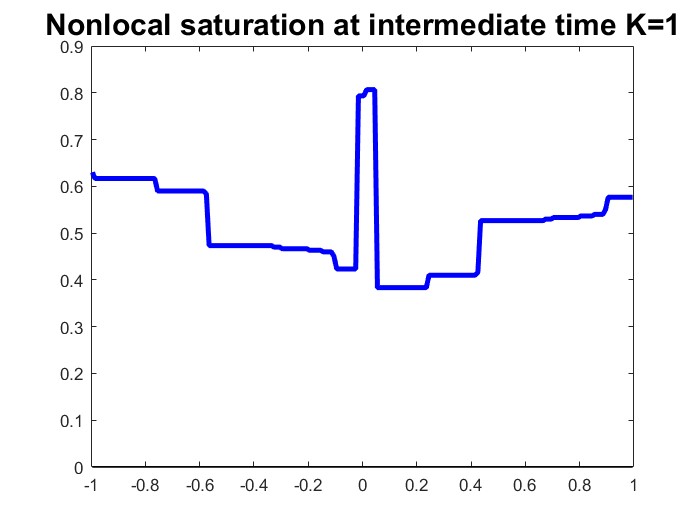}\!\!\includegraphics[scale=0.26]{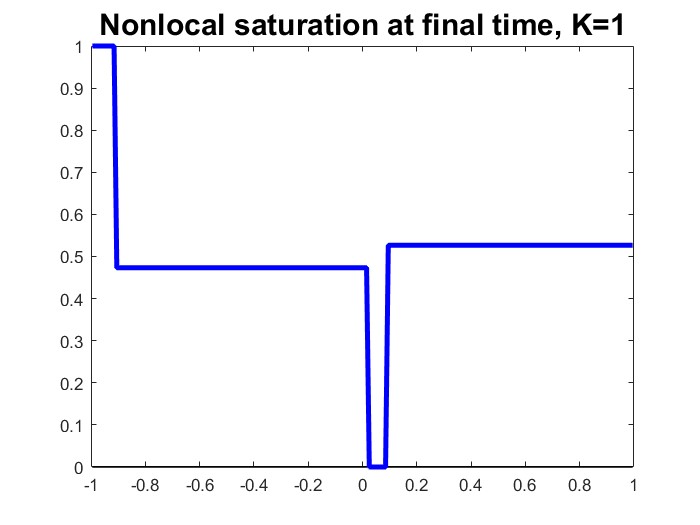}
\captionof{figure}{Drifts and saturations for nonlocal saturation model, $K=1$. Drifts in first and second rows for  PDE model (\S 3.5.3) and SDEs (\S 4.2), respectively. Saturations in third and fourth rows for  PDE model (\S 3.5.3) and SDEs (\S 4.2), respectively. All of them at intermediate time on the left and at final time on the right.} 
\label{fig:NL_K1_driftsatu}
\end{minipage}


     \begin{minipage}{\textwidth}
\centering
\includegraphics[scale=0.24]{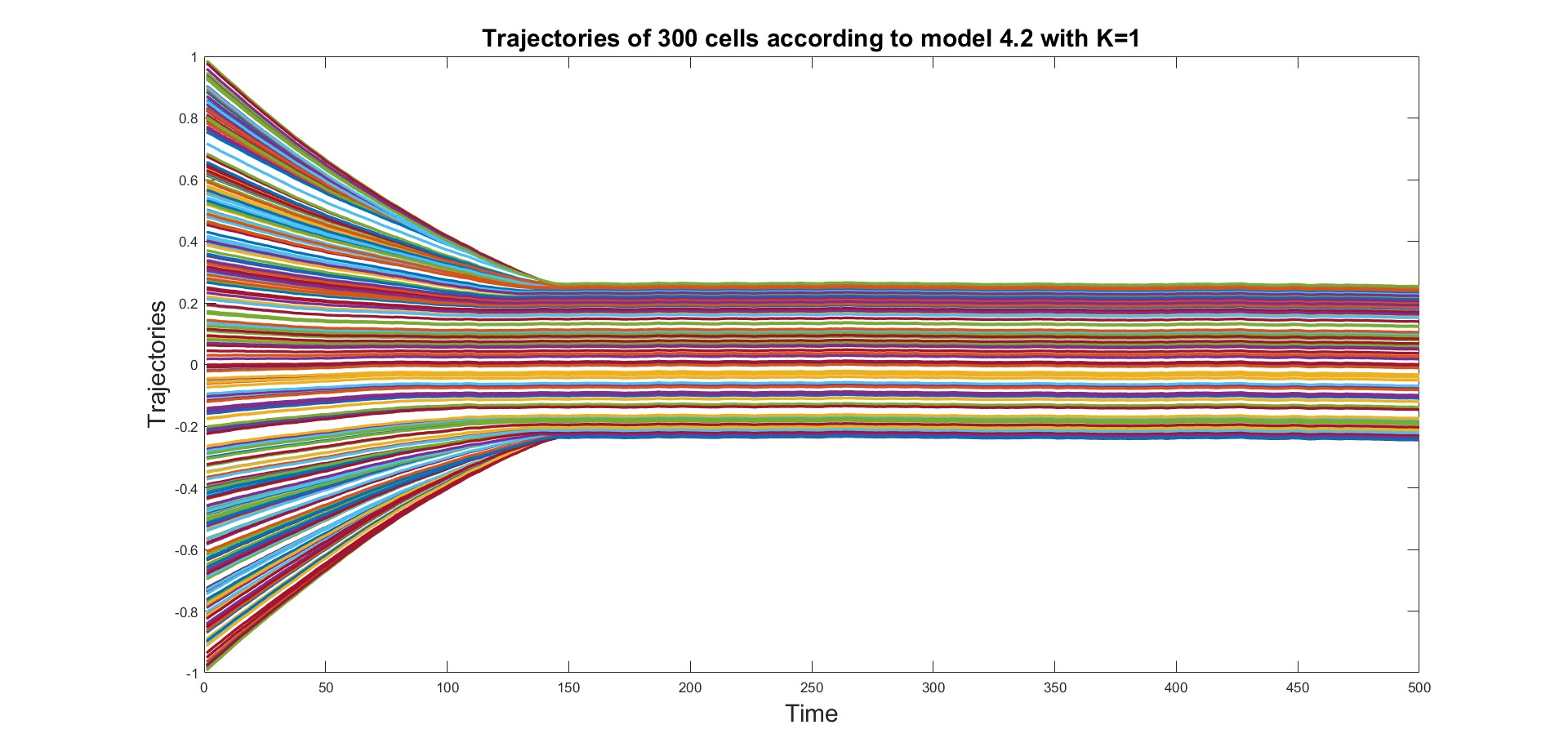}\\

\includegraphics[scale=0.18]{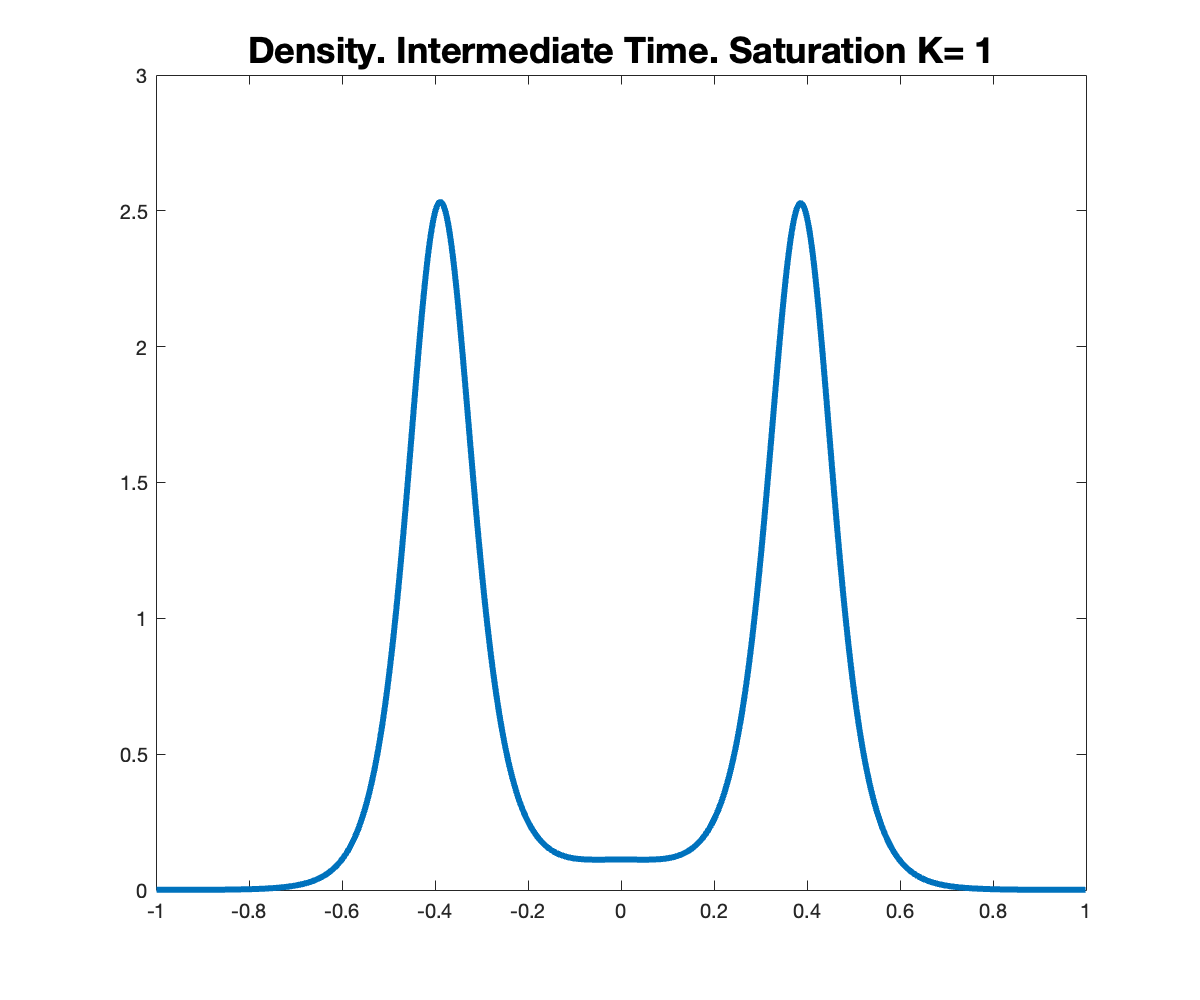}\includegraphics[scale=0.18]{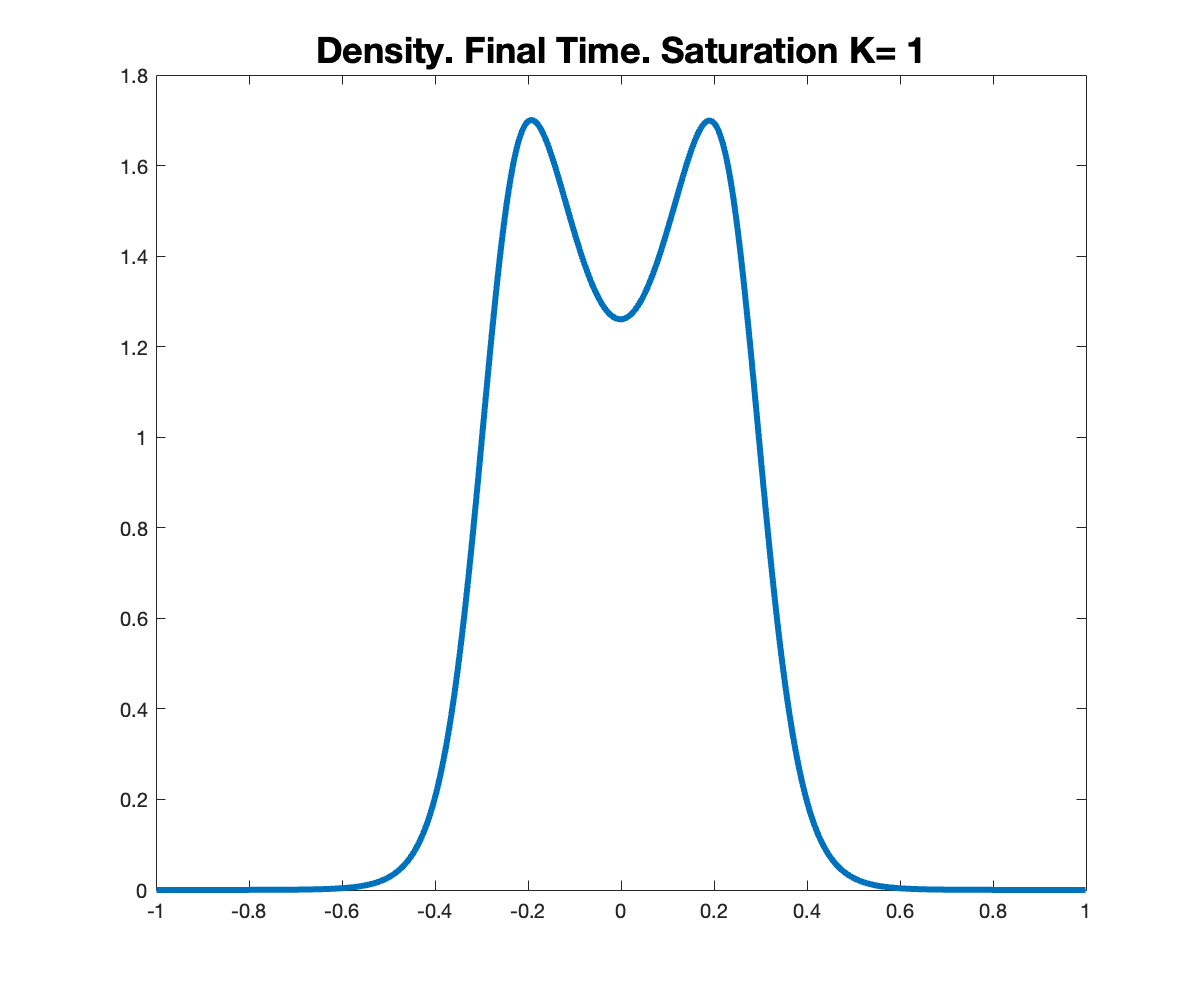}\\
\includegraphics[scale=0.3]{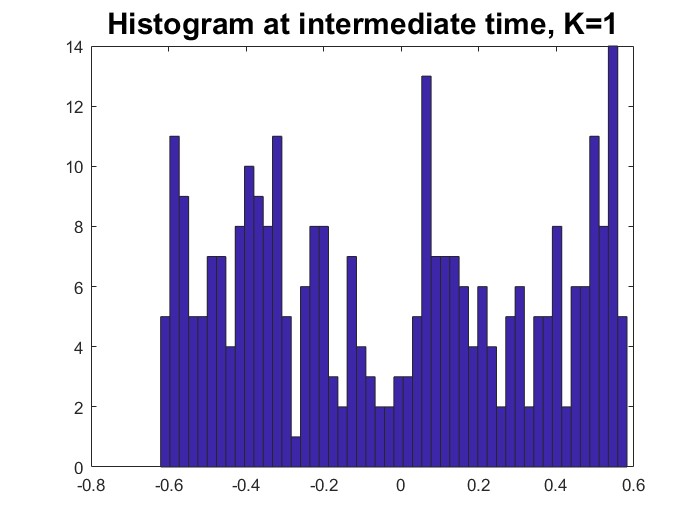}\!\!\includegraphics[scale=0.3]{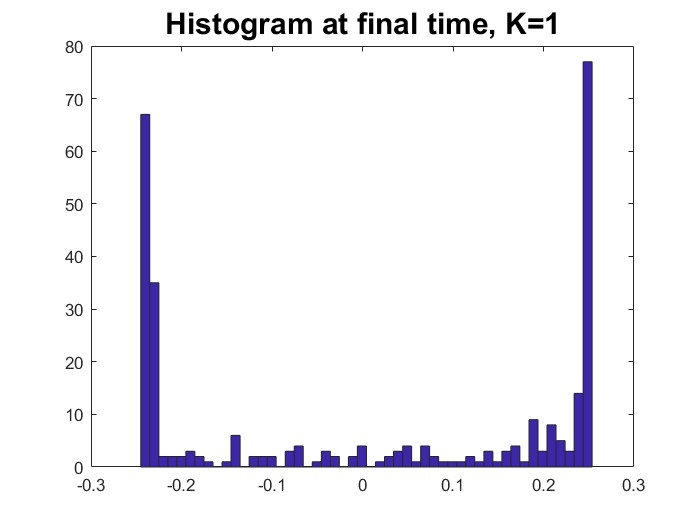}
\captionof{figure}{Nonlocal saturation model, $K=1$ and constant weight. First row trajectories with $N=300$ cells, second row densities of the PDE model (\S 3.5.3) and third one histograms of SDEs (\S 4.2). Both at intermediate time on the left and at final time on the right.} 
\label{fig:NL_K1_cw_histogramas}
\end{minipage}

 \begin{minipage}{\textwidth}
\centering
\includegraphics[scale=0.15]{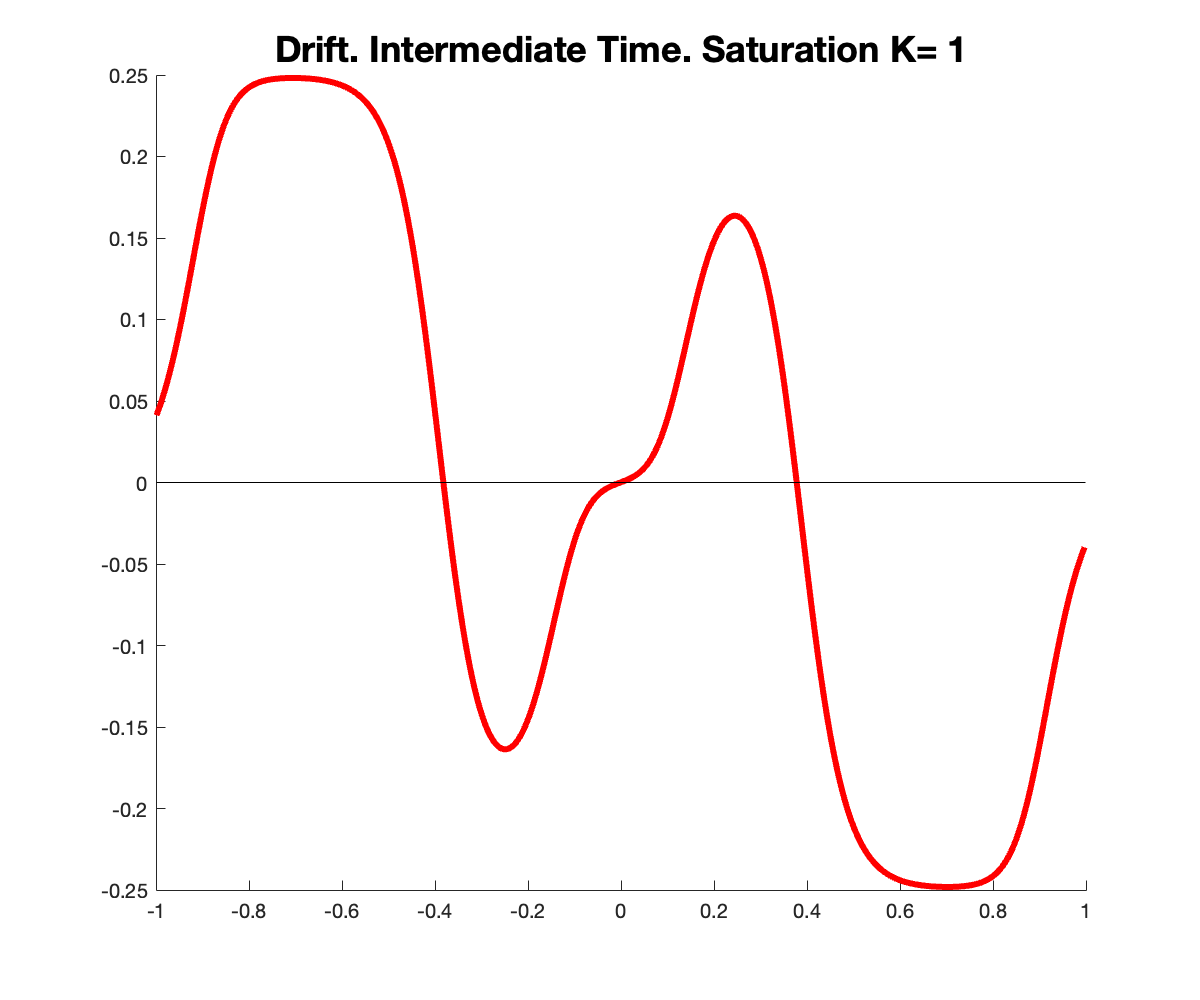}\!\!\!\!\includegraphics[scale=0.15]{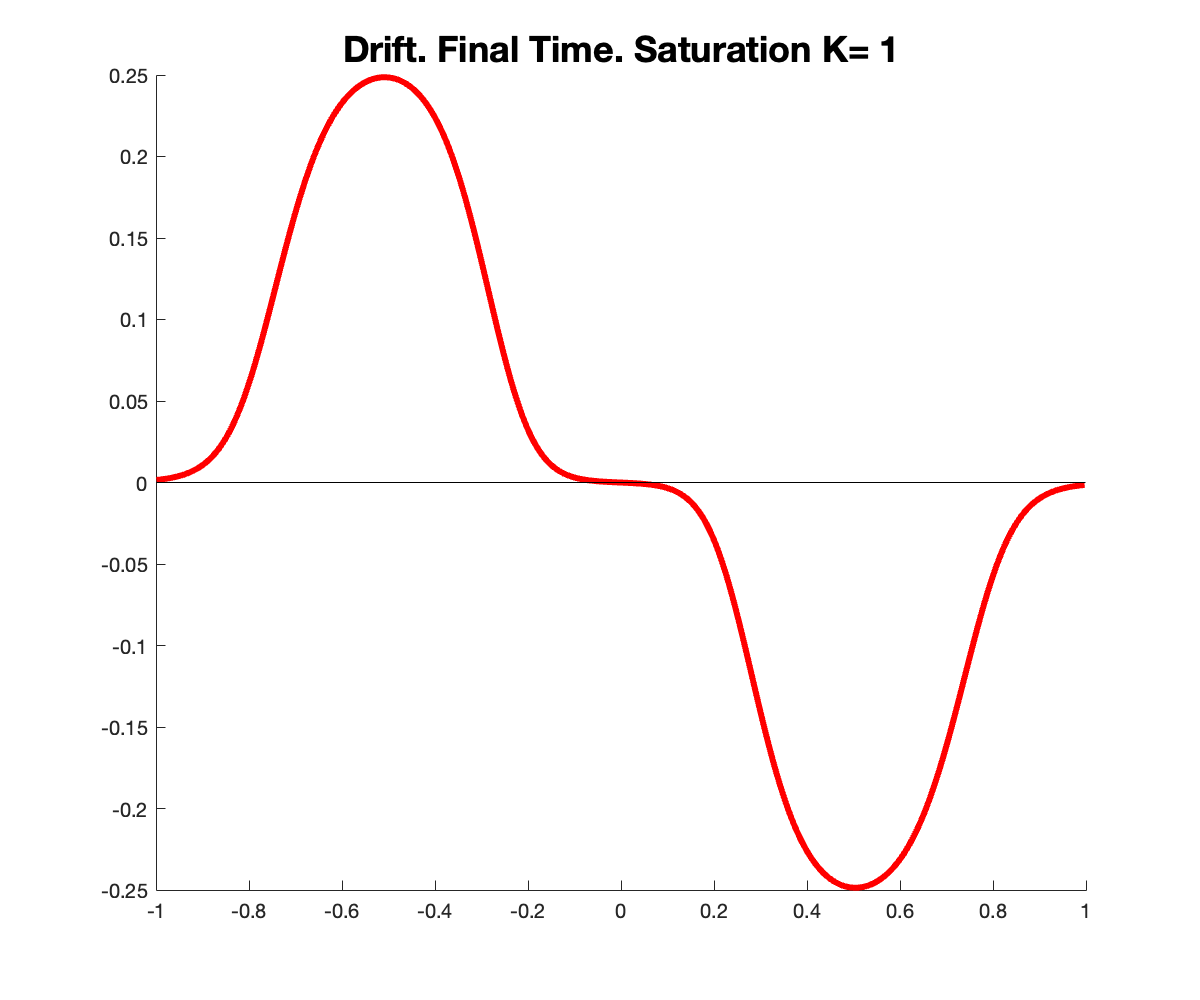}\\

\includegraphics[scale=0.26]{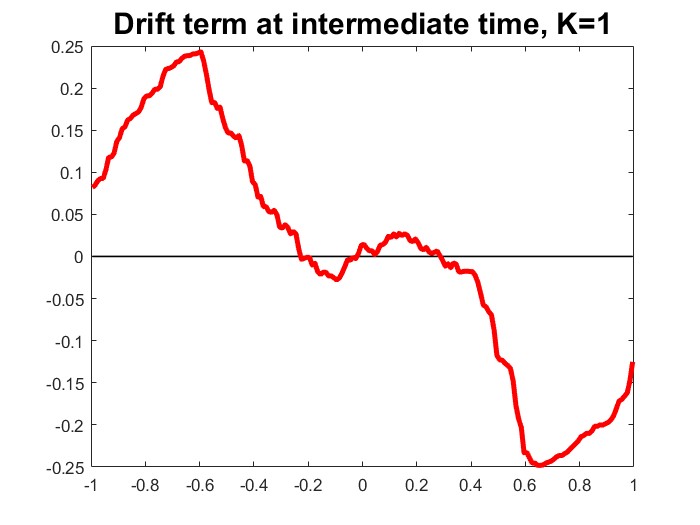}\!\!\includegraphics[scale=0.26]{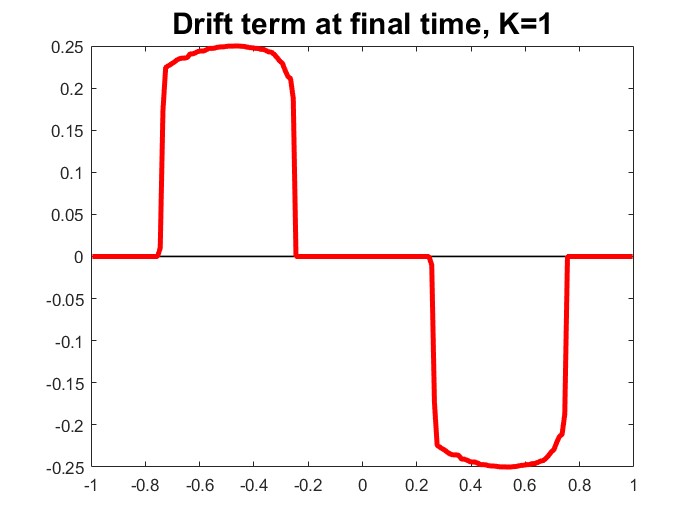}
\includegraphics[scale=0.15]{DEFINITIVAS/Modelo3K1/Saturation_Intermedia_Time_K=1__u0=.5+.1rand.png}\!\!\!\!\includegraphics[scale=0.15]{DEFINITIVAS/Modelo3K1/Saturation_Final_Time_K=1__u0=.5+.1rand.png}\\

\includegraphics[scale=0.26]{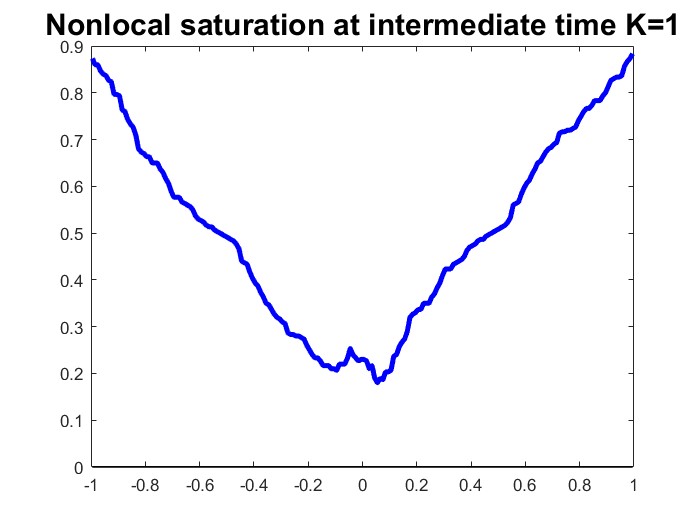}\!\!\includegraphics[scale=0.26]{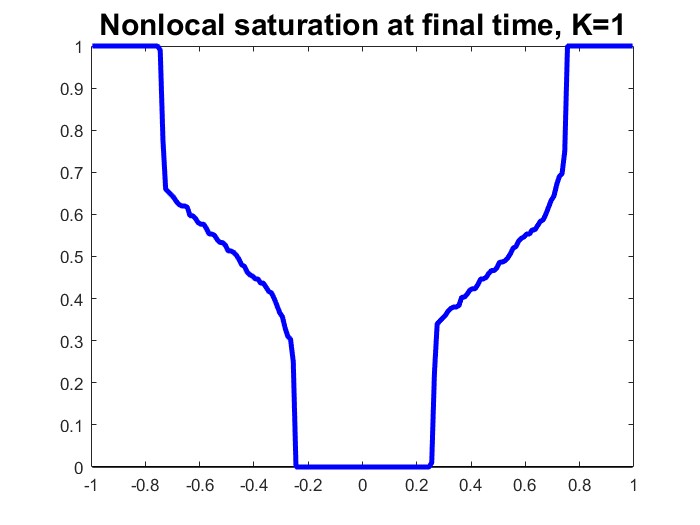}
\captionof{figure}{Drifts and saturations for nonlocal saturation model, $K=1$ and constant weight. Drifts in first and second rows for  PDE model (\S 3.5.3) and SDEs (\S 4.2), respectively. Saturations in third and fourth rows for  PDE model (\S 3.5.3) and SDEs (\S 4.2), respectively. All of them at intermediate time on the left and at final time on the right.} 
\label{fig:NL_K1_driftsatu}
\end{minipage}

\subsection{Weak repulsion $K\sim0.6$}\label{k=06}

\

 \begin{minipage}{\textwidth}
\centering
\includegraphics[scale=0.22]{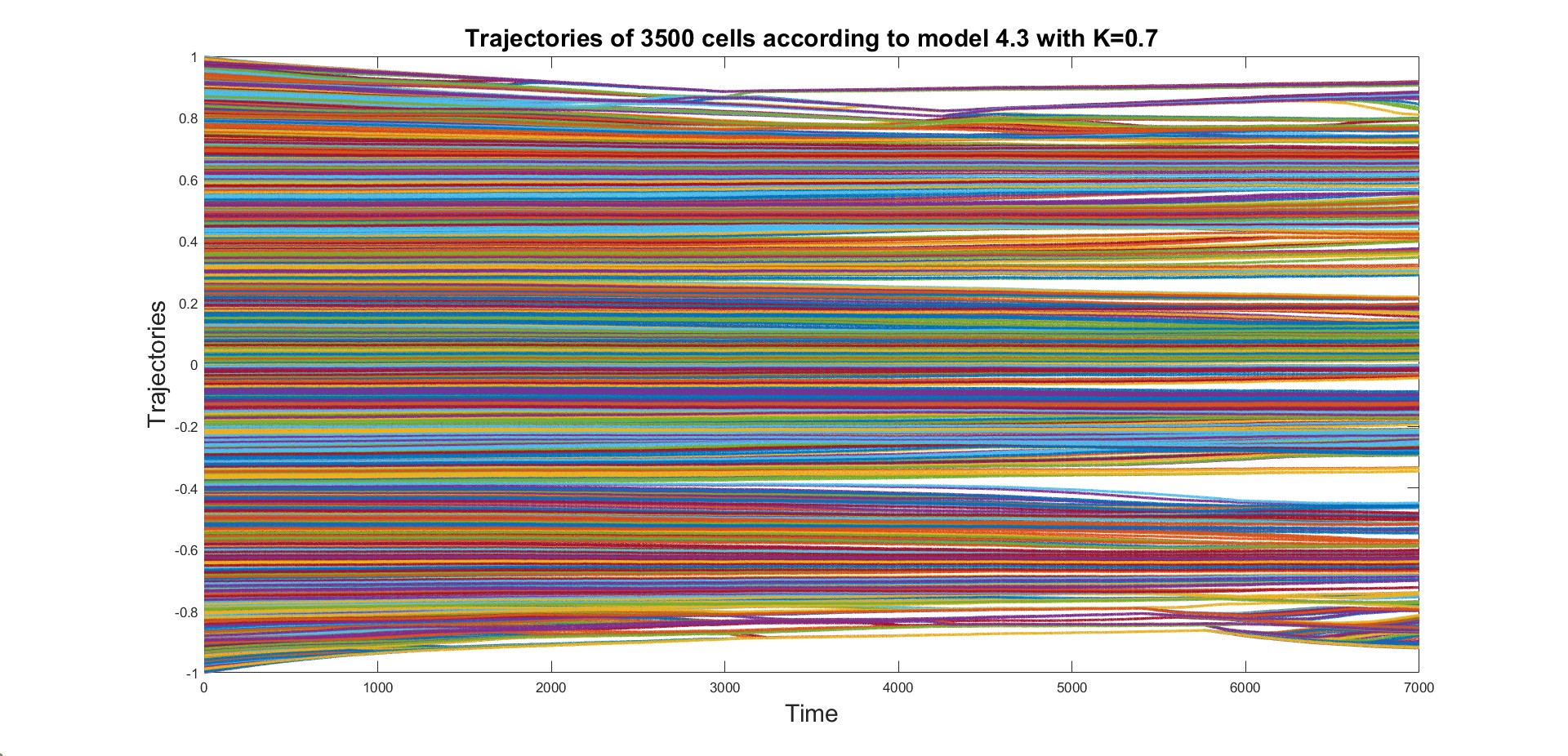}\\

\includegraphics[scale=0.17]{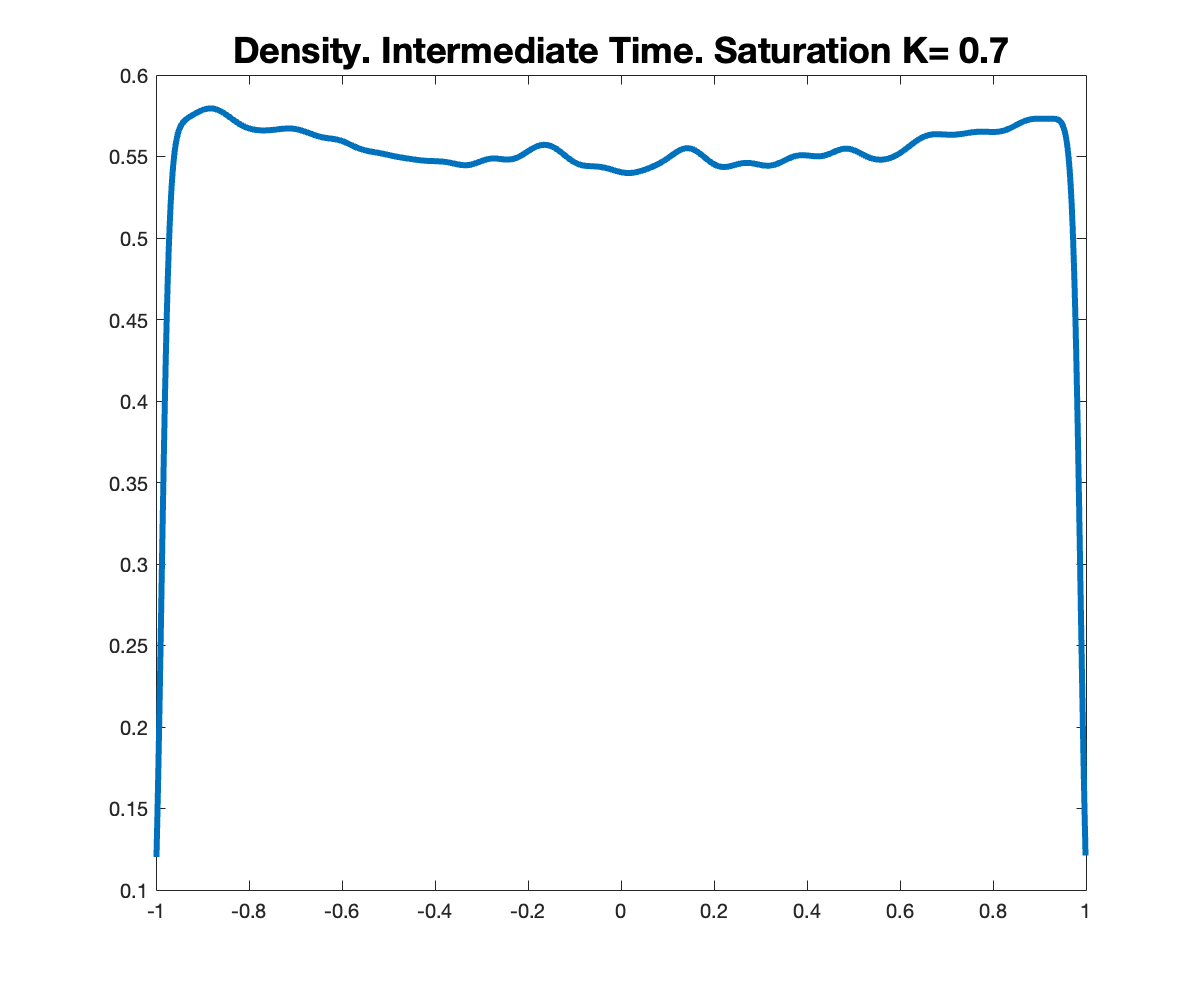}\!\!\includegraphics[scale=0.17]{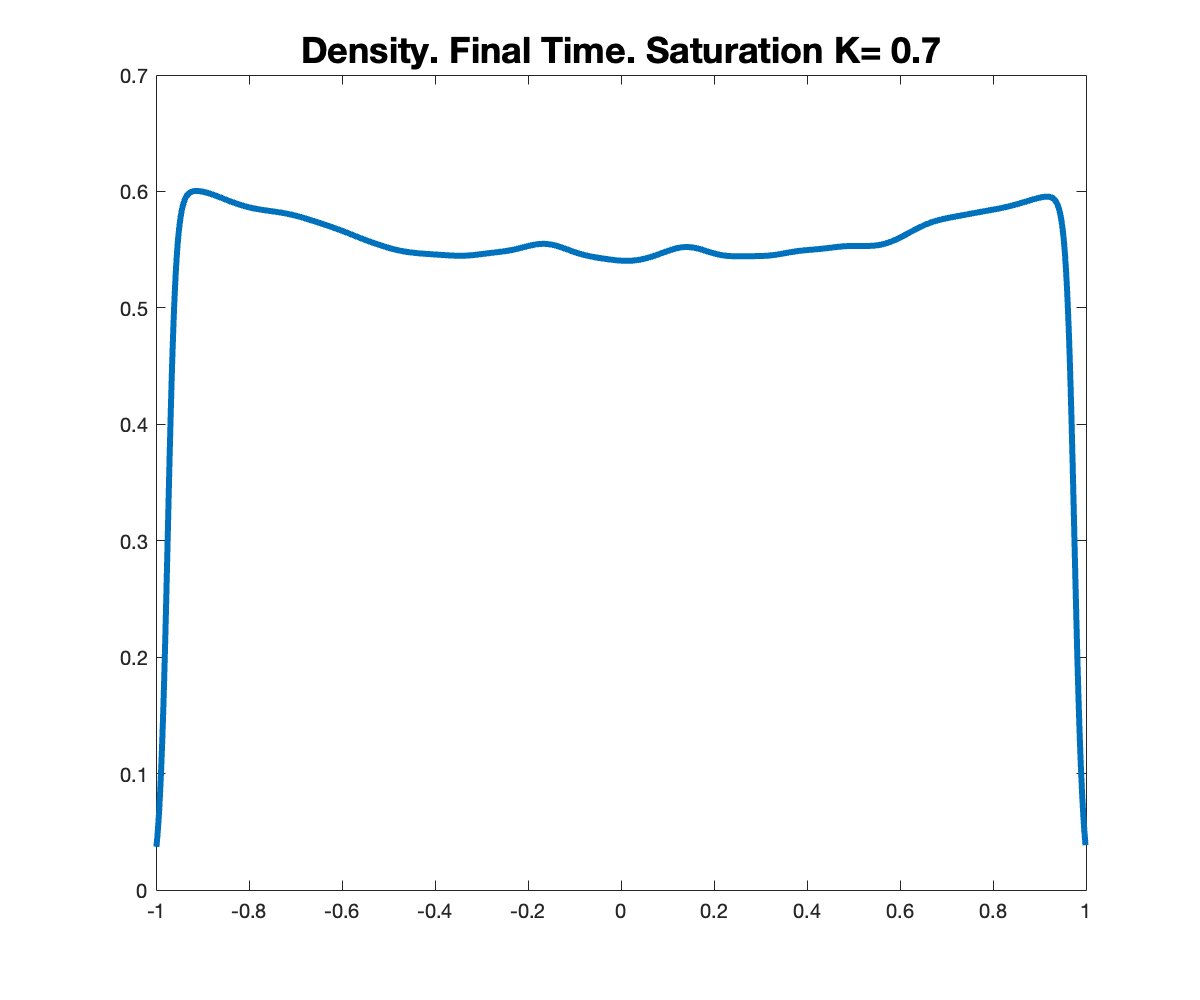}\\
\includegraphics[scale=0.28]{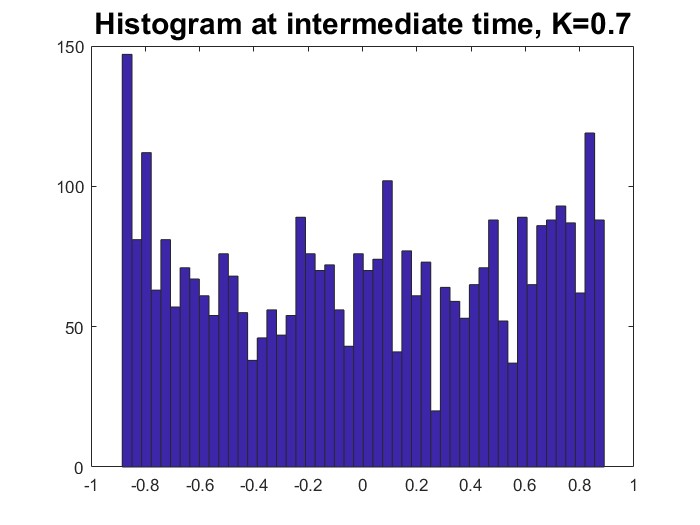}\!\!\includegraphics[scale=0.28]{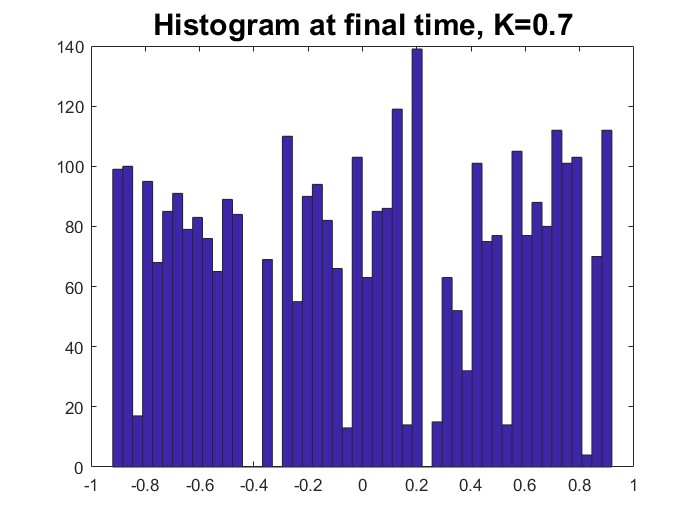}
\captionof{figure}{Local saturation model, $K=0.7$. First row trajectories with $N=3500$ cells, second row densities of the PDE model (\S 3.5.2) and third one histograms of SDEs (\S 4.3). Both at intermediate time on the left and at final time on the right.
} 
\label{fig:Carrillo__K07_histogramas}
\end{minipage}

  \begin{minipage}{\textwidth}
\centering
\includegraphics[scale=0.17]{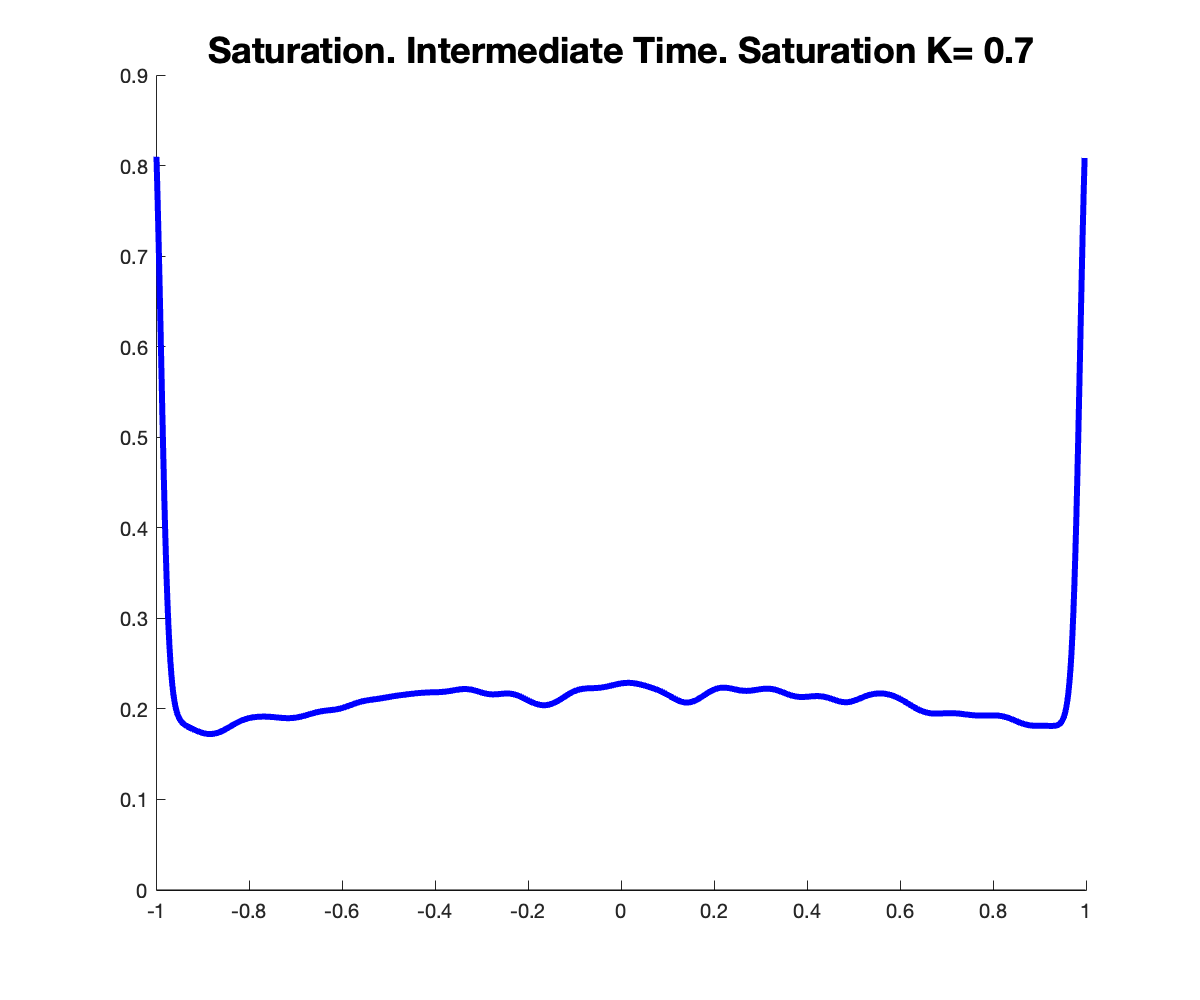}\!\!\!\!\!\!\!\!\!\includegraphics[scale=0.17]{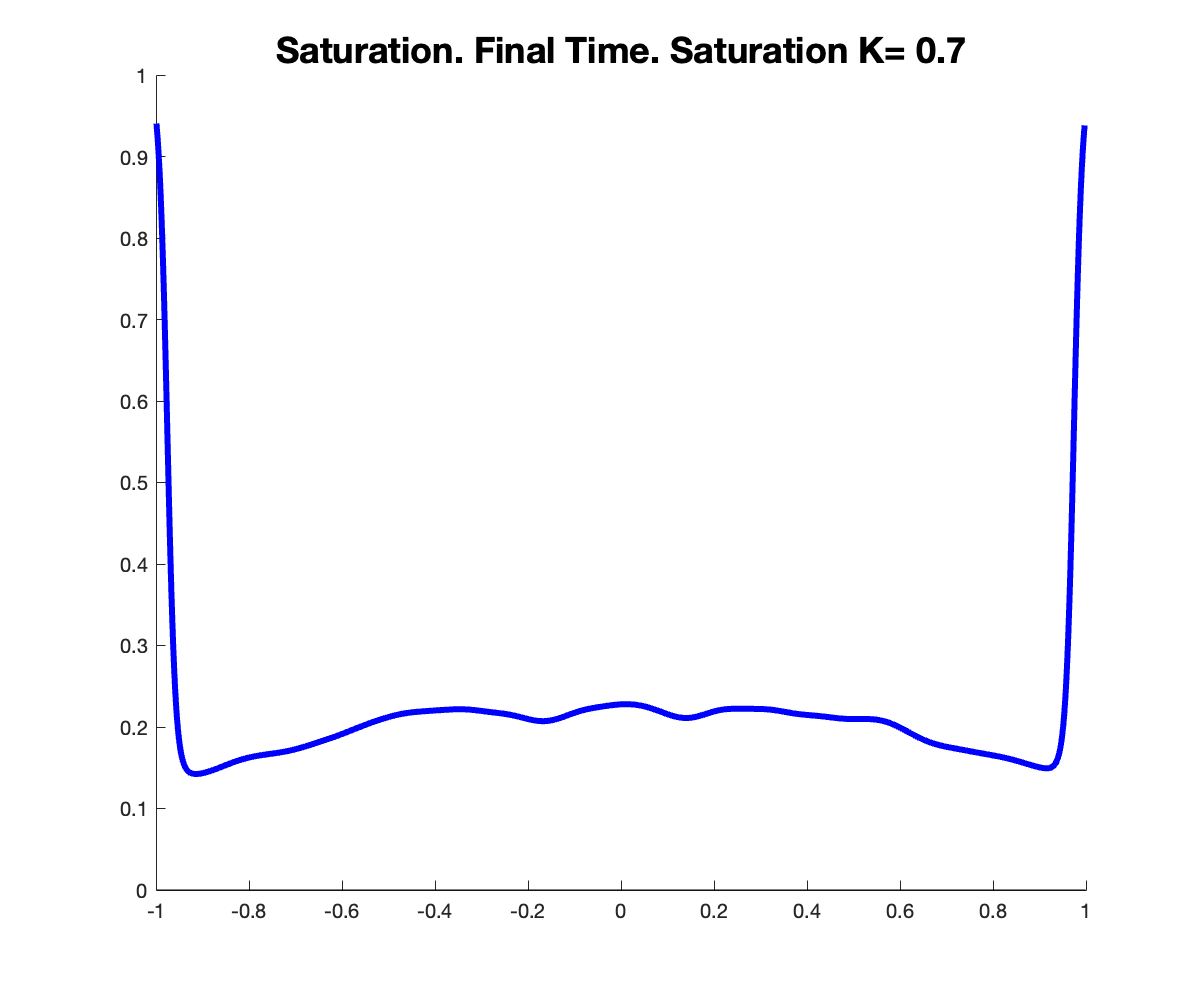}\\
\hspace{-1cm}\includegraphics[scale=0.25]{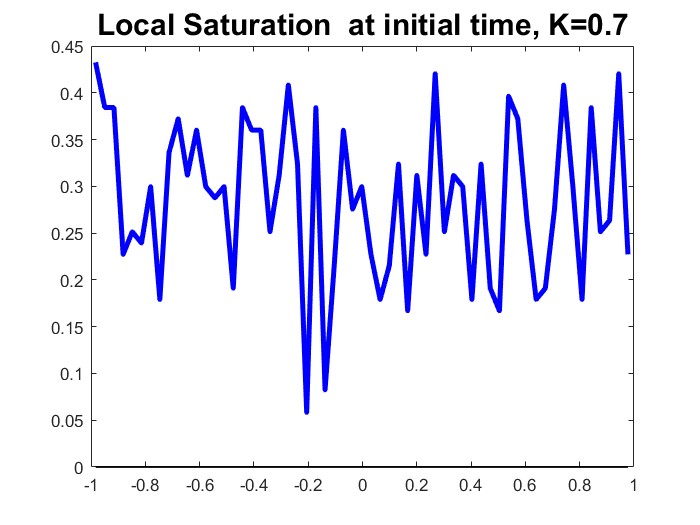}\!\!\!\!\!\!\includegraphics[scale=0.25]{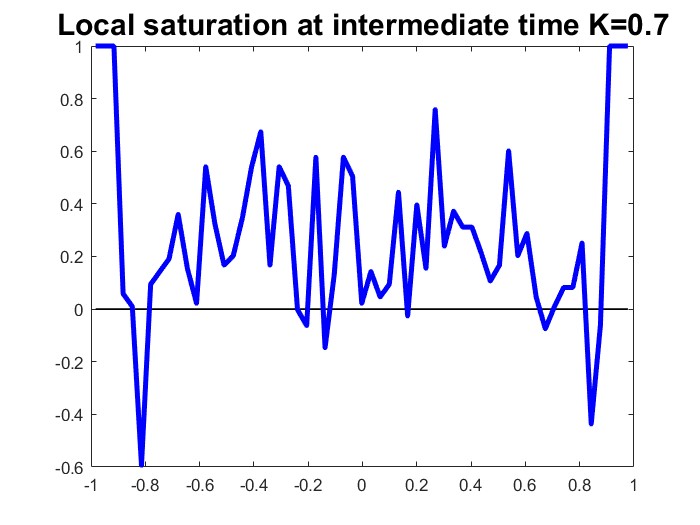}
\!\!\!\!\!\!\includegraphics[scale=0.25]{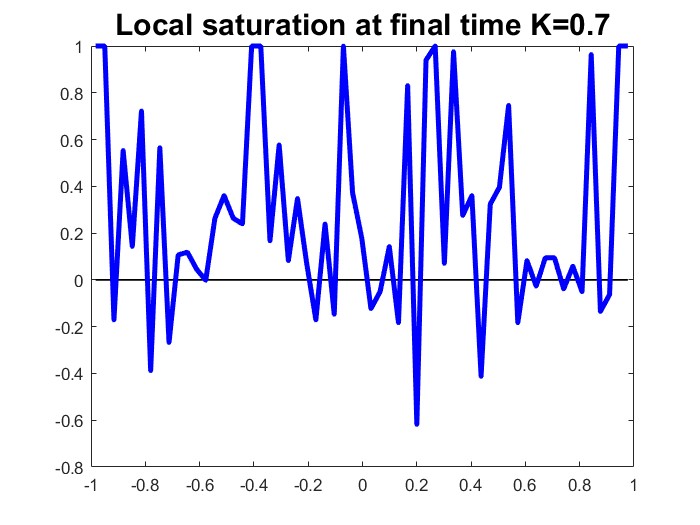}\\
\includegraphics[scale=0.17]{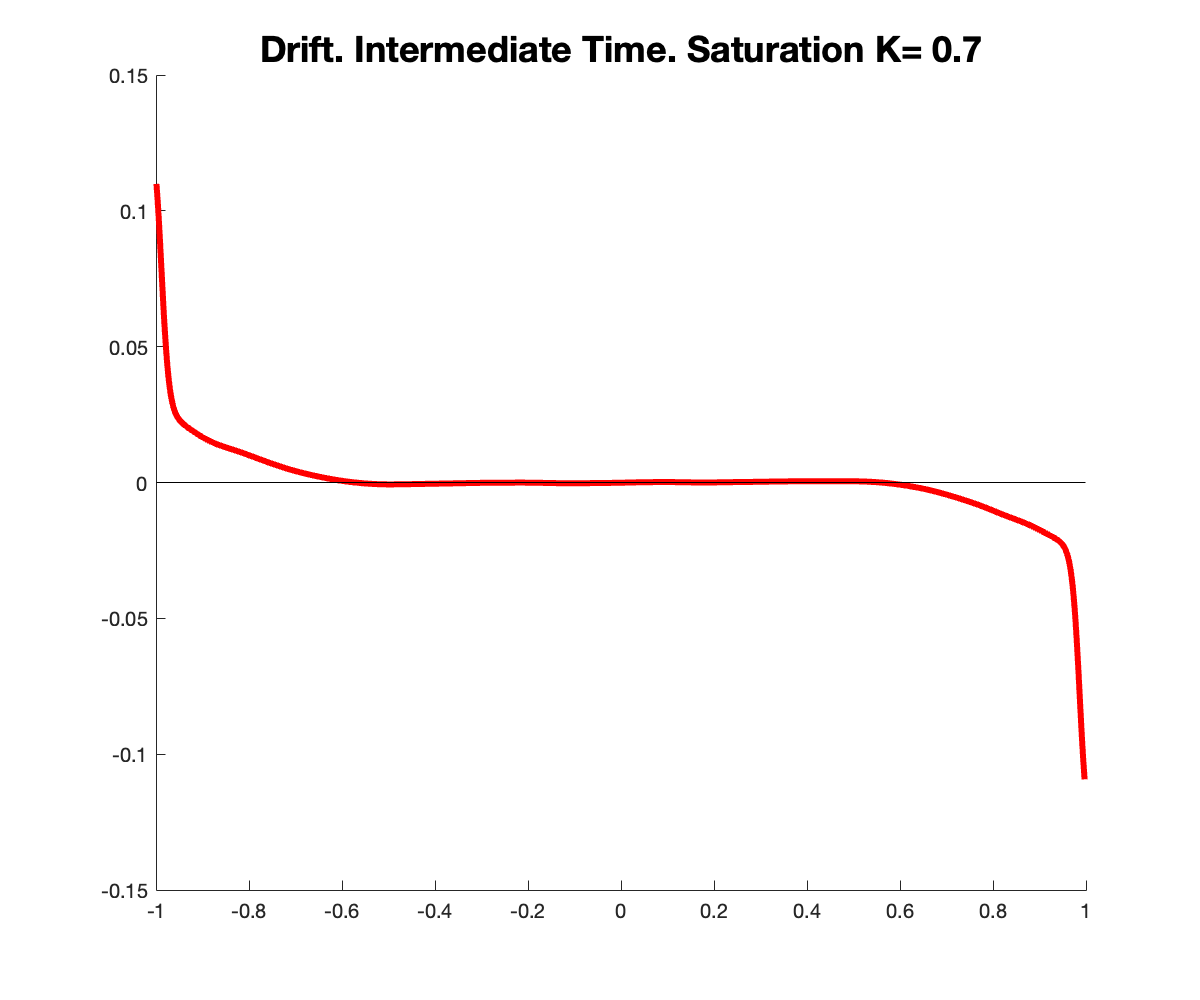}\includegraphics[scale=0.17]{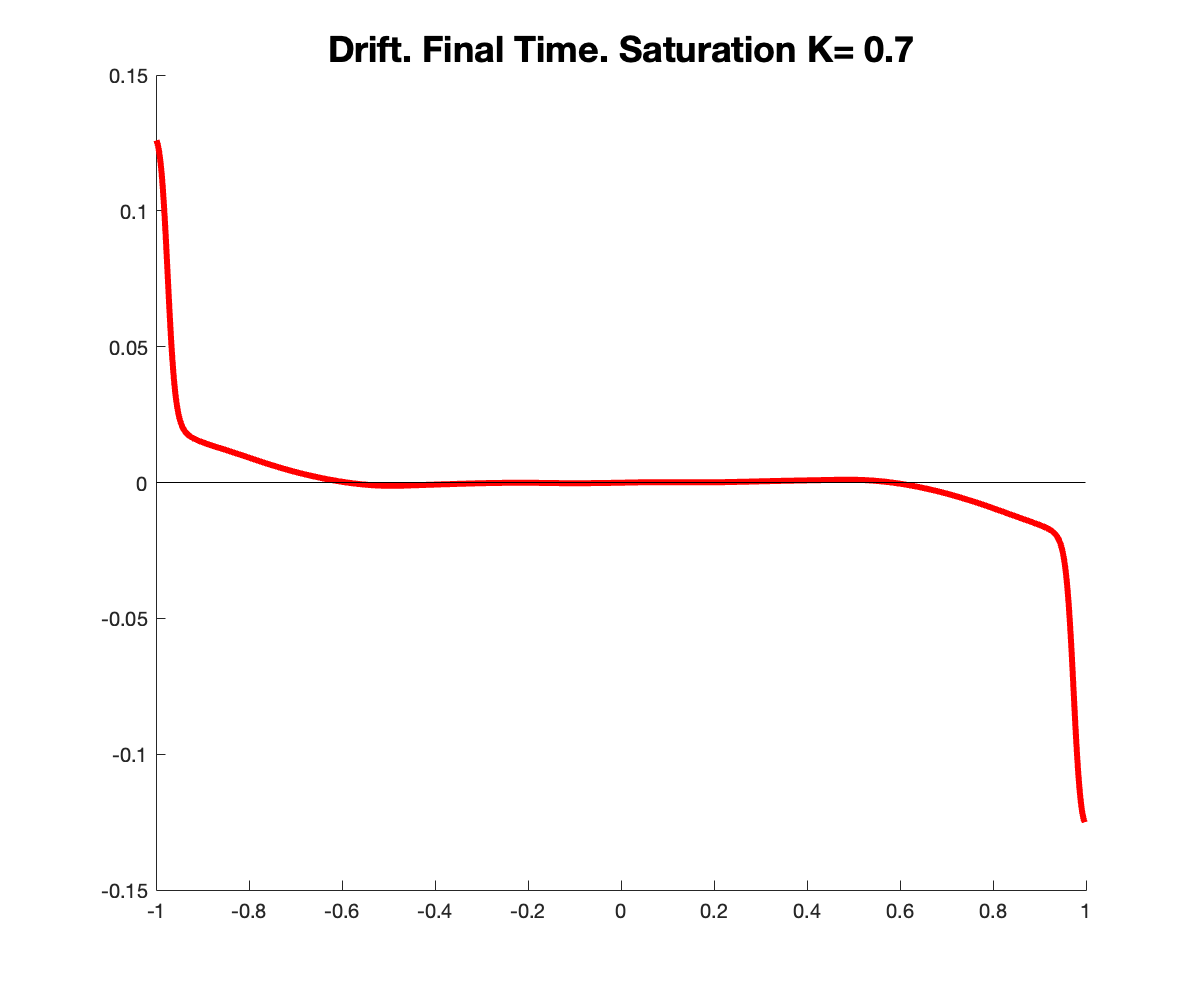}

\includegraphics[scale=0.28]{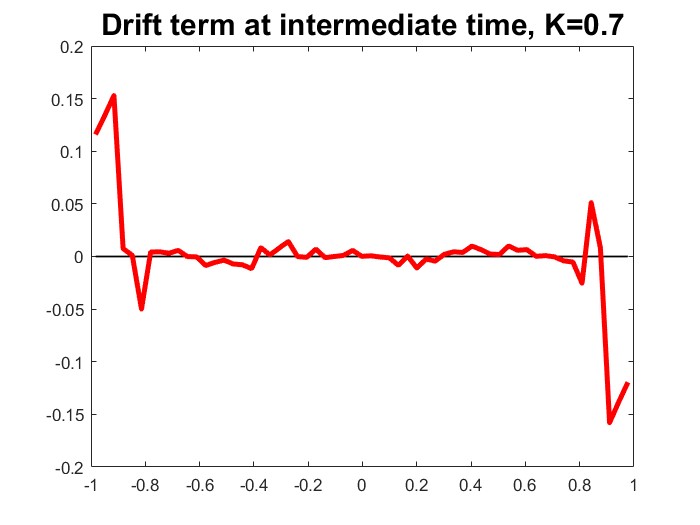}\!\!\includegraphics[scale=0.28]{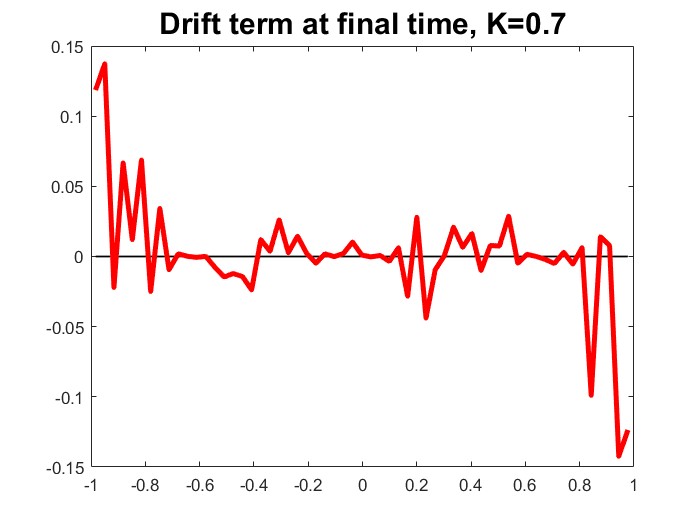}
\captionof{figure}{Local saturation coefficients and drifts for  PDE model (\S 3.5.3) and SDEs (\S 4.2), respectively  $K=0.7$. } 
\label{fig:Carrillo__K07_Drift}

\end{minipage}


    \begin{minipage}{\textwidth}
\centering
\includegraphics[scale=0.24]{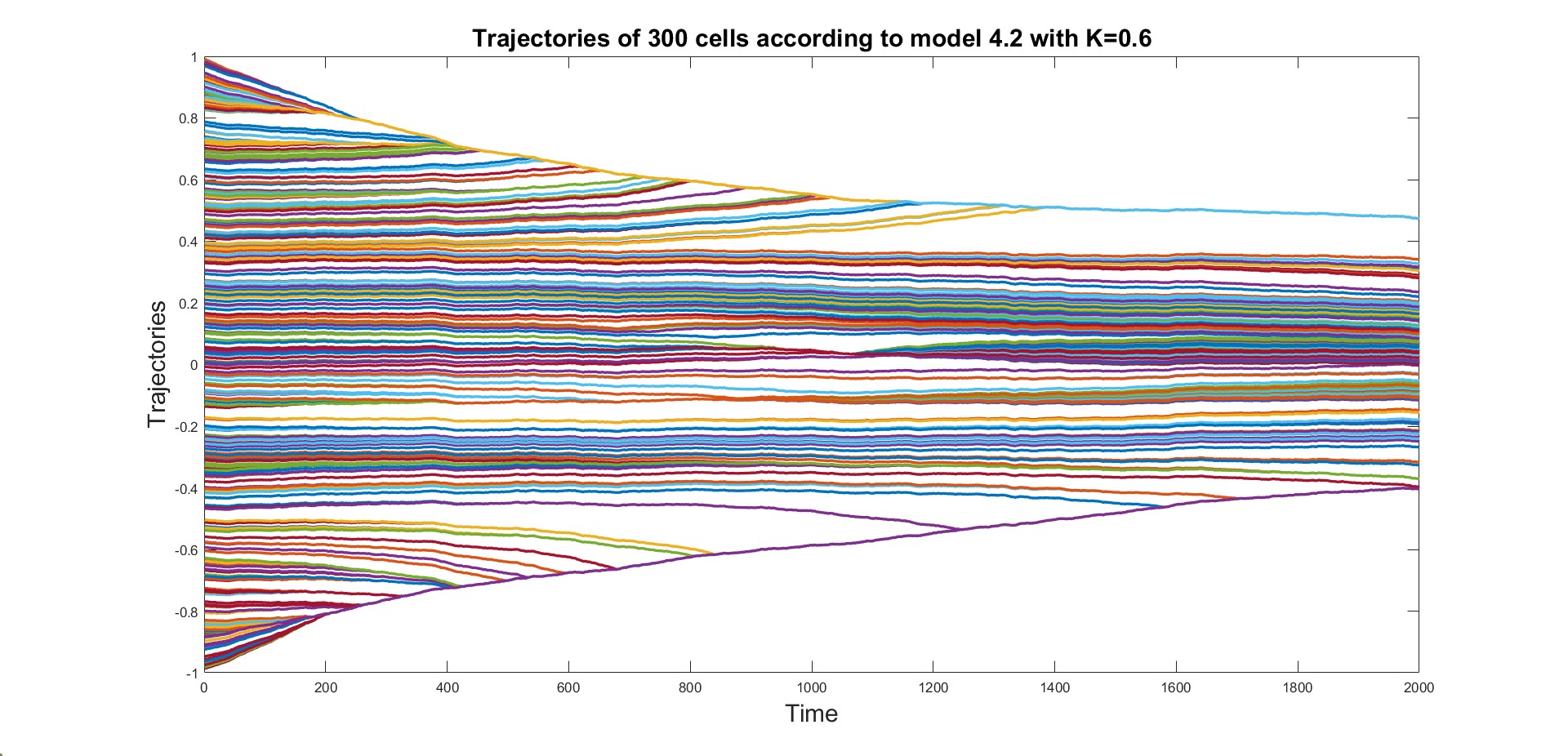}\\
\includegraphics[scale=0.18]{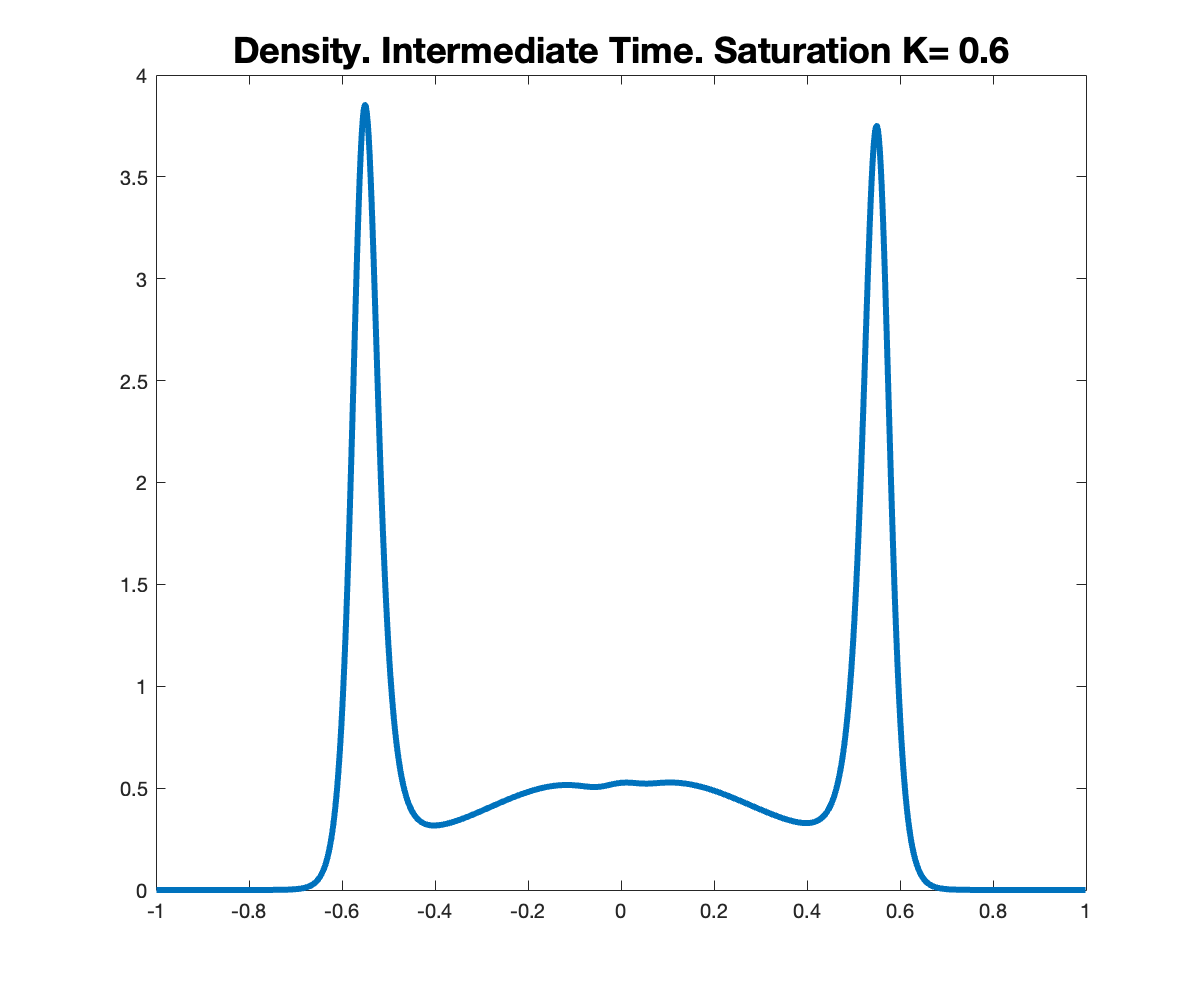}\!\!\includegraphics[scale=0.18]{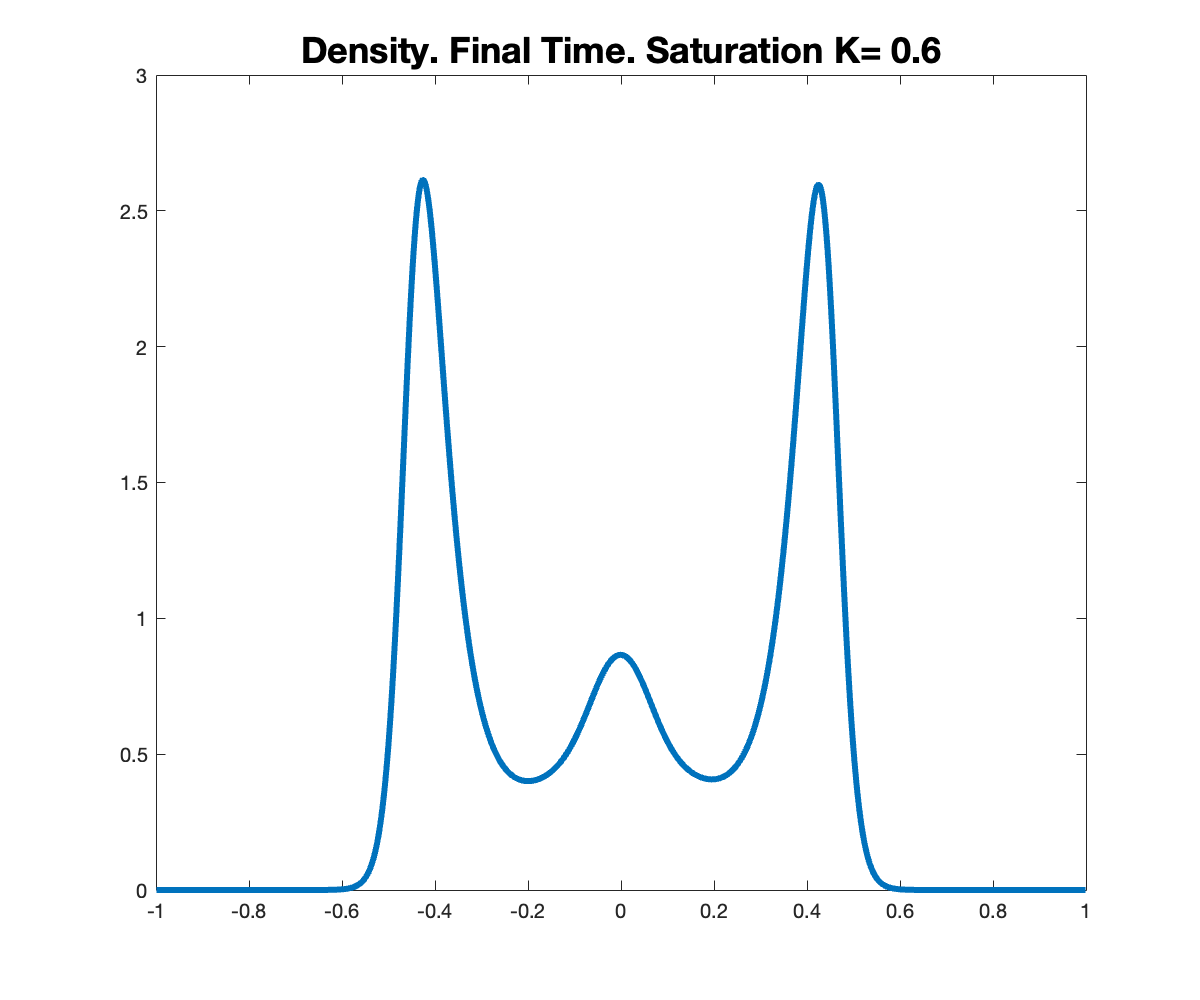}\\
\includegraphics[scale=0.3]{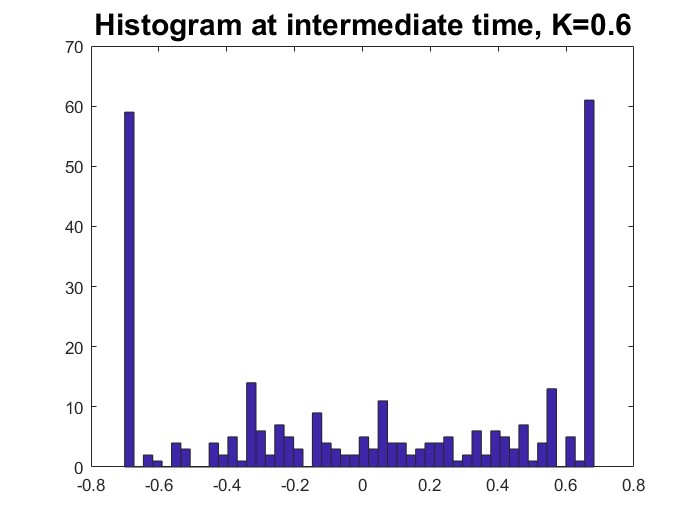}\!\!\includegraphics[scale=0.3]{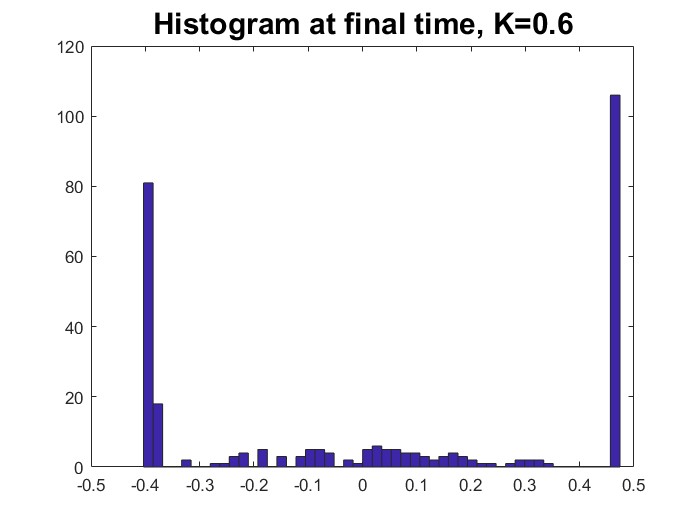}
\captionof{figure}{Nonlocal saturation model, $K=0.6$. First row trajectories with $N=300$ cells, second row densities of the PDE model (\S 3.5.3) and third one histograms of SDEs (\S 4.2). Both at intermediate time on the left and at final time on the right.} 
\label{fig:NL_K06_histogramas}
\end{minipage}

 \begin{minipage}{\textwidth}
\centering
\includegraphics[scale=0.15]{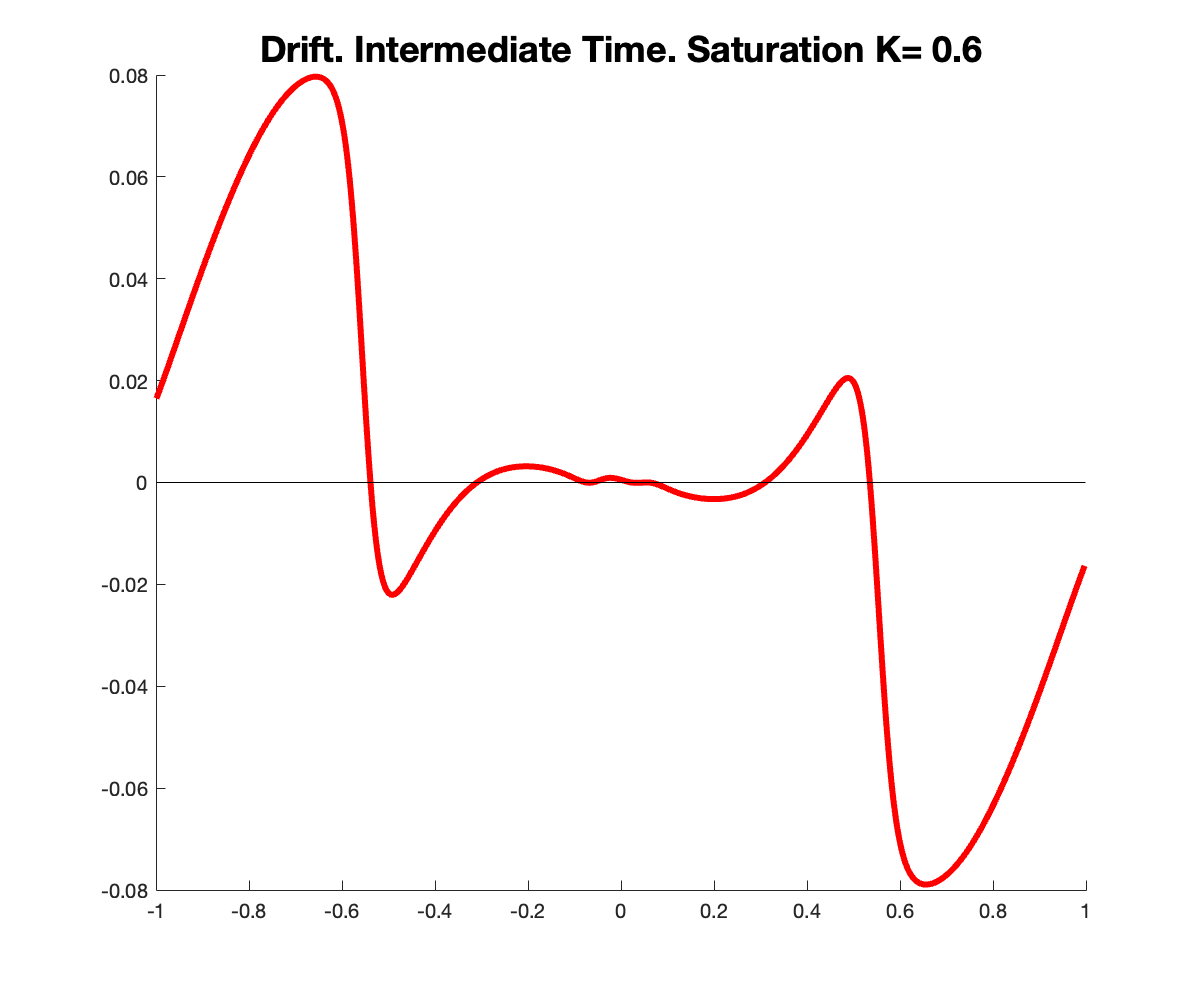}\!\!\!\!\includegraphics[scale=0.15]{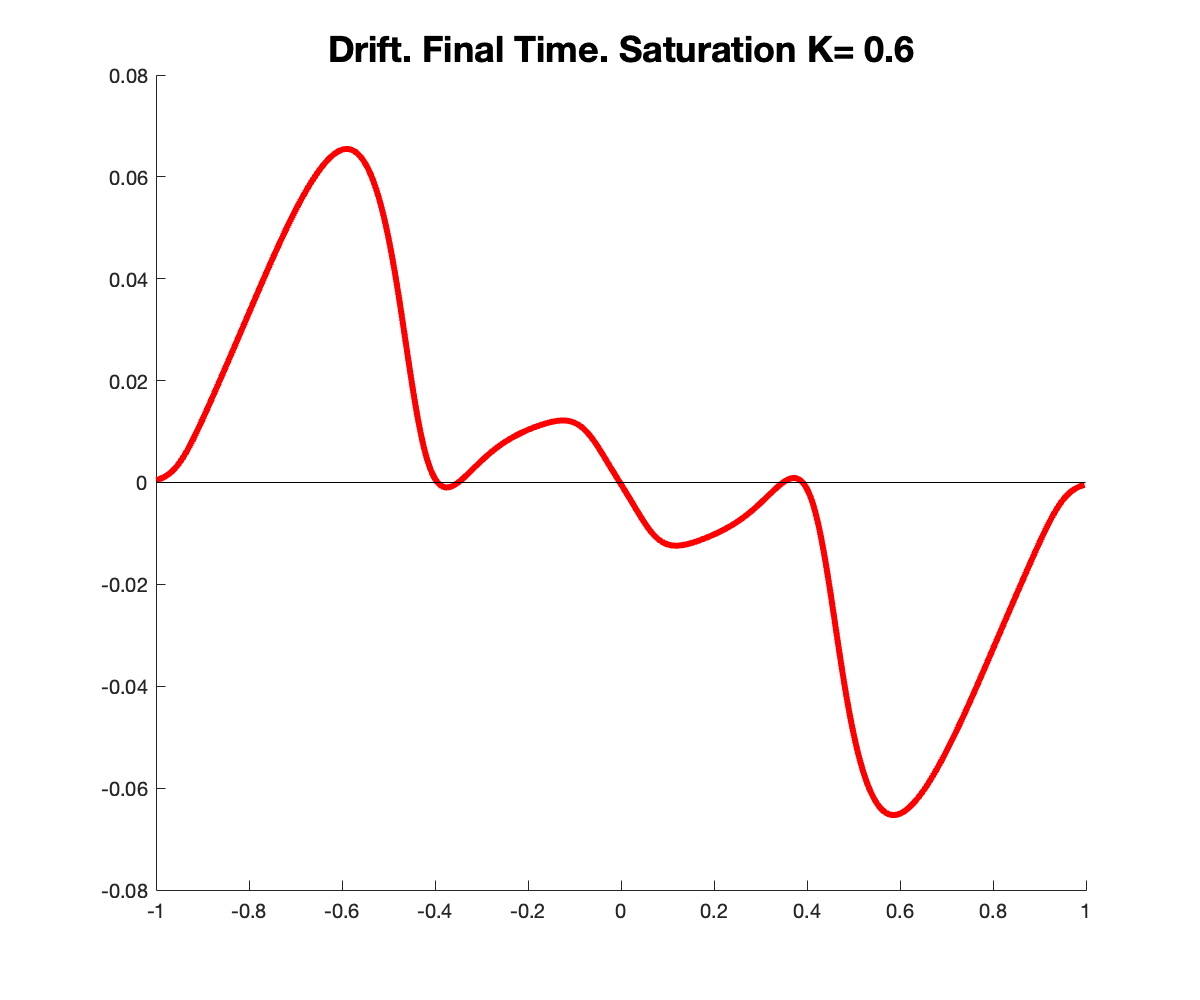}\\

\includegraphics[scale=0.26]{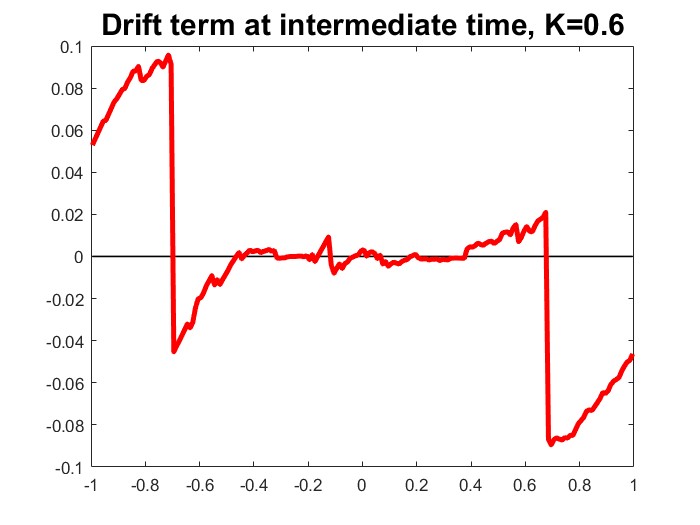}\!\!\includegraphics[scale=0.26]{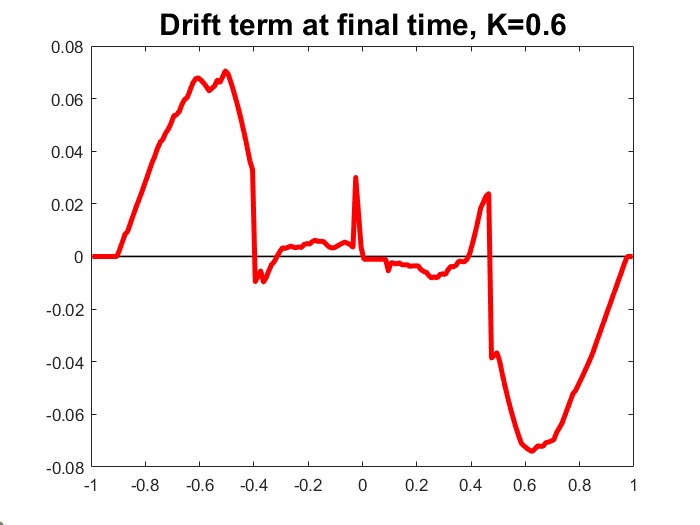}
\includegraphics[scale=0.15]{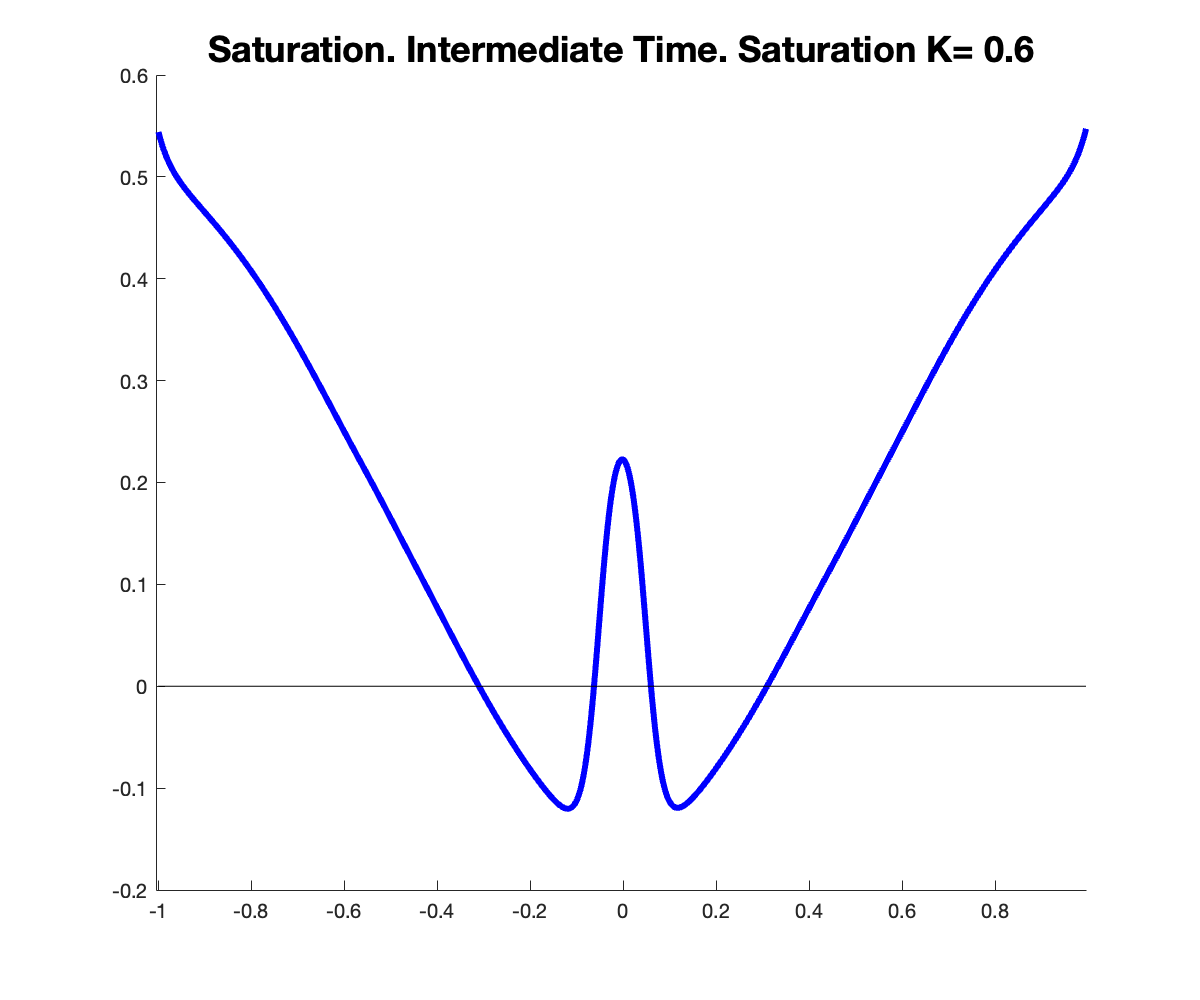}\!\!\!\!\includegraphics[scale=0.15]{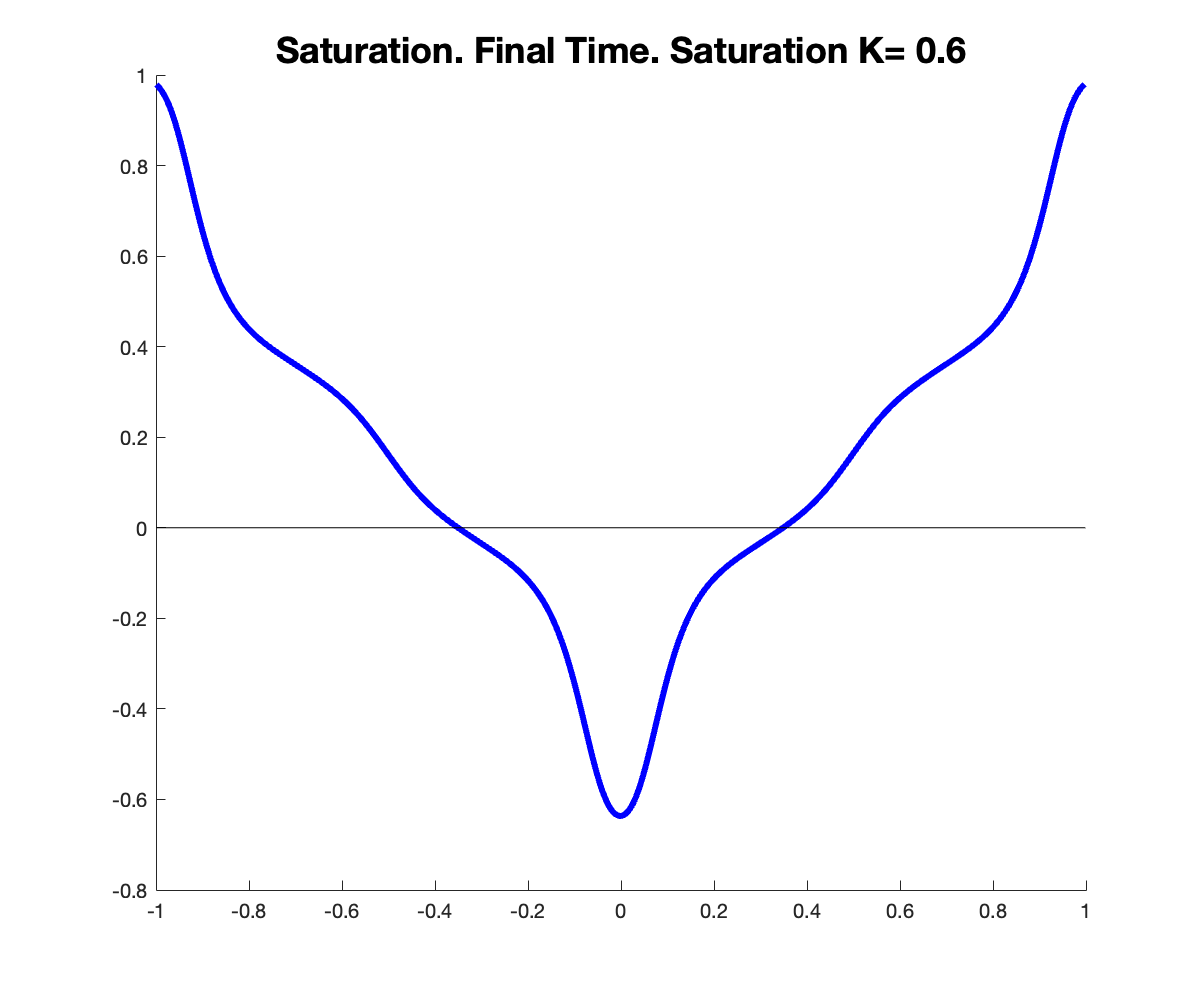}\\

\includegraphics[scale=0.26]{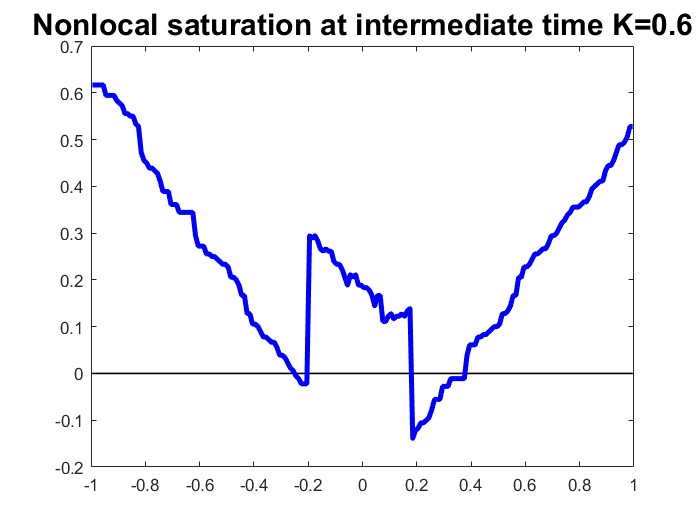}\includegraphics[scale=0.26]{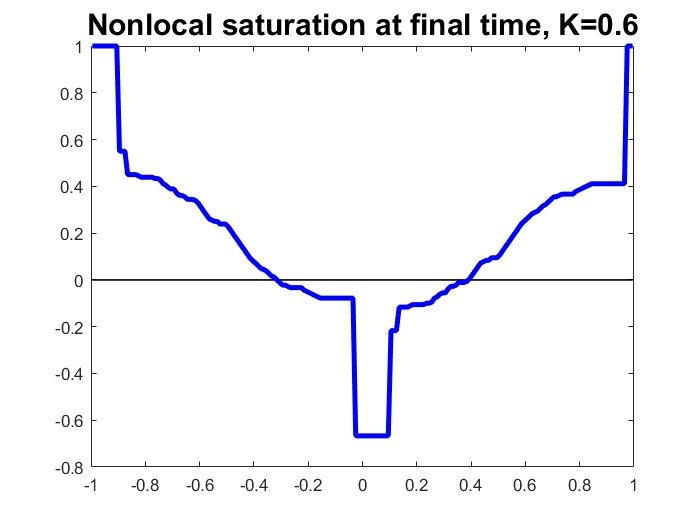}
\captionof{figure}{Drifts and saturations for nonlocal saturation model, $K=0.6$. Drifts in first and second rows for  PDE model (\S 3.5.3) and SDEs (\S 4.2), respectively. Saturations in third and fourth rows for  PDE model (\S 3.5.3) and SDEs (\S 4.2), respectively. All of them at intermediate time on the left and at final time on the right.} 
\label{fig:NL_K06_driftsatu}
\end{minipage}

\newpage
    \begin{minipage}{\textwidth}
\centering
\includegraphics[scale=0.24]{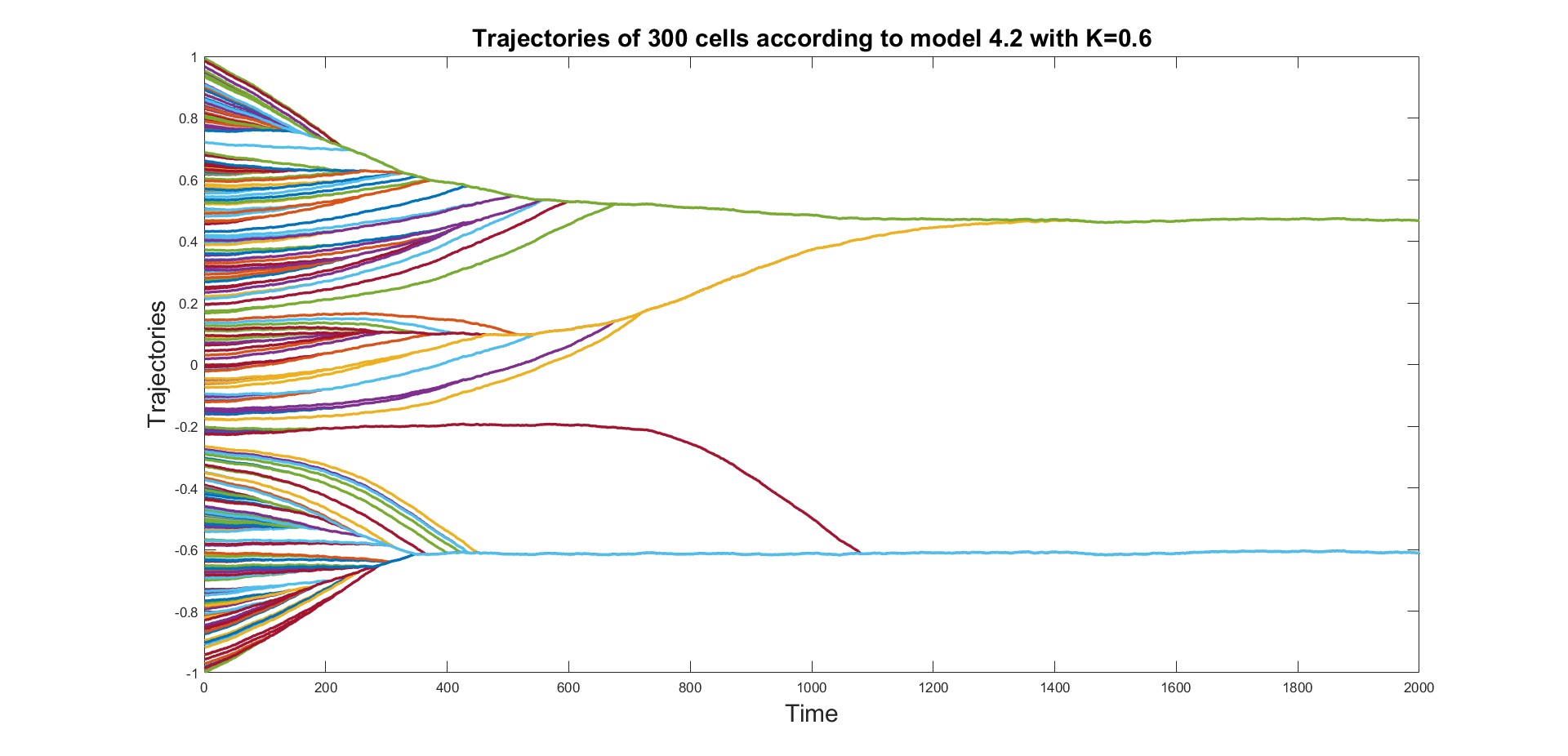}\\

\includegraphics[scale=0.18]{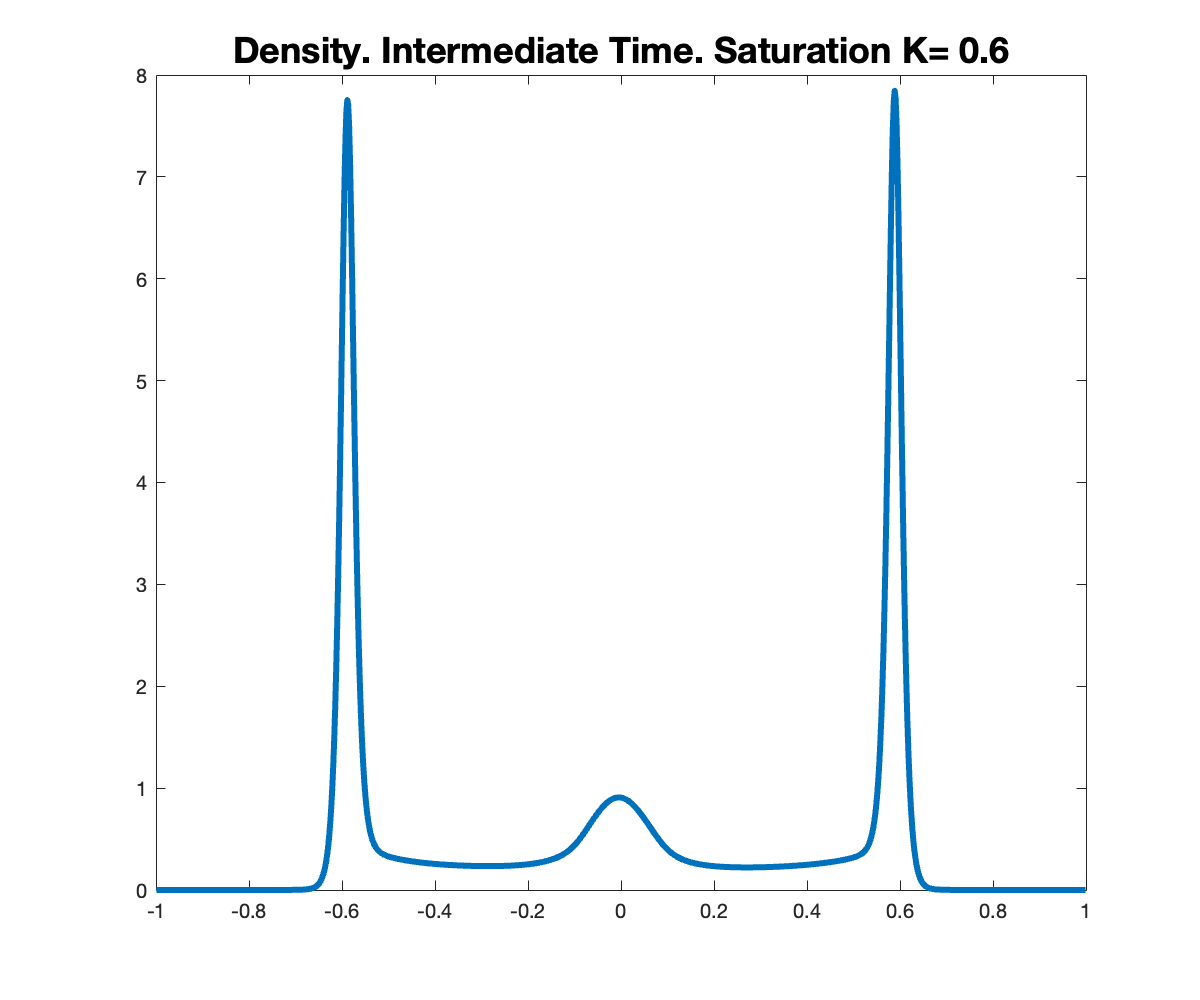}\!\!\includegraphics[scale=0.18]{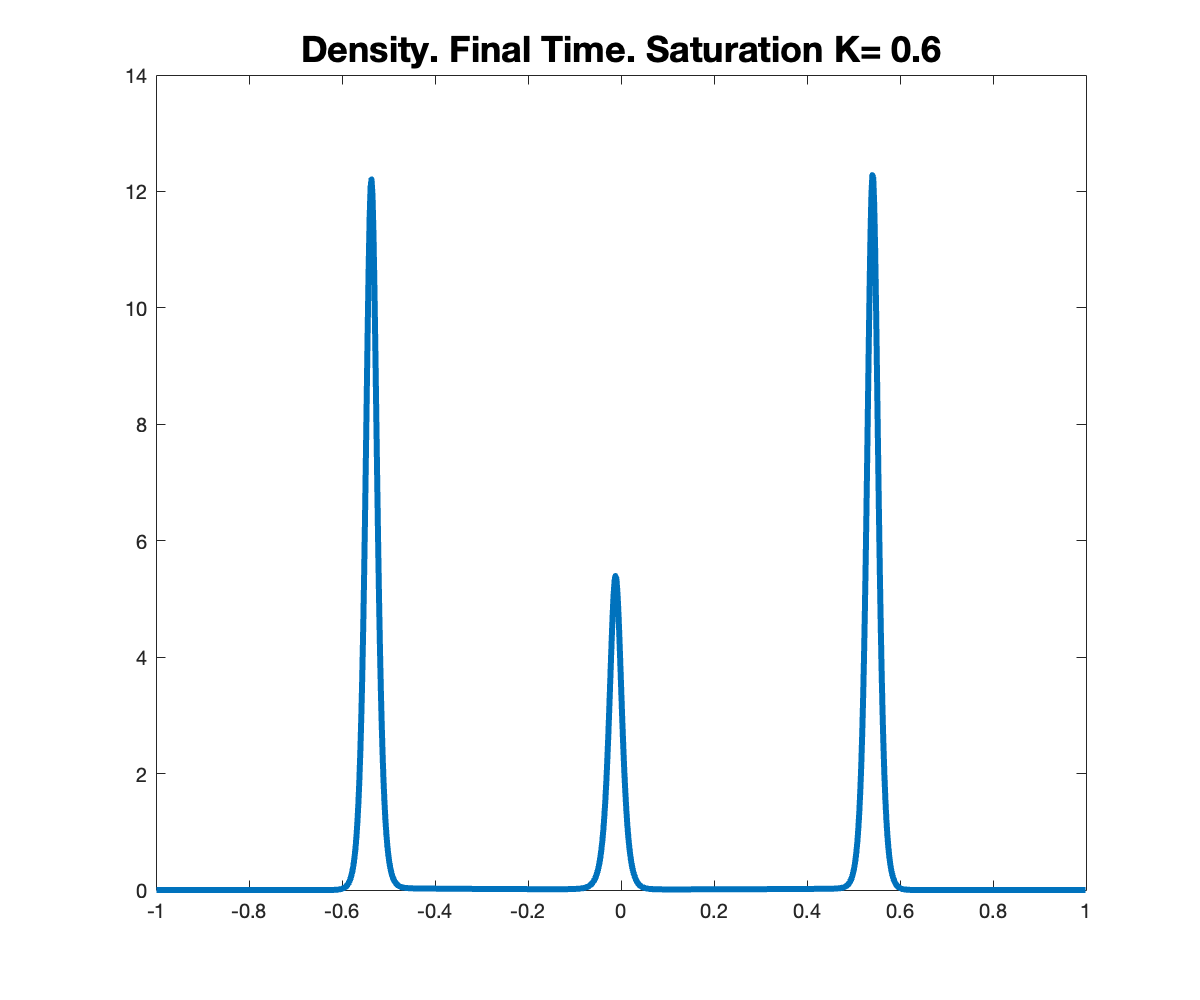}\\
\includegraphics[scale=0.3]{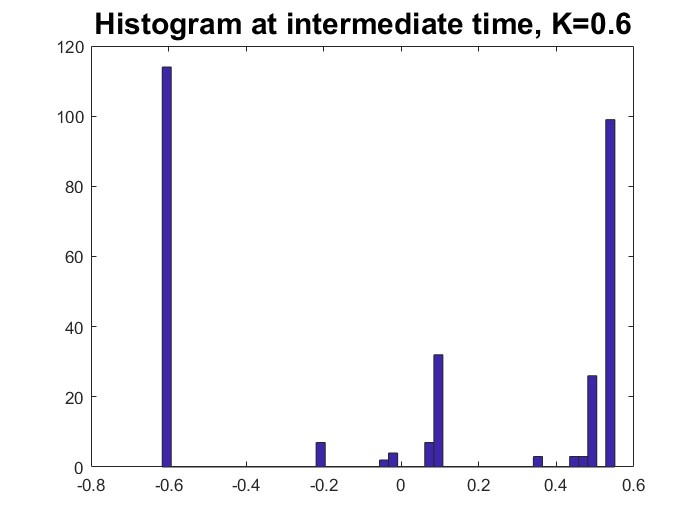}\!\!\includegraphics[scale=0.3]{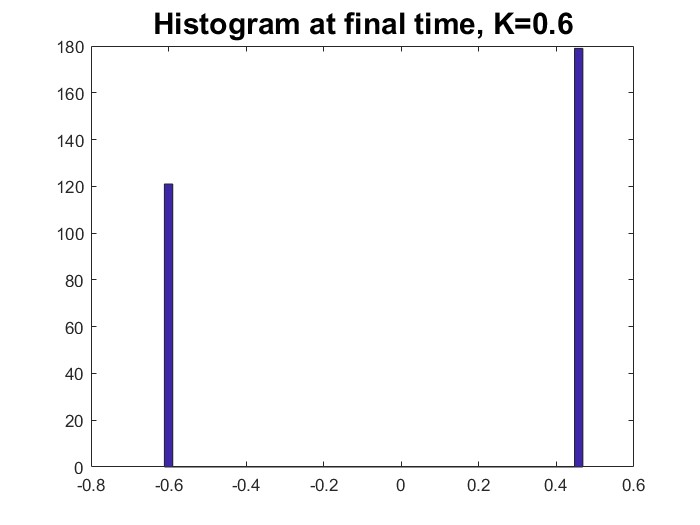}
\captionof{figure}{Nonlocal saturation model, $K=0.6$, linear weight in saturation. First row trajectories with $N=300$ cells, second row densities of the PDE model (\S 3.5.3) and third one histograms of SDEs (\S 4.2). Both at intermediate time on the left and at final time on the right.} 
\label{fig:NL_K06pesolineal_histogramas} 
\end{minipage}

 \begin{minipage}{\textwidth}
\centering
\includegraphics[scale=0.15]{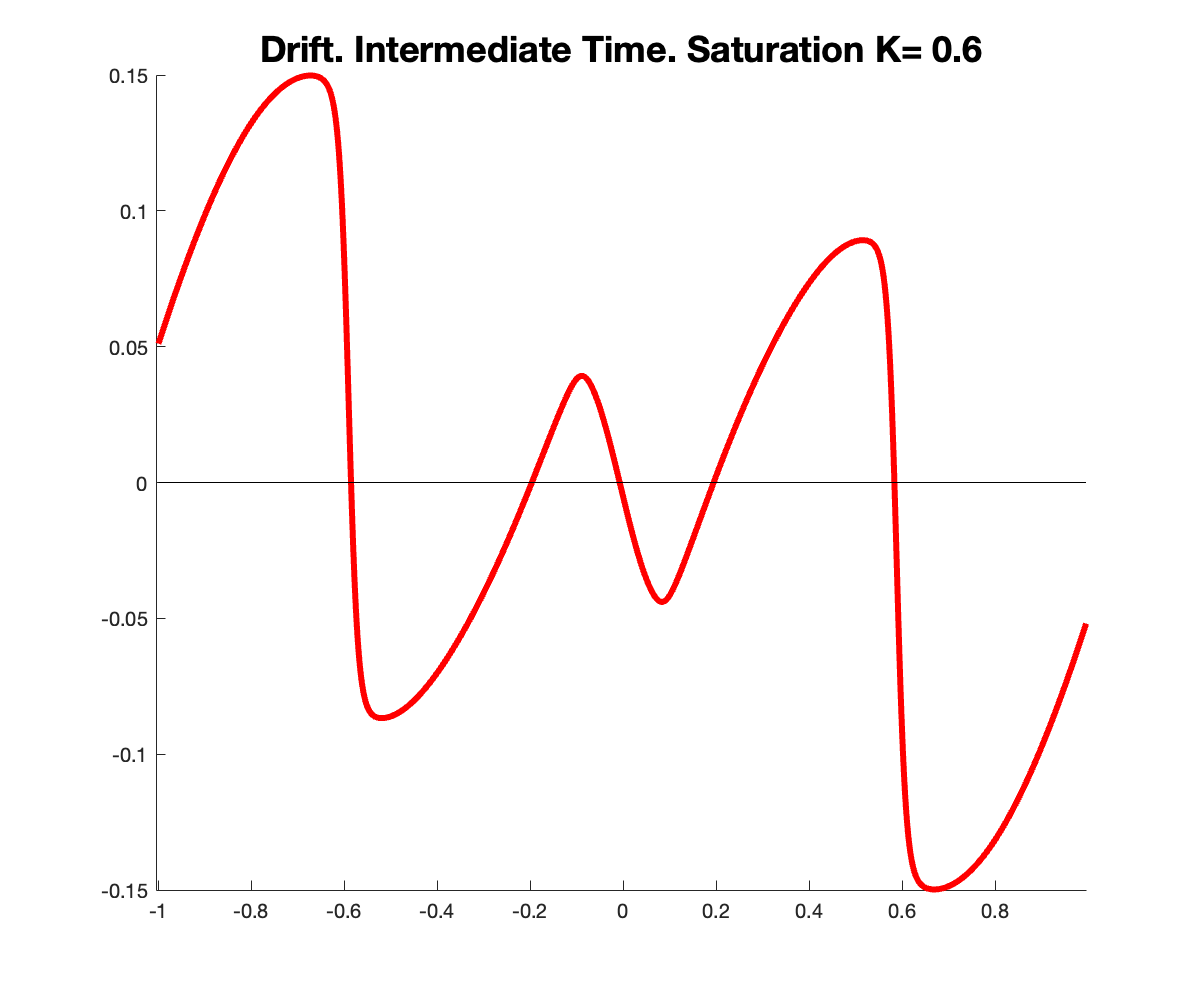}\!\!\!\!\includegraphics[scale=0.15]{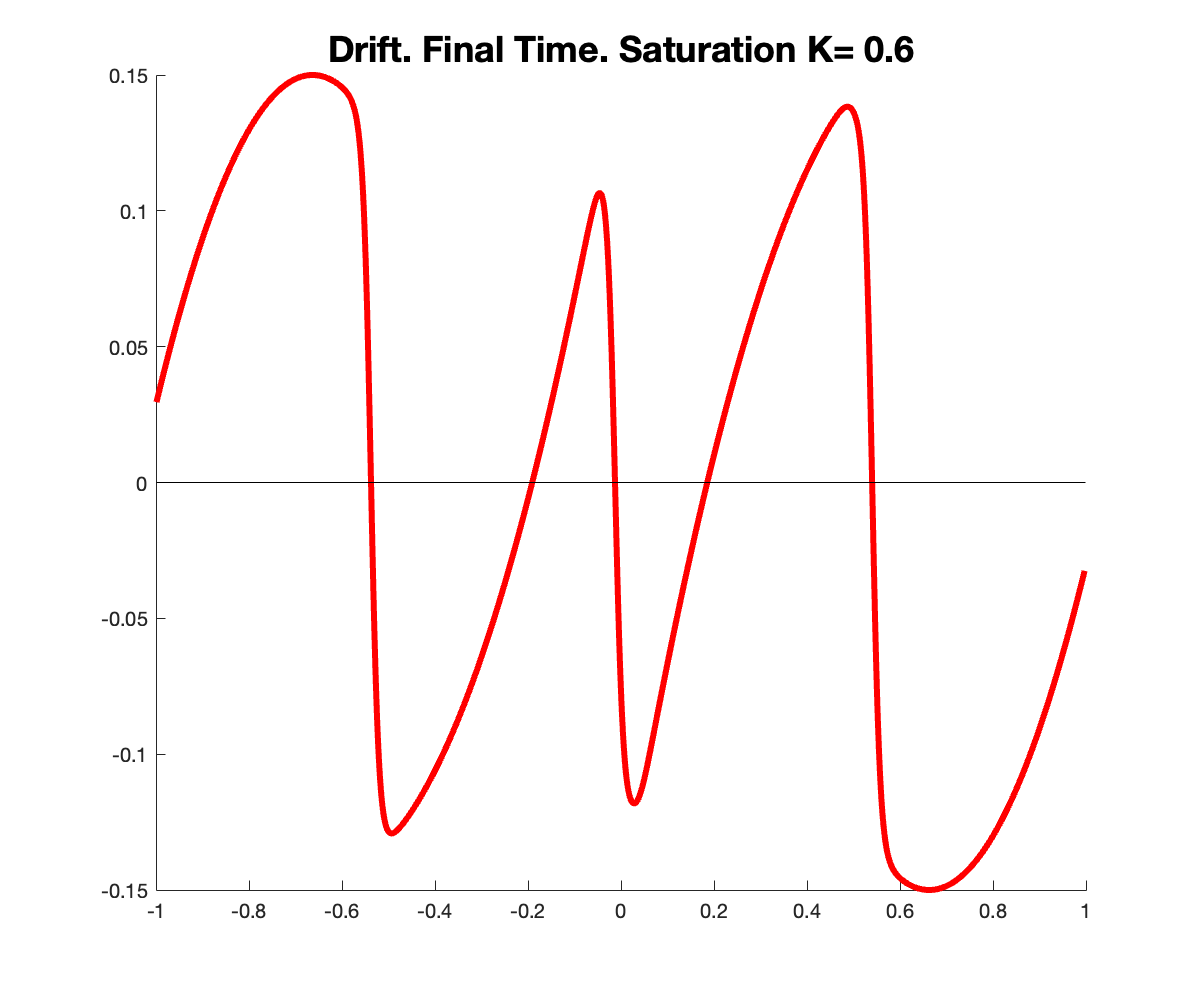}\\
\includegraphics[scale=0.26]{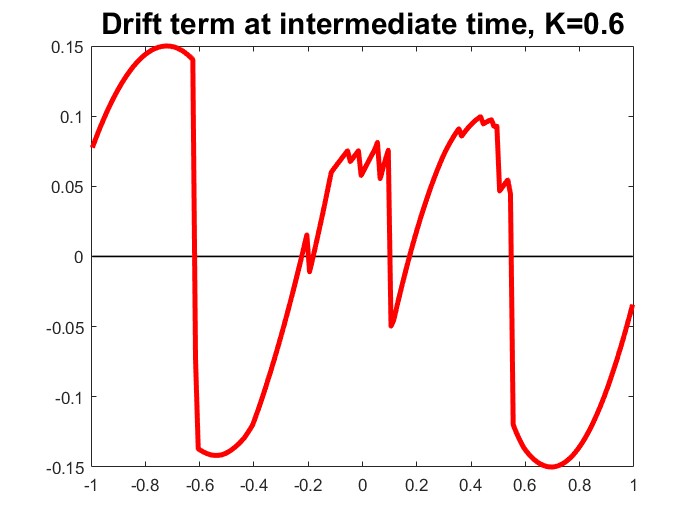}\!\!\includegraphics[scale=0.26]{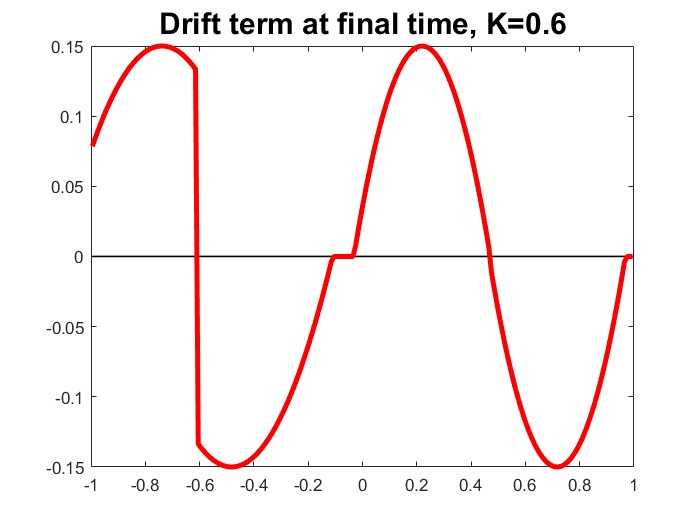}\\
\includegraphics[scale=0.15]{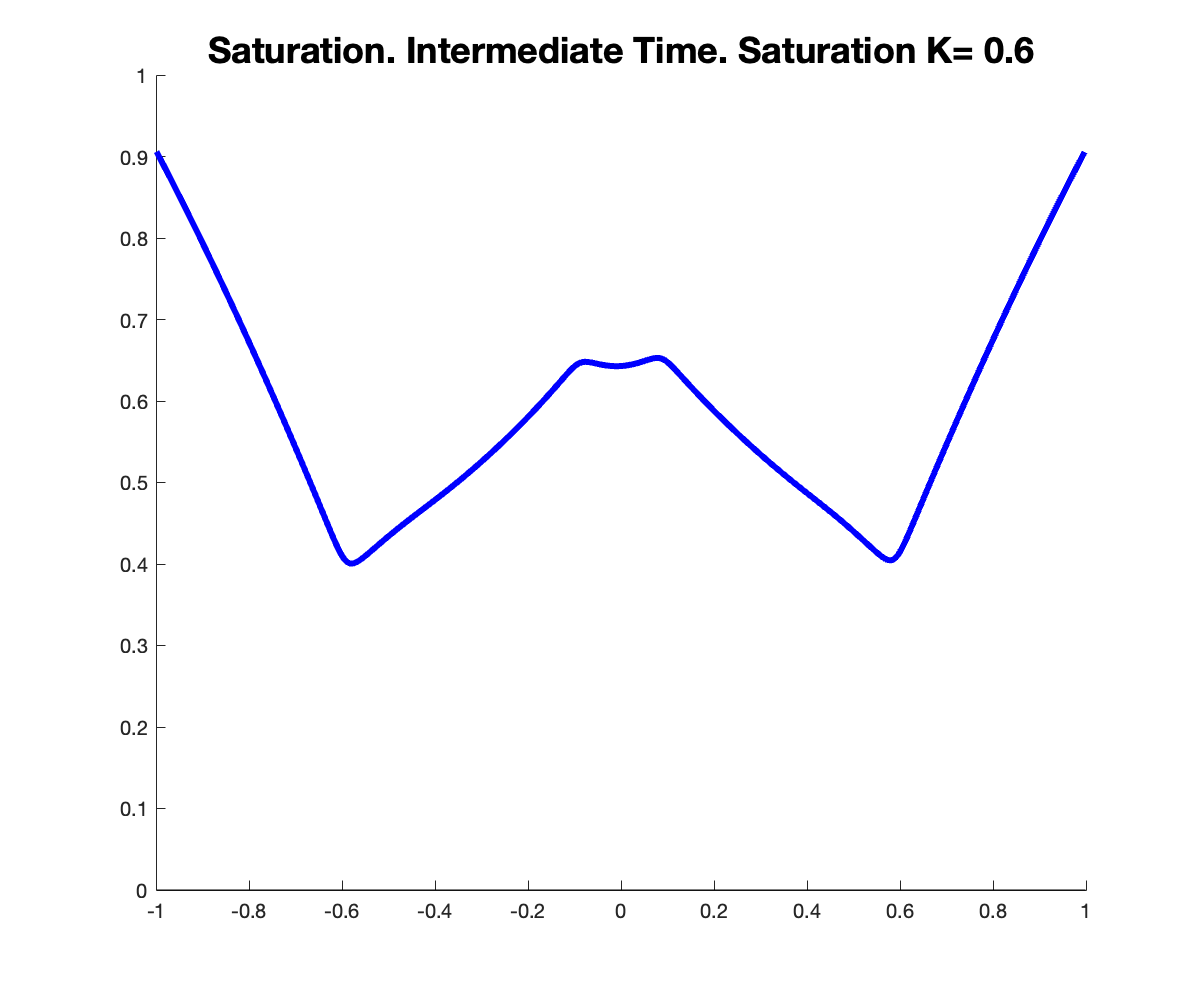}\!\!\!\!\includegraphics[scale=0.15]{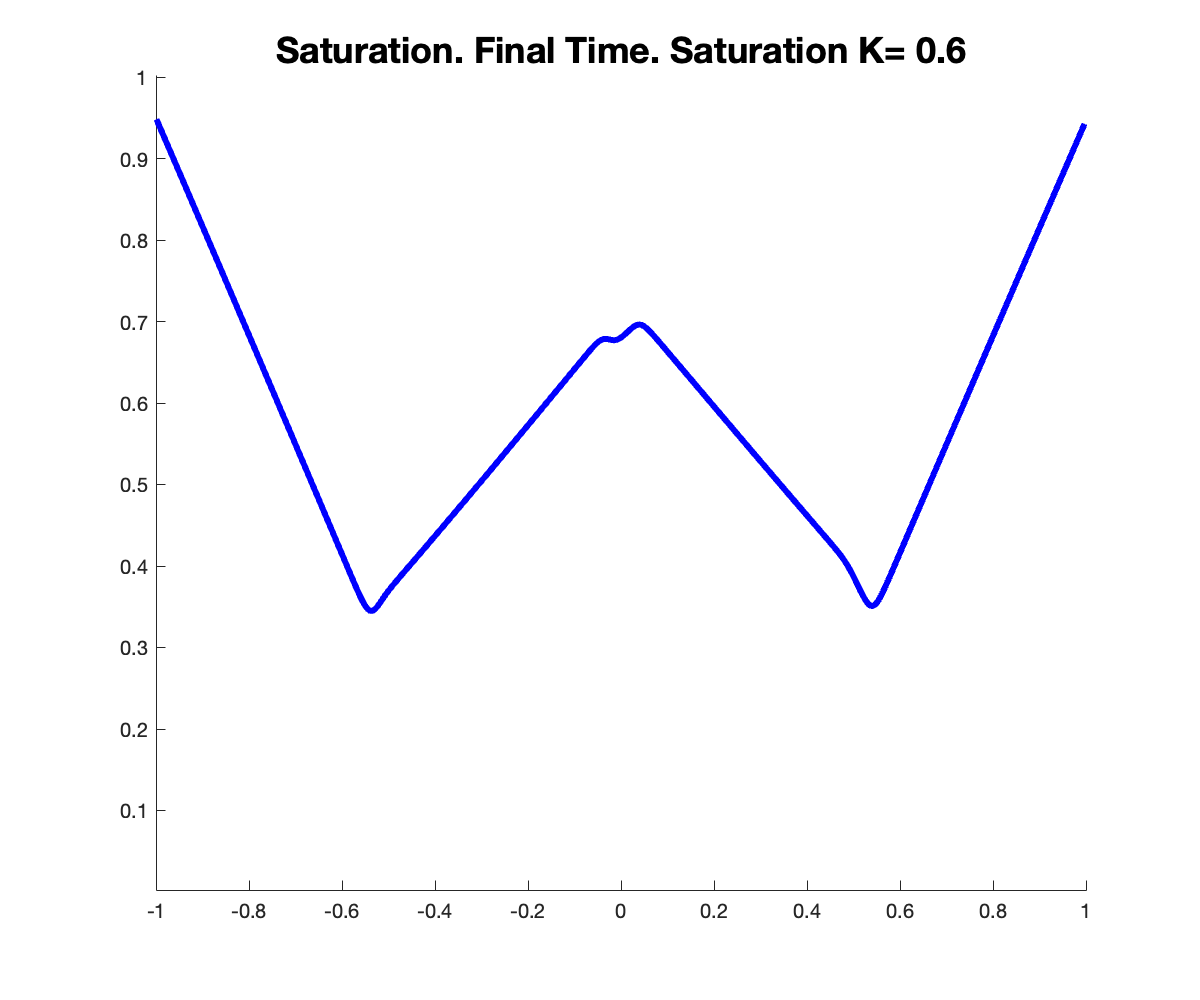}\\
\includegraphics[scale=0.26]{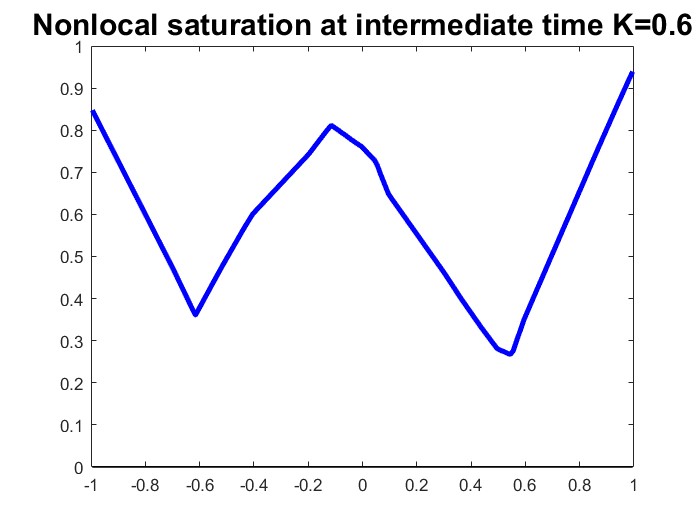}\includegraphics[scale=0.26]{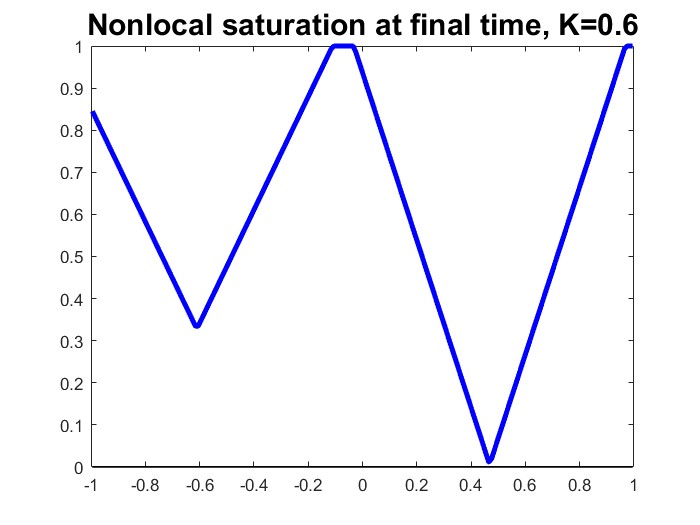}
\captionof{figure}{Drifts and saturations for nonlocal saturation model, $K=0.6$, linear weight in saturation. Drifts in first and second rows for  PDE model (\S 3.5.3) and SDEs (\S 4.2), respectively. Saturations in third and fourth rows for  PDE model (\S 3.5.3) and SDEs (\S 4.2), respectively. All of them at intermediate time on the left and at final time on the right.} 
\label{fig:NL_K06pesolineal_driftsatu}
\end{minipage}


   \begin{minipage}{\textwidth}
\centering
\includegraphics[scale=0.24]{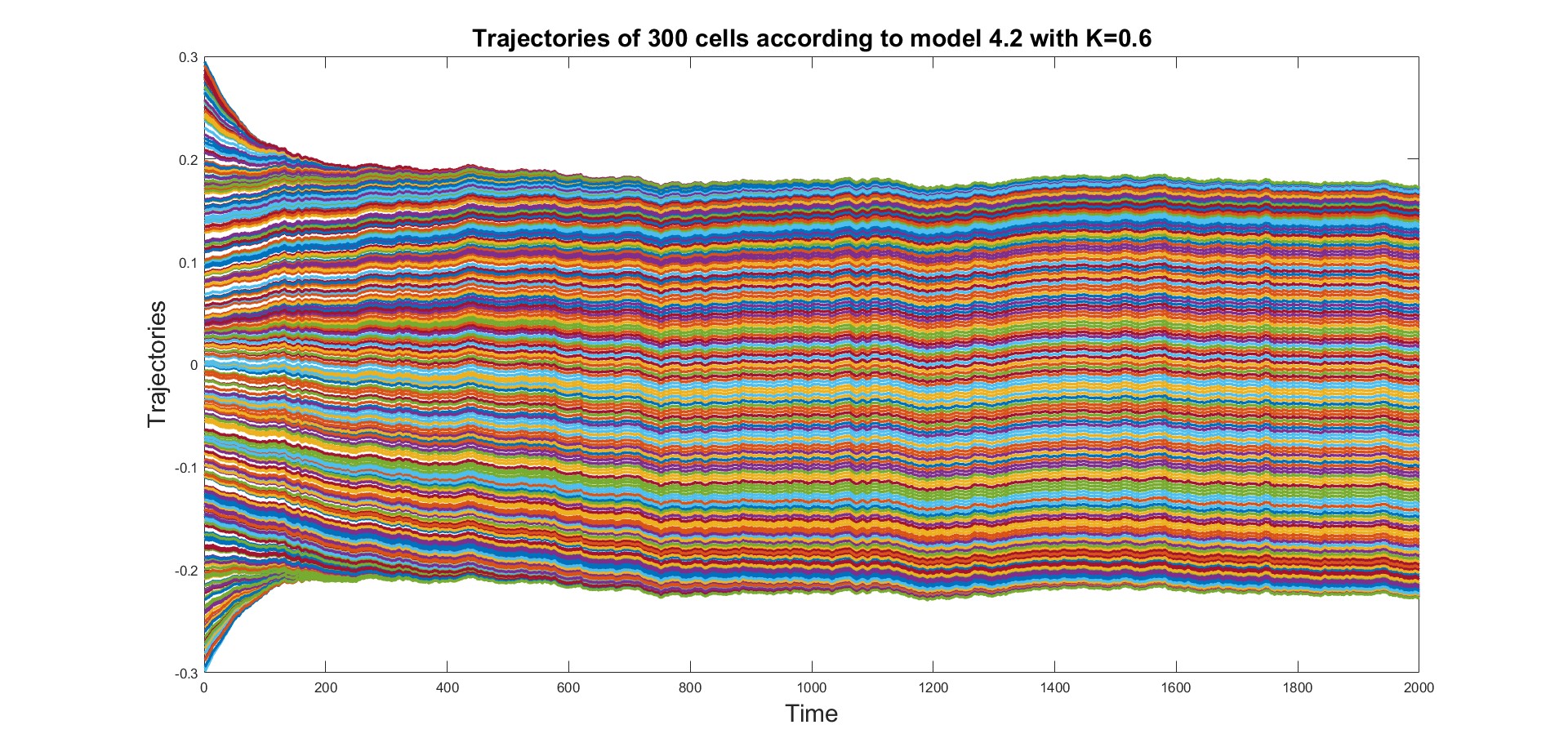}\\

\includegraphics[scale=0.18]{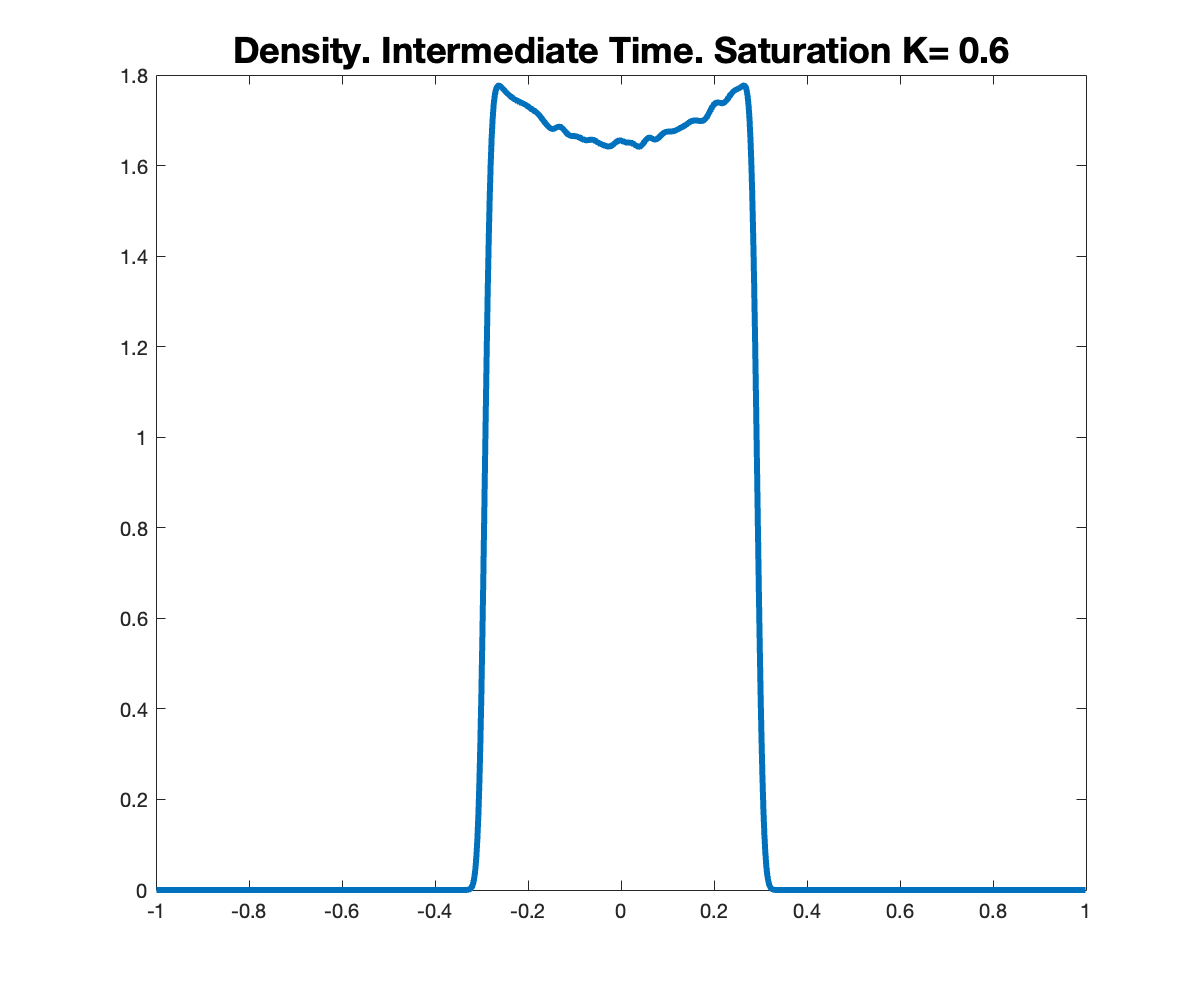}\!\!\includegraphics[scale=0.18]{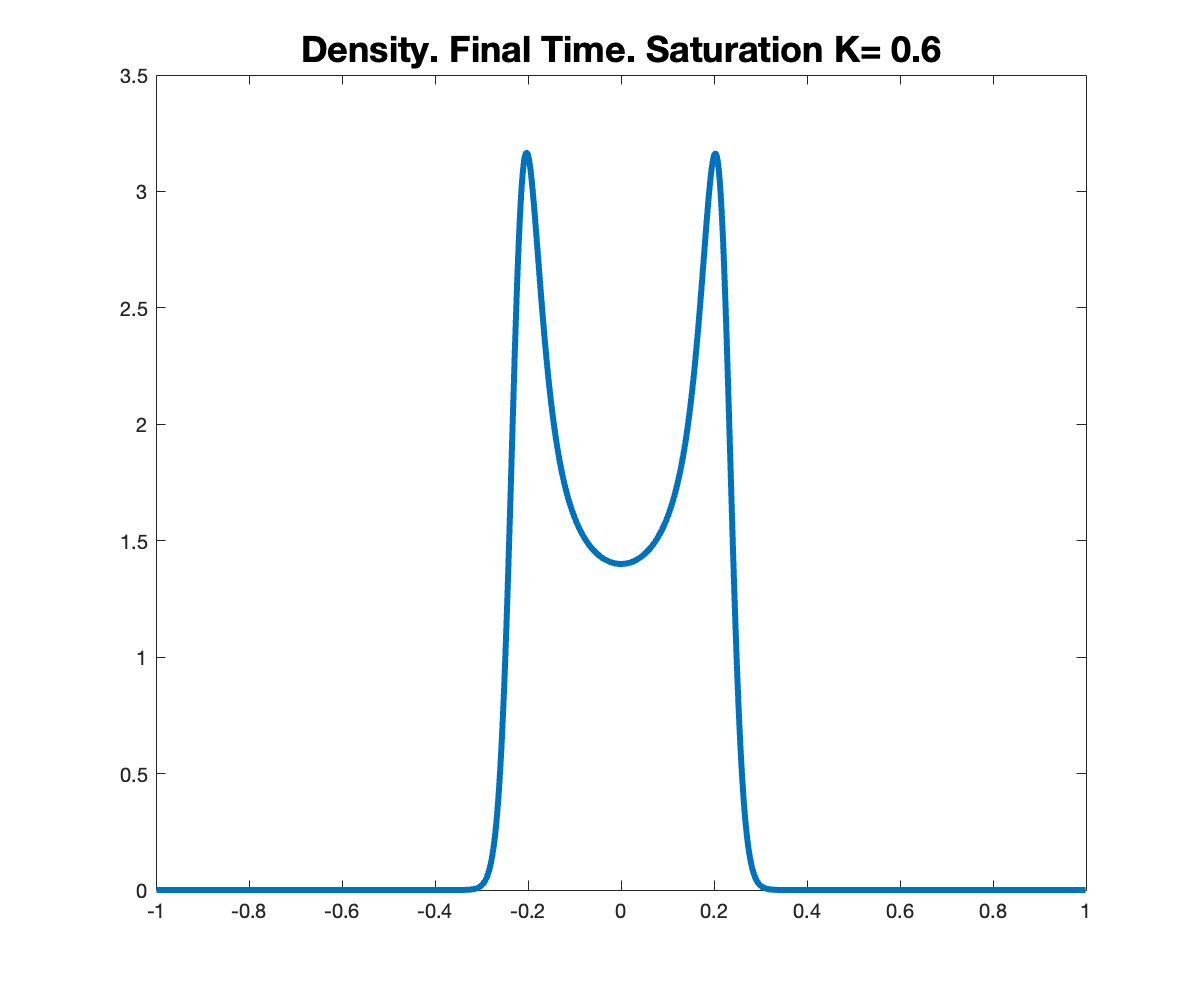}\\
\includegraphics[scale=0.3]{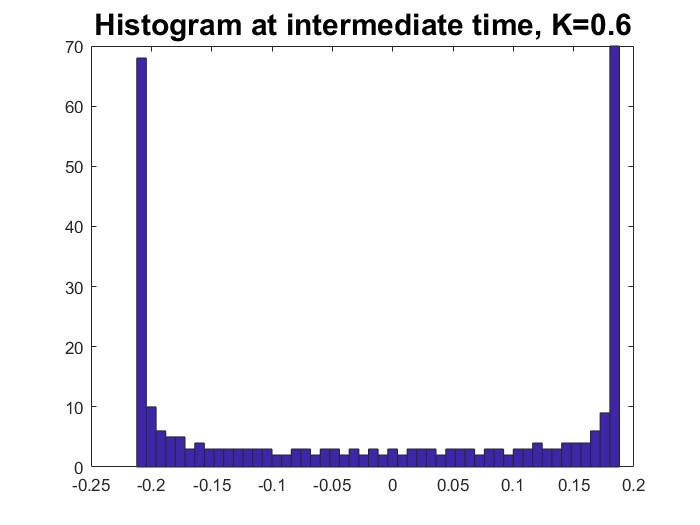}\!\!\includegraphics[scale=0.3]{DEFINITIVAS/Modelo3K06concentrado/mod3-hist-final-k06-con.jpg}
\captionof{figure}{Nonlocal saturation model, $K=0.6$, Concentrated initial data within $[-0.3,0.3]$. First row trajectories with $N=300$ cells, second row densities of the PDE model (\S 3.5.3) and third one histograms of SDEs (\S 4.2). Both at intermediate time on the left and at final time on the right.} 
\label{fig:NL_K06concentrado_histogramas}
\end{minipage}

 \begin{minipage}{\textwidth}
\centering
\includegraphics[scale=0.15]{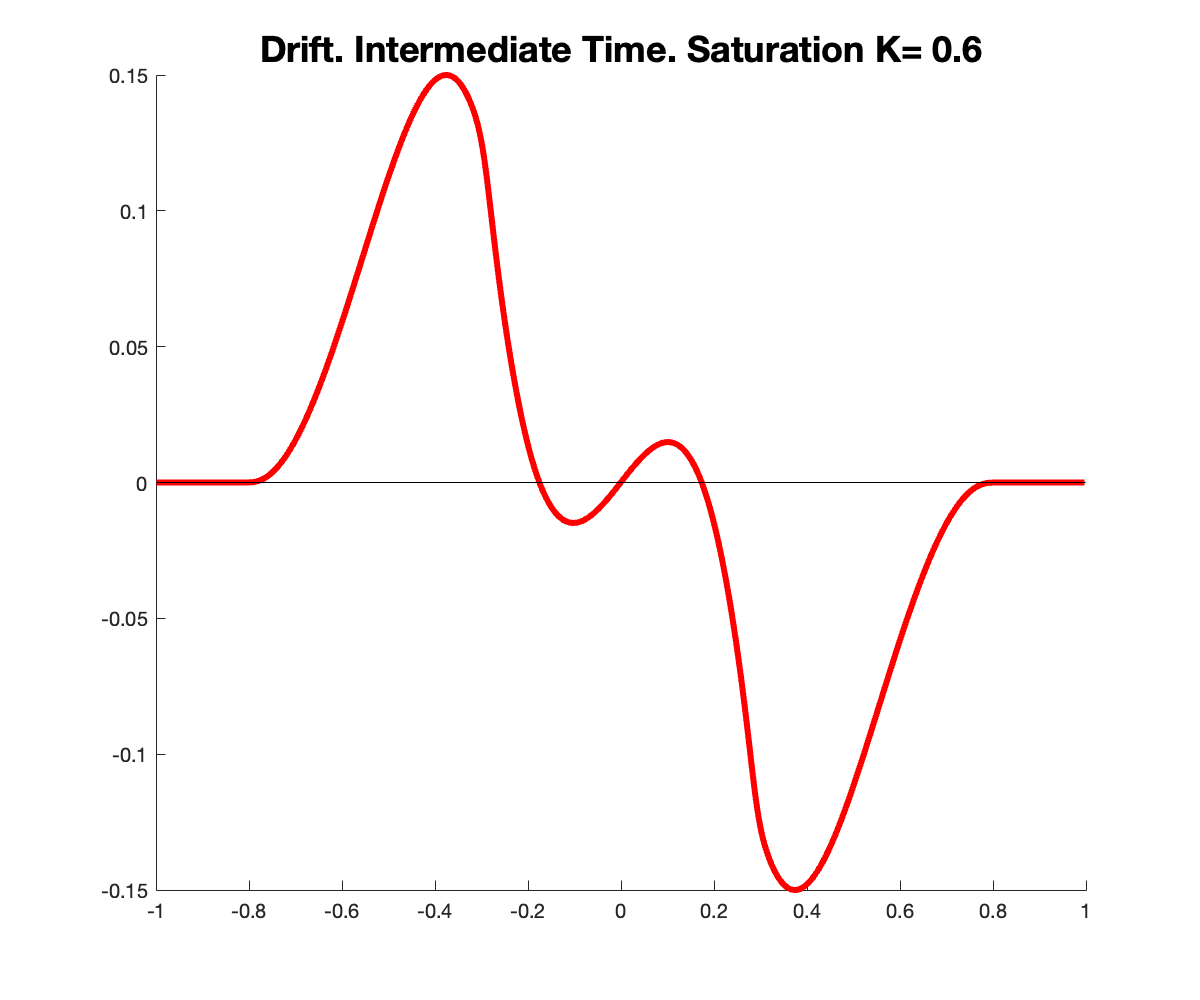}\!\!\!\!\includegraphics[scale=0.15]{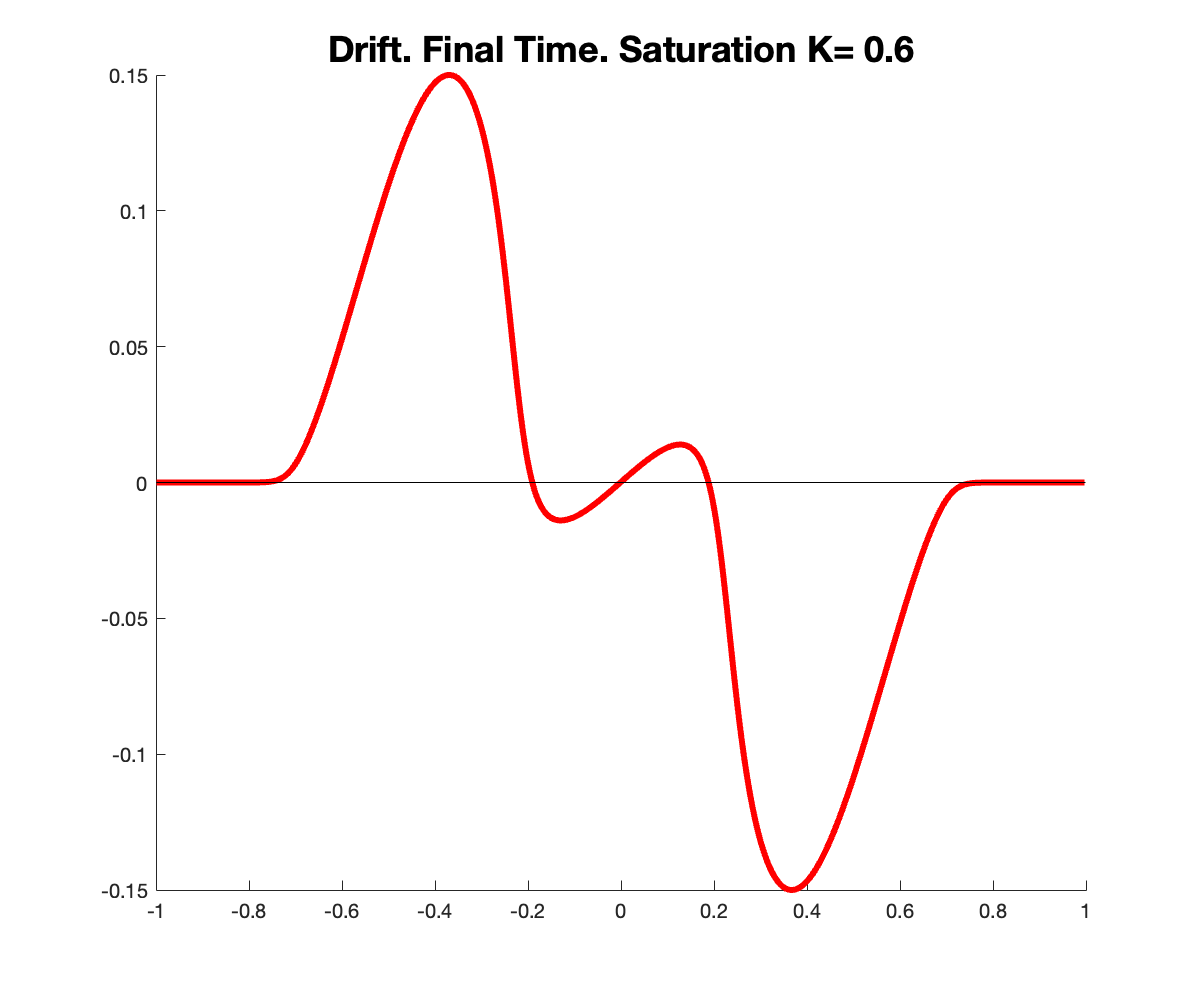}\\
\includegraphics[scale=0.26]{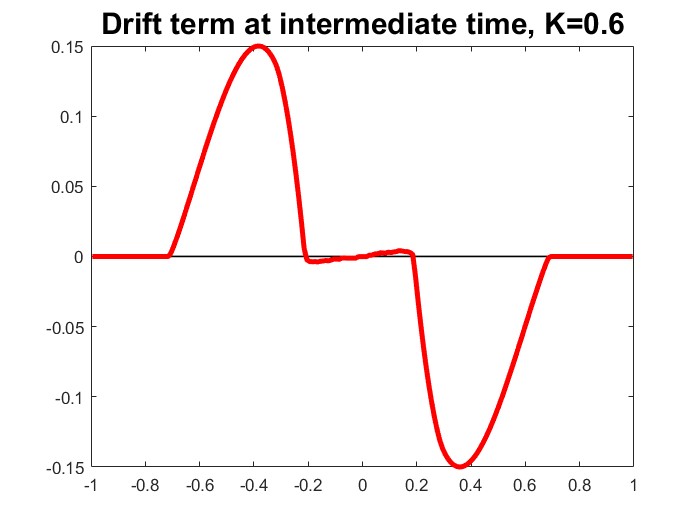}\!\!\includegraphics[scale=0.26]{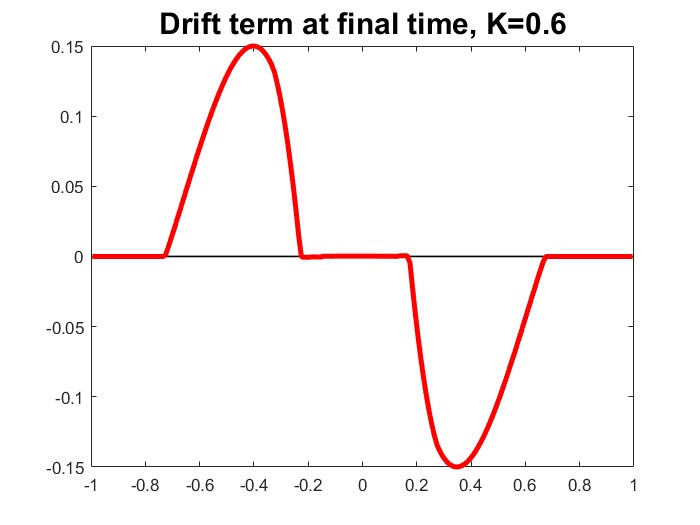}\\
\includegraphics[scale=0.15]{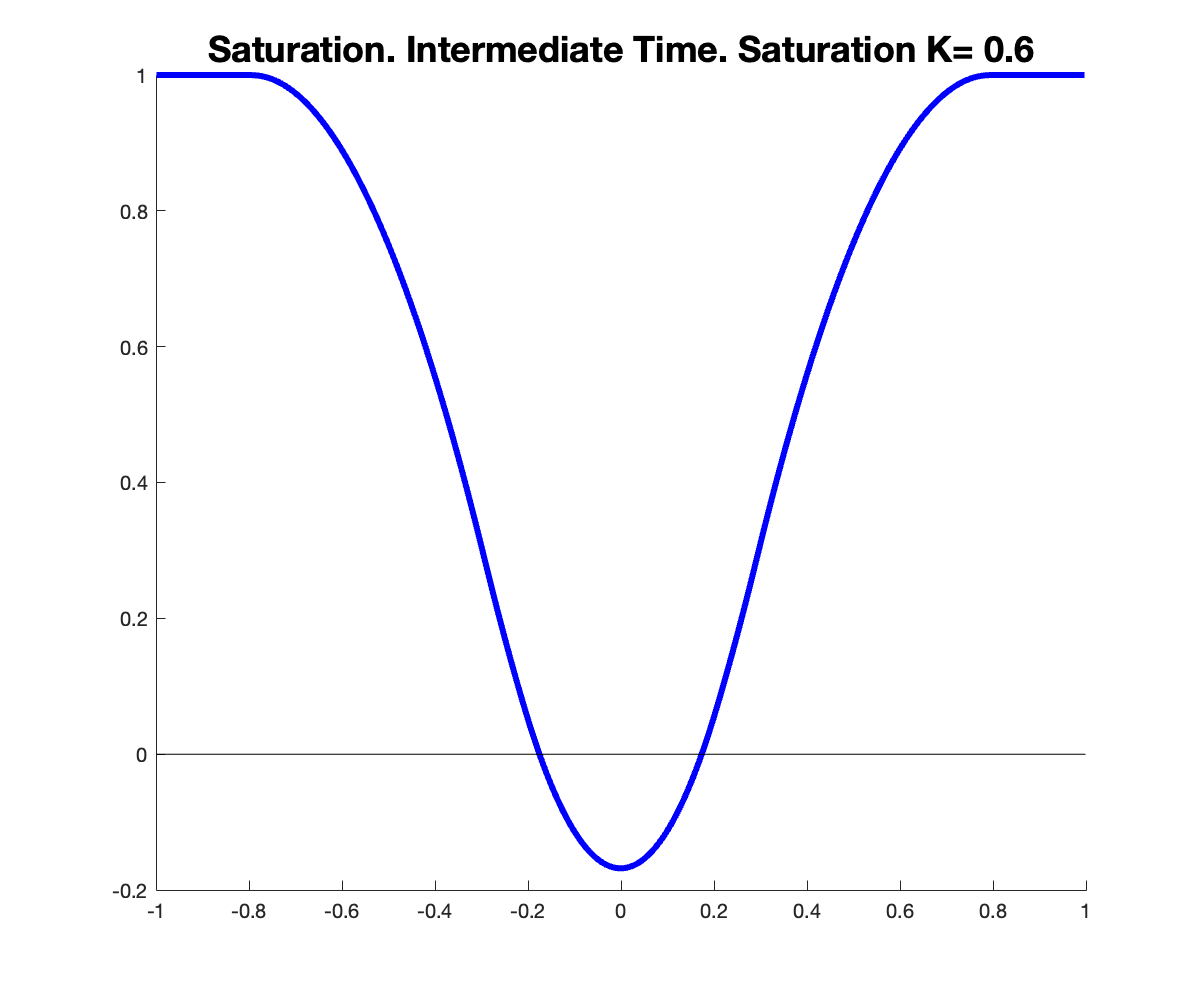}\!\!\!\!\includegraphics[scale=0.15]{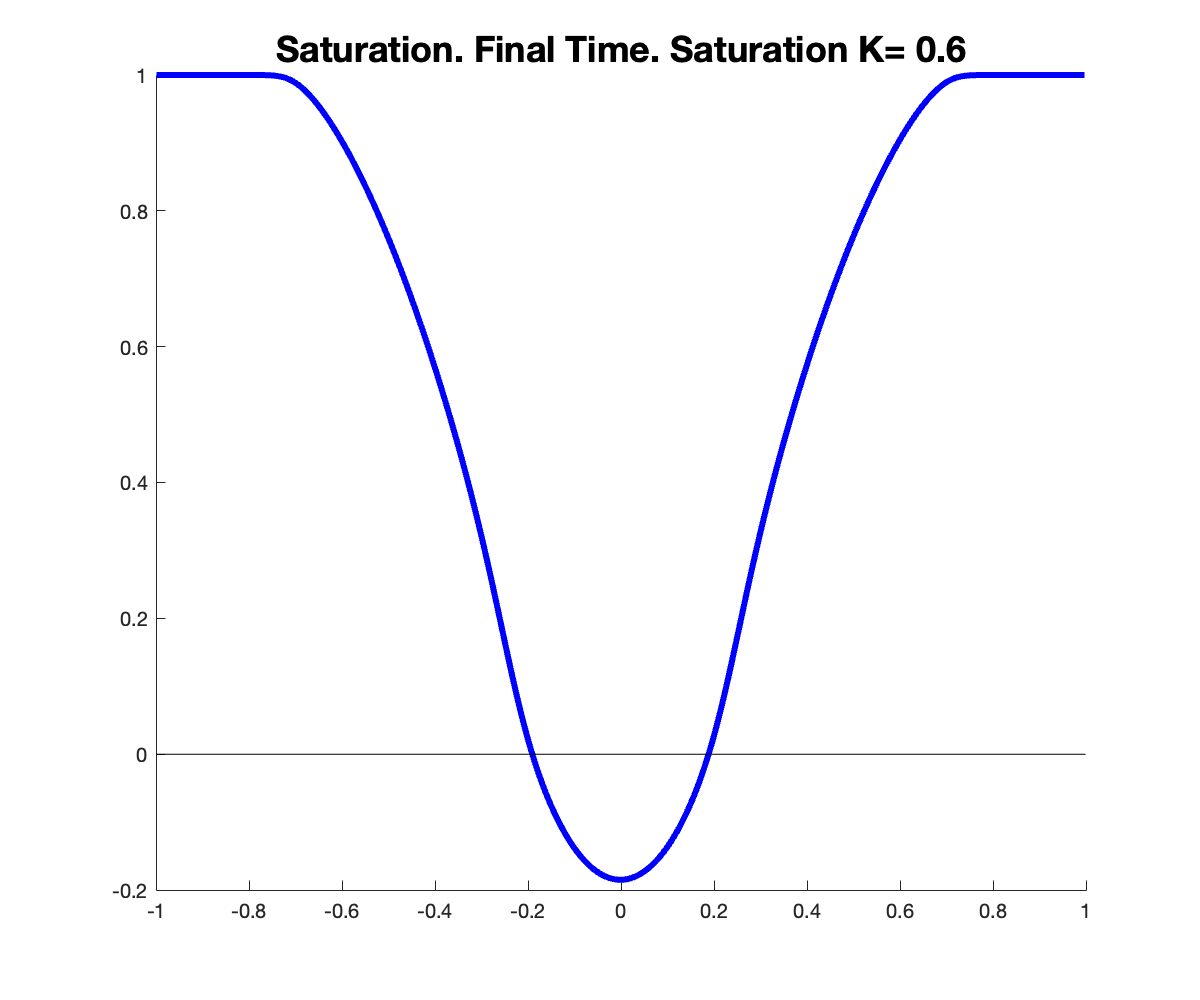}\\
\includegraphics[scale=0.26]{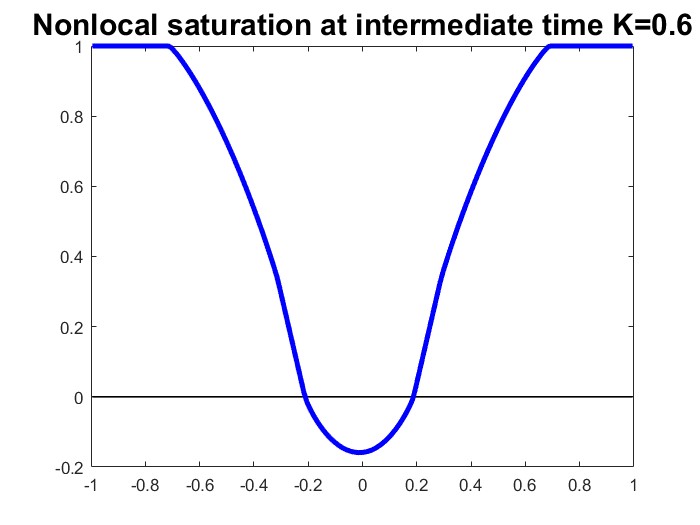}\includegraphics[scale=0.26]{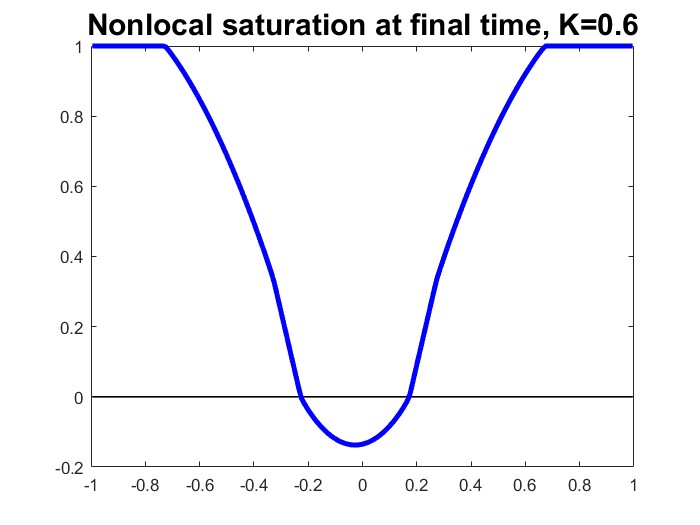}
\captionof{figure}{Drifts and saturations for nonlocal saturation model, $K=0.6$, Concentrated initial data within $[-0.3,0.3]$. Drifts in first and second rows for  PDE model (\S 3.5.3) and SDEs (\S 4.2), respectively. Saturations in third and fourth rows for  PDE model (\S 3.5.3) and SDEs (\S 4.2), respectively. All of them at intermediate time on the left and at final time on the right.} 
\label{fig:NL_K06concentrado_driftsatu}
\end{minipage}


\subsection{Intermediate repulsion $K=0.4$}\label{k=04}

\


 \begin{minipage}{\textwidth}
\centering
\includegraphics[scale=0.22]{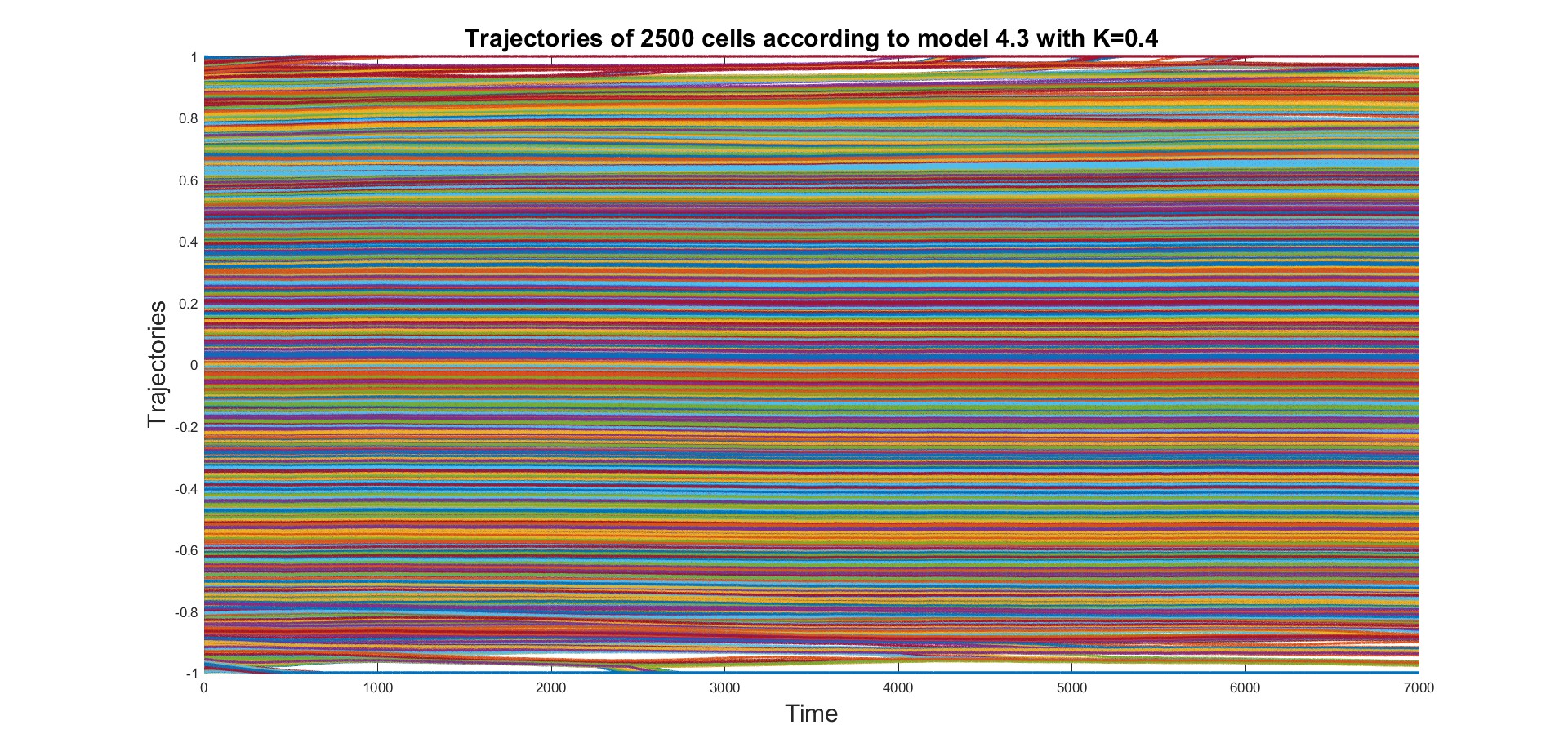}\\

\includegraphics[scale=0.17]{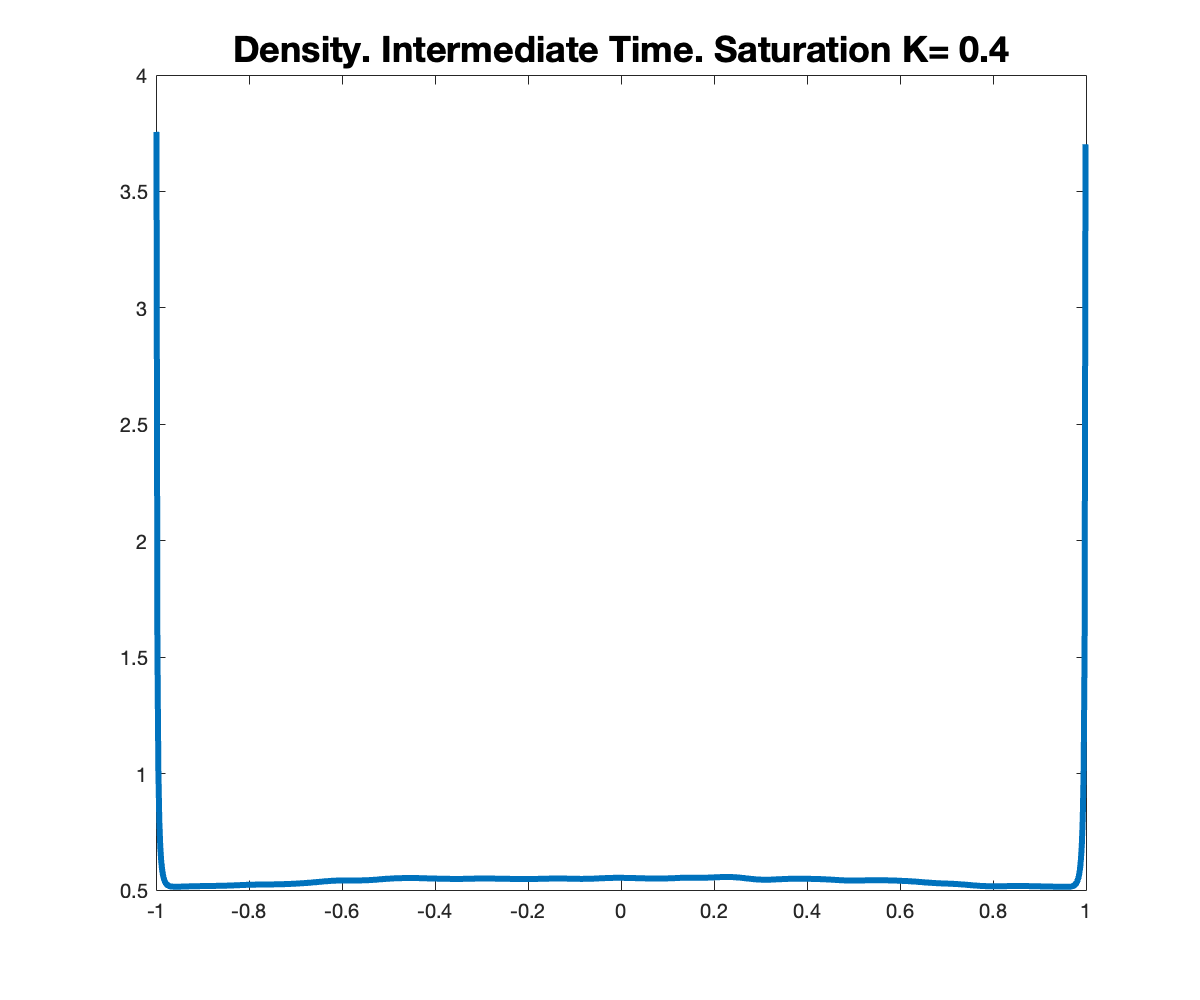}\!\!\includegraphics[scale=0.17]{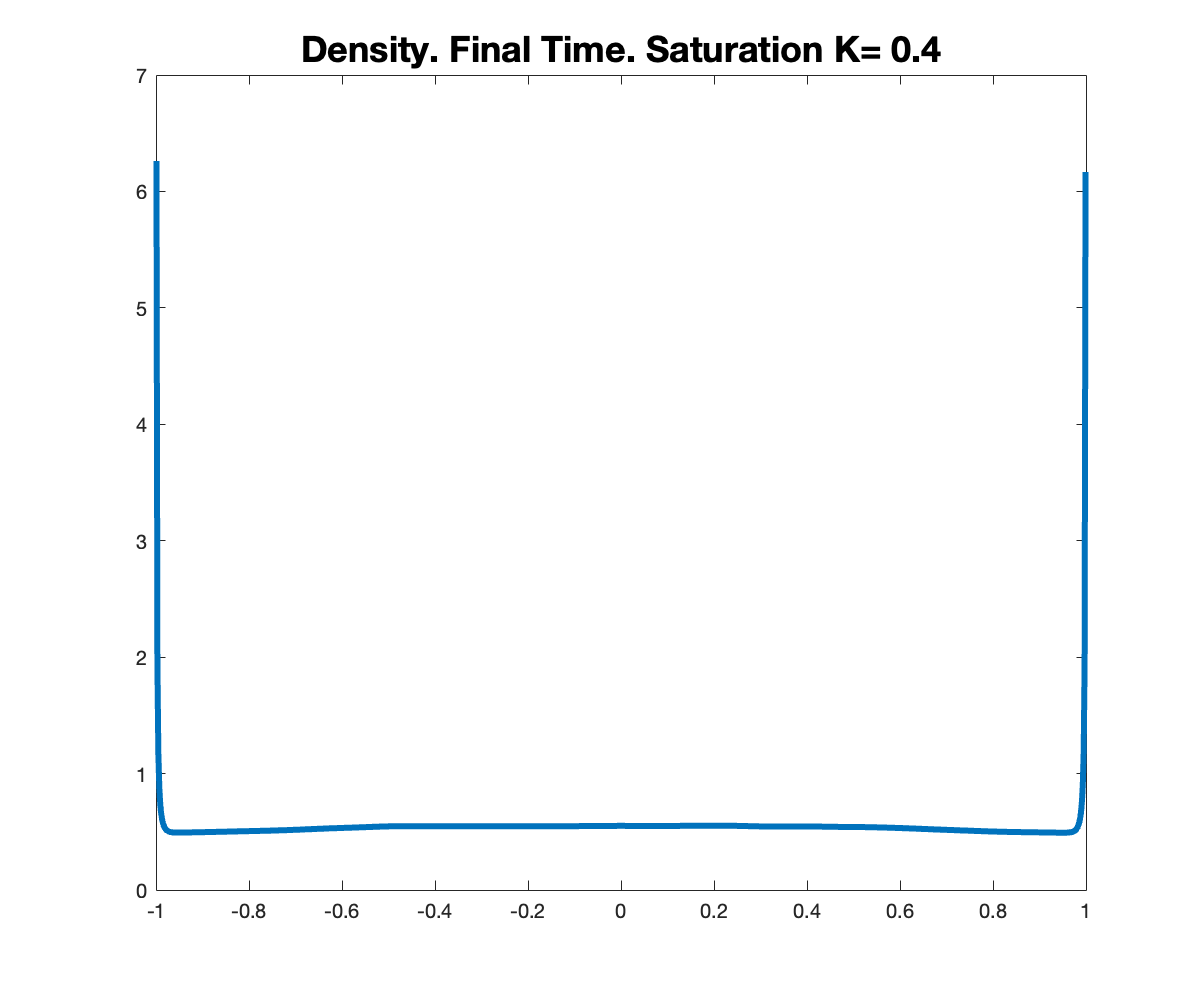}\\
\includegraphics[scale=0.28]{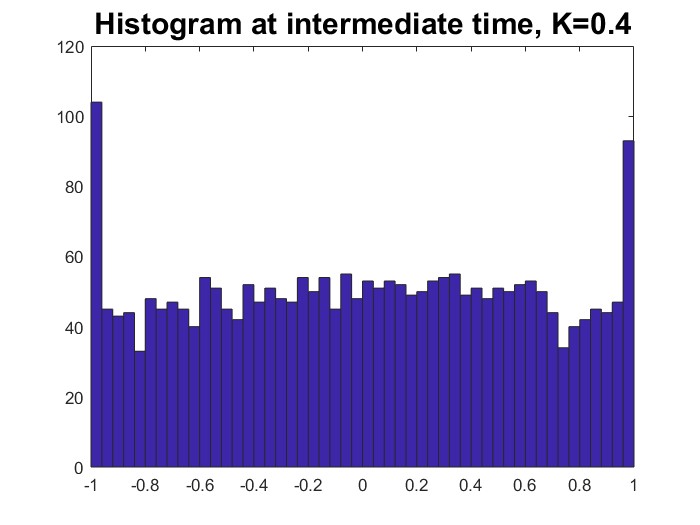}\!\!\includegraphics[scale=0.28]{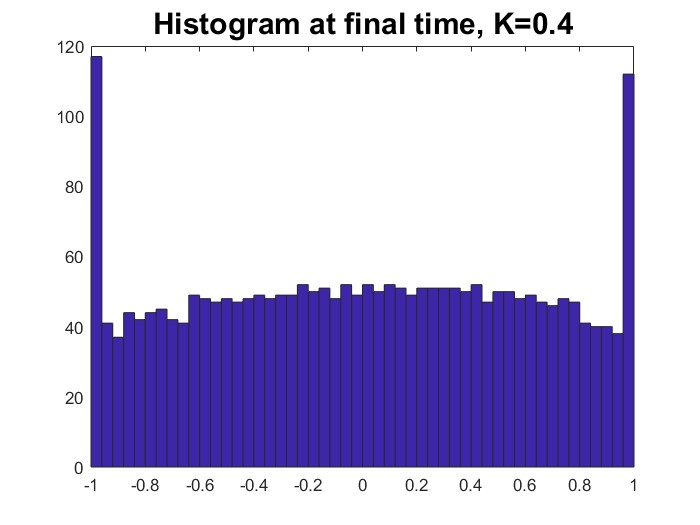}
\captionof{figure}{Local saturation model, $K=0.4$. First row trajectories with $N=2500$ cells, second row densities of the PDE model (\S 3.5.2) and third one histograms of SDEs (\S 4.3). Both at intermediate time on the left and at final time on the right.
} 
\label{fig:Carrillo_K04_histogramas}
\end{minipage}

\newpage

 \begin{minipage}{\textwidth}
\centering
\includegraphics[scale=0.17]{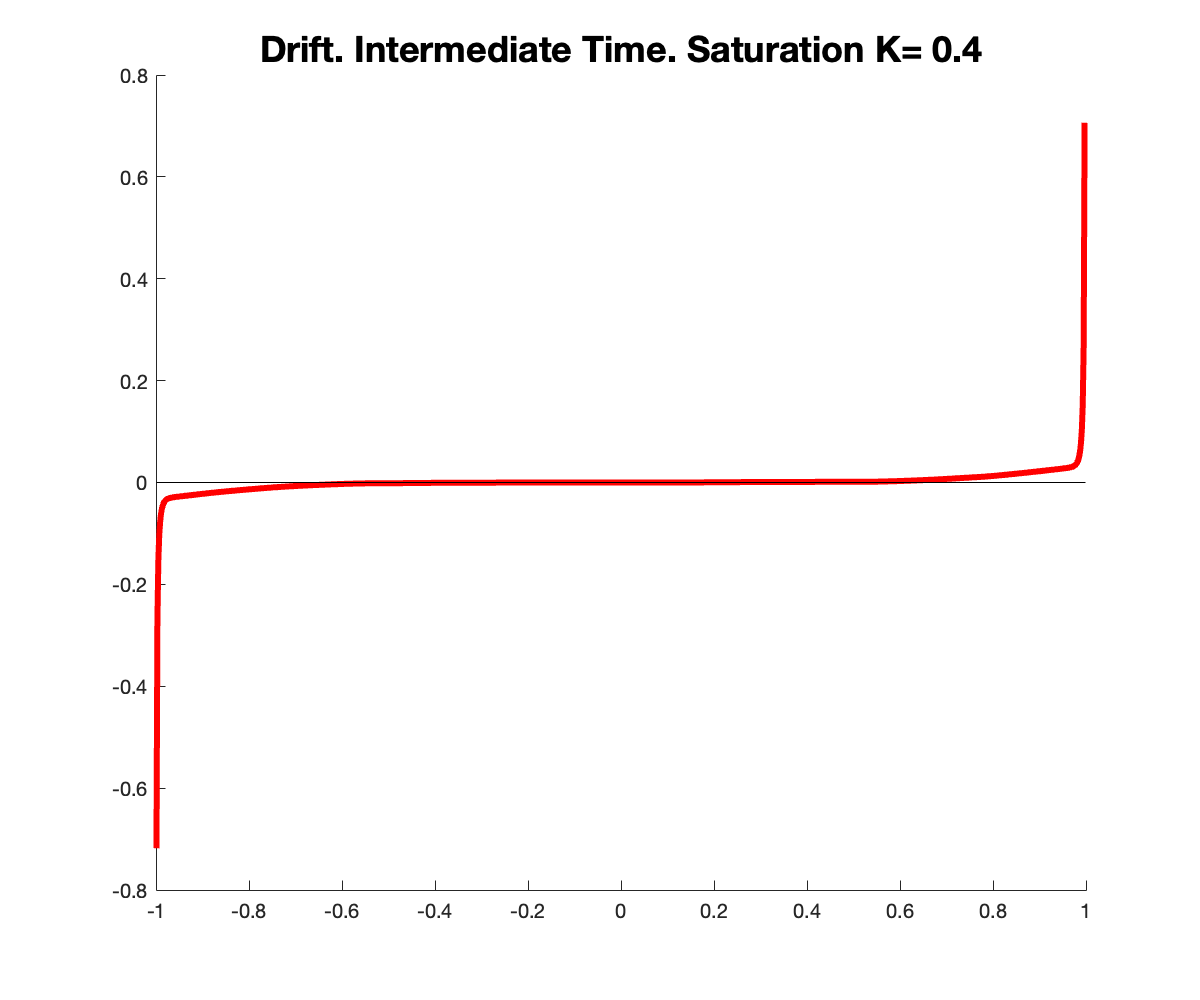}\!\!\!\!\includegraphics[scale=0.17]{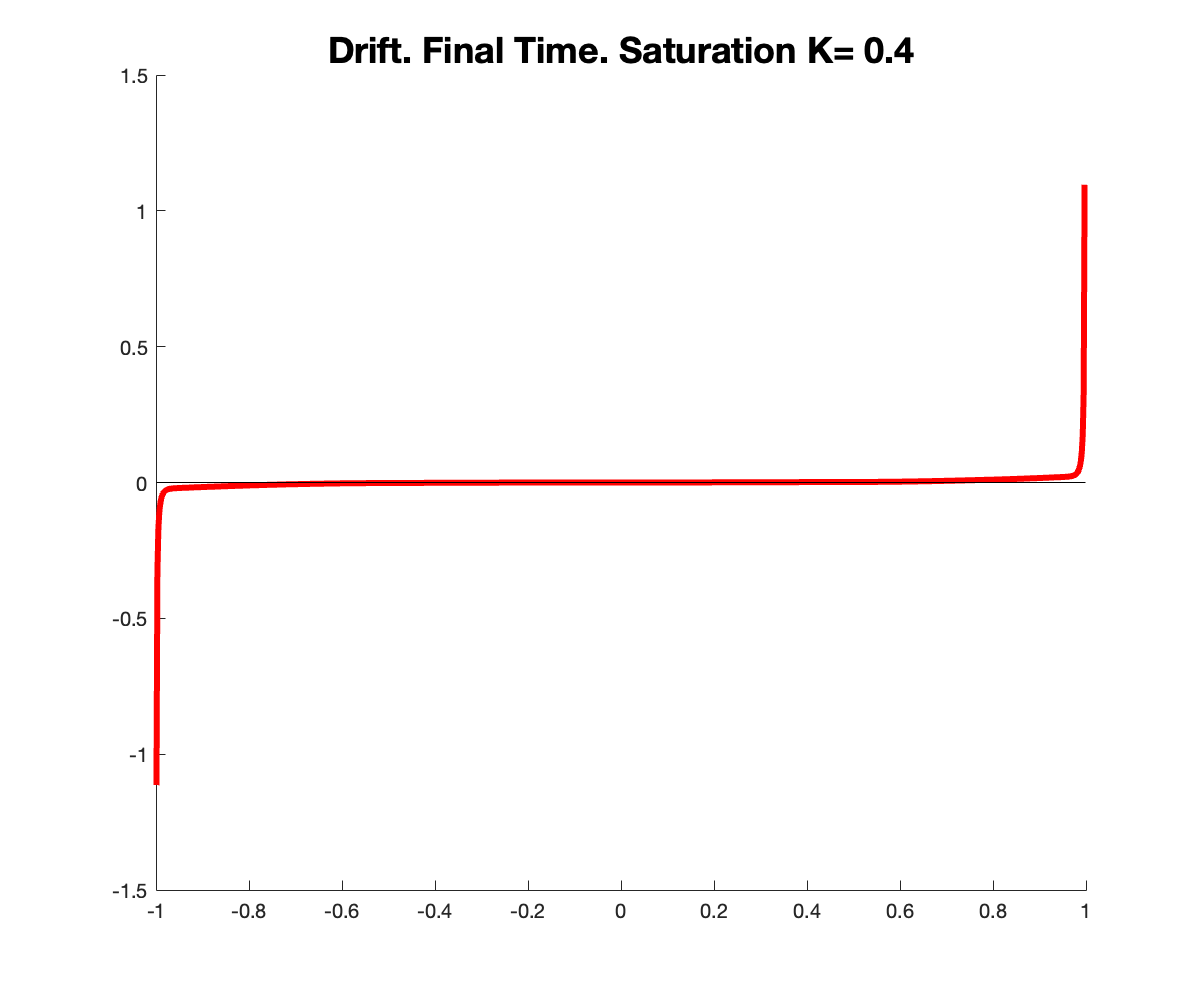}\\

\includegraphics[scale=0.3]{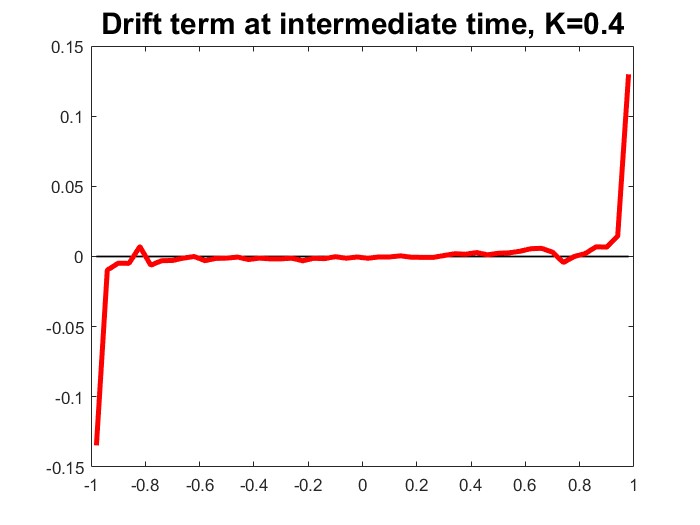}\!\!\includegraphics[scale=0.3]{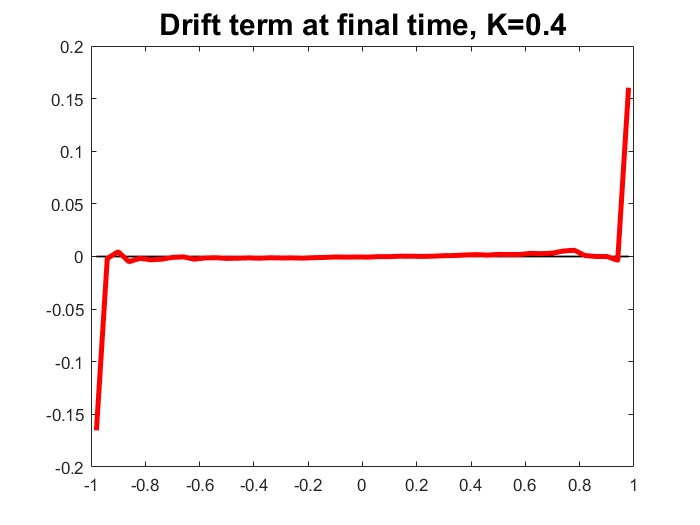}
\captionof{figure}{Drifts for local saturation model, $K=0.4$. First row PDE model (\S 3.5.2) and second row SDEs (\S 4.3). Both at intermediate time on the left and at final time on the right.} 
\label{fig:Carrillo_K04_Drift}
\end{minipage}



     \begin{minipage}{\textwidth}
\centering
\includegraphics[scale=0.24]{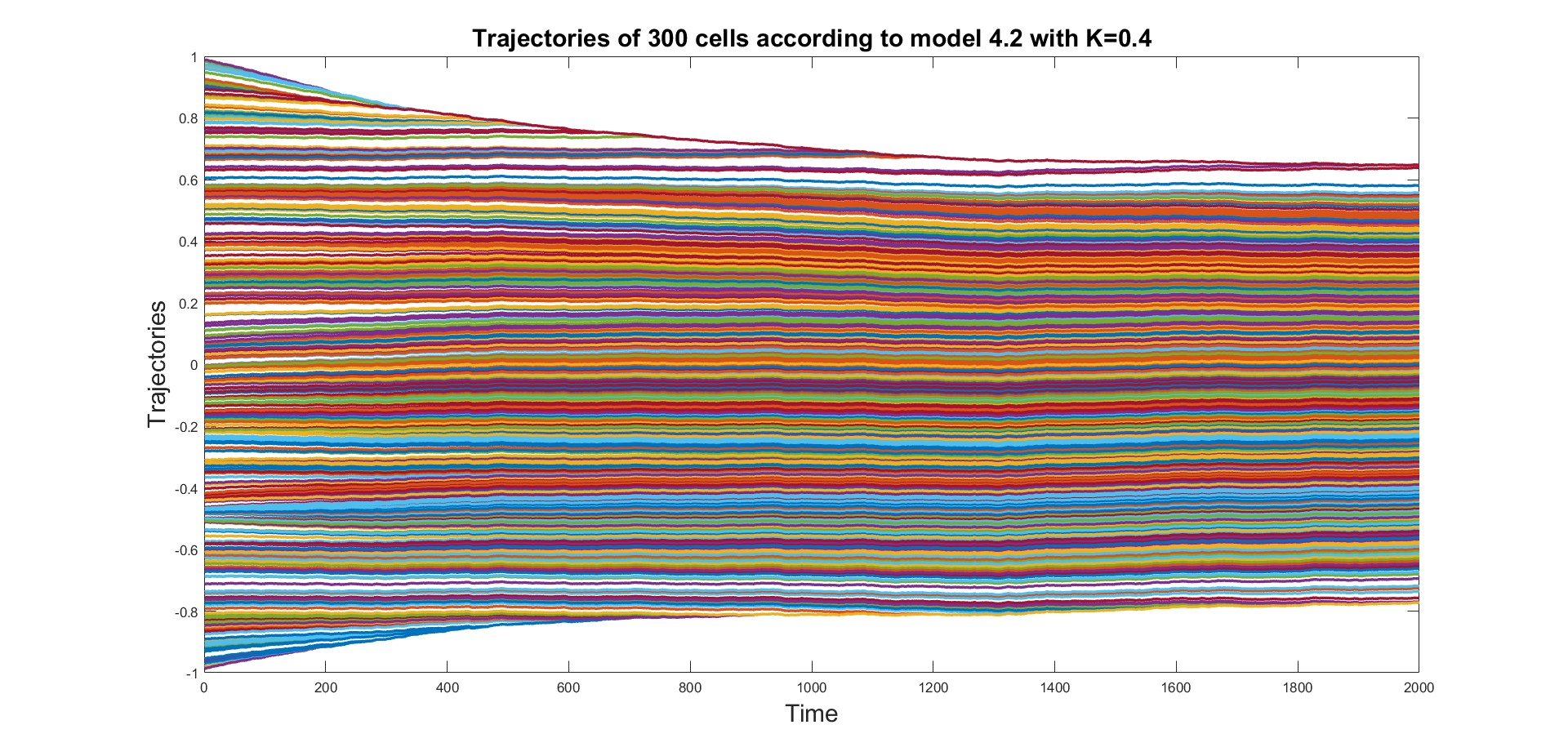}\\

\includegraphics[scale=0.35]{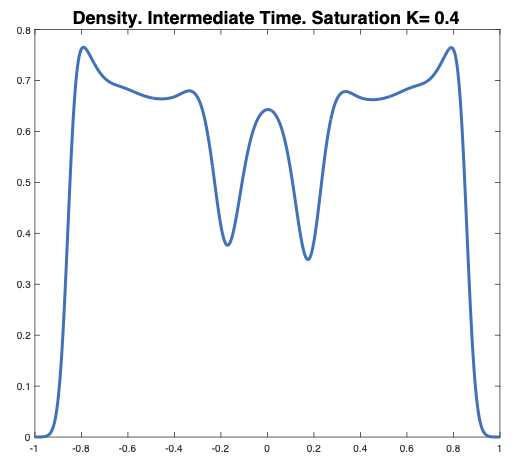}\,\,\includegraphics[scale=0.35]{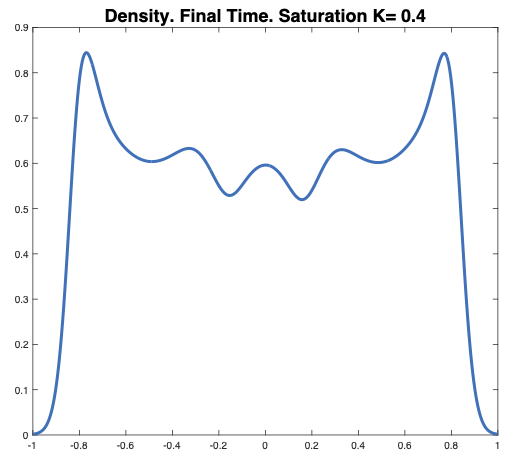}\\
\includegraphics[scale=0.3]{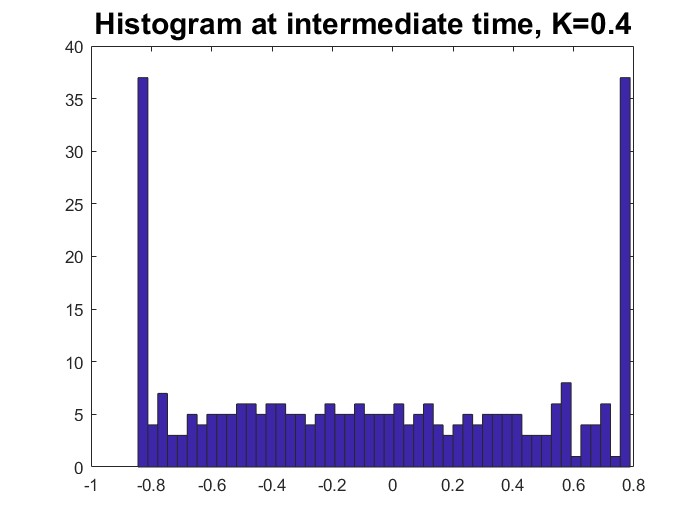}\!\!\includegraphics[scale=0.3]{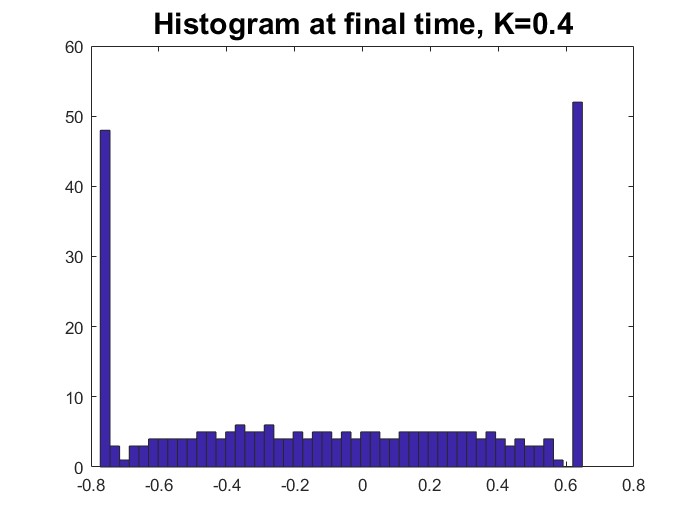}
\captionof{figure}{Nonlocal saturation model, $K=0.4$. First row trajectories with $N=300$ cells, second row densities of the PDE model (\S 3.5.3) and third one histograms of SDEs (\S 4.2). Both at intermediate time on the left and at final time on the right.} 
\label{fig:NL_K04_histogramas}
\end{minipage}

 \begin{minipage}{\textwidth}
\centering
\includegraphics[scale=0.32]{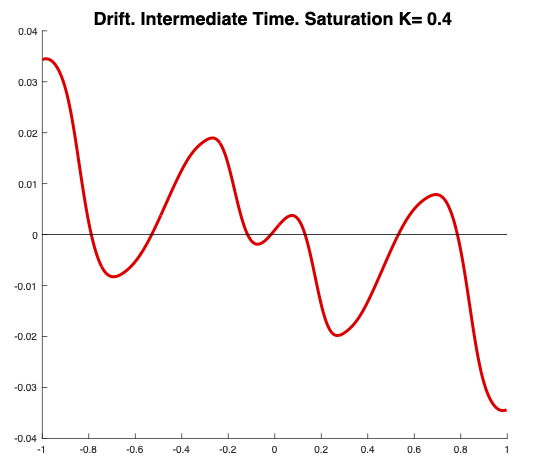}\,\,\,\,\includegraphics[scale=0.32]{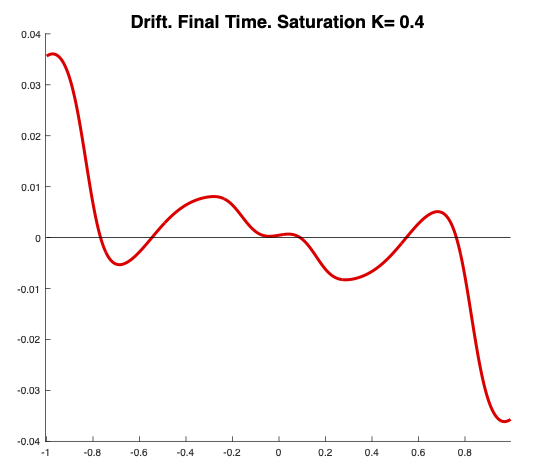}\\

\includegraphics[scale=0.27]{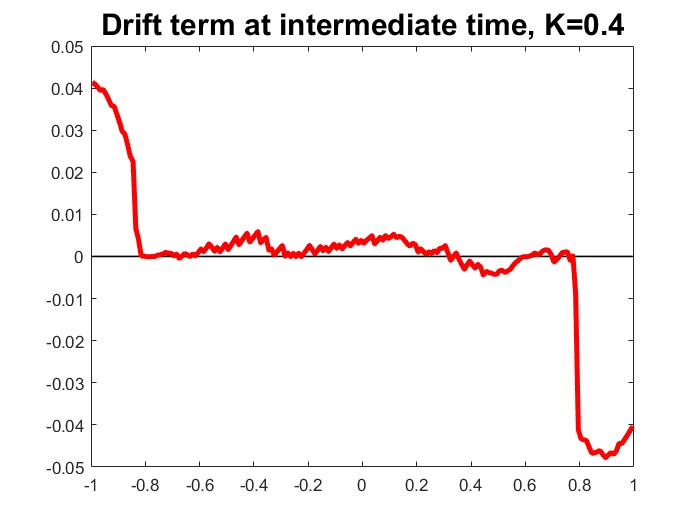}\!\!\includegraphics[scale=0.27]{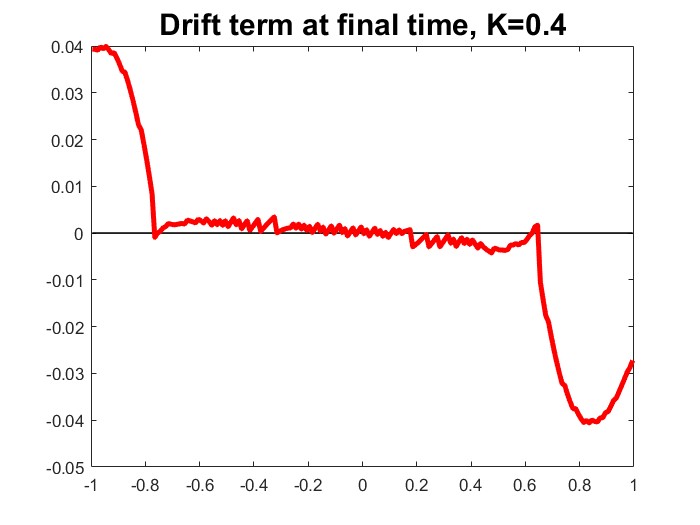}

\includegraphics[scale=0.32]{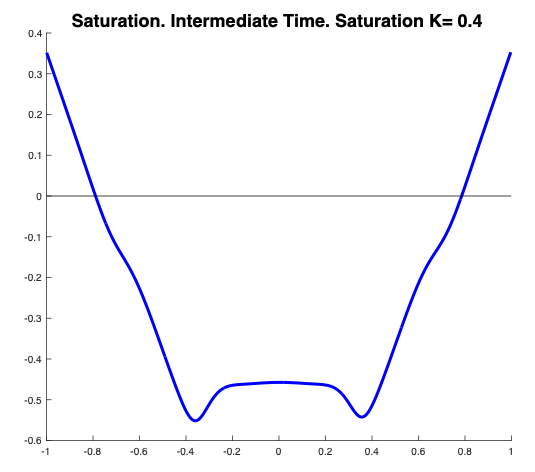}\,\,\,\includegraphics[scale=0.32]{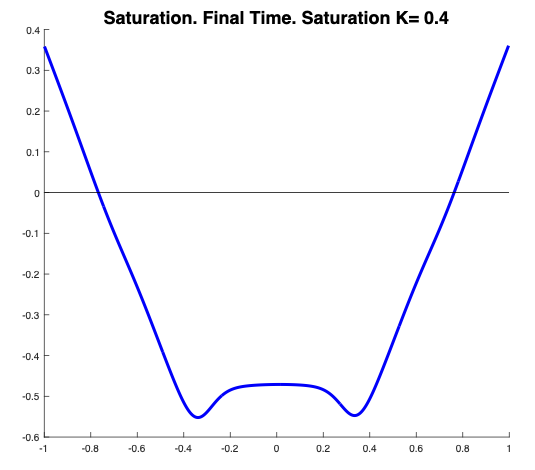}\\

\includegraphics[scale=0.27]{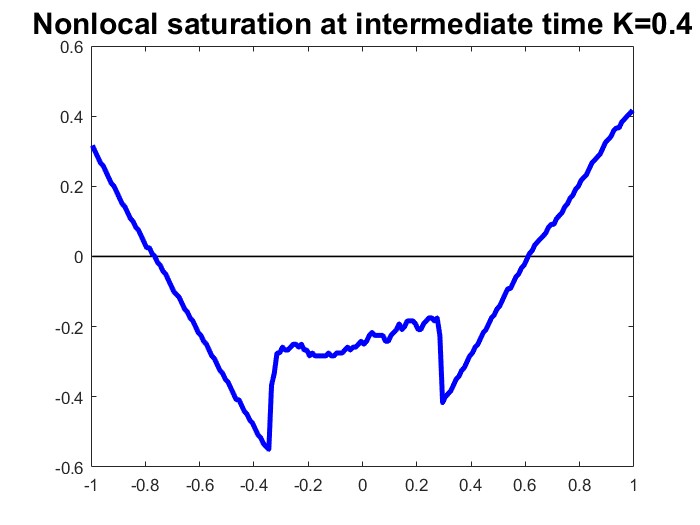}\includegraphics[scale=0.27]{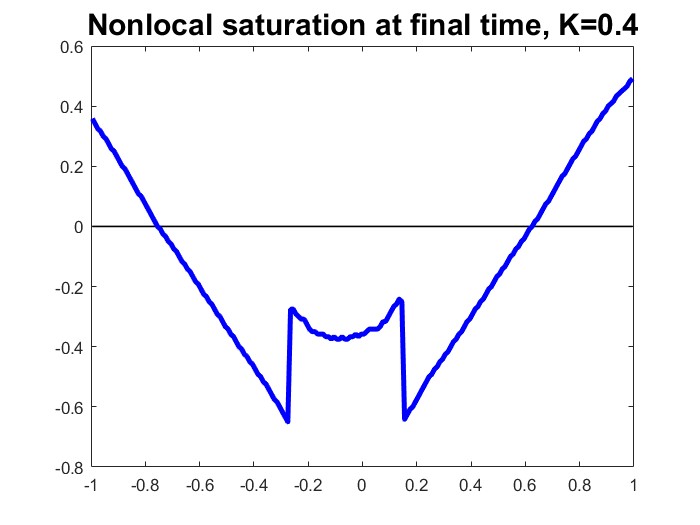}
\captionof{figure}{Drifts and saturations for nonlocal saturation model, $K=0.4$. Drifts in first and second rows for  PDE model (\S 3.5.3) and SDEs (\S 4.2), respectively. Saturations in third and fourth rows for  PDE model (\S 3.5.3) and SDEs (\S 4.2), respectively. All of them at intermediate time on the left and at final time on the right.} 
\label{fig:NL_K04_driftsatu}
\end{minipage}


     \begin{minipage}{\textwidth}
\centering
\includegraphics[scale=0.24]{DEFINITIVAS/Modelo3K04tubos/mod3-tray-k04.jpg}\\

\includegraphics[scale=0.17]{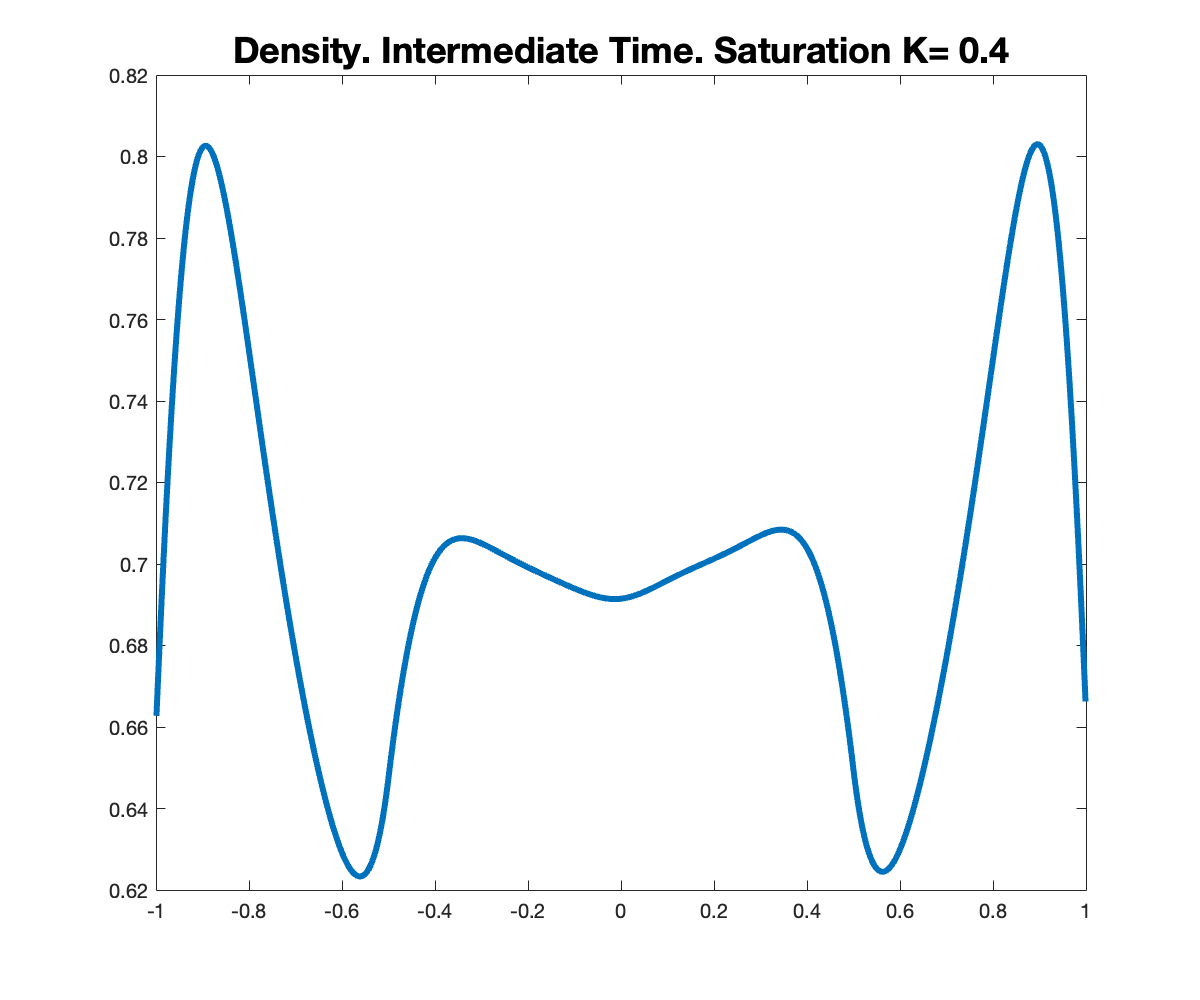}\,\,\includegraphics[scale=0.17]{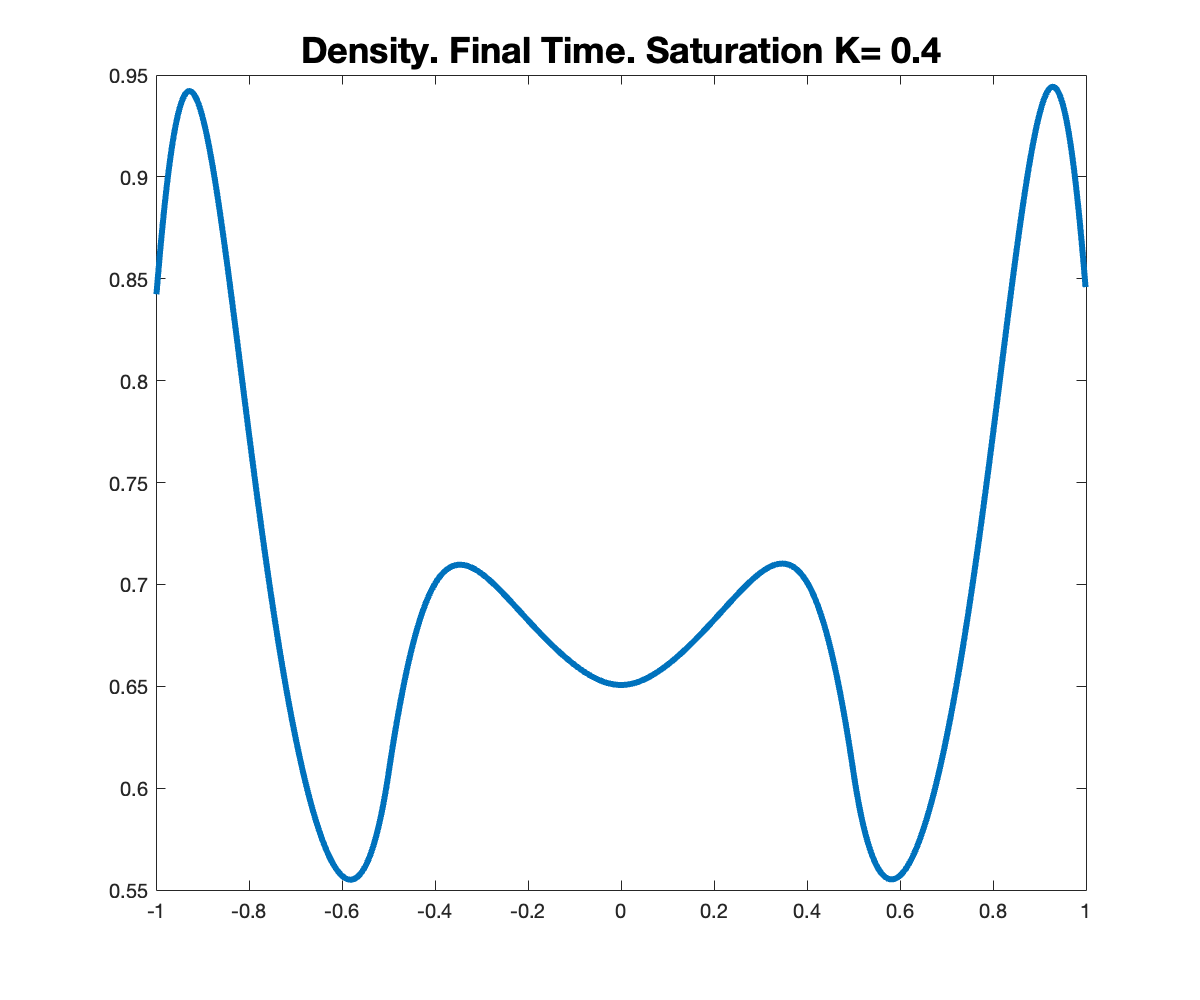}\\
\includegraphics[scale=0.3]{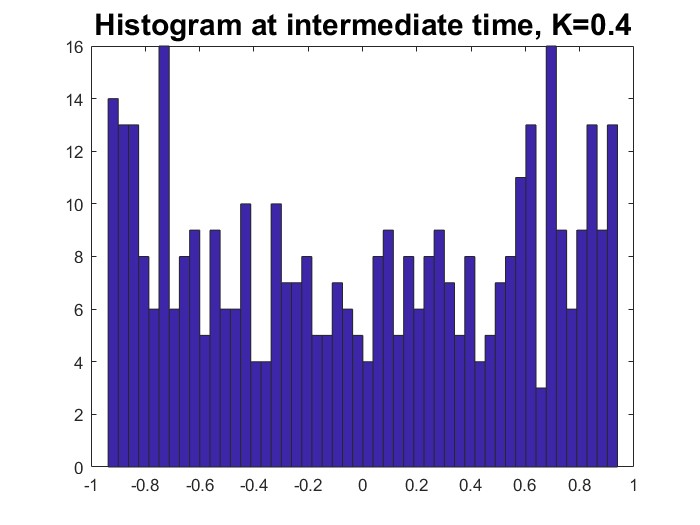}\!\!\includegraphics[scale=0.3]{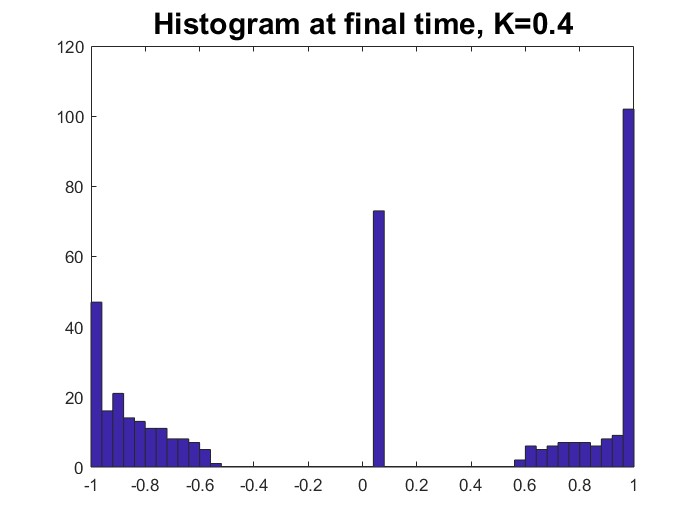}
\captionof{figure}{Nonlocal saturation model, $K=0.4$ and constant weight. First row trajectories with $N=300$ cells, second row densities of the PDE model (\S 3.5.3) and third one histograms of SDEs (\S 4.2). Both at intermediate time on the left and at final time on the right.} 
\label{fig:NL_K04_cw_histogramas}
\end{minipage}

 \begin{minipage}{\textwidth}
\centering
\includegraphics[scale=0.15]{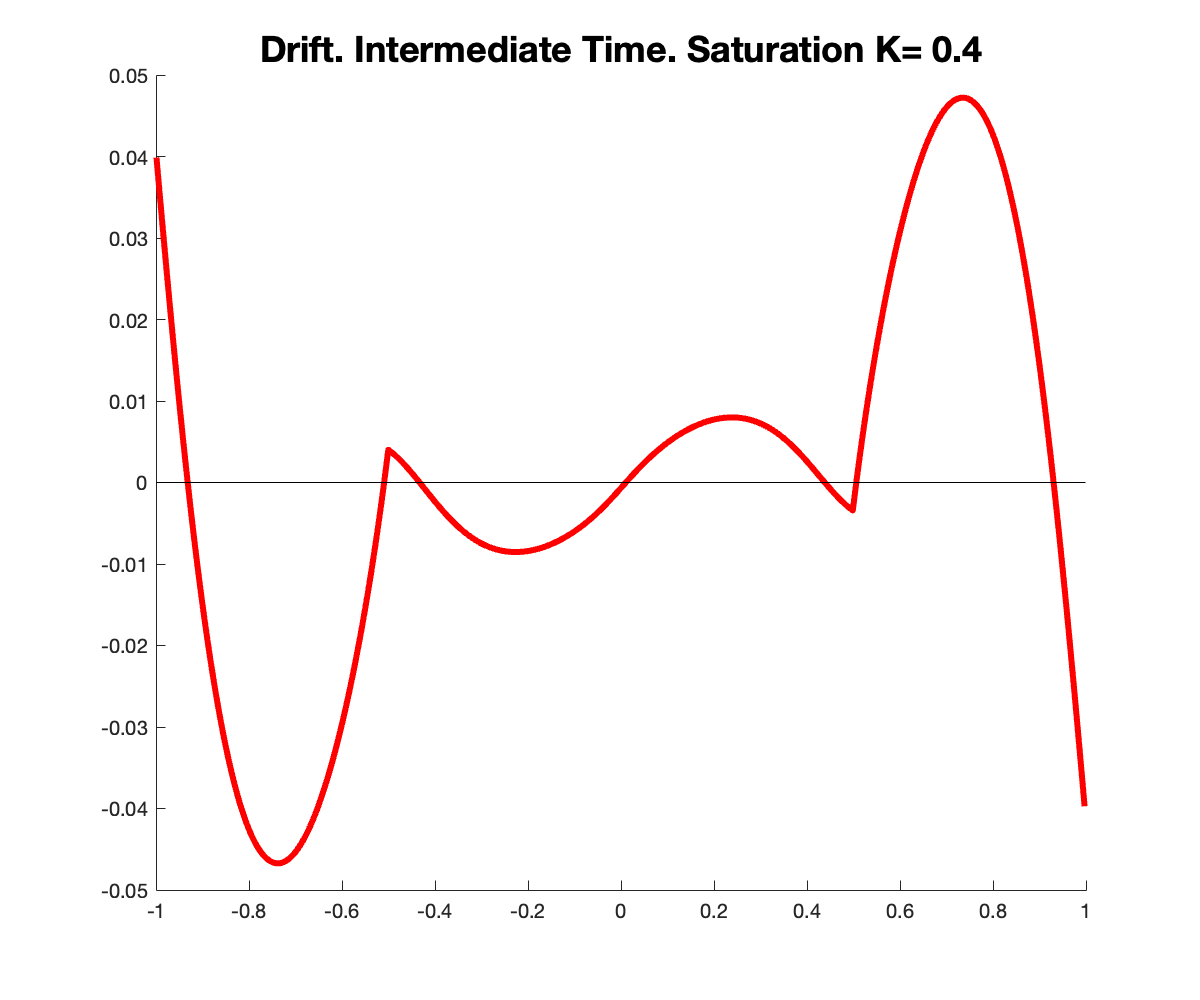}\,\,\,\,\includegraphics[scale=0.15]{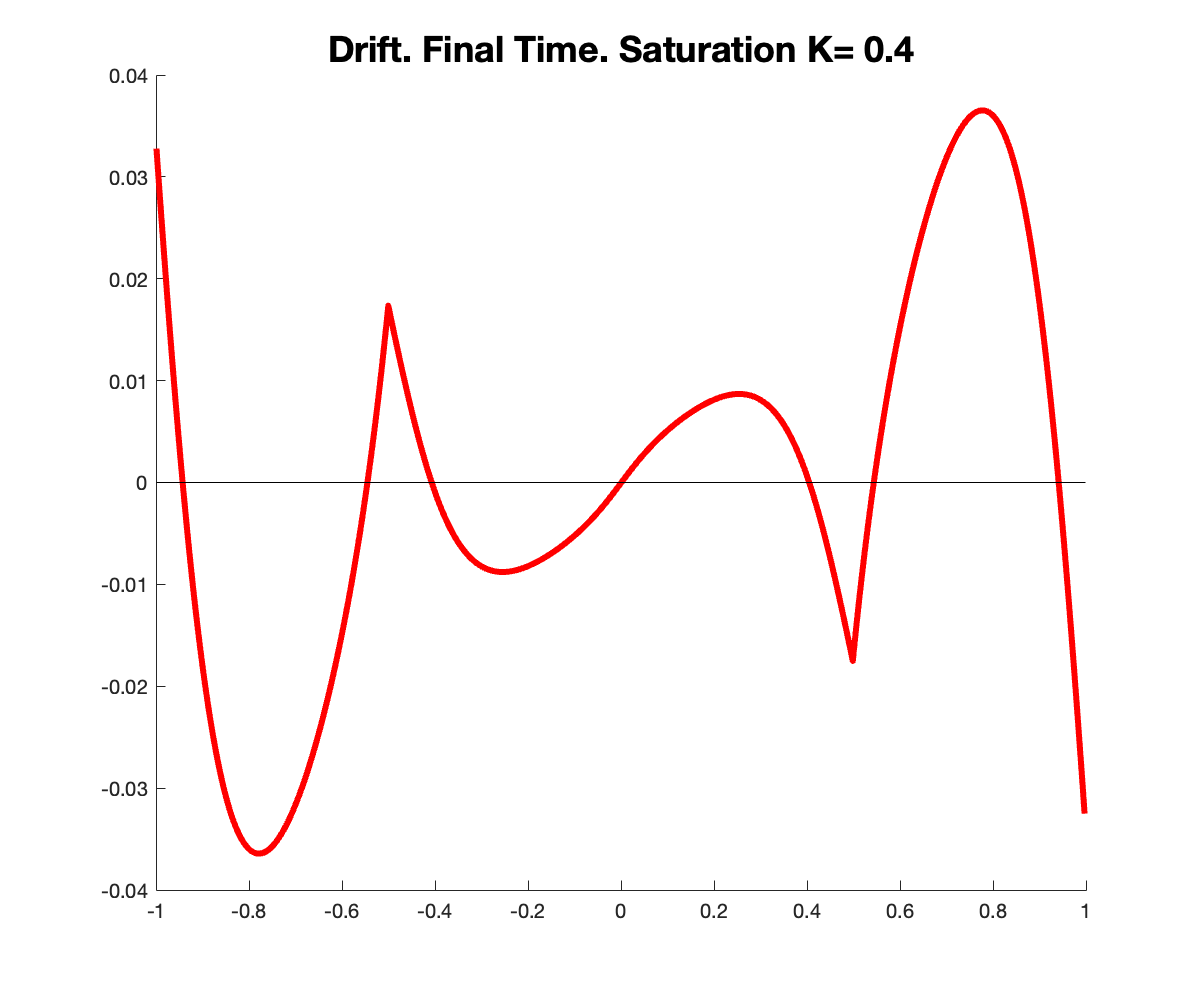}\\

\includegraphics[scale=0.27]{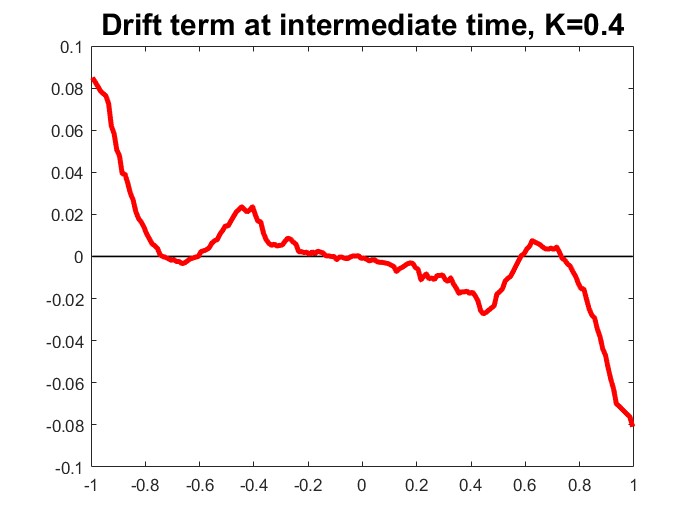}\!\!\includegraphics[scale=0.27]{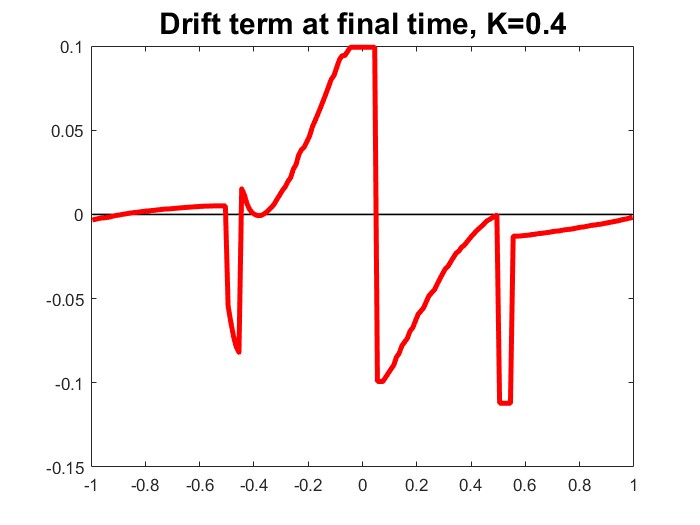}

\includegraphics[scale=0.15]{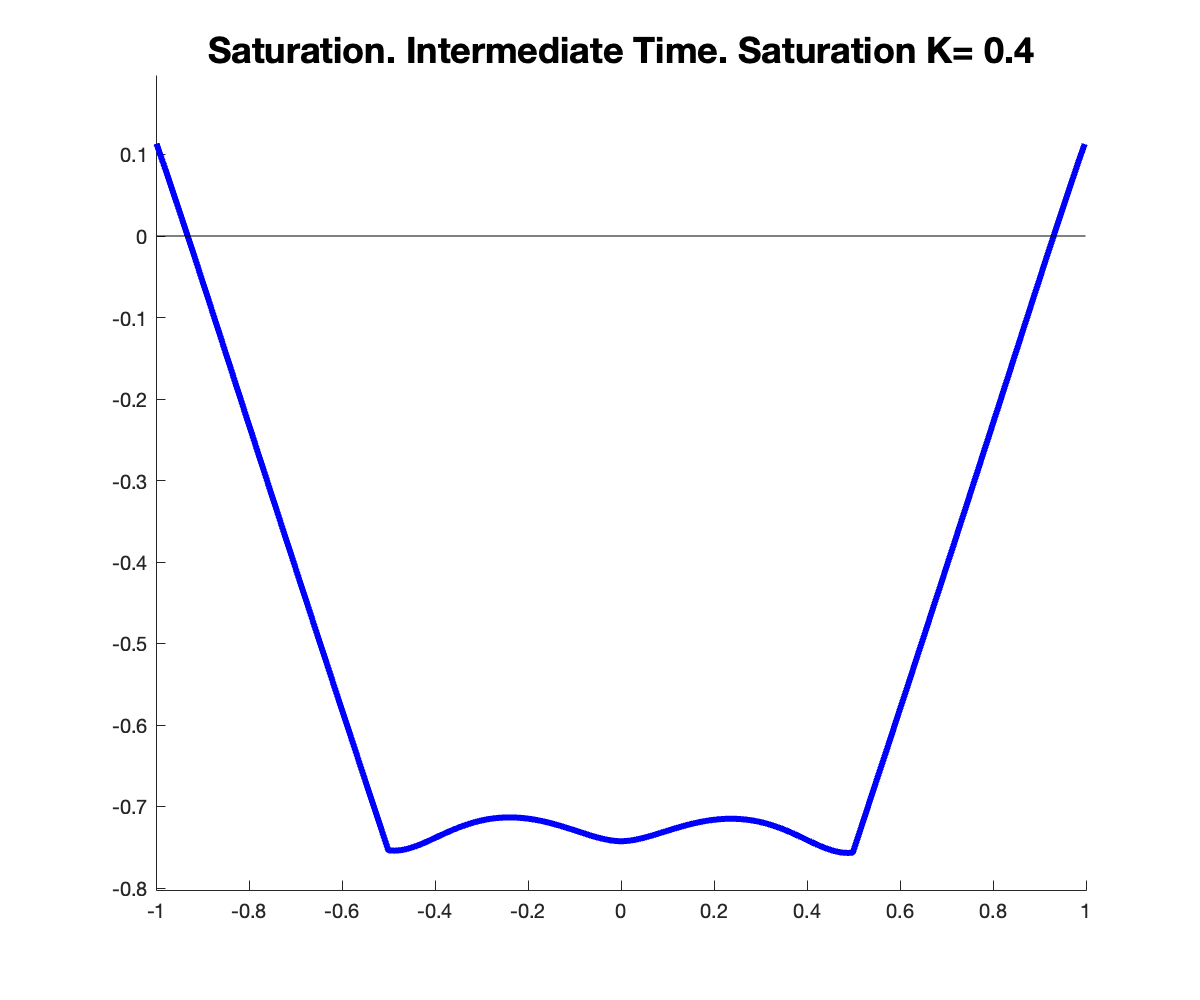}\,\,\,\includegraphics[scale=0.15]{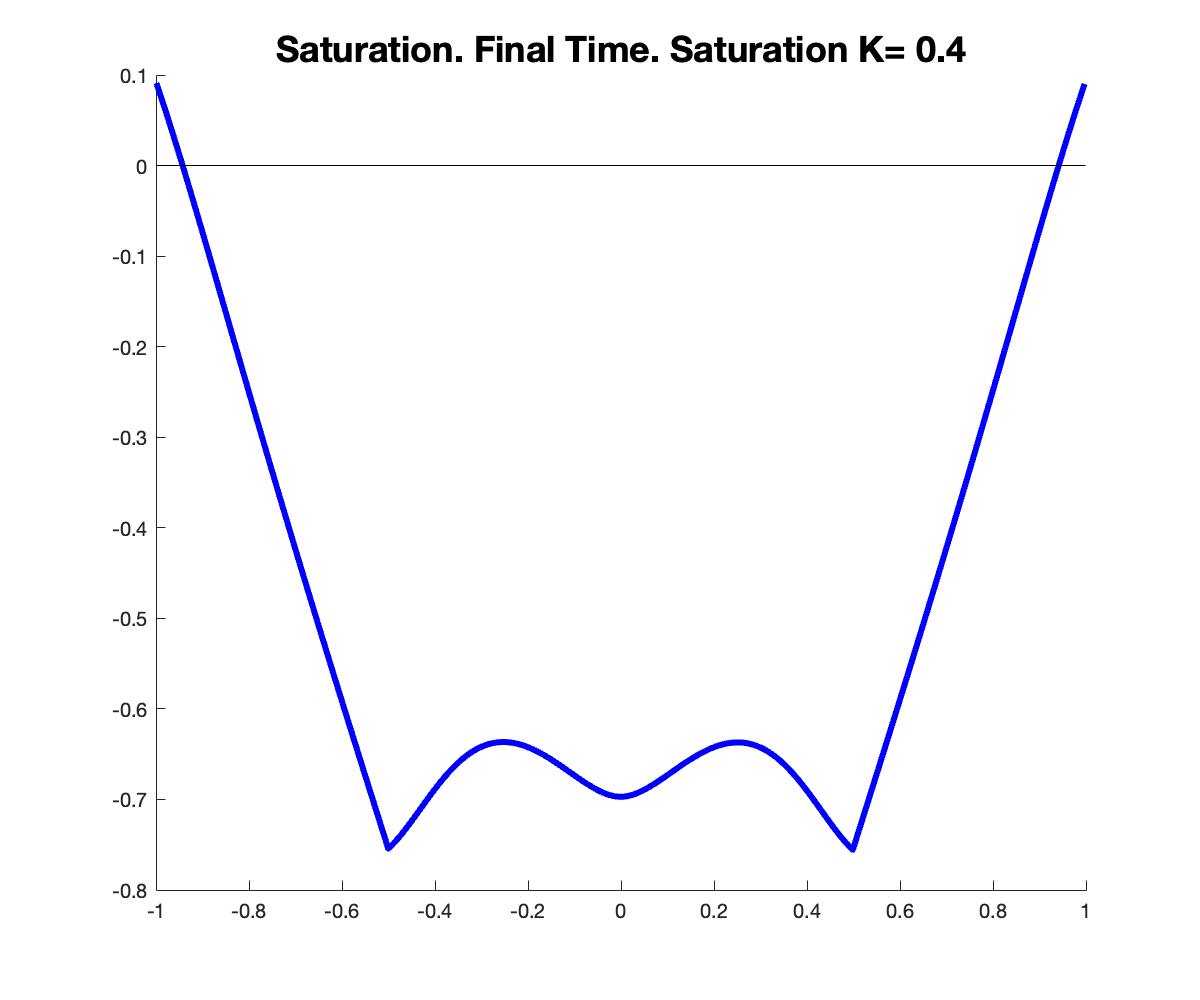}\\

\includegraphics[scale=0.27]{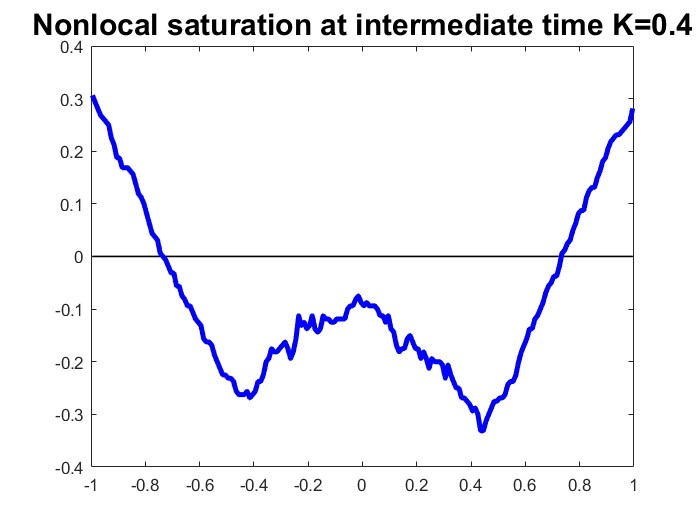}\includegraphics[scale=0.27]{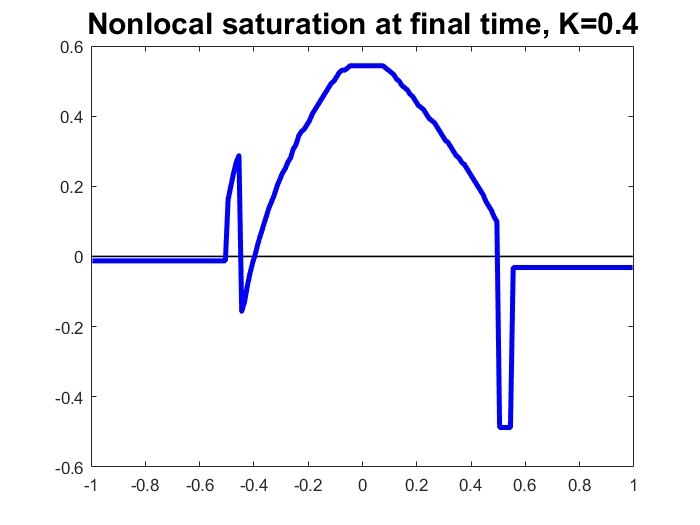}
\captionof{figure}{Drifts and saturations for nonlocal saturation model, $K=0.4$ and constant weight. Drifts in first and second rows for  PDE model (\S 3.5.3) and SDEs (\S 4.2), respectively. Saturations in third and fourth rows for  PDE model (\S 3.5.3) and SDEs (\S 4.2), respectively. All of them at intermediate time on the left and at final time on the right.} 
\label{fig:NL_K04_cw_driftsatu}
\end{minipage}

\subsection{Strong repulsion  $K=0.2$}\label{k=02}

\



    \begin{minipage}{\textwidth}
\centering
\includegraphics[scale=0.24]{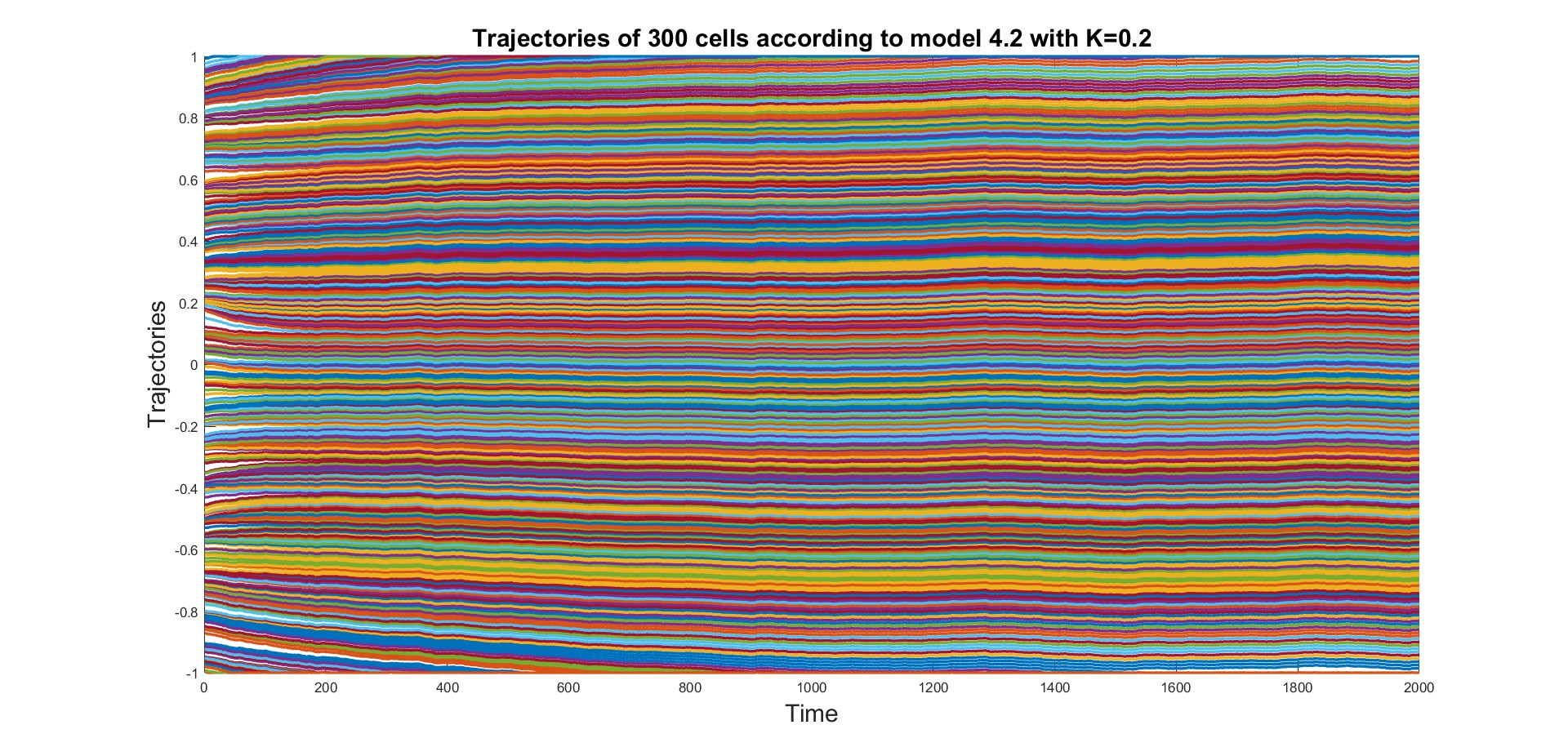}\\
\includegraphics[scale=0.18]{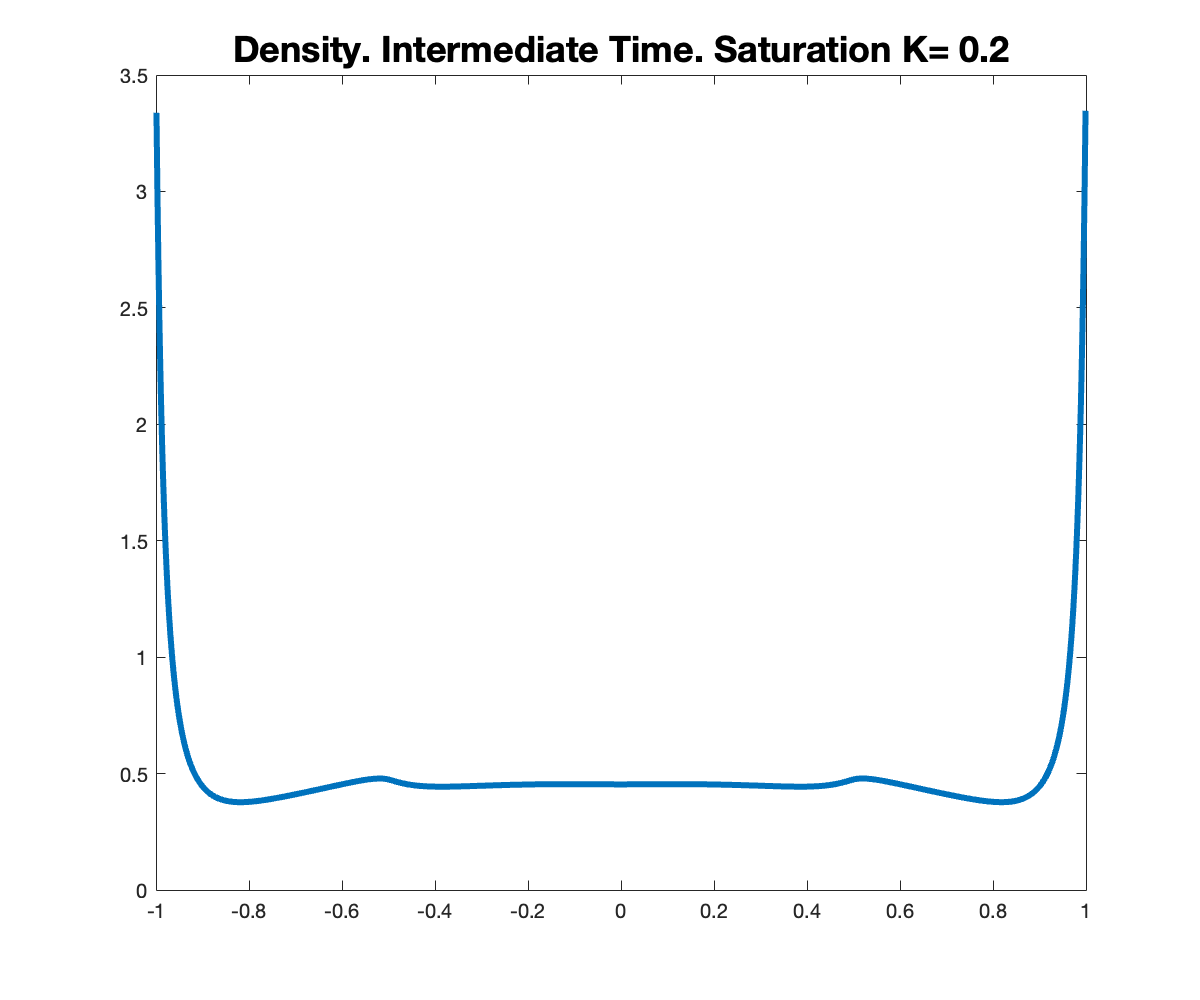}\!\!\includegraphics[scale=0.18]{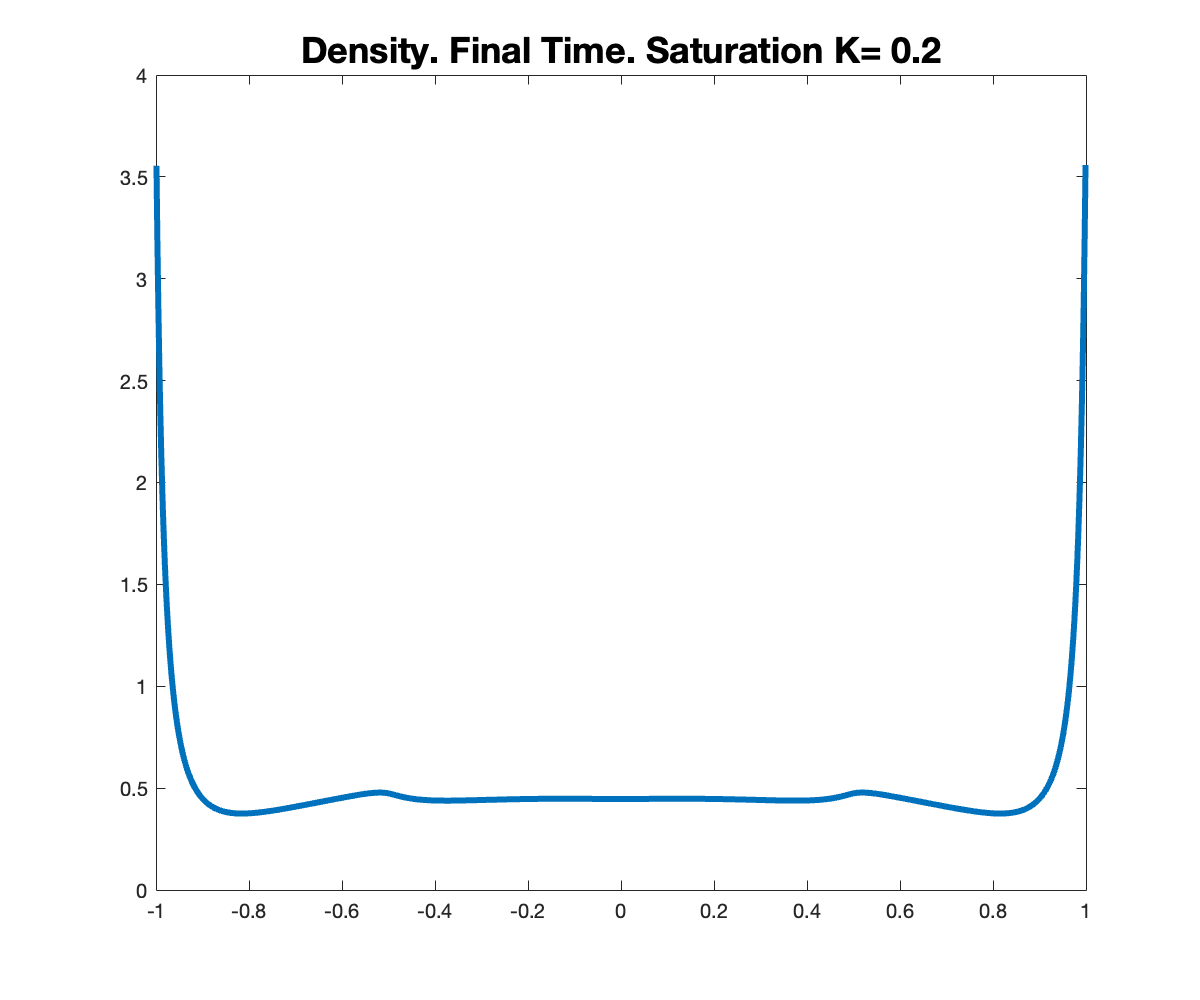}\\
\includegraphics[scale=0.3]{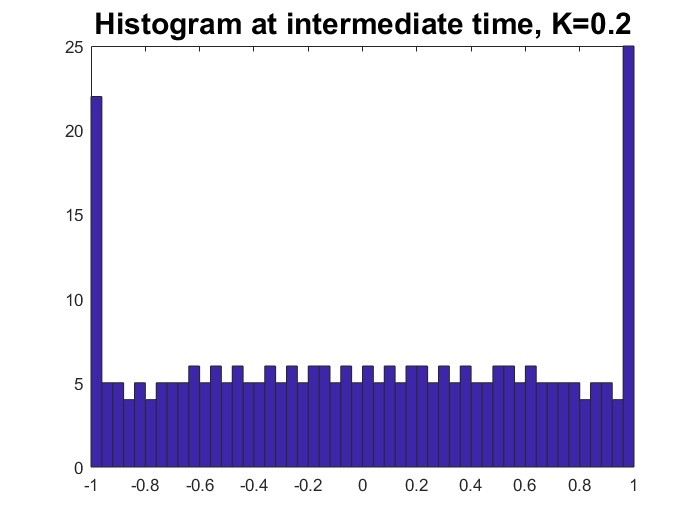}\!\!\includegraphics[scale=0.3]{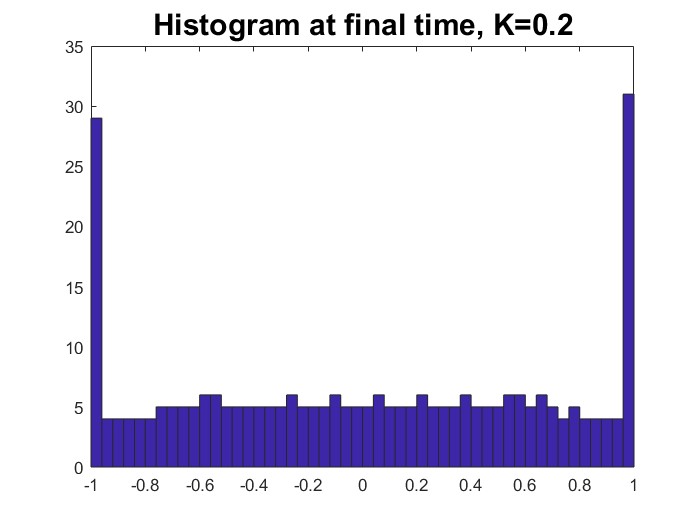}
\captionof{figure}{Nonlocal saturation model, $K=0.2$. First row trajectories with $N=300$ cells, second row densities of the PDE model (\S 3.5.3) and third one histograms of SDEs (\S 4.2). Both at intermediate time on the left and at final time on the right.} 
\label{fig:NL_K02_histogramas}
\end{minipage}

 \begin{minipage}{\textwidth}
\centering
\includegraphics[scale=0.15]{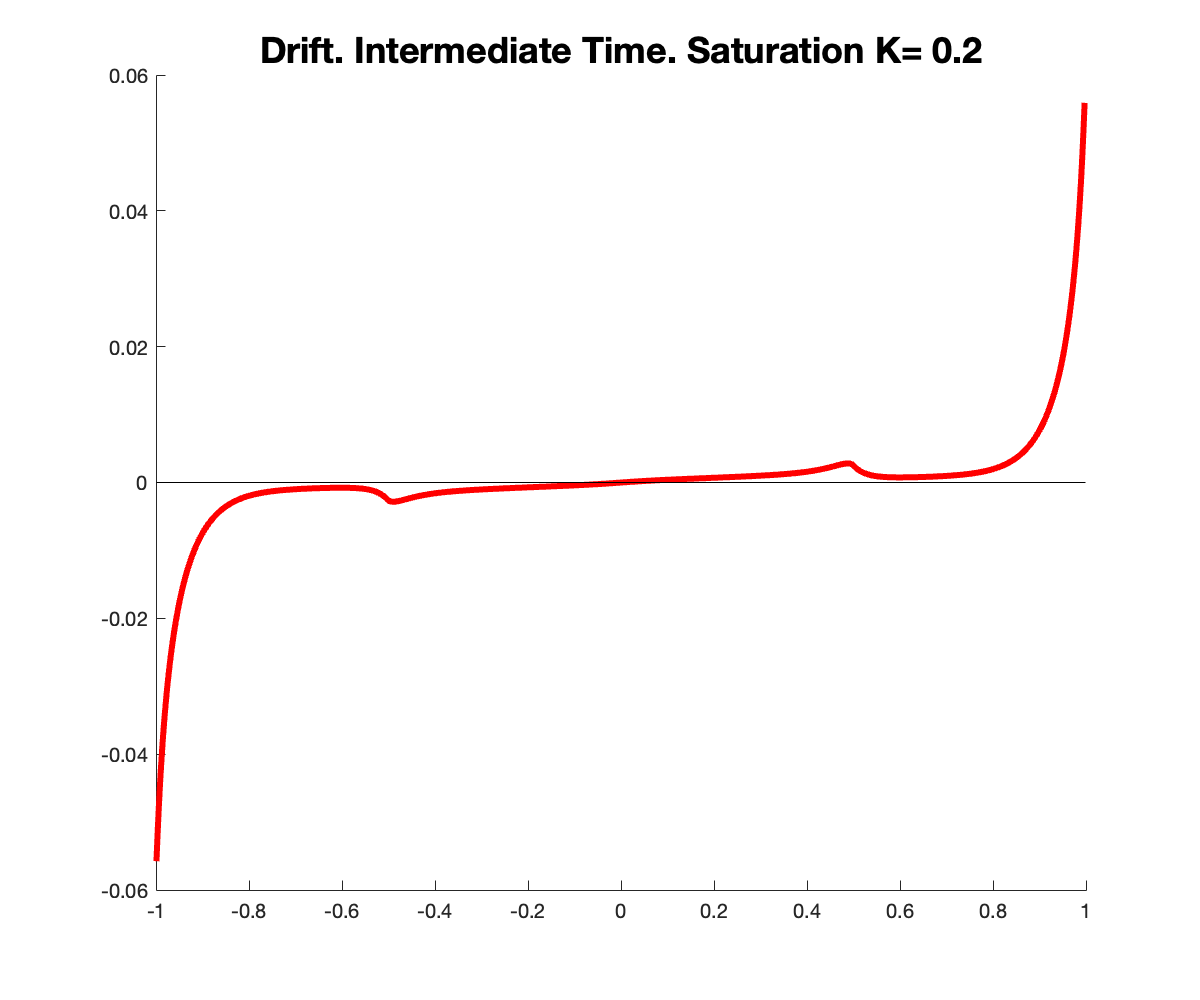}\!\!\!\!\includegraphics[scale=0.15]{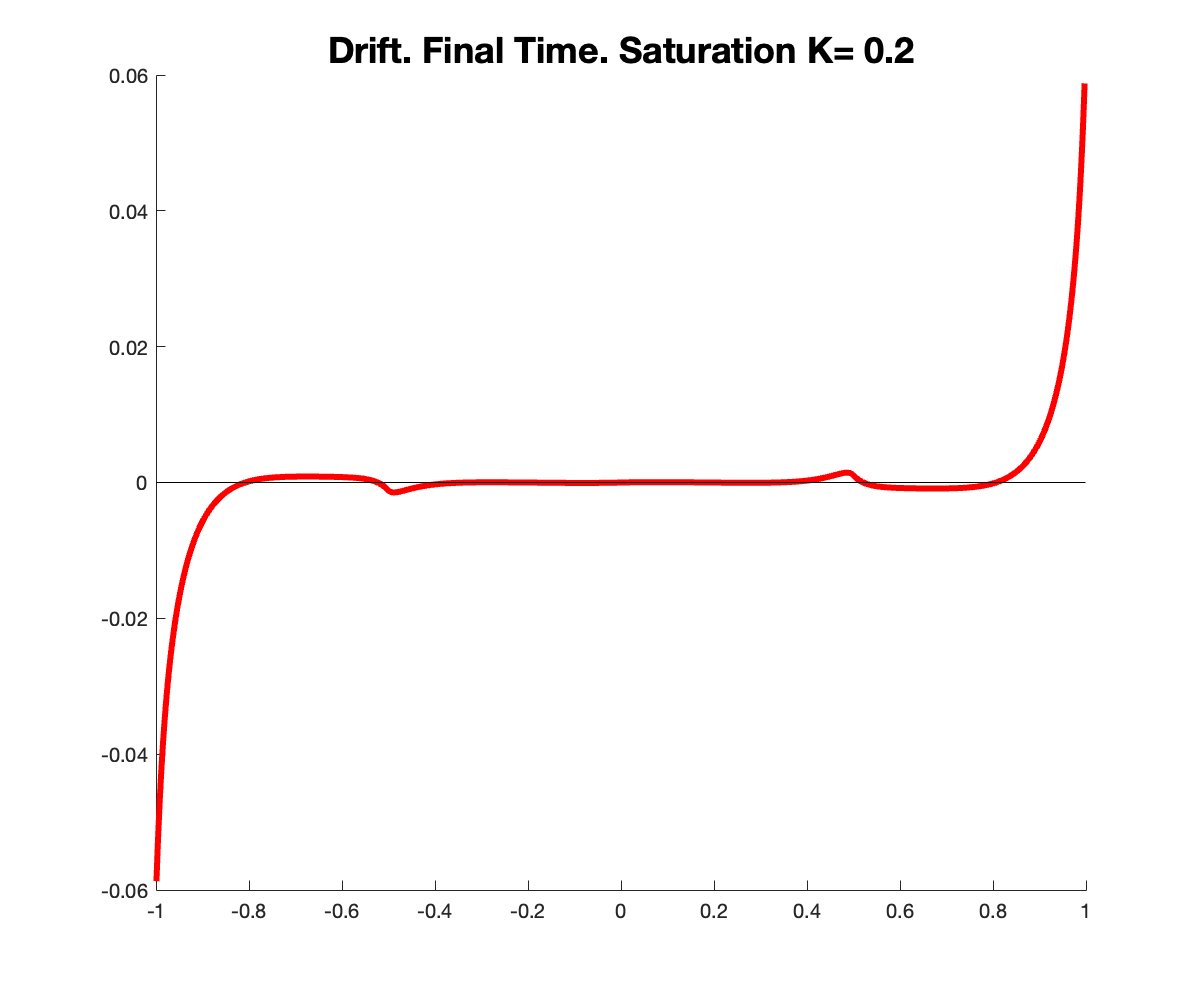}\\

\includegraphics[scale=0.26]{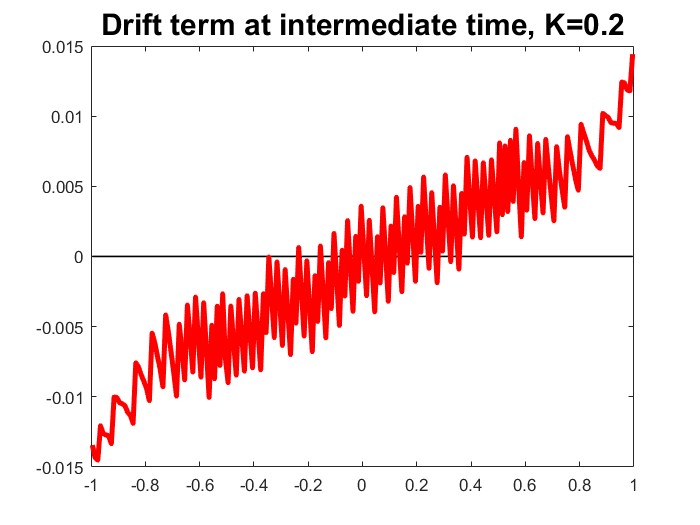}\!\!\includegraphics[scale=0.26]{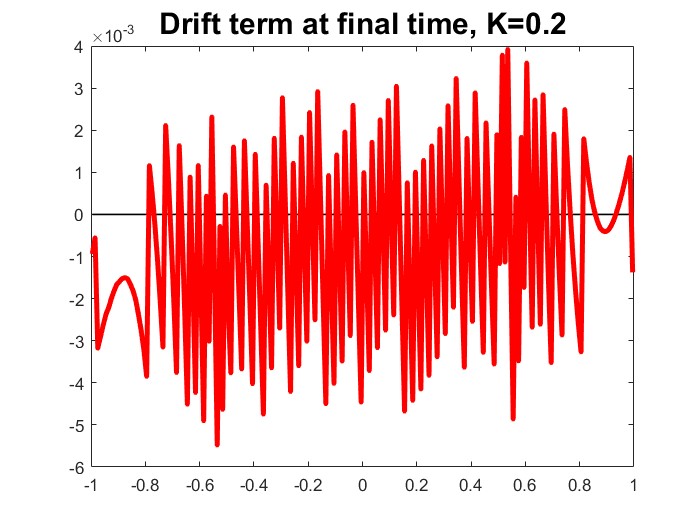}
\includegraphics[scale=0.15]{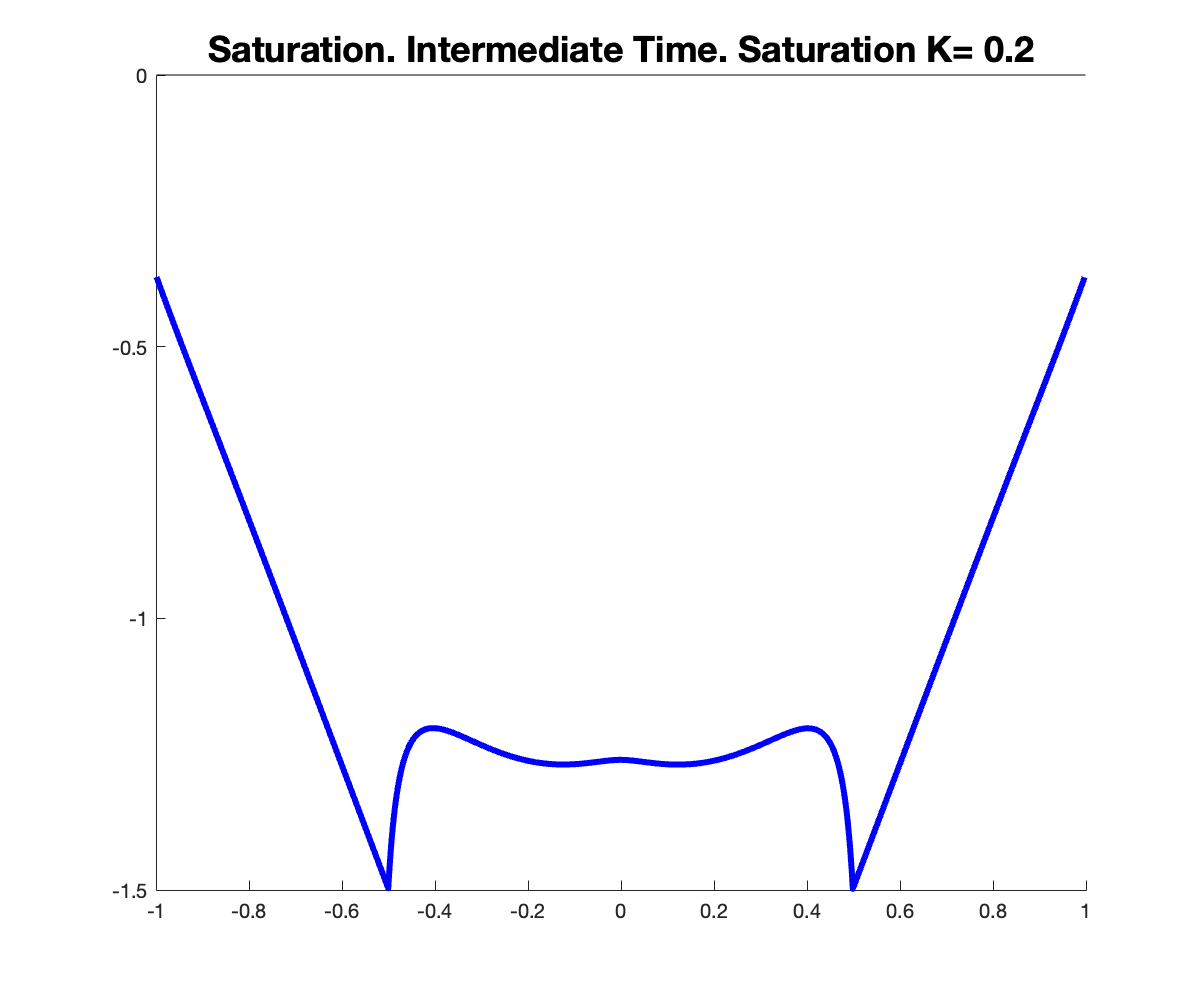}\!\!\!\!\includegraphics[scale=0.15]{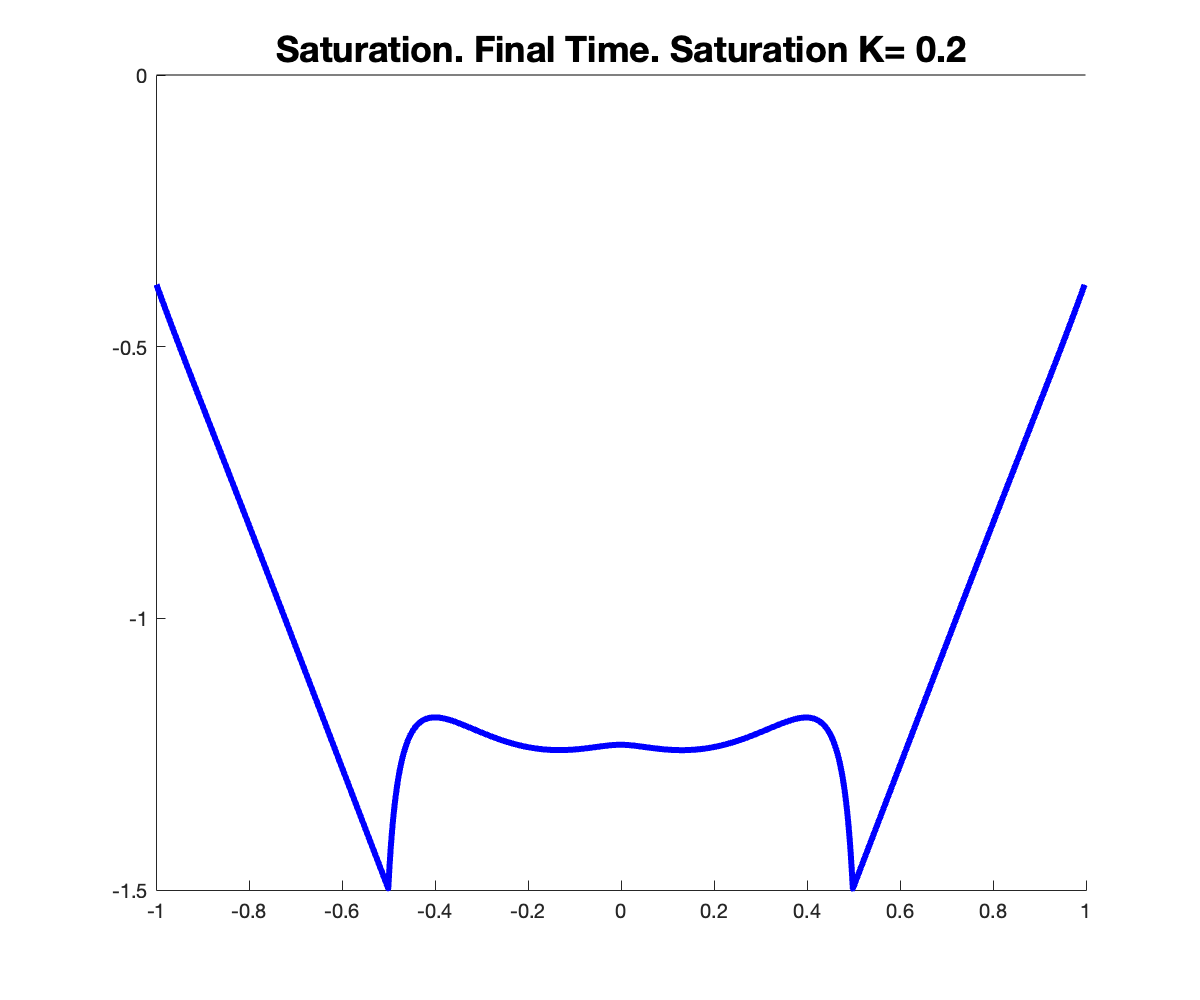}\\

\includegraphics[scale=0.26] {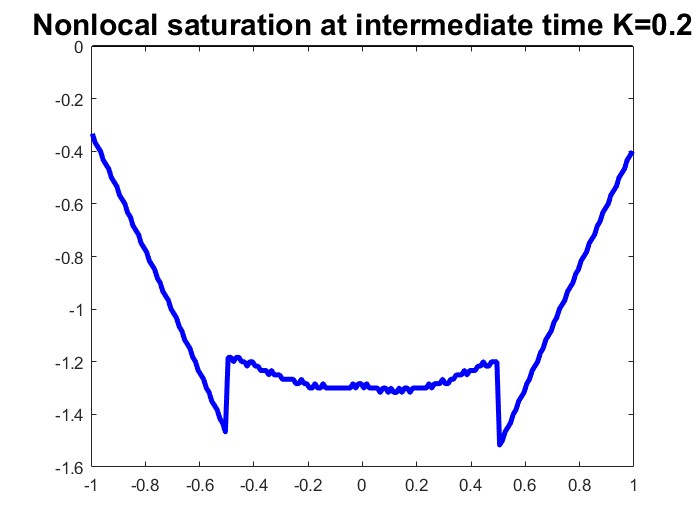}\includegraphics[scale=0.26]{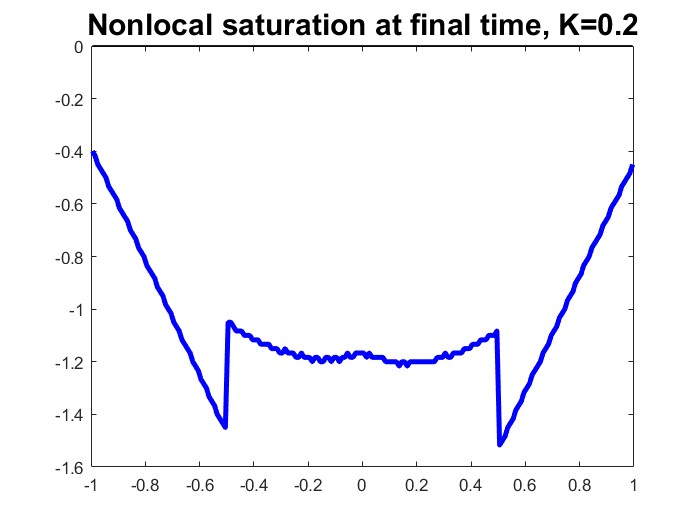}
\captionof{figure}{Drifts and saturations for nonlocal saturation model, $K=0.2$. Drifts in first and second rows for  PDE model (\S 3.5.3) and SDEs (\S 4.2), respectively. Saturations in third and fourth rows for  PDE model (\S 3.5.3) and SDEs (\S 4.2), respectively. All of them at intermediate time on the left and at final time on the right.} 
\label{fig:NL_K02_driftsatu}
\end{minipage}


   \begin{minipage}{\textwidth}
\centering
\includegraphics[scale=0.24]{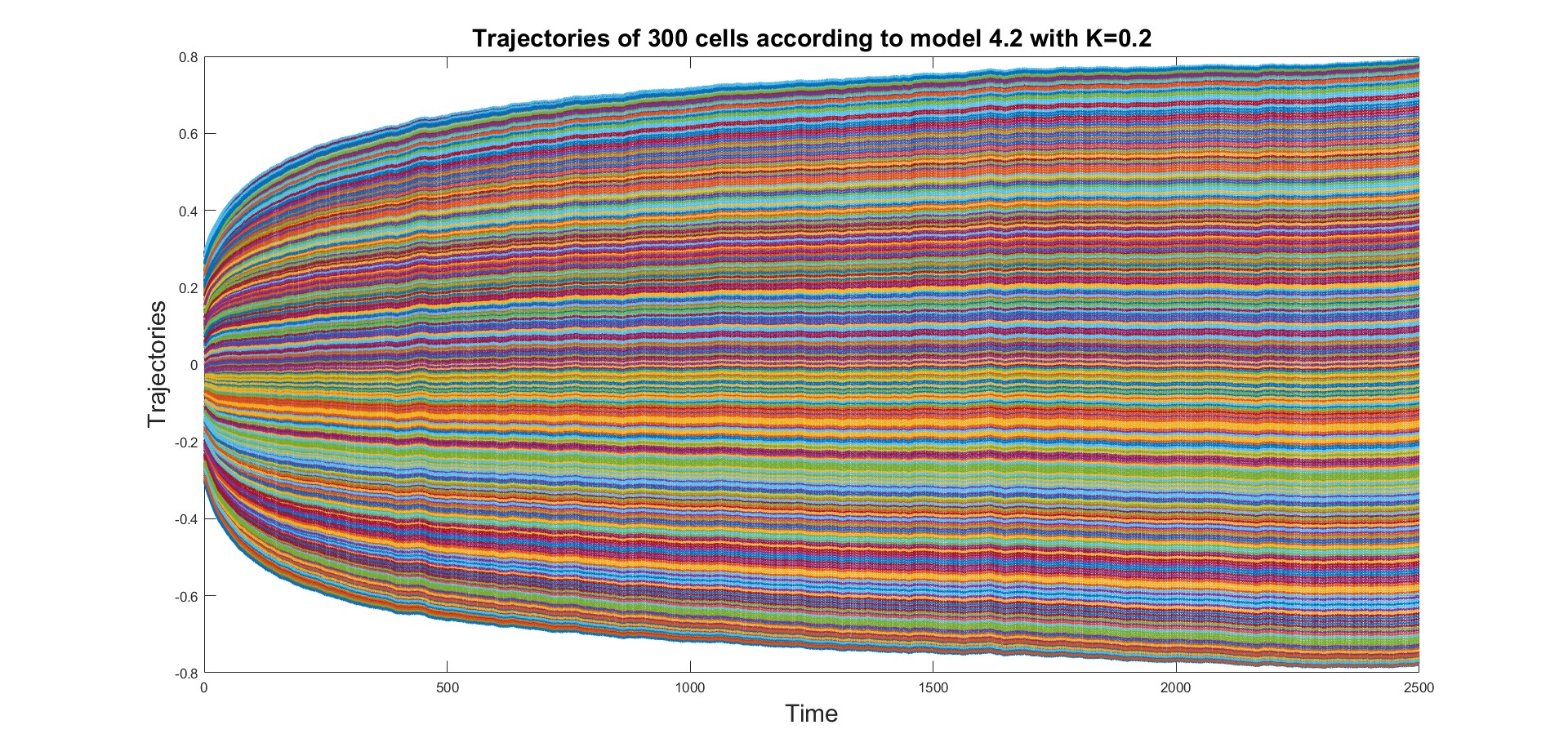}\\

\includegraphics[scale=0.18]{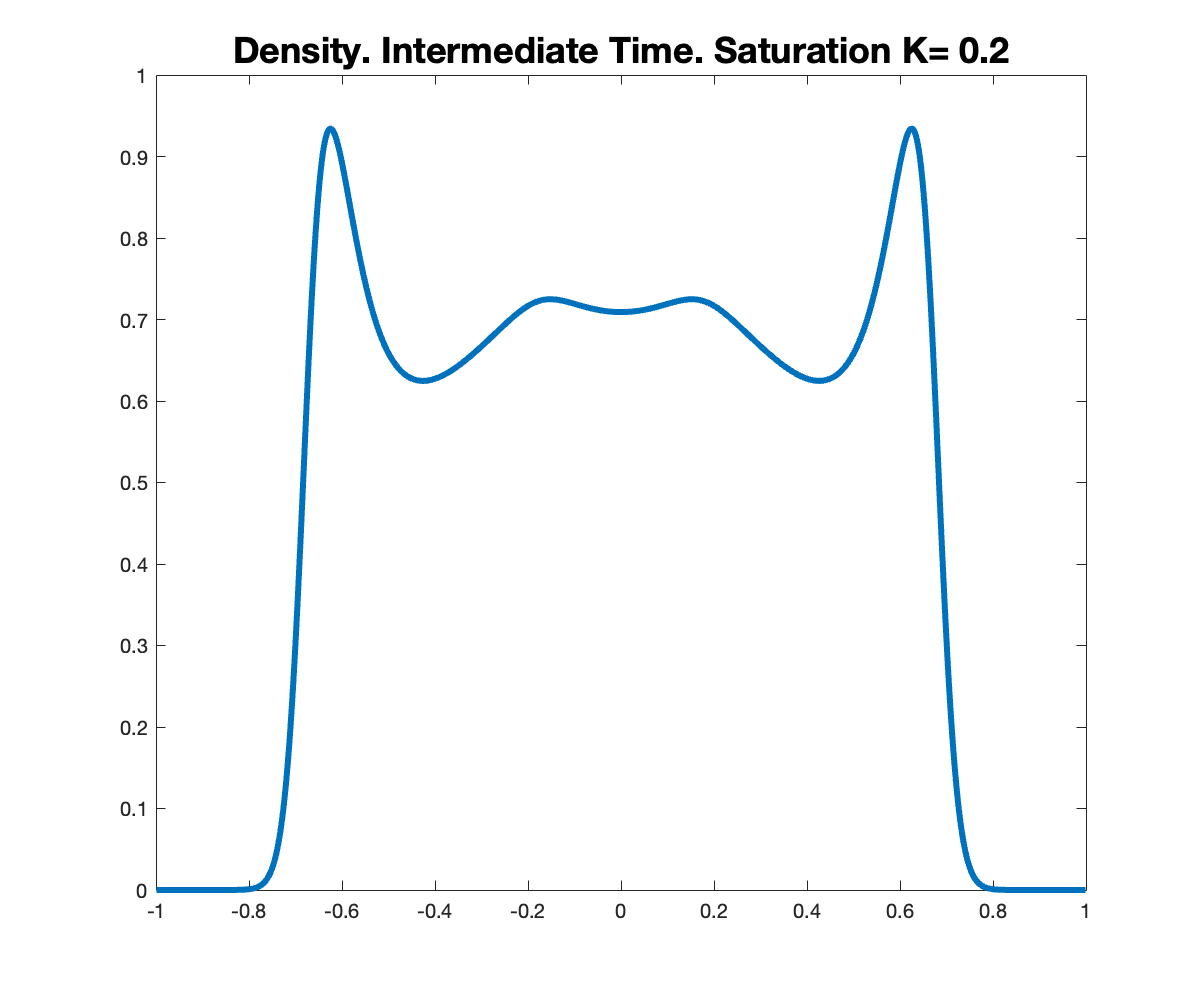}\!\!\includegraphics[scale=0.18]{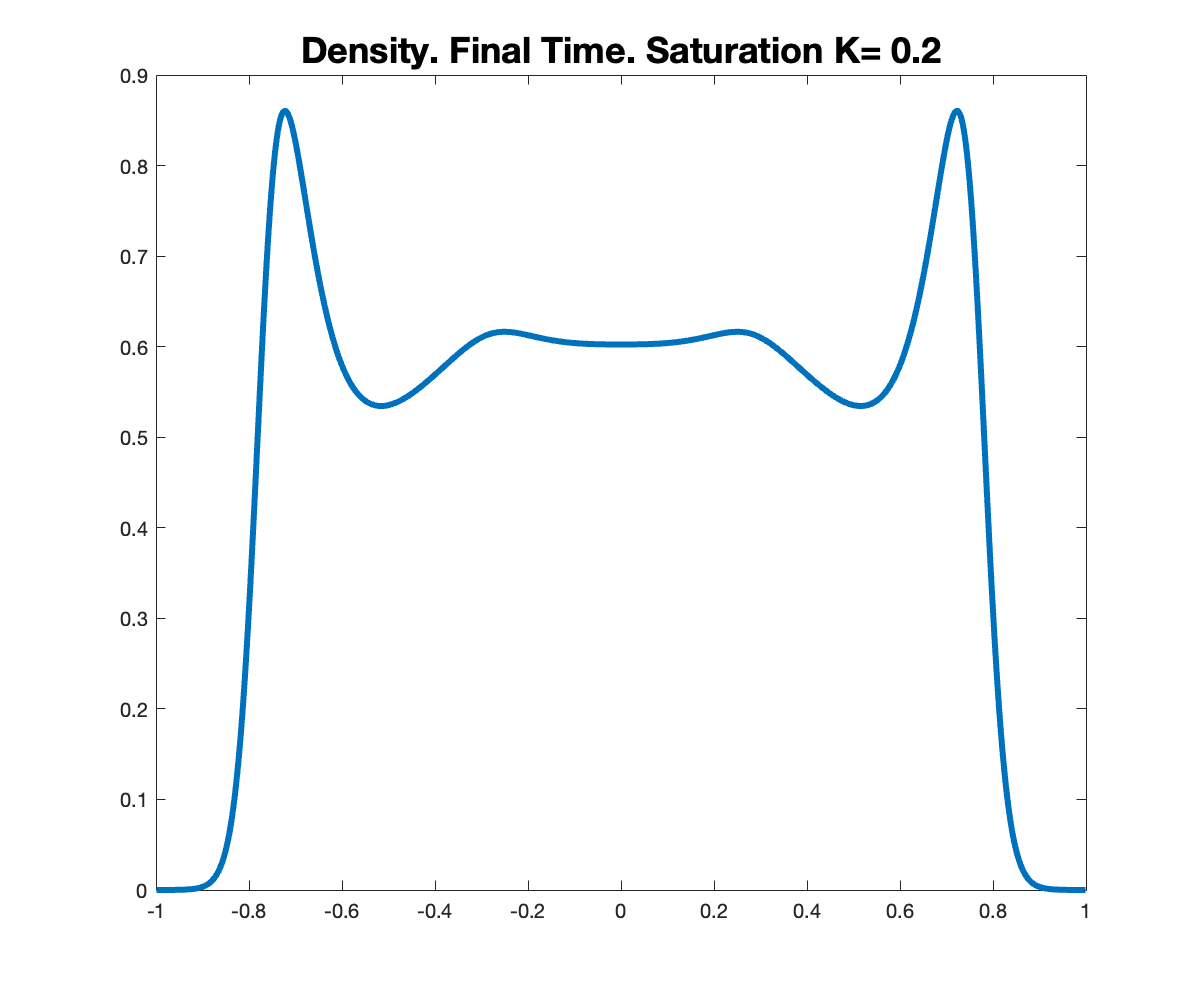}\\
\includegraphics[scale=0.3]{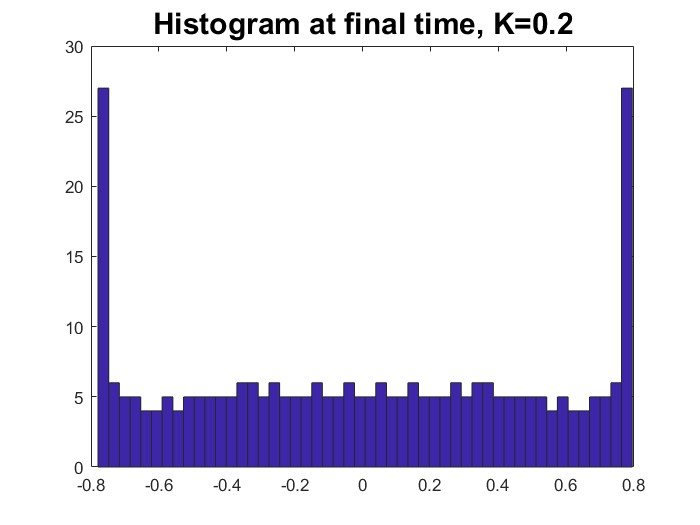}\!\!\includegraphics[scale=0.3]{DEFINITIVAS/Modelo3K02concentrado/mod3-hist-final-k02-con.jpg}
\captionof{figure}{Nonlocal saturation model, $K=0.2$, Concentrated initial data within $[-0.3,0.3]$. First row trajectories with $N=300$ cells, second row densities of the PDE model (\S 3.5.3) and third one histograms of SDEs (\S 4.2). Both at intermediate time on the left and at final time on the right.} 
\label{fig:NL_K02concentrado_histogramas}
\end{minipage}

 \begin{minipage}{\textwidth}
\centering
\includegraphics[scale=0.15]{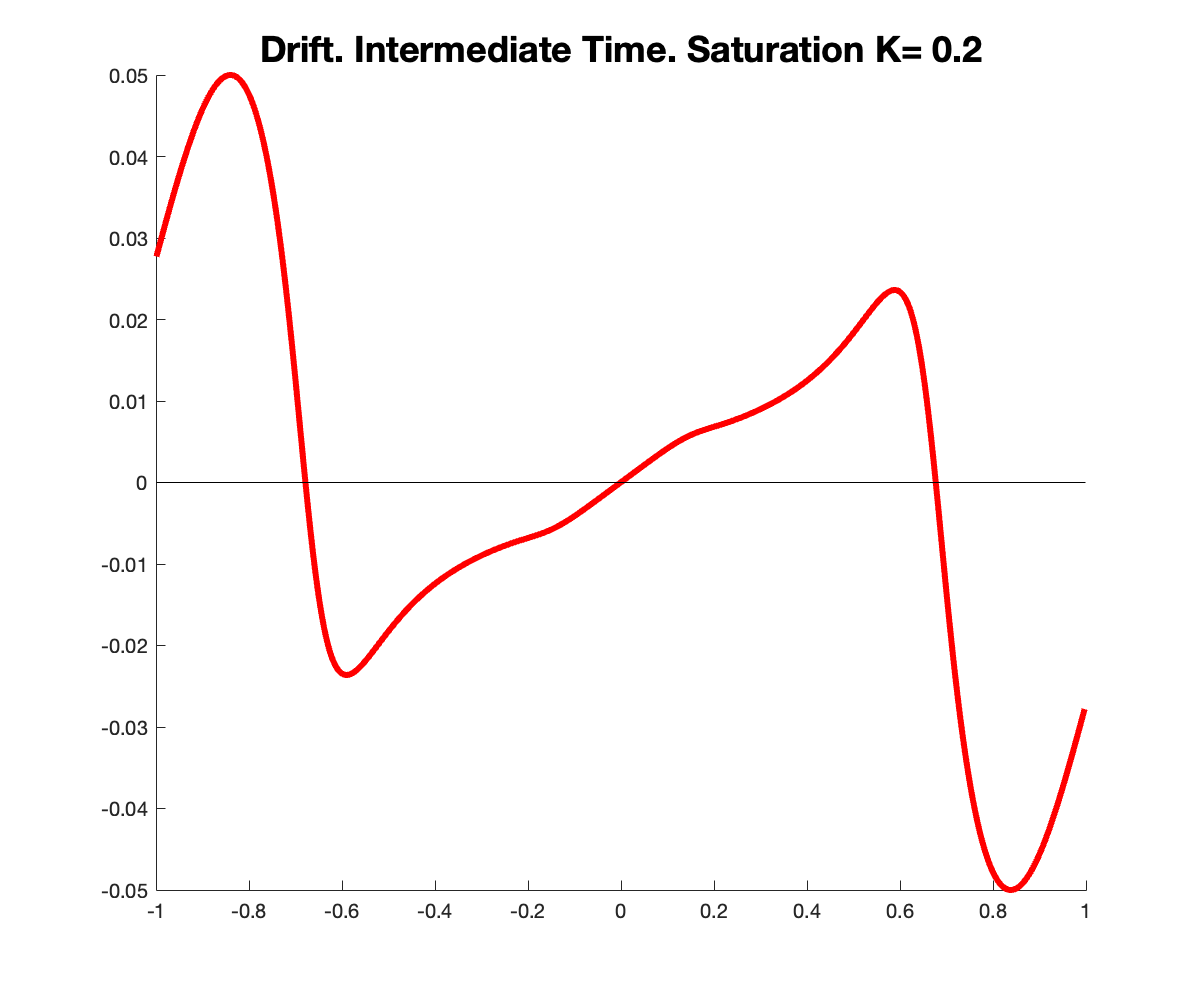}\!\!\!\!\includegraphics[scale=0.15]{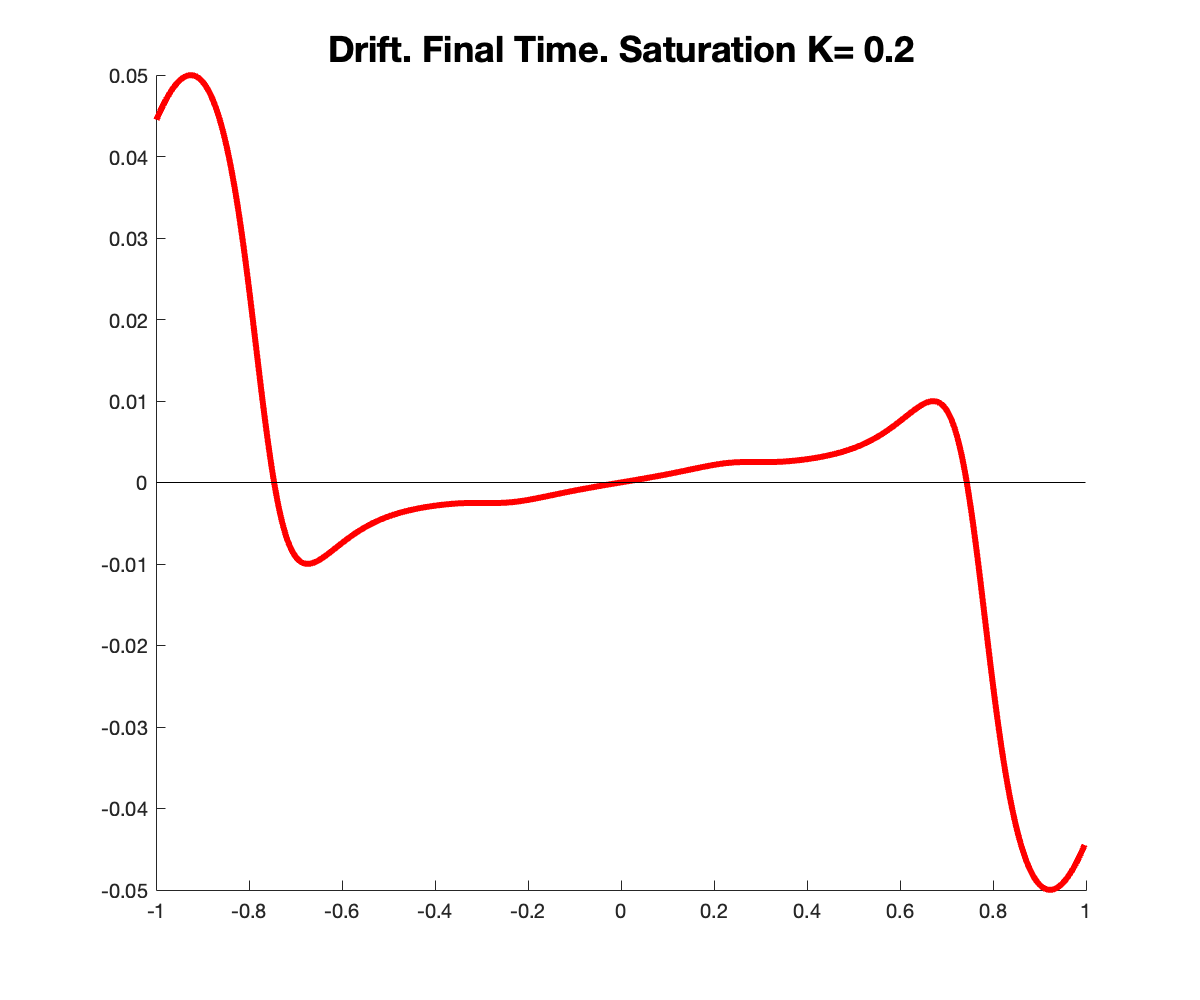}\\
\includegraphics[scale=0.26]{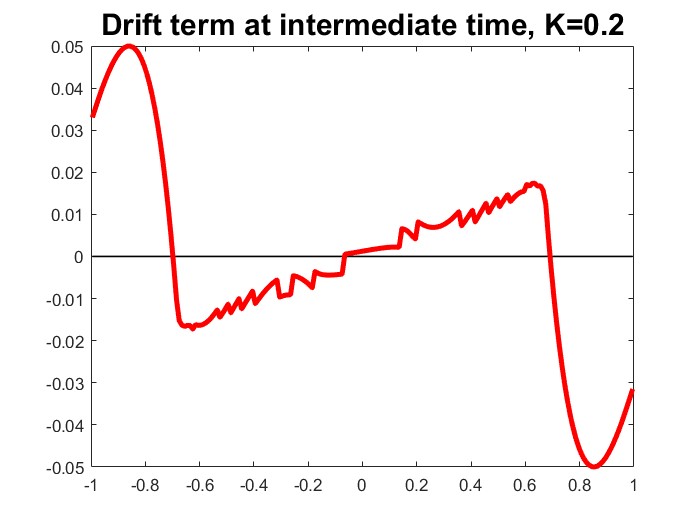}\!\!\includegraphics[scale=0.26]{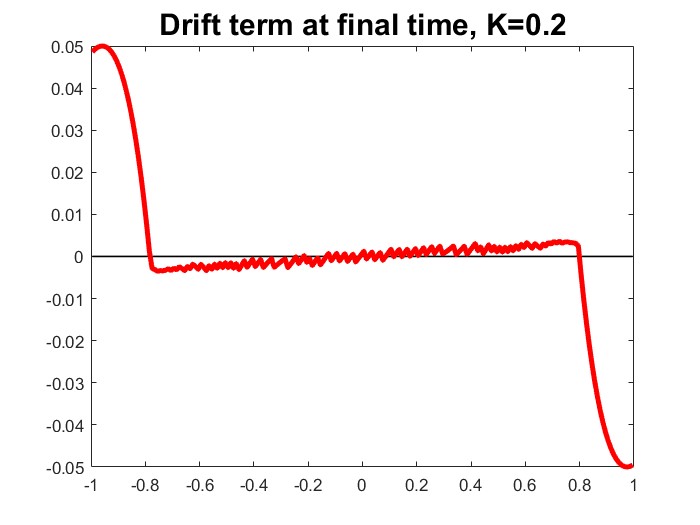}\\
\includegraphics[scale=0.15]{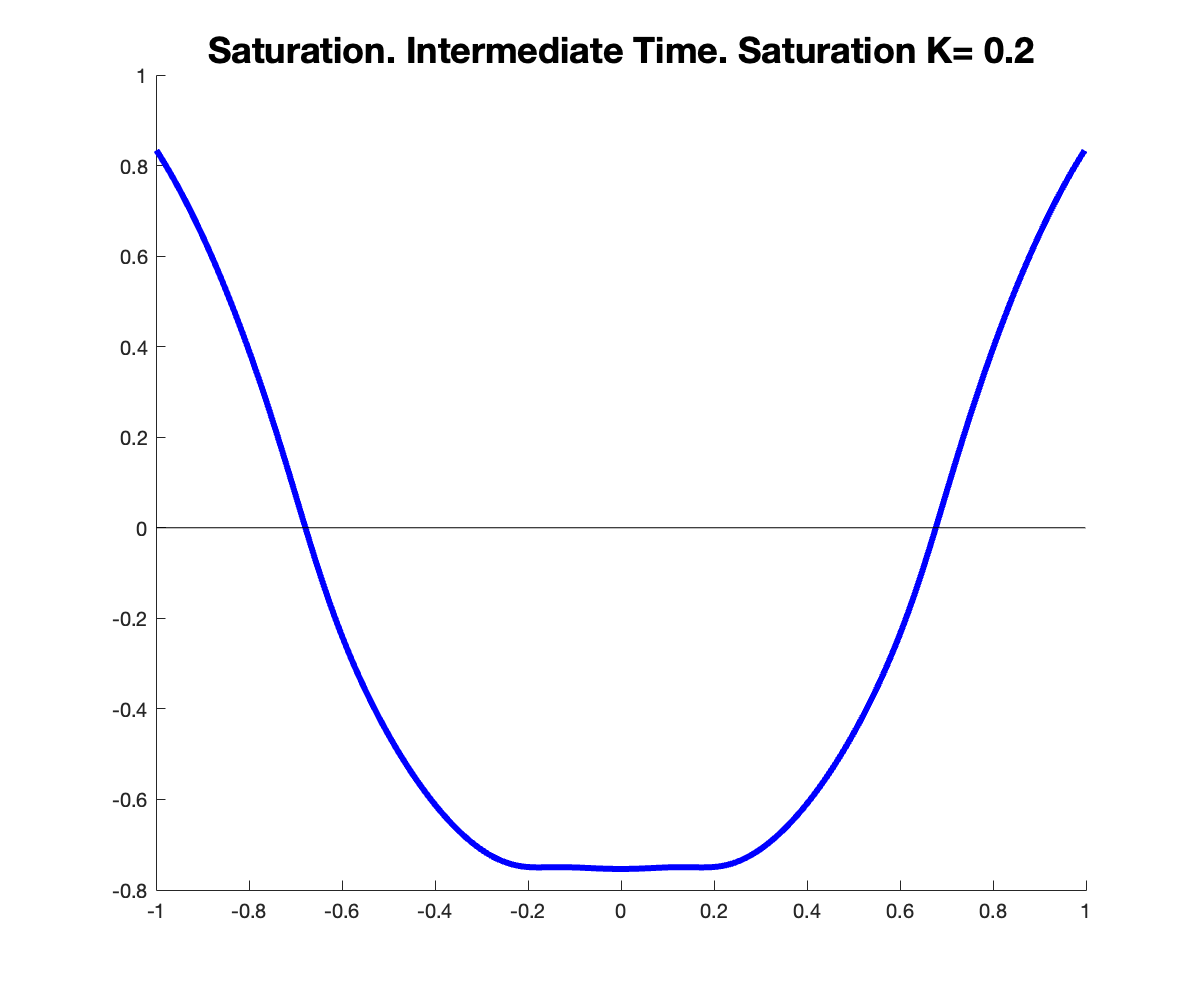}\!\!\!\!\includegraphics[scale=0.15]{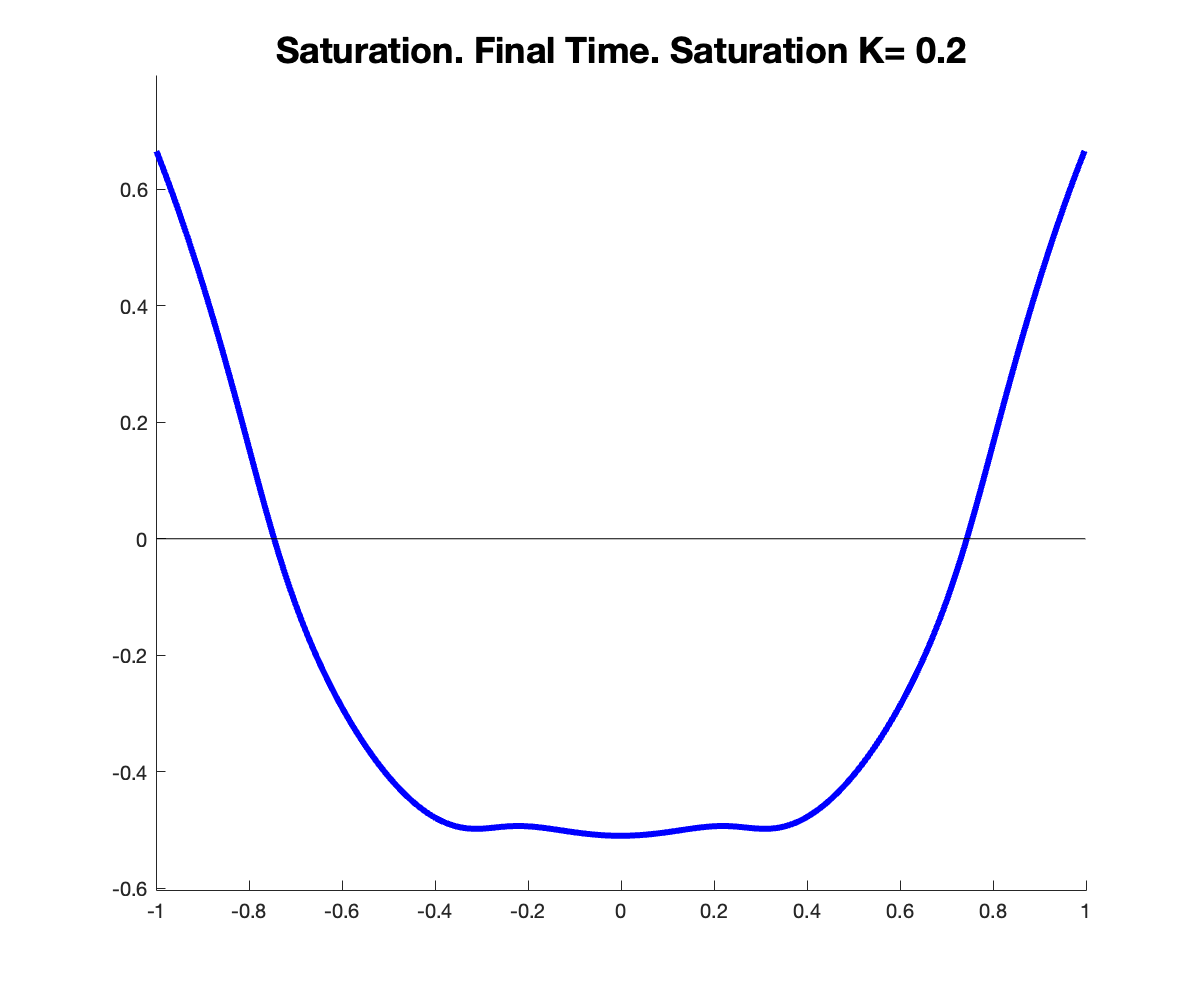}\\
\includegraphics[scale=0.26]{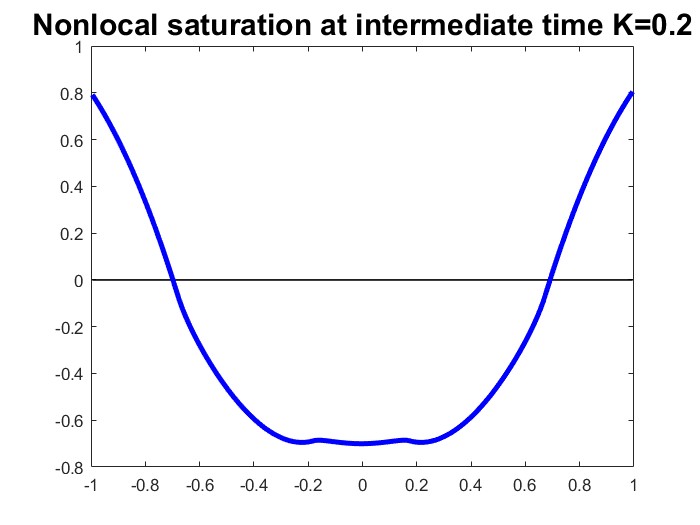}\includegraphics[scale=0.26]{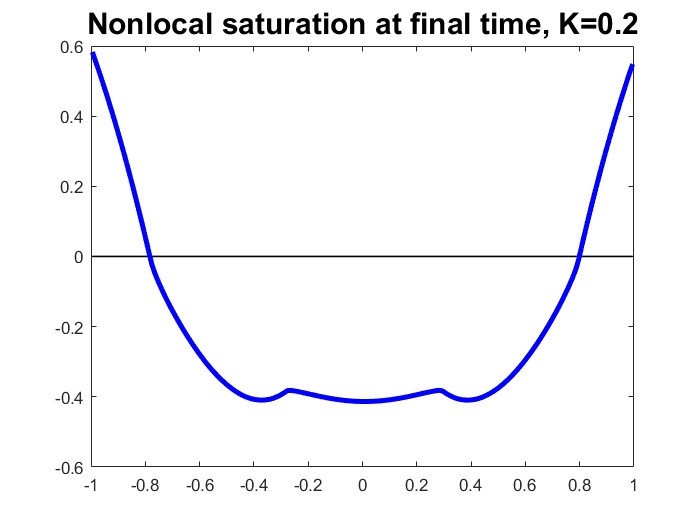}
\captionof{figure}{Drifts and saturations for nonlocal saturation model, $K=0.2$, Concentrated initial data within $[-0.3,0.3]$. Drifts in first and second rows for  PDE model (\S 3.5.3) and SDEs (\S 4.2), respectively. Saturations in third and fourth rows for  PDE model (\S 3.5.3) and SDEs (\S 4.2), respectively. All of them at intermediate time on the left and at final time on the right.} 
\label{fig:NL_K02concentrado_driftsatu}
\end{minipage}




\section{Conclusions and perspectives}\label{conclusiones}

In this work we propose two different discrete approaches to approximate cell-cell adhesion models. One the one hand, from Eulerian perspective as a balance of gain and loss terms deduced  in Section~3, on the other hand following particles trajectories in Section~4. In both cases, they describe identical dynamics when applied to the specific models: APS in \cite{APS} and the local saturation model studied in \cite{CMSTT}, but considering homogeneous boundary conditions, instead of periodic boundary conditions as in the mentioned papers. 

In addition, we also introduce a new cell-cell adhesion model, accounting for saturation effects  by multiplying the drift term by a nonlocal coefficient. Eulerian and Lagrangian approach also show similar dynamics. However, the behaviour of the solutions obeying nonlocal saturation effect clearly differs   with respect to the  models in \cite{APS} and \cite{CMSTT}. Indeed, our model  is suitable  to describe very different dynamics: it is capable to disseminate the cells along the whole domain, as well as it can recover  high aggregation peaks as in the APS model, passing by intermediate scenarios of moderate aggregation, just varying the values of $K$ and the weight $w$.

Moreover, our model includes saturation or even repulsion from too crowded areas as part of the interaction terms, in contrast to \cite{CMSTT} where this repulsion is driven by a porous media type diffusion. 

 As a continuation of this work, we wish to model cell motion due to external cues, that could modify their natural trajectories. The action of certain signals or  irregularities in the media could drive the cells to prescribed areas and produce long range motion. Some nonlocal operators can diffuse the particles in these anisotropic or prescribed way, instead of  arranging them averaged, and  produce nonlocal jumps. This anisotropic long range diffusion would compete with the interaction drift term, that should also accounts for saturation or repulsion from too populated areas.

\

\bigskip

\bigskip

\noindent{\bf Acknowledgements:} This research was initially conducted while CGF visited University of Seville thanks to the Feder project entitled {\it Diferentes Perspectivas para Modelos
Biomatem\'aticos: Modelizaci\'on, An\'alisis y Aproximaci\'on} (US-1381261). At present FGG and MPL  are supported by project PID2023-149182NB-I00, and AS is supported by project PID2023-149509NB-I00. The authors wish to thank deeply to both  referees  for their numerous and valuable comments and suggestions, which improved unquestionably the content and presentation of a previous version of the manuscript.  

\bibliographystyle{amsplain}
\bibliography{biblio}

\end{document}